\documentclass{amsart}

\usepackage[a4paper,margin=1.2in]{geometry}

\usepackage[english]{babel}
\usepackage[utf8]{inputenc}
\usepackage[T1]{fontenc}
\usepackage{amsmath}
\usepackage{amsfonts}
\usepackage{mathrsfs}
\usepackage{mathtools}
\mathtoolsset{showonlyrefs}

\usepackage{esint}
\usepackage[colorlinks=true, allcolors=red, linktocpage=true]{hyperref}
\usepackage{amssymb}
\usepackage{graphicx}
\usepackage{commath}
\usepackage{cases}
\usepackage{bm}
\usepackage{dsfont}
\usepackage{xcolor}
\usepackage{verbatim}
\usepackage{todonotes}

\usepackage{amsthm}

\theoremstyle{plain}

\newtheorem{thm}{Theorem}[section]
\newtheorem{prop}[thm]{Proposition}
\newtheorem{lem}[thm]{Lemma}

\newtheorem{cor}[thm]{Corollary}

\numberwithin{equation}{section}

\theoremstyle{definition}

\newtheorem{defn}[thm]{Definition}

\newtheorem{rem}[thm]{Remark}

\newtheorem{notation}[thm]{Notation}

\makeatletter
\newcommand*{\rom}[1]{\expandafter\@slowromancap\romannumeral #1@}
\makeatother

\newcommand{\ve}{\varepsilon}

\newcommand{\N}{\mathbb{N}}

\newcommand{\Pee}{\mathbb{P}}

\newcommand{\R}{\mathbb{R}}

\newcommand{\Z}{\mathbb{Z}}

\newcommand{\calB}{\mathcal{B}}

\newcommand{\calF}{\mathcal{F}}
\newcommand{\calH}{\mathcal{H}}

\newcommand{\calL}{\mathcal{L}}
\newcommand{\calM}{\mathcal{M}}

\newcommand{\calS}{\mathcal{S}}
\newcommand{\calT}{\mathcal{T}}

\newcommand{\diam}{\mathrm{diam}}
\newcommand{\spann}{{\mathrm{span}}}
\newcommand{\eucl}{{\mathrm{eucl}}}

\def\dist{{\mathop\mathrm{dist}}}
\def\supp{\mathop\mathrm{supp}}

\newcommand{\cubes}{\mathcal{D}}

\newcommand{\cuberatio}{A_0}

\newcommand{\Cseparation}{C_3}

\newcommand{\pdiam}{\diam_p}

\newcommand{\adim}{n}

\newcommand{\Tree}{\mathrm{Tree}}

\newcommand{\udensity}{A}
\newcommand{\ldensity}{\tau} 
\newcommand{\compdoub}{C_1}
\newcommand{\cubefactor}{C_0}

\def\X{\mathbf{X}}
\def\Y{\mathbf{Y}}
\def\bZ{\mathbf{Z}}
\def\W{\mathbf{W}}

\newcommand{\CA}{C_1}
\newcommand{\CF}{C_2}
\newcommand{\CFprime}{C_4}
\newcommand{\Csigma}{C_3}

\newcommand{\nocontentsline}[3]{}
\newcommand{\tocless}[2]{\bgroup\let\addcontentsline=\nocontentsline#1{#2}\egroup}

\usepackage{microtype}
\usepackage{tikz}
\usetikzlibrary{arrows.meta,calc,patterns,decorations.pathreplacing}
\usepackage{float}
\usepackage{hyperref}
\DeclareMathOperator{\Lip}{Lip}

\title{Characterization of parabolic rectifiability in terms of Jones' square function}

\author{Simon Bortz}

\author{Aakash Gurung}

\address{Simon Bortz \& Aakash Gurung, Department of Mathematics, University of Alabama, Tuscaloosa, AL 35487, USA}
\email{sbortz@ua.edu, agurung3@crimson.ua.edu}

\author{Pablo Hidalgo-Palencia}

\address{Pablo Hidalgo-Palencia, Departament de Matemàtiques i Informàtica, Universitat de Barcelona, Gran Via de les Corts Catalanes 585, 08007 Barcelona, Spain}
\email{pablo.hidalgo@ub.edu}

\author{Matthew Hyde}

\address{Matthew Hyde, University of Jyväskylä, Department of Mathematics and Statistics, P.O. Box 35 (MaD), FI-40014 University of Jyväskylä, Finland}

\email{matthew.j.hyde@jyu.fi}

\begin{document}

\begin{abstract}
    We show that a measure $\mu$ in $\R^n \times \R$ equipped with the parabolic metric is rectifiable (in the parabolic sense) if and only if the Jones' square function is finite $\mu$-almost everywhere. This square function is constructed using the parabolic version of the Jones' beta numbers, which control how close $\mu$ is to being supported in time-independent planes (these planes are the paradigmatic example of parabolic rectifiable set). Our only background assumption is that $\mu$ has positive and finite densities (qualitatively) almost everywhere.

    A key feature of our proof is the use of a parabolic version of the Reifenberg parametrization theorem, which in the Euclidean setting is due to David and Toro (MAMS, 2012), and has recently been extended to the parabolic setting by Hoffman, Mourgoglou and the fourth author. In the latter work, the authors construct a bi-H\"older map. Here, under an additional summability condition, we improve this construction to obtain a parabolic regular bi-Lipschitz map. This improvement in regularity plays a crucial role in our arguments.  
\end{abstract}

\maketitle

\vspace{-.6cm}

\tocless{\section*}{Disclaimer for this version}

The results in this manuscript are parabolic extensions of the Euclidean papers of Tolsa \cite{tolsa2015characterization} and Azzam and Tolsa \cite{azzam2015characterization}. In the advanced stages of writing this article, we learned of the preprint \cite{koirala2026parabolic} which announces a proof of the main result here (and actually, many other results). While the article \cite{koirala2026parabolic} does provide a Lean certificate for these results (based on 176000 lines of Lean code, according to the author), there are no real details about the strategy or the proof itself. 

\vspace{.1cm}

This article is a polished version of the parts of the manuscript which were already written at the time of the announcement \cite{koirala2026parabolic}. There are several new difficulties that arise only in the parabolic setting, and these are the ones we treat in detail in this manuscript, including detailed proofs of all the truly parabolic challenges. However, in order to offer a fully self-contained proof of our main result, we should include an abstract stopping-time arguments from \cite{azzam2015characterization} that carries over to the parabolic setting with essentially no change. These arguments have not yet been written out in full and are therefore only presented as a sketch.  

\vspace{.1cm}

While we believe this article offers genuine new insights into the main problem, and parabolic rectifiability in general, it was not our original intention to make the manuscript public before these remaining details had been completed. Unfortunately, we decided that we could not delay the upload any further; the missing details will appear in a revised version as soon as possible. We apologize for the possible inconveniences that this may offer to the reader, this is just one of the many unusual scenarios in which we find ourselves as we transition from the pre-AI to the post-AI world.

\setcounter{tocdepth}{1} % hide subsections in table of contents
\tableofcontents

\section{Introduction}

Rectifiability is a central concept in modern geometric measure theory and, more generally, analysis. It has deep connections with different areas, like harmonic analysis: for instance, in \cite{nazarov2014riesz}, it is shown that the if Riesz transforms are bounded in $L^2$ weighted by the surface measure of a set, then such set is rectifiable. This, in turn, has implications in PDE theory: in \cite{azzam2016rectifiability} it was shown that the subset of the boundary of a domain in which the harmonic and surface measure are absolutely continuous, is in fact rectifiable. Similarly, in free boundary problems, it is often the case that lower dimensional strata of the singular set of the free boundary are rectifiable, and this is very informative about their structure; see e.g. \cite{edelen2019quantitative, kriventsov2025rectifiability}. 

\vspace{.1cm}

In the Euclidean setting, a finite Borel measure $\mu$ is $(n-1)$-rectifiable if it is absolutely condition with respect to Hausdorff measure $(n-1)$-dimensional Hausdorff measure and if we can cover $\mu$-a.e. its support with a countable number of Lipschitz images of $\R^{n-1}$. Tolsa \cite{tolsa2015characterization} and Azzam-Tolsa \cite{azzam2015characterization} showed that (under background density assumptions) rectifiability  of $\mu$ is equivalent to the square function associated to the Jones' beta numbers of $\mu$ being finite $\mu$-a.e. (these coefficients measure how close is $\mu$ to being supported over planes, at all scales; we will present them more rigorously soon).  

\vspace{.1cm}

In connection with parabolic PDE, rectifiability necessarily takes a different form. Indeed, the natural metric of space-time $\R^n \times \R$ that respects the natural homogeneity/invariance of the heat equation (the equation is invariant under the transformation $(X, t) \mapsto (\lambda X, \lambda^2 t)$ for $\lambda > 0$), is not just the Euclidean metric in $\R^{n+1}$. Under this metric, parabolic Lipschitz functions are Lipschitz in space, but only Hölder-$1/2$ in time. And it turns out that covering a set (like $\supp \mu$) by images of these parabolic Lipschitz functions is not enough to provide meaningful applications (see e.g. the counterexample of \cite{kaufman1980singularity}). 

\vspace{.1cm}

After the quantitative results of Lewis and Murray \cite{lewis1995method}, and Hofmann, Martell, Nyström and the first author \cite{bortz2025solvability}, and the qualitative result of Egert, Ferris, Saari and the first author \cite{bortz2026parabolic} (among several other relevant related works), it is now understood what is the correct notion of parabolic rectifiability in order to obtain the analogous PDE applications as in the Euclidean case. Indeed, not only parabolic Lipschitzianity is needed, but an additional \textit{regularity} in time (certain control on the oscillations of the half-order time derivatives, see Section~\ref{sec:flatness}), and concretely, parabolic flatness means being close to vertical (that is, time-independent) planes (as opposed to arbitrary planes in the Euclidean setting). This imposes real new difficulties when trying to show that a set or a measure is parabolic rectifiable. In this paper, we overcome these difficulties to show that a measure is parabolic rectifiable if and only if the associated parabolic Jones square function is finite $\mu$-a.e. (under background density hypotheses), thus, adapting \cite{tolsa2015characterization, azzam2015characterization} to the parabolic setting.

\subsection{Main result}

We quantify how flat a measure $\mu$ in space-time $\R^n \times \R$ is by using the parabolic version of the beta numbers originally introduced by Jones in \cite{jones1990rectifiable}. Beta numbers measure, at a concrete scale and location, and in a scale-invariant fashion, how close the set is to a \textit{vertical} plane:
\begin{equation*}
    \beta_\mu((X, t), r)
    :=
    \inf_P \left( \frac{1}{r^{n+1}} \int_{B_p((X, t),r)} \left( \frac{\dist_p((Y, s), P)}{r} \right)^2 \, d\mu(Y, s) \right)^{1/2},
\end{equation*}
where $P$ runs over the family of affine planes in $\R^{n+1}$ which contain a line in the time direction, the so-called \textit{vertical} or \textit{$t$-independent} planes. Vertical planes are the prototypical example of a parabolic rectifiable set. For this and more notation, like parabolic distance $\dist_p$ or cylinders $B_p$, we direct the reader to Section~\ref{sec:setting}. 

In light of this definition, our main result reads as follows.

\begin{thm}[Summability of beta numbers is equivalent to rectifiability] \label{thm:main}
    Let $\mu$ be a finite Borel measure in space-time $\R^n \times \R$ endowed with the parabolic metric, and assume that 
    \begin{equation*}
        0 
        < 
        \Theta^{n+1, *}((X, t), \mu)
        :=
        \limsup_{r \to 0} \frac{\mu(B_p((X, t),r))}{r^{n+1}} 
        < 
        + \infty, 
        \quad 
        \text{for $\mu$-a.e. $(X, t) \in \R^n \times \R$.}
    \end{equation*}
    Then, the square function of the beta numbers is finite $\mu$-a.e., that is, 
    \begin{equation}\label{e:finite-jones-beta-main-thm}
        \int_0^1 \beta_\mu((X, t), r)^2 \frac{dr}{r} < + \infty, 
        \qquad 
        \text{for $\mu$-a.e. $(X, t) \in \R^n \times \R$,}
    \end{equation}
    if and only if $\mu$ is parabolic rectifiable, that is, $\mu \ll \calH^{n+1}_p$ and there exist countably many regular Lip(1,1/2) graphs $\Gamma_j$ such that
    \[
    \mu \left( \mathbb R^{n+1} \setminus \bigcup_{j=1}^{\infty} \Gamma_j \right) 
    =
    0. 
    \]
\end{thm}

For the definition of regular Lip(1,1/2) graphs, we refer to Section~\ref{sec:flatness}.

\begin{rem}[Homogeneities and dimension]
    As a quick note to understand the homogeneities in our main result: the prototypical parabolic rectifiable set is a vertical plane, such as $P := \{ (X, t) \in \R^n \times \R : X_n = 0\}$. A good example of measure to test our theorem would be $\mu := \mathcal{H}^{n+1}_p \lfloor P$. Indeed, although $P$ has topological dimension $n$, the fact that we have endowed space-time with the parabolic metric (see Section~\ref{sec:setting}) makes $P$ have homogeneous dimension $n+1$. In the same way, the ambient space-time $\R^n \times \R$ has homogeneous dimension $n+2$.
\end{rem}

\begin{rem}[Background assumptions]
    The implication that rectifiability yields $\mu$-a.e. finiteness of the beta square function does not require the background assumption on the density, as was that case in the Euclidean counterpart, see \cite{tolsa2015characterization}.
\end{rem}

%%%%%%%%%%%%%%%%%%%%%%%%%%%%%%%%%%%%%%%%%%%%%%%%%%%%%%%%%%%%%%%%%%%%%%
%%%%%%%%%%%%%%%%%%%%%%%%%%%%%%%%%%%%%%%%%%%%%%%%%%%%%%%%%%%%%%%%%%%%%%
%%%%%%%%%%%%%%%%%%%%%%%%%%%%%%%%%%%%%%%%%%%%%%%%%%%%%%%%%%%%%%%%%%%%%%
%%%%%%%%%%%%%%%%%%%%%%%%%%%%%%%%%%%%%%%%%%%%%%%%%%%%%%%%%%%%%%%%%%%%%%
%%%%%%%%%%%%%%%%%%%%%%%%%%%%%%%%%%%%%%%%%%%%%%%%%%%%%%%%%%%%%%%%%%%%%%

\subsection{Sketch of the proof and organization of the article}

After a section giving the necessary definitions and preliminaries (Section~\ref{sec:preliminaries}), the article is divided into two parts. In the first part, we prove the implication that finiteness of the square function of the beta numbers implies rectifiability (that is, the parabolic version of the theorem of Azzam and Tolsa \cite{azzam2015characterization}): the proof is given across Sections~\ref{s:Reif}, \ref{sec:regular} and \ref{sec:Jones-implies-rectifiability}. The second part deals with the converse implication that rectifiability implies finiteness of the Jones' square function (that is, the parabolic version of the theorem of Tolsa \cite{tolsa2015characterization}): the proof is given in Sections~\ref{sec:preliminaries_2}, \ref{sec:graph} and \ref{proof-of-the-main-theorem}.\\

\noindent\textbf{First implication.} The proof that finiteness of the beta square function implies parabolic rectifiability, requires harder work to be adapted from the Euclidean analogue of Azzam and Tolsa \cite{azzam2015characterization}. Indeed, the main problem when extending Azzam-Tolsa's proof to the parabolic setting is that we need to prove that the graphs $\Gamma_j$ are not only (parabolic) Lipschitz, but also \textit{regular} (in the sense of Definition~\ref{def:regular}), which is crucial to prove parabolic rectifiability. 

To prove that the approximating graphs are regular Lip(1, 1/2), one needs to strengthen the construction in \cite[Subsection 6.2]{azzam2015characterization}, in which these graphs are constructed using a theorem of David-Toro \cite{david2012reifenberg}, inspired by earlier work of Reifenberg. Given a grid of nested balls (associated to a stopping-time regime in which the densities of $\mu$ and the $\beta$ numbers are under control) and the associated best approximating planes (which form controlled angles with their neighbors because of the control on the $\beta$ numbers), the David-Toro argument constructs a bi-H\"older map by carefully \textit{gluing together} these planes across all scales. Under an additional summability condition on the angles between consecutive planes, one can guarantee that this map is bi-Lipschitz, which is what is used in \cite{azzam2015characterization}. The bi-H\"older version of this theorem has recently been implemented in the parabolic setting by Hoffman, Mourgoglou and the fourth author \cite{hoffman2026parabolic}. In this paper we deal with the bi-Lipschitz result. Moreover, we show how to obtain a map which is \textbf{regular} parabolic bi-Lipschitz. 

More precisely, in Section~\ref{s:Reif}, we show that the construction from \cite{hoffman2026parabolic} can be imported to our setting (see Theorem~\ref{t:DT}, and Subsection~\ref{sec:lip}), and actually can be enhanced to a bi-Lipschitz map (see Theorem~\ref{t:bi-Lip}). The key idea is that, as happened in the Euclidean setting, the control on the beta square function gives a control over how much the angles can change and how those changes accumulate over scales. So under the assumption of smallness of the beta square function, the constructed parametrization can be chosen not to oscillate too much, which ultimately gives bi-Lipschitzianity. This argument is carried out along Subsection~\ref{sec:bilip}. We note that, as well as this bi-Lipschitz construction, we show that the sequence of approximating set i.e., the images of the functions in the approximation, can be written locally as Lipschitz graphs with precise estimates on the first-order time and second order spatial derivatives (see Lemma \ref{l:regular-graphs}). Although, we found a way \textit{in this case} to prove parabolic regularity without this, we believe that such properties will be essential in future works. 

In Section~\ref{sec:regular}, we show that the constructed graphs are actually \textbf{regular} Lip(1, 1/2), at least in appropriate stopping-time regimes; see Proposition~\ref{p:construct-regular-images}. After some technical work adapting the arguments of \cite{azzam2015characterization} to the parabolic setting, we arrive at Subsection~\ref{subsec:regular}, where regularity of the map is shown (Proposition~\ref{p:f-regular-Lipschitz}), which is the main novelty. We need to control how the beta numbers (or even better, the Omega numbers) accumulate over scales: after separating several easier cases, like scales which are much larger (Lemma~\ref{l:L}), smaller (Lemma~\ref{l:far}) or very far (Lemma~\ref{l:F}) from the scales associated to the stopping-time regime under consideration, we are left to estimate what happens precisely at those relevant scales. Upon establishing a precise relationship between the stopping-time cubes and the geometrical scales in Lemma~\ref{l:correspondence-CCBP}, we exploit further properties of the construction and the stopping-time regime to show the required bound for the beta numbers in Lemma~\ref{e:omega-estimate-global}.

Finally, in Section~\ref{sec:Jones-implies-rectifiability}, once we know that the graphs constructed via the Reifenberg parametrization scheme are regular parabolic bi-Lipschitz, we can simply follow the rest of the abstract arguments from \cite{azzam2015characterization} to finish the proof of the implication in Theorem~\ref{thm:main}.\\

\noindent\textbf{Second implication.} 
The proof of the implication of Theorem~\ref{thm:main} in which parabolic rectifiability implies finiteness of the square function of the beta numbers is an adaptation to the parabolic setting of the Euclidean work of Tolsa in \cite{tolsa2015characterization}. We start by some preliminary notation (in Section~\ref{sec:preliminaries_2}) which sometimes differs from the previous sections.

In Section~\ref{sec:graph}, we prove a related version of the required implication, which instead of considering the square function of the beta numbers, uses the alpha numbers instead. The main result is Theorem~\ref{thm:graph-alpha}, in which the square function of the $\alpha_\mu$ numbers (measuring flatness of $\mu$) are seen to be $\sigma$-a.e. finite, where $\sigma$ is the surface measure of a regular Lip(1,1/2) graph $\Gamma$. In order to prove this, one first restricts to the case where $\mu$ is not just an arbitrary measure, but it supported over $\Gamma$ (Subsection~\ref{the-supported-case}): the proof in this case is very harmonic analytical in nature, using a Calderón-Zygmund decomposition of $\mu$ with respect to $\sigma$ to handle cancellations. For the general case where $\mu$ is not necessarily supported over $\Gamma$, we construct (in Subsection~\ref{the-approximating-measure}) a measure $\eta$ supported over $\Gamma$ which somehow resembles $\mu$ (we use Whitney cubes to encapsulate the behavior at different scales, and then project them towards $\Gamma$), and then show that the alpha numbers of $\mu$ are controlled with respect to those for $\eta$ up to some error terms (in Subsection~\ref{the-good-cubes}). Since $\eta$ is supported over $\Gamma$, Subsection~\ref{the-supported-case} provides its control; and we handle the rest of error terms in Subsection~\ref{the-main-lemma-proof}.

In Section~\ref{proof-of-the-main-theorem}, we prove that, whenever the alpha-square function is finite $\mu$-a.e., then so is the beta-square function, which jointly with the previous section, gives the full implication of Theorem~\ref{thm:main}. To prove the main result of this section, namely Proposition~\ref{prop:alpha-beta}, we perform a stopping-time argument, restricting to cubes with controlled densities and small alpha numbers (Subsection~\ref{s:stopping-time}). In Subsection~\ref{s:angles}, we control angles between approximating planes at neighboring scales, and iterating these along chains of cubes in the stopping-time regimes that finish in a good cube. Then, in Subsection~\ref{s:proof_alpha-beta}, we use the results from the last section and our choice for the stopping-time regimes, to first show that they imply bounds for the beta-square function restricted to the good cubes, and later to all the cubes in the regime, which finishes the proof of Proposition~\ref{prop:alpha-beta}.

%%%%%%%%%%%%%%%%%%%%%%%%%%%%%%%%%%%%%%%%%%%%%%%%%%%%%%%%%%%%%
%%%%%%%%%%%%%%%%%%%%%%%%%%%%%%%%%%%%%%%%%%%%%%%%%%%%%%%%%%%%%
%%%%%%%%%%%%%%%%%%%%%%%%%%%%%%%%%%%%%%%%%%%%%%%%%%%%%%%%%%%%%
%%%%%%%%%%%%%%%%%%%%%%%%%%%%%%%%%%%%%%%%%%%%%%%%%%%%%%%%%%%%%
%%%%%%%%%%%%%%%%%%%%%%%%%%%%%%%%%%%%%%%%%%%%%%%%%%%%%%%%%%%%%
%%%%%%%%%%%%%%%%%%%%%%%%%%%%%%%%%%%%%%%%%%%%%%%%%%%%%%%%%%%%%
%%%%%%%%%%%%%%%%%%%%%%%%%%%%%%%%%%%%%%%%%%%%%%%%%%%%%%%%%%%%%
%%%%%%%%%%%%%%%%%%%%%%%%%%%%%%%%%%%%%%%%%%%%%%%%%%%%%%%%%%%%%
%%%%%%%%%%%%%%%%%%%%%%%%%%%%%%%%%%%%%%%%%%%%%%%%%%%%%%%%%%%%%

\newpage 

\section{Preliminaries} \label{sec:preliminaries}

We will utilize the customary notation that $C,c > 0$ stands for a constant that may change from line to line, and $X \lesssim Y$ is an equivalent notation to $X \leq C Y$. We sometimes make explicit what the implicit constant depends on using the notation $X \lesssim_A Y$. In turn, $X \approx Y$ mean $X \lesssim Y \lesssim X$, and $X \approx_{\kappa} Y$ is defined in the obvious way. 

\subsection{Parabolic setting} \label{sec:setting}

\begin{itemize}
    \item (Points) We denote points in space-time $\R^{n+1}$ by $(X,t) \in \R^\adim \times \R$.

    \item (Parabolic norm) Let $|\cdot|$ denote the standard Euclidean norm on $\R^\adim$ and define the parabolic norm on $\R^{\adim+1}$ by
    \[\|(X,t)\| \coloneqq |X| + |t|^\frac{1}{2} .\]

    \item (Parabolic distance) For $(X,t),(Y,s) \in \R^{\adim+1},$ define 
    \[ d_p(X,t,Y,s) := d_p((X,t),(Y,s)) := \|(X,t) - (Y,s)\|.\]
    Diameters of sets with respect to this metric are denoted by $\diam_p$.

    \item (Parabolic cylinders) For $(X,t) \in \R^{\adim+1}$ and $r > 0$, let 
    \[ B_p(X,t,r) := B_p((X, t), r) := \big\{ (Y,s) \in \R^{\adim+1} \colon |X-Y| \leq r \mbox{ and } |t-s| \leq r^2 \big\}. \] 
    We sometimes denote this cylinder by $B_r(X, t)$, too. We write $r(B)$ for the radius of the cylinder, and  $\theta B$ (for some $\theta > 0$) is the cylinder with the same center as $B$, and radius $\theta r(B)$.

    \item (Vertical planes) An affine plane $V \subseteq \R^{\adim+1}$ is \textit{$t$-independent} if $V$ contains a line in the $t$-direction. Such a plane is also called \textit{vertical}.

    \item (Parabolic Hausdorff measure and dimension) We define parabolic Hausdorff measure and dimension as usual, but using the parabolic metric. Indeed, for $\eta,\delta > 0$ and $E \subseteq \R^{\adim+1}$, define 
	\[ 
    \mathcal{H}^\eta_{p,\delta}(E) 
    := 
    \inf \left\{ \sum_j \pdiam(E_j)^\eta \colon E \subseteq \bigcup_j E_j \mbox{ and } \pdiam(E_j) \leq \delta\right\}.
    \]
	Then, define the \textit{parabolic Hausdorff measure of dimension} $\eta$ by 
	\[ \mathcal{H}^\eta_p(E) := \lim_{\delta \rightarrow 0} \mathcal{H}^\eta_{p,\delta}(E). \]
    We also define parabolic Hausdorff dimension by 
    \[
    \text{dim}_{\mathcal{H}_p}(E) 
    :=
    \inf \big\{ \eta > 0 : \mathcal{H}^\eta_p(E) = 0 \big\}
    =
    \sup \big\{ \eta > 0 : \mathcal{H}^\eta_p(E) = \infty \big\}.
    \]

    \item (Bilateral distance between sets) Given two closed sets $E, F \subset \R^{n+1}$, $(X, t) \in \R^{n+1}$ and $r > 0$, we denote 
    \begin{equation*}
        d_{(X, t), r}(E, F)
        :=
        \frac{1}{r} \max \left\{ \sup_{(Y, s) \in E \cap B_p(X, t)} \dist_p((Y, s), F), \sup_{(Y, s) \in F \cap B_p(X, t,r)} \dist_p((Y, s), E) \right\}.
    \end{equation*}
    Sometimes we also denote this quantity $d_{B_p(X, t ,r)}(E, F)$.
\end{itemize}

\subsection{Flatness coefficients and parabolic rectifiability} \label{sec:flatness}

As already said in the introduction, we quantify how flat a measure is by using the parabolic version of the beta numbers originally introduced by Jones in \cite{jones1990rectifiable}:

\begin{defn}[Beta numbers] \label{def:beta}
For a given $(X, t) \in \R^n \times \R$ and $r > 0$, we define
    \begin{equation*}
        \beta_\mu((X, t), r)
        :=
        \inf \left( \frac{1}{r^{n+1}} \int_{B_p(X, t,r)} \left( \frac{d_p((Y, s), P)}{r} \right)^2 \, d\mu(Y, s) \right)^{1/2},
    \end{equation*}
    where the infimum is taken for $P$ running over the family of vertical planes.
\end{defn}

For the implication \textit{rectifiability implies bounded Jones function}, we first show how to estimate the $\alpha$-numbers of Tolsa. Unlike $\beta$-numbers, the $\alpha$-number is a two-sided condition (asking that both the measure is close to the plane and the plane is close to the measure) which also detects variations in density. 

\begin{defn}\label{d:alpha-numbers}[Alpha numbers]
    Given two locally finite Borel measures $\mu,\nu$ and a
    parabolic ball $B$, we set
    \[
    \dist_B(\mu,\nu)
    :=\sup\left\{\left|\int f\,d\mu-\int f\,d\nu\right|:
    \Lip_{d_{p}}f\le1,\ \supp f\subset B\right\}.
    \]
    Then, for a point $(X,t) \in \supp(\mu)$ and a radius $r > 0$, we set 
    \[
    \alpha_\mu((X,t),r)
    :=\frac1{r^{s+1}}\inf_{a\ge0,V}
    \dist_{B_{p}((X,t),r)}
      \bigl(\mu,a\mathcal H_{p}^{s}\lfloor V\bigr),
    \]
    where the infimum is over $a\ge0$ and vertical hyperplanes $V$.
    We call a measure of the form $a\mathcal H_{p}^{s}\lfloor V$,
    with $a\ge0$ and $V$ a vertical hyperplane, a \textit{flat measure}.
    The condition $\supp f\subset B$ localizes the comparison to $B$, so
    the two measures may have different total masses.
\end{defn}

Similarly, we can measure how flat a function is at a concrete scale by finding its deviation from the best-approximating vertical plane. In this case, we take into account the parametrization of both the function and the plane.

\begin{defn}[Omega numbers]
    Given $f: \R^n \longrightarrow \R^{n+1}$ ($\R^n$ here means $n-1$ spatial coordinates and the time one), and $I \subset \R^n$ is a parabolic cylinder, we define its $\Omega$ coefficient by
    \begin{equation} \label{def:Omega_number}
        \Omega_f(I)
        :=
        \inf_{\substack{A \in \R^{n-1} \\ b \in \R}} 
        \left( \frac{1}{\abs{I}} \iint_{I} \left( \frac{\abs{f(X, t) - (A \cdot X + b, t)}}{\ell(I)} \right)^2 \, dX dt  \right)^{1/2}
    \end{equation}
    If $g : \R^n \to \R$ denotes just the height coordinate over the supporting plane, we use the previous definition with $f(X, t) := (X, g(X, t), t)$.
\end{defn}

\begin{defn}[Parabolic Lipschitz functions]\label{regfngraph.def}
    A function $g:\mathbb{R}^{n} \to \mathbb{R}$ on $n$-dimensional space-time is \textit{parabolic Lipschitz} (or Lip(1,1/2)) if there exists a constant $L$ such that
    \[|g(x,t) - g(y,s)| \le Ld_p((x,t),(y,s)), 
    \qquad 
    \forall \, x \in \R^{n-1}, t \in \R\]
\end{defn}

\begin{defn}[Gamma numbers] \label{def:gamma}
    For a Lip(1,1/2) function $g:\R^n \to \R$, we define the $\gamma$-numbers to be 
    \[\gamma_g((X, t), r) \coloneqq \inf_{a, \vec{b}}\left(\frac{1}{|B_p(X, t ,r)|}  \iint_{B_p(X, t,r)} \frac{|g(y,s) - (a + \vec{b}\cdot y)|^2}{r^2} \, dy \,ds\right)^\frac{1}{2}.\]
\end{defn}

\begin{defn}[Parabolic uniform rectifiability] \label{def:UR}
    Let $\Sigma \subset \R^{n+1}$ be closed, and define its surface measure by
    \begin{equation*}
        \nu:=\mathcal H_p^{n+1}\lfloor\Sigma.
    \end{equation*}
    We say that $\Sigma$ is \textit{parabolic Ahlfors–David regular} if
    \begin{equation*}
    C^{-1}r^{n+1}\le \nu(B_p(X,t,r))\le Cr^{n+1},
    \qquad \forall \, (X, t) \in\Sigma,\; 0<r<\diam_p(\Sigma).
    \end{equation*}
    If $\Sigma$ is parabolic Ahlfors-David regular, we say that $\Sigma$ is \textit{parabolic uniformly rectifiable} if there exists $M<\infty$ such that, if we define $\beta_\Sigma :=\beta_\nu$ as in Definition~\ref{def:beta}, it holds 
    \begin{equation} \label{eq:UR_via_betas}
        \sup_{\substack{(X, t) \in \Sigma \\ R > 0}} \frac{1}{R^{n+1}} \int_{\Sigma \cap B_R(X, t)} \int_0^R \beta_\Sigma((Y, s),r)^2 \,\frac{dr}{r} d\nu(Y, s) 
        \leq 
        M.
    \end{equation}    
\end{defn}

\begin{defn}[Regular parabolic Lipschitz functions] \label{def:regular}
    A scalar function $g:\R^{n-1} \times \R \to \R$ is \textbf{regular} parabolic Lipschitz (regular Lip(1,1/2)) if it is parabolic Lipschitz and its coefficients from Definition~\ref{def:gamma} satisfy
    \begin{equation*}
        \sup_{\substack{(X, t) \in \R^{n-1} \times \R \\ R > 0}} \frac{1}{R^{n+1}} \int_{B_p(X,t,R)} \int_0^R \gamma_g((Y,s),r)^2 \,\frac{dr}{r} \,dYds
        \leq
        M
    \end{equation*}
\end{defn}

\begin{rem} \label{rem:UR_with_betas}
    For a Lip(1,1/2) function $g: \R^{n-1} \times \R \to \R$, being regular Lip(1,1/2) is actually equivalent to the fact that its graph 
    \begin{equation*}
        \Gamma_g := \big\{(X,g(X,t),t):(X,t)\in\mathbb R^{n-1}\times\mathbb R\big\}
    \end{equation*}
    is parabolic uniformly rectifiable. Moreover, this equivalence is quantitative, that is, the constants in the definition of $g$ being regular Lip(1,1/2) is uniform on the Lip(1,1/2) character of $g$ and the uniform rectifiability constant of $\Gamma_g$. See \cite[Subsection 2.7]{bortz2025solvability} for further details. 
\end{rem}

\begin{rem}
	Regular parabolic Lipschitz can be alternatively defined by asking that $D_t^{1/2}f$ belongs to the parabolic BMO space\footnote{This means replacing standard Euclidean cubes by parabolic cylinders.}, where
	\[D_t^{1/2}f(x,t): = c\  \text{p.v.}\int_{\mathbb{R}} \frac{f(x,s) - f(x,t)}{|s-t|} \, ds\]
    is understood distributionally,	and $c$ is an absolute constant. For more details on the equivalence, see \cite[Remark 2.11]{BHS}. 
    In this paper, we prefer to focus on the more geometrical version involving beta/gamma numbers.
\end{rem}

%%%%%%%%%%%%%%%%%%%%%%%%%%%%%%%%%%%%%%%%%%%%%%%%%%%%%%%

\subsection{Dyadic cubes adapted to the support of an arbitrary measure} \label{sec:cubes}

We shall use a system of dyadic-like ``cubes'' on the set $\supp \mu$, even if we have very little a priori information about its structure. For that purpose, we use the construction of David and Mattila \cite[Theorem 3.2]{DM}, which moreover gives us useful information about which cubes in the construction are doubling, and estimates illustrating how the measure decays in the presence of non-doubling cubes. The construction of \cite{DM} is conducted in the Euclidean setting, but the proof works without changes in doubling spaces such as parabolic space-time (as already pointed out in other works like \cite{VZ}).\footnote{Note the difference: the construction requires that the ambient space (in our case, space-time $\R^n \times \R$ endowed with the parabolic metric) is doubling, not that the set under study (in our case, $E := \supp \mu$, where $\mu$ is a possibly non-doubling measure) is doubling.}

\begin{lem}[Dyadic cubes adapted to the support of $\mu$] \label{lem:cubes}
    Let $C_0>1$ and $A_0>5000\,C_0$. Then, for any Radon measure $\mu$ on space-time $\R^n \times \R$ (endowed with the parabolic metric), if define $E := \supp\mu$, there exists a family of ``cubes'' $Q$, indexed by generations $\cubes_k$, $k \in \Z$, such that the following hold:

    \begin{itemize}
        \item (Each generation partitions $E$) For each $k \in \N$, $E$ is the disjoint union of the ``cubes'' $Q$, $Q\in\cubes_k$.

        \item (Generations are nested) If $Q\in\cubes_\ell$ and $R\in\cubes_k$, for $k<\ell$, then either $Q \cap R=\varnothing$ or $Q \subset R$.

        \item (Cubes contain and are contained in cylinders of similar radius) For each $k \in \N$ and $Q\in\cubes_k$, there is a cylinder $B(Q)=B(z_Q,r(Q))$ such that
        \begin{equation*}
            z_Q\in E, \qquad A_0^{-k}\leq r(Q)\leq C_0\,A_0^{-k},
        \end{equation*}
        \begin{equation*}
            E\cap B(Q) \; \subset \; Q \; \subset \; E\cap 28\,B(Q),
        \end{equation*}
        and 
        \begin{equation*}
            \text{the cylinders $5B(Q)$, $Q\in\cubes_k$, are disjoint.}
        \end{equation*}

        \item (Doubling cubes) Denote by $\cubes_k^{db}$ the family of cubes $Q\in\cubes_k$ for which
        \begin{equation}\label{eqdob22}
            \mu(100B(Q))\leq C_0\,\mu(B(Q)).
        \end{equation}

        \item (Small measure of non-doubling cubes) When $Q\in\cubes_k\setminus \cubes_k^{db}$, it holds $r(Q)=A_0^{-k}$ and 
        \begin{equation}\label{eq:non_doub}
            \mu(100B(Q))\leq C_0^{-\ell}\,\mu(100^{\ell+1}B(Q))\quad
            \mbox{for all $\ell\geq1$ such that $100^\ell\leq C_0$.}
        \end{equation}
    \end{itemize}
\end{lem}

\begin{notation}[More notation associated to the dyadic grid] \label{not:david-mattila}
    We set
    \begin{equation*}
        \cubes := \bigcup_{k\geq0}\cubes_k, 
        \qquad 
        \cubes^{db} := \bigcup_{k\geq0}\cubes_k^{db}.
    \end{equation*}
    For $Q\in\cubes$, we denote the cubes below $Q$ by 
    \begin{equation*}
        \cubes(Q) := \{P\in\cubes:P\subset Q\},
        \qquad 
        \cubes^{db}(Q) := \cubes^{db}\cap \cubes(Q),
    \end{equation*} 
    and its \textit{generation} and \textit{side-length} by 
    \begin{equation*}
        J(Q) := k, 
        \qquad 
        \ell(Q) := 56\,C_0\,A_0^{-k}.
    \end{equation*}
    We note that 
    \begin{equation*}
        \frac1{28}\,C_0^{-1}\ell(Q)\leq \diam(Q)\leq\ell(Q), 
        \qquad 
        r(Q)\approx\diam(Q)\approx\ell(Q)
    \end{equation*}
    We also call $z_Q$ the \textit{center} of $Q$, and the cube $Q'\in \cubes_{k-1}$ such that $Q \subset Q'$ the \textit{parent} of $Q$.
\end{notation}

Along the manuscript, we will very often utilize the so-called stopping-time regions (sometimes also called \textit{trees}). With this, we refer to a family of cubes in $\cubes$ which are connected (via the natural relation of parent). Let us introduce some notation for these regions, that will be useful later.

\begin{notation}[Stopping-time regions]
    Let $\calS \subset \cubes$ be a family of cubes. We call $\calS$ a stopping-time region if the following hold:
    \begin{enumerate}
        \item There is a single maximal/top cube, namely there exists $Q(\calS)$ such that $Q \subset Q(\calS)$ for all $Q \in \calS$. 

        \item The family is semi-coherent, that is, if $Q_1, Q_2 \in \calS$ and $Q_1 \subset Q_2$, then all the intermediate cubes also lie in $\calS$, namely 
        \begin{equation*}
            Q \in \cubes, 
            \quad 
            Q_1 \subset Q \subset Q_2 
            \quad \implies \quad 
            Q \in \calS.
        \end{equation*}
    \end{enumerate}
    Given such a stopping-time region $\calS$, we denote its minimal cubes by 
    \begin{equation*}
        \min(\calS) 
        :=
        \big\{ Q \in \calS : \text{there is no } Q' \in \calS \text{ such that } Q' \subsetneq Q \big\},
    \end{equation*}
    and we also call its contact set the points which are covered by infinitely many cubes in $\calS$, of shrinking radii, namely
    \begin{equation*}
        z(\calS) := \big\{ (X, t) \in E : \text{ there exist $Q_i \in \calS$ with $\ell(Q_i) \to 0$ and $(X, t) \in Q_i$} \big\}.
    \end{equation*}
\end{notation}

Although the construction of the David-Mattila cubes works for any $C_0 > 1$ and $A_0 > 5000C_0$, in order to draw relevant conclusions from the property \eqref{eq:non_doub} by iterating it along chains of consecutive non-doubling cubes, it is more useful to choose 
\begin{equation} \label{eq:choiceC_0A_0}
    \text{$C_0$ is large enough, and $A_0 = C C_0^{100}$,}
\end{equation}
for an absolute constant $C$, as is done in \cite[(4.1)]{DM}. Therefore, we fix this choice in what follows.

With this choice of constants, the following consequences follow:

\begin{lem}[{Cubes can be a.e. covered by doubling cubes, \cite[Lemma 5.28]{DM}}]\label{lem:cover_doubling}
    Let $R\in\cubes$. Then there exists a family of doubling cubes $\{Q_i\}_{i\in I}\subset \cubes^{db}$, with $Q_i\subset R$ for all $i$, such that their union covers $\mu$-almost all $R$.
\end{lem}

\begin{lem}[{Decay of measure along non-doubling chains, \cite[Lemma 5.31]{DM}}]  \label{lemcad22}
    Let $R\in\cubes$ and let $Q\subset R$ be a cube such that all the intermediate cubes $S$,
    $Q\subsetneq S\subsetneq R$ are non-doubling (i.e.\ belong to $\cubes \setminus \cubes^{db}$).
    Then
    \begin{equation}\label{eqdk88}
    \mu(100B(Q))\leq A_0^{-10(n+1)(J(Q)-J(R)-1)}\mu(100B(R)).
    \end{equation}
\end{lem}

\begin{rem}
    As explained before \cite[Lemma 5.31]{DM}, the result is a consequence of the fact that one can apply \eqref{eq:non_doub} with $\ell$ going up to $\ell_0 = \lfloor \log_{100} C_0 - 1\rfloor$, which yields 
    \begin{equation*}
        \mu(100B(Q))
        \leq 
        C_0^{-\ell_0} \mu(C_0 B(Q))
        \leq 
        A_0^{-10} \mu(\widehat{Q}),
    \end{equation*}
    where $\widehat{Q}$ is the parent of $Q$. This is a consequence of the fact that $C_0^{-\ell_0} \approx C_0^{-\log C_0}$ decays faster than any power of $C_0$. Therefore, choosing $C_0$ large enough, $C_0^{-\ell_0}$ is smaller than any power of $A_0 \approx C_0^{100}$ (recall \eqref{eq:choiceC_0A_0}), and one can simply choose the power 10.
\end{rem}    

\begin{notation}[Density]
    Given a ball $B\subset \R^d$, we consider its density:
    $$\Theta_\mu(B) := \frac{\mu(B)}{r(B)^{n+1}}.$$
    We will also write $\Theta_\mu(x,r)$ instead of $\Theta_\mu(B(x,r))$.
\end{notation}

From the preceding lemma we deduce, as in \cite[Lemma 2.4]{azzam2015characterization}:

\begin{lem}[Decay of density along non-doubling chains] \label{lemcad23}
Let $Q,R\in\cubes$ be as in Lemma \ref{lemcad22}.
Then
$$\Theta_\mu(100B(Q))\leq C_0\,A_0^{-9(n+1)(J(Q)-J(R)-1)}\,\Theta_\mu(100B(R))$$
and
$$\sum_{\substack{S\in\cubes \\ Q\subset S\subset R}} \Theta_\mu(100B(S)) \lesssim_{C_0, A_0} \,\Theta_\mu(100B(R)).$$
\end{lem}

%%%%%%%%%%%%%%%%%%%%%%%%%%%%%%%%%%%%%%%%%%%%%%%%%%%%%%%%%
%%%%%%%%%%%%%%%%%%%%%%%%%%%%%%%%%%%%%%%%%%%%%%%%%%%%%%%%%

\subsection{Auxiliary results}

\begin{lem}[Dorronsoro estimates under second-derivative bounds] \label{l:Dorronsoro-for-C^2-maps}
	Let $g \colon \R^{n} \to \R$ ($\R^n$ is here seen as $n-1$ spatial variables plus the time variable) be Lipschitz with  $C^2$ with derivative estimate
			\begin{align}
			|\partial_t g| \lesssim r_{k}^{-1} \quad \mbox{ and } \quad  |D^2_X g| \lesssim r_{k}^{-1}. 
		\end{align}
	Then, 
	\begin{align}
		\sum_{I \subseteq J} \Omega_g(3I)^2|I| \lesssim |J|
	\end{align}
	for all parabolic cubes $J \subseteq \R^n$ with $\ell(J) \lesssim r_k$. In this sum, $I$ runs over the standard grid of dyadic (parabolic) cylinders in space time $\R^n = \R^{n-1} \times \R$.
\end{lem}
\begin{proof}
	Although the result is standard, let us include a proof for completeness. Fix some $J \subseteq \R^{n}$ with $\ell(J) \lesssim r_k^{-1}$. For each $I \subseteq J$, let $g_I$ be the Taylor polynomial defined by
	\begin{align}
		g_I(X,t) := g(X_I,t_I) + D_Xg(X_I,t_I) \cdot (X - X_I) \quad \mbox{for each } (X,t) \in \R^n = \R^{n-1} \times \R. 
	\end{align}
	Then, if $(X,t) \in \R^n$, we have, by the mean value theorem and out hypothesis on $g$,
	\begin{multline*}
		\abs{g(X, t) - g_I(X, t)} 
        \leq 
        \abs{g(X, t) - g(X, t_I)}
        + \abs{g(X, t_I) - g(X_I, t_I) - D_X g(X_I, t_I) \cdot (X - X_I)}
        \\ \lesssim 
        r_k^{-1} \abs{t-t_I} 
        + r_k^{-1} \abs{X - X_I}^2
        \lesssim 
        \ell(I)^2 r_k^{-1},
	\end{multline*}
    where the last estimate is due to the definition of parabolic cubes, using the parabolic metric. This implies that, if we choose $g_I$ in the definition of $\Omega_g(3I)$, we can estimate 
    \begin{multline*}
        \sum_{I \subseteq J} \Omega_g(3I)^2|I|
        \lesssim 
        \sum_{I \subseteq J} \iint_{3I} \left( \frac{\ell(I)^2 r_k^{-1}}{\ell(I)} \right)^2 \, dX dt
        \approx 
        \sum_{I \subseteq J} \ell(I)^2 r_k^{-2} |I|
        \\ =
        r_k^{-2} \sum_{m \geq 0} \sum_{\substack{I \subseteq J \\ \ell(I) = 2^{-m}\ell(J)}} \ell(I)^2 \, |I|
        =
        r_k^{-2} \ell(J)^2 \sum_{m \geq 0} 2^{-2m} \!\!\! \sum_{\substack{I \subseteq J \\ \ell(I) = 2^{-m}\ell(J)}}  \!\!\! |I|
        \\ \leq 
        r_k^{-2} \ell(J)^2 \sum_{m \geq 0} 2^{-2m} |J|
        \lesssim 
        r_k^{-2} \ell(J)^2 |J| 
        \lesssim 
        |J|.
    \end{multline*}
\end{proof}

\subsection{Balanced balls} 

We are going to define balanced balls following \cite[Section 3]{azzam2015characterization} as those who identify $n$ independent directions in which $\mu$ has mass at a particular scale, and hence it will be prone to resemble a plane. The correct adaptation of the notion to the parabolic setting takes into account that the desired geometric property is that $\mu$ is close to be supported over a vertical plane, so really we should only care about the $n-1$ spatial directions that are spanned (for that, we need $n$ spatial points). 

First, let us recall the following elementary fact about linearly independent vectors, if measured in a quantitative way. The proof follows just by using appropriately the triangle inequality.

\begin{lem}[Stability of the notion of quantitative linear independence] \label{lem:li}
    Let $\gamma > 0$, and assume that the points $X_0, \ldots, X_n \in \R^n$ satisfy, for certain $r > 0$,
    \begin{equation*}
        \dist_\eucl (X_i, \spann \{X_0, \ldots, X_{i-1}\}) 
        \geq 
        \gamma r, 
        \qquad 
        \forall \, 1 \leq i \leq n.
    \end{equation*}
    Then, there exists $t_0 = t_0(\gamma) > 0$ small enough such that if $Y_i \in B_{t_0r}(X_i)$ for certain points $Y_i \in \R^n$, then it holds 
    \begin{equation*}
        \dist_\eucl (Y_i, \spann \{Y_0, \ldots, Y_{i-1}\}) 
        \geq 
        \frac12 \gamma r, 
        \qquad 
        \forall \, 1 \leq i \leq n.
    \end{equation*}
\end{lem}

Next, we prove the analogous result of \cite[Lemma 3.1]{azzam2015characterization} in the parabolic setting, giving an alternative for the behavior of $\mu$ at each scale. The desirable behavior (in order for it to reflect flatness) would be that the mass of $\mu$ is distributed close to a vertical plane, that is, a geometric structure which has $n-1$ spatial dimension and is essentially constant in time (which adds also the spatial dimension, which has homogeneity 2 according to the parabolic homogeneous Hausdorff dimension). If this does not happen, the idea is to show that then $\mu$ has to concentrate close to a lower dimensional structure, and therefore there will be a steep increase of the density of $\mu$. More concretely, we prove the following result.

\begin{lem}[Alternative: closeness to a vertical plane, or increase in density] \label{lem:balanced}
    Let $0 < t, \gamma < 1$. Let also $\mu$ be a Radon measure in $\R^n \times \R$ with the parabolic metric, $B$ a (parabolic) ball of radius $r$. Then, there exists $\varepsilon = \varepsilon(t) > 0$ small enough such that one of the following two alternatives hold:
    \begin{enumerate}
        \item ($\mu$ has mass over $n-1$ spatial directions in $B$) There exist points $\{(X_i, t_i)\}_{i=0, \ldots, n-1}$, all contained in $B$, such that 
        \begin{equation*}
            \mu(B_p(X_i, t_i ,tr) \cap B) \geq \varepsilon \mu(B), 
            \qquad 
            \forall \, 0 \leq i \leq n-1,
        \end{equation*}
        and also 
        \begin{equation*}
            \dist_\eucl(X_i, \spann\{X_0, \ldots, X_{i-1}\}) \geq \gamma r, 
            \qquad 
            \forall \, 1 \leq i \leq n-1.
        \end{equation*}

        \item (Density increases significantly at a smaller scale) There exist a family $\{B_i\}_{i \in \mathcal{I}_B}$ of radii $4\gamma r$, centered inside $B$, such that $\{10B_i\}_{i \in \mathcal{I}_B}$ are pairwise disjoint, they satisfy 
        \begin{equation*}
            \sum_{i \in \mathcal{I}_B} \mu(B_i) 
            \gtrsim 
            \mu(B),
        \end{equation*}
        and also 
        \begin{equation*}
            \Theta_\mu(B_i) \gtrsim \gamma^{-1} \Theta_\mu(B), 
            \qquad 
            \forall \, i \in \mathcal{I}_B.
        \end{equation*}
    \end{enumerate}
\end{lem}
\begin{proof}
    Let us choose the points inductively in a greedy fashion, trying to verify the first alternative when possible. Indeed, choose $(X_0, t_0)$ such that 
    \begin{equation*}
        \mu(B_p(X_0, t_0,r) \cap B) \geq \frac12 \sup_{(X, t) \in B} \mu(B_p(X, t,r) \cap B).
    \end{equation*}
    Note that, since $B$ can be covered by at most $\approx t^{-(n+2)}$ cylinders of radius $tr$, it must be 
    \begin{equation*}
        \mu(B_p(X_0, t_0 ,r) \cap B) 
        \geq 
        c t^{n+2} \mu(B)
        \geq 
        \varepsilon \mu(B)
    \end{equation*}
    if we choose $\varepsilon = \varepsilon(t)$ small enough. This gives the first step of the construction. 

    Now let us define the subsequent points by induction. Assuming that we have constructed $(X_0, t_0), \ldots, (X_{i-1}, t_{i-1})$ satisfying the first alternative (concretely, $\mu(B_{p}(X_i, t_i ,r) \cap B) \geq \varepsilon \mu(B)$), we choose the next point $(X_i, t_i)$ satisfying 
    \begin{equation*}
        \mu(B_p(X_0, t_0 ,r ) \cap B) 
        \geq 
        \frac12 \sup_{\substack{(X, t) \in B \\ \dist_\eucl(X_i, \spann\{X_0, \ldots, X_{i-1}\}) \geq \gamma r}} \mu(B_p(X, t,r) \cap B),
    \end{equation*}
    which guarantees that the second condition of the first alternative is met. At this point, there are two options: either $\mu(B_p(X_0, t_0,r) \cap B) \geq \varepsilon \mu(B)$ and we can continue the construction, or we stop the construction because the condition does not hold. 
    In case this condition is always met and we can successfully complete the construction of the $n$ points, the first alternative of the lemma holds, and the proof finishes. Otherwise, let us call $i_0$ the first index for which this condition is not satisfied, and forces us to stop the construction. Clearly, $i_0 \in \{1, \ldots, n-1\}$.

    Let us now verify the second alternative, whenever the construction above stopped at the stage $i_0$. The idea is that the density should increase because the mass of $\mu$ will be accumulated very close to the vertical plane spanned by the first $i_0-1$ points (which is a lower dimensional set), namely the following slab:
    \begin{equation*}
        S
        :=
        \Big\{ (X, t) \in B :  \dist_\eucl(X, \spann\{X_0, \ldots, X_{i_0-1}\}) \leq \gamma r \Big\},
    \end{equation*}
    which is nothing more than an $\gamma r$-neighborhood (within $B$) of the spatial plane spanned by the first points, crossed with time (so that it becomes a vertical plane). 
    
    Note that, by definition of $i_0$, every ball of radius $tr$ contained in $B \setminus S$ satisfies 
    \begin{equation*}
        \mu(B_{p}(X, t ,tr) \cap B) \leq 2\varepsilon \mu(B).
    \end{equation*}
    Since $B \setminus S$ can trivially be covered by $\lesssim t^{-(n+2)}$ such balls, it holds 
    \begin{equation*}
        \mu(B \setminus S) 
        \leq 
        C t^{-(n+2)} \varepsilon \mu(B)
        \leq 
        \frac12 \mu(B)
    \end{equation*}
    if we choose $\varepsilon = \varepsilon(t)$ small enough. This implies that 
    \begin{equation} \label{eq:measure_slab}
        \mu(S \cap B) \geq \frac12 \mu(B).
    \end{equation}

    Now, by definition of $S$, it follows $i_0-1$ fixed spatial directions (those of $X_0, \ldots, X_{i_0-1}$) and the time direction (it is close to vertical), and we have restricted the rest of coordinates ($n-(i_0-1)$ spatial ones) to a $\gamma r$-neighborhood, so it is elementary that we can cover $S \subset \bigcup_{i \in \mathcal{I}_B^{\star \star}} B_i$ for a family of cylinders of radius $4\gamma r$, with 
    \begin{equation*}
        \abs{\mathcal{I}_B^{\star \star}}
        \leq 
        N
        :=
        C \frac{r^{i_0-1} r^2 (\gamma r)^{n-(i_0-1)}}{(4\gamma r)^{n+2}}
        =
        C \gamma^{-i_0-1}.
    \end{equation*}
    Let us retain only the balls carrying more mass, namely 
    \begin{equation*}
        \mathcal{I}_B^{\star} 
        :=
        \left\{ i \in \mathcal{I}_B^{\star \star} : \mu(B_i) \geq \frac{\mu(B)}{4N} \right\}, 
    \end{equation*}
    so that 
    \begin{equation*}
        \mu \left( \bigcup_{i \in \mathcal{I}_B^{\star \star} \setminus \mathcal{I}_B^{\star}} B_i \right)
        \leq 
        \sum_{i \in \mathcal{I}_B^{\star \star} \setminus \mathcal{I}_B^{\star}} \mu(B_i)
        \leq 
        N \frac{\mu(B)}{4N} 
        =
        \frac{\mu(B)}{4},
    \end{equation*}
    and since the balls with indices in $\mathcal{I}_B^{\star \star}$ covered $S \cap B$, \eqref{eq:measure_slab} implies 
    \begin{equation*}
        \mu \left( \bigcup_{i \in \mathcal{I}_B^{\star}} B_i \right)
        \geq 
        \frac{\mu(B)}{4}.
    \end{equation*}
    We continue by an elementary covering argument, which allows to restrict to a subset of the indices of $\mathcal{I}_B^{\star}$, which we denote by $\mathcal{I}_B$, so that the balls $\{10B_i\}_{i \in \mathcal{I}_B}$ are pairwise disjoint, and it still holds 
    \begin{equation*}
        \mu \left( \bigcup_{i \in \mathcal{I}_B} B_i \right)
        \gtrsim 
        \mu \left( \bigcup_{i \in \mathcal{I}_B^{\star}} B_i \right)
        \gtrsim 
        \mu(B).
    \end{equation*}
    Lastly, we confirm that density has increased for this family of balls. Indeed, for each $i \in \mathcal{I}_B$, since then $i \in \mathcal{I}_B^{\star}$, we can compute 
    \begin{equation*}
        \Theta_\mu(B_i)
        =
        \frac{\mu(B_i)}{r(B_i)^{n+1}}
        =
        \frac{\frac{\mu(B)}{4N}}{(4\gamma r)^{n+1}}
        =
        \frac{\mu(B)}{r^{n+1}} \frac{1}{4^{n+2} N \gamma^{n+1}}
        \geq 
        \Theta_\mu(B) \frac{1}{C\gamma^{-i_0-1}\gamma^{n+1}}
        \geq 
        \frac{1}{C\gamma} \Theta_\mu(B)
    \end{equation*}
    because $i_0 \leq n-1$ and $\gamma \in (0, 1)$. This finishes the proof of the second alternative.
\end{proof} 

With the above lemmas in mind, we define the concept of \textit{balanced balls} (following the notation in \cite{azzam2015characterization}), as those cylinders which contain enough $\mu$-mass close to a vertical plane at that scale (which is of course something desirable in our context, for it means that $\mu$ is somehow flat, in the sense of parabolic rectifiability, at the given scale). We do so by identifying $n-1$ linearly independent spatial directions (through giving $n$ points) that contain mass of $\mu$, in a quantitative way. We also make sure that the parameters are chosen correctly so that we can even perturb the points a little and they will still be linearly independent.

\begin{defn}[Balanced balls] 
    Given $\gamma > 0$, we say that the family of balls $\{B_r(X_i, t_i)\}_{i=0, \ldots, n}$, where $(X_i, t_i) \in B$ ($B$ is a fixed ball) is $\gamma$-balanced if 
    \begin{enumerate}
        \item (Each ball contains mass of $\mu$) For each $0 \leq i \leq n$ it holds 
        \begin{equation*}
            \mu(B_{p}(X_i, t_i ,r) \cap B) 
            \geq 
            \varepsilon 
        \end{equation*}
        where $t_0 = t_0(\gamma) > 0$ comes from Lemma~\ref{lem:li}, and then $\varepsilon = \varepsilon(t_0) > 0$ comes from Lemma~\ref{lem:balanced}.

        \item (The centers are quantitatively linearly independent) For each $1 \leq n$ it holds 
        \begin{equation*}
            \dist_\eucl (X_i, \spann \{X_0, \ldots, X_{i-1}\}) 
            \geq 
            \gamma r.
        \end{equation*}
    \end{enumerate}
\end{defn}

Since the concept of parabolic balanced balls is actually only really looking at the spatial directions (the time direction is always left free, for the model planes are vertical/$t$-independent), the following purely euclidean lemma of \cite{azzam2015characterization} will be useful in the sequel (concretely, in Lemma~\ref{l:angle-control}, where we compare the angles of the best approximating planes associated to contiguous scales).

\begin{lem}[{Angle control for spatial planes which are mutually close to a given well-separated set of points, \cite[Lemma 6.4]{azzam2015characterization}}] \label{lem:angle_planes}
    Suppose that $P_1$ and $P_2$ are hyperplanes in $\R^n$ (only spatial coordinates are considered here), and $\mathcal{X} := \{X_0, X_1, \ldots, X_n\}$ satisfy 
    \begin{enumerate}
        \item The points in $\mathcal{X}$ are quantitatively linearly independent, in the sense that 
        \begin{equation*}
            \eta(\mathcal{X}) 
            := 
            \min_{0 \leq i \leq n} \frac{\dist_{\eucl}(X_i, \spann(\mathcal{X} \setminus \{X_i\})}{\diam_{\eucl}(\mathcal{X})}
        \end{equation*}
        satisfies $\eta(\mathcal{X}) \in (0, 1)$, and 

        \item Both planes are quantitatively close to $\mathcal{X}$, in the sense that 
        \begin{equation*}
            \max_{\substack{0 \leq i \leq n \\ 1 \leq j \leq 2}} \frac{\dist_{\eucl}(X_i, P_j)}{\diam_{\eucl}(\mathcal{X})} 
            < 
            \varepsilon(\mathcal{X})
        \end{equation*}        
        for some $\varepsilon(\mathcal{X}) \in (0, \frac{\eta(\mathcal{X})}{2n})$.
    \end{enumerate}
    Then, it holds 
    \begin{equation*}
        \dist_{\eucl} (Y, P_1)
        \leq 
        \varepsilon \left( \frac{2n}{\eta(\mathcal{X})} \dist_{\eucl}(Y, \mathcal{X}) + \diam_\eucl(\mathcal{X}) \right), 
        \qquad 
        \forall \, Y \in P_2.
    \end{equation*}
\end{lem}

%%%%%%%%%%%%%%%%%%%%%%%%%%%%%%%%%%%%%%%%%%%%%%%%%%%%%%%%%%
%%%%%%%%%%%%%%%%%%%%%%%%%%%%%%%%%%%%%%%%%%%%%%%%%%%%%%%%%%

\newpage 
    
\section{Reifenberg parametrizations}\label{s:Reif}

	In this section, we discuss the parabolic version of David and Toro's version of the Reifenberg parametrization theorem \cite{david2012reifenberg} that was recently obtained in \cite{hoffman2026parabolic}. In \cite{hoffman2026parabolic}, the authors' main interest was in obtaining a map whose image was Reifenberg flat, i.e., the regularity of this map did not play such an important role. In this paper, we will need to obtain a parametrising map which is regular parabolic bi-Lipschitz. This is be possible under some additional assumption, as we will see.

	\subsection{Statement of the H\"older Reifenberg Theorem} \label{sec:lip}

	Following \cite{david2012reifenberg}, the parabolic Reifenberg theorem in \cite{hoffman2026parabolic} is stated in terms of \textit{coherent collection of balls and planes}, as defined right below. The idea is to have a mesh of balls of different sizes, organized by generations, such that every ball has a plane associated to it, and these planes oscillate slowly when we vary the scale and the location. The model case should be that of covering a smooth set of codimension one with many balls at different scales, using as planes the approximate tangents at the given ball: smoothness of the set will make these planes vary smoothly.
	
	\begin{defn}[$\varepsilon$-CCBP, Coherent Collection of Balls and Planes]  \label{d:CCBP}
		Let $\ve > 0$, $k_0 \in \Z$ and $P_0$ be a vertical plane. For $k \geq k_0$, let 
        \begin{equation*}
            r_k := 10^{-k},
        \end{equation*} 
        and let $\{ (X_{j,k},t_{j,k}) \}_{j \in J_k}$ be a collection of points in $\R^{\adim+1}$ such that 
		\begin{align}\label{e:net}
			d_p( (X_{i,k},t_{i,k}) , (X_{j,k},t_{j,k}) ) \geq r_k \quad \mbox{ for all } i,j \in J_k, \ i\not=j.
		\end{align}
		Additionally, let $B_{j,k} := B_{r_k}(X_{j,k},t_{j,k})$ and suppose $P_{j,k}$ are vertical planes such that $(X_{j,k},t_{j,k}) \in P_{j,k}.$ We say that 
        \begin{align}
            ( P_{k_0} , \{B_{j,k}\},\{P_{j,k}\})_{k \geq k_0, \ j \in J_k} \quad  \mbox{ is an \textit{$\ve$-coherent collection of balls and planes}},
        \end{align}
        $\ve$-CCBP for short, if it satisfies the following conditions:
		\begin{enumerate}
			\item \label{CCBP1} (Generations are somehow nested) If $k \geq k_0$ and $j \in J_k$, then
			\begin{align}\label{e:cinV^2}
				(X_{j,k},t_{j,k}) \in V^2_{k-1}, 
			\end{align}
			where we define the following enlarged regions
			\begin{align}\label{e:V_k}
				V_{k}^\lambda \coloneqq \bigcup_{j \in J_k} \lambda B_{j,k} 
                \qquad 
                \text{for $\lambda > 0$.}
			\end{align}
            
			\item \label{CCBP2} (Planes in the first generation are almost parallel to $P_0$) If $j \in J_{k_0}$,  then
			\begin{align}\label{e:cdist0}
				d_{100B_{j,k_0}}(P_{j,k_0},P_0) \leq \ve.
			\end{align}
            
			\item \label{CCBP3} (Small angle of close-by planes in the same generation) If $k \geq k_0$ and $i,j \in J_k$, then
			\begin{align}\label{e:c100rk}
                \norm{(X_{i,k},t_{i,k}) - (X_{j,k},t_{j,k})} \leq 100 r_k
                \quad \implies \quad 
				d_{100B_{i,k}}(P_{i,k},P_{j,k}) \leq \ve.
			\end{align}
            
			\item \label{CCBP4} (Small angle of close-by planes of consecutive generations) If $k \geq k_0 $, $i \in J_k$ and $j \in J_{k+1}$, then
			\begin{align}\label{e:c2r_k}
                \norm{(X_{i,k},t_{i,k}) - (X_{j,k+1},t_{j,k+1})} \leq 4r_k
                \quad \implies \quad 
				d_{20B_{i,k}}(P_{i,k},P_{j,k+1}) \leq \ve.
			\end{align}
		\end{enumerate}
	\end{defn}
	
	The Reifenberg map will be obtained as limit of a sequence of functions $\{f_k\}_{k \geq k_0}$, each of which is a projection onto the nearby balls at the generation $k$ of the CCBP. To capture this notion of which balls are actually nearby, we utilize a smooth partition of unity adapted to the CCBP.

    \begin{notation}[Partition of the unity associated to the CCBP]
        For each $k \geq k_0$, let $\{\theta_{j,k}\}$ be a collection of functions such that 
        \begin{equation} \label{e:supp_theta}
            \supp(\theta_{j,k}) \subseteq 10B_{j,k}
        \end{equation}
        such that $0 \leq \theta_{j,k} \leq 1,$  
		\begin{align}\label{e:=1}
			\sum_{j \in J_k} \theta_{j,k}(X,t) = 1 \quad \mbox{ for all } (X,t) \in V_k^9;
		\end{align}
		\begin{align}\label{e:lip-theta}
			|\theta_{j,k}(X,t) - \theta_{j,k}(Y,s)| \lesssim r_k^{-1} d_p((X,t),(Y,s)), 
            \qquad \text{for all } (X, t), (Y, s) \in \R^{n+1},
		\end{align}
		\begin{align}\label{e:deriv-theta}
			|\partial_t^m D^\alpha_X\theta_{j,k}(X,t)| \lesssim_m r_k^{-2m - |\alpha|} \quad \mbox{ for all } (X,t) \in \R^{n+1}, m \in \N, \alpha \in \N^n.
		\end{align}
    \end{notation}

    The associated projections are the following. 

    \begin{notation}[Projections associated to the CCBP]
        For $k \geq k_0$ and $j \in J_k$, let $\pi_{j,k} : \R^{n+1} \to \R^{n+1}$ denote the orthogonal projection onto $P_{j,k}$ and $\pi_{j,k}^\perp$ denote the orthogonal projection onto the linear subspace orthogonal to $P_{j,k}$, i.e., $P_{j,k}^\perp$. Their differentials $D\pi_{j,k}$ and $D\pi_{j,k}^\perp$ are, respectively, the orthogonal projection onto the linear subspace parallel to $P_{j,k}$, denoted by $P_{j,k}'$, and the orthogonal projection onto $P_{j,k}^\perp$, i.e., $D\pi_{j,k}^\perp = \pi_{j,k}^\perp$. 
    \end{notation}

    Having settled the notation for the partition of unity and the projections we are ready to define the Reifenberg map as an iteration of these projections.
    
	\begin{notation}[The iterative maps associated to the CCBP]
		For $k \geq k_0$ and $j \in J_k$, define
		\begin{align}\label{e:sigma_k} 
			\sigma_k(X,t) :=  (X,t) + \sum_{j \in J_k} \theta_{j,k}(X,t) \big[ \pi_{j,k}(X,t) - (X,t) \big]. 
		\end{align}
		Then, inductively, define a sequence of functions $\{f_k\}_{k \geq k_0}$ such that $f_k : P_0 \rightarrow \R^{\adim+1}$ by setting
		\begin{align}\label{e:deff}
			f_{k_0}(X,t) := (X,t) \quad \text{and} \quad f_{k+1}= \sigma_k \circ f_k. 
		\end{align}
        \end{notation}
	
	We can now state the main Reifenberg theorem in \cite{hoffman2026parabolic}. There are some extra conditions in the statement there that we will not need in this paper. For that reason, we have chosen to omit them. 
	
	\begin{thm}[bi-H\"older Reifenberg]  \label{t:DT}
		There is $\ve_0 >0$ such that the following holds. Let $0 < \ve < \ve_0$, fix $k_0 \in \mathbb Z$ and suppose $(P_0,\{B_{j,k}\},\{P_{j,k}\})_{k \geq k_0, \ j \in J_k}$ is an $\ve$-CCBP. Then, the sequence $\{f_k\}_{k \geq k_0}$ converges uniformly to a bi-H\"older map $f: P_{k_0} \rightarrow \R^{\adim+1}$ which satisfies the following properties. 
		\begin{enumerate}
			\item \label{DT1} We have 
			\begin{equation}
				E_\infty \coloneqq \bigcap_{K=k_0+1}^\infty \overline{\bigcup_{k=K}^\infty \{(X_{j,k},t_{j,k})\}_{j \in J_k}} 
                \;\; \subseteq \;\;
                f(P_0) =: \Sigma.
			\end{equation}
			\item \label{DT3} We have $\|f(X,t)-(X,t)\| \lesssim \ve$ for all $(X,t) \in P_0.$
			\item \label{DT4} For $k \geq k_0,$
			\begin{equation}\label{e:10}
				\sigma_k(X,t) = (X,t) \ \text{and} \ D\sigma_k(X,t)=I \ \quad \text{for} \ (X,t) \in \R^{\adim+1} \setminus V^{10}_k.
			\end{equation}
			\item \label{DT5} Let $\Sigma_{k_0} := P_0$ and 
			\begin{equation}
				\Sigma_k := \sigma_k(\Sigma_{k-1}), \,\,\, k >k_0.
			\end{equation}
			For $k \geq k_0$, there is a $\emph{Lip}(1,\tfrac{1}{2})$ function $A_{j,k}:P_{j,k} \cap 49B_{j,k} \rightarrow P_{j,k}^\perp$, which is $C^2$ such that 
            \begin{equation*}
                |A_{j,k}(X_{j,k},t_{j,k})| \ \lesssim \ve r_k 
                \quad \text{and} \quad 
                \| A_{j,k}\|_{\emph{Lip}(1,1/2)} \lesssim \ve,
            \end{equation*}
            and if $\Gamma_{j,k}$ is its graph over $P_{j,k}$, then 
			\begin{equation}
				\Sigma_k \cap D(X_{j,k},t_{j,k},P_{j,k},49r_k) = \Gamma_{j,k} \cap D(X_{j,k},t_{j,k},P_{j,k},49r_k)
			\end{equation}
			where
			\begin{equation}
				D(X,t,V,r) := \big\{z+w : z \in V \cap B_p(X,t,r), \ w \in V^\perp \cap B_p(0,0,r)\big\}.
			\end{equation}

			\item \label{DT8} For all $k \geq k_0$,
			\begin{align}\label{Sigma_ksigma_k}
				\|\sigma_k(X,t) - (X,t)\| \lesssim \ve r_k, 
                \qquad 
                \text{for all $(X,t) \in \Sigma_k$,}
			\end{align}
			and moreover,
			\begin{align}\label{Sigma_kSigma}
				{\dist}_p((X,t),\Sigma) \lesssim \ve r_k  \qquad \mbox{ for all } \ (X,t) \in \Sigma_k.
			\end{align}		
		\end{enumerate}
	\end{thm}

	\subsection{Bi-Lipschitz Reifenberg Theorem} \label{sec:bilip}
	
	In this section, we show how to obtain bi-Lipschitz regularity of the map provided via Theorem \ref{t:DT}, under the additional assumption that the angles between consecutive planes are square summable. To quantify this, it will be convenient to introduce the following quantity. For $k \geq k_0$ and $(X,t) \in \Sigma_k$, let 
	\begin{align}
		\ve_k(X,t) 
        \coloneqq 
        \sup 
        \left\{ d_{(X,t),10^4r_l}(P_{i,k} , P_{j,l}) \colon
        \begin{array}{c}
        i \in J_k, \ j \in J_l, \ |l - k| \leq   2,\\
        (X,t) \in 100B_{i,k} \cap 100B_{j,l}
        \end{array}
        \right\}.
	\end{align}
	The main result is then the following.

	\begin{thm}[Bi-Lipschitz Reifenberg] \label{t:bi-Lip}
		Let $\ve_0$ be as in Theorem \ref{t:DT} and assume that 
		\[(P_0,\{B_{j,k}\},\{P_{j,k}\})_{k \geq k_0, \ j \in J_k} \text{   is an $\ve$-CCBP} \]
		 for some $0 < \ve < \ve_0$. Suppose further that the maps $\{f_k\}_{k \geq k_0}$ defined in \eqref{e:deff} satisfy the condition 
		\begin{align}\label{e:ve-summable}
			\sum_{k \geq k_0} \ve_k(f_k(X,t))^2 \leq M \quad \mbox{ for all } (X,t) \in P_0. 
		\end{align}
		for some $M \geq 0$. Then, the map $f$ from Theorem \ref{t:DT} is bi-Lipschitz with constants depending only on $M$. 
	\end{thm}

    Our goal for the remainder of the section is to prove Theorem \ref{t:bi-Lip}. Thus, we shall fix for the remainder of the section a collection
	\[(P_0,\{B_{j,k}\},\{P_{j,k}\})_{k \geq k_0, \ j \in J_k}\]
which is an $\ve$-CCBP for some $0 < \ve < \ve_0$.

\subsection{Extra properties of $\sigma_k$}
	 Before dealing directly with the bi-Lipschitz estimates, we will need some preliminary observation regarding the construction of the maps $f,f_k$ and $\sigma_k$. These observations are not present (or needed) in \cite{hoffman2026parabolic}, but are essential for us here. At various points, we will be interested in how the mapping $\sigma_k$ (and those mappings which define $\sigma_k$) interact with the tangent vectors of $\Sigma_k$. Thus, we will let 
     \begin{align}
         T\Sigma_k(X,t) \mbox{ denote the tangent plane to } \Sigma_k \mbox{ at } (X,t). 
     \end{align}
     We note that this is well-defined since $\Sigma_k$ is a $C^2$-manifold by Theorem \ref{t:DT}.

\begin{lem}\label{l:tangent-small-angle-to-plane}
	For each $j \in J_k$, $(X,t) \in \Sigma_k \cap 10B_{j,k}$ and $v = (v_t,v_X) \in T\Sigma_k(X,t)$, we have
	\begin{align}\label{e:tangent-small-angle-to-plane}
		| D\pi_{j,k}(X,t) \cdot v - v | \lesssim  \ve | v_X |. 
	\end{align}
\end{lem}

\begin{proof}
	Fix $j \in J_k$, $(X,t)$ and $v$ as in the statement. Since $P_{j,k}$ is vertical, the map $\pi_{j,k}$ can be written as 
	\begin{align}
		\pi_{j,k}(X,t) = (\tilde{\pi}_{j,k}(X), t) \quad \mbox{ for all } (X,t) \in \R^{n+1}. 
	\end{align}
	Thus, $D\pi_{j,k}$ only modifies the spatial variables, and the left-hand side of \eqref{e:tangent-small-angle-to-plane} can be written as $| D_X \tilde{\pi}_{j,k} \cdot v_X - v_X |$. Recalling from Theorem \ref{t:DT} that $\Sigma_k$ coincide with a $\textrm{Lip}(1,1/2,C\ve)$ graph over $P_{j,k}$ in $10B_{j,k}$, and $D_X \tilde{\pi}_{j,k}$ is the orthogonal projection onto $\tilde{P}_{j,k}$, we have 
	\begin{align}
		| D_X \tilde{\pi}_{j,k} \cdot v_X - v_X | \lesssim \ve |v_X|, 
	\end{align}
	as required. 
\end{proof}

\begin{lem}\label{l:sigma_k-bi-lip}
	The map $\sigma_k \colon \Sigma_k \to \Sigma_{k+1}$ is $(1+C\ve)$-bi-Lipschitz. 
\end{lem}

\begin{proof}
	Fix $(X,t) , (Y,s) \in \Sigma_k$. If $(X,t) , (Y,s) \in (V_k^{10})^c$, then the map is the identity and there is nothing to show. Suppose instead that one of the points lies in $V_k^{10}$. Without loss of generality, we can assume that $(X,t) \in 10B_{j,k}$ for some $j \in J_k$. If $(Y,s) \in 45B_{j,k}$, then the bi-Lipschitz estimates follow from Lemma 3.19 in \cite{hoffman2026parabolic}. If $(Y,s) \not\in 45B_{j,k}$, then the bi-Lipschitz estimate follows from the fact that $d_p(\sigma_k(X,t),(X,t)) \lesssim \ve \, r_k \lesssim \ve \, d_p((X,t),(Y,s))$, see Theorem \ref{t:DT} \eqref{DT8}. 
\end{proof}

\

\begin{lem}\label{l:change-in-tangent-super}
    For each $(X,t) \in \Sigma_k$ and $v = (v_X,v_t) \in T\Sigma_k(X,t)$, we have 
    \begin{align}
        |D \sigma_k(X,t) \cdot v - v|  \lesssim \ve \min\left\{ \| v\| , |v_X| + r_k^{-1}|v_t| \right\}. 
    \end{align}
\end{lem}

\begin{proof}
    We start by proving the upper bound in terms of $\ve \|v\|$, which follows directly from Lemma \ref{l:sigma_k-bi-lip}. Indeed, since $\Sigma_k$ is $C^2$, we can find a curve $\gamma \colon [-1,1] \to \Sigma_k$ such that $\gamma'(0) = v$. Unraveling definitions and using the bi-Lipschitz estimate for the map $\sigma_k-\text{Id}$ from Lemma \ref{l:sigma_k-bi-lip}, we get 
	\begin{align}
		|D\sigma_k(X,t) \cdot v - v| &= \| D\sigma_k(X,t) \cdot v - v  \| 
        = \| (\sigma_k \circ \gamma)'(0) - \gamma'(0) \| \\
		&= \lim_{h \to 0} \frac{\norm{\sigma_k(\gamma(h)) - \sigma_k(\gamma(0)) - (\gamma(h) - \gamma(0))}}{\abs{h}} \\
		&\lesssim \ve  \lim_{h \to 0} \frac{ \norm{\gamma(h) - \gamma(0)}}{\abs{h}} 
        = \ve \| v\|,
	\end{align}
    as required. We now prove the upper bound in terms of $|v_X| + r_k^{-1}|v_t|$. Recalling that 
	\begin{align}
		\sigma_k(X,t) = (X,t) + \sum_{j \in J_k} \theta_{j,k}(X,t) \left( \pi_{j,k}(X,t) - (X,t) \right) ,
	\end{align}
	we have  
	\begin{align}
    \begin{split}\label{e:Dsigma.v-v}
        |D\sigma_k(X,t) \cdot v - v | &\leq \sum_{j \in J_k} \theta_{j,k}(X,t)| D\pi_{j,k}\cdot v  - v | \\
		&\hspace{2em}  \sum_{j \in J_k} |D\theta_{j,k}(X,t) \cdot v| |\pi_{j,k}(X,t) - (X,t) |. 
    \end{split}
	\end{align}
	By applying Lemma \ref{l:tangent-small-angle-to-plane} and recalling that $\sum_{j \in J_k} \theta_{j,k}(X,t) \lesssim 1$, the first term on the right-hand side of \eqref{e:Dsigma.v-v} is at most $C\ve |v_X|$. Recalling the derivative estimates on $\theta_{j,k}$ and the fact that $\{10B_{j,k}\}$ have bounded overlap, the second term is at most 
	\begin{align}
		C ( r_k^{-1}|v_X| + r_k^{-2} |v_t| )\, \ve \ r_k = C \ve (|v_X| + r_k^{-1}|v_t|) , 
	\end{align}
	as required.
\end{proof}

As well as the coefficient $\ve_k$ which measure the maximum angles between nearby planes in consecutive generations, it will be convenient to introduce the quantity $\ve_k'$ which measure the maximum angles between nearby planes in the \textit{same} generation. Indeed, define 
\begin{align}
    \ve_k'(X,t) \coloneqq \sup \left\{ d_{(X_{i,k},t_{i,k}),100r_k}(P_{i,k} , P_{j,k}) \colon i,j \in J_k, \ (X,t) \in 10B_{i,k} \cap 10B_{j,k} \right\} . 
\end{align}
This is also present in \cite{david2012reifenberg}, however, we caution that reader that in \cite{david2012reifenberg}, $\ve_k$ is the coefficient which controls angles in a single generation and $\ve_k'$ is the coefficient which controls angles across consecutive generation, i.e., the roles are reversed there. 

It is clear from the definitions that 
\begin{align}\label{e:ve-prime-smaller-than-ve}
    \ve_k'(X,t) \leq \ve_k(X,t) \quad \mbox{ for all } \quad (X,t) \in \R^{n+1}. 
\end{align}
Moreover, by the same proof as in \cite[Lemma 7.20]{david2012reifenberg}, one has 
\begin{align}\label{e:ve-prime-k-1-smaller-than-ve-k}
    \ve_{k-1}'(\sigma_{k-1}^{-1}(X,t)) \lesssim \ve_k(X,t) \quad \mbox{ for all } \quad (X,t) \in \Sigma_k \cap V_k^8.   
\end{align}

The following is consequence of \cite[Lemmas 3.22 \& 3.23]{hoffman2026parabolic}. We note that the upper bound there is in terms of $\ve$, rather than $\ve_k'(X,t)$. Examining the proof of these results, the estimates hold with any quantity which controls angle between planes $P_{i,k}$ and $P_{j,k}$ such that $(X,t) \in 10B_{i,k} \cap 10B_{j,k}$. In particular, it holds with the quantity $\ve_k'(X,t)$. 

\begin{lem}\label{l:parallel-to-plane}
    For each $k \geq 0$, $(X,t) \in V_k^8$ and $j \in J_k$ such that $(X,t) \in 10B_{j,k}$, we have 
    \begin{align}
        |\sigma_k(X,t) - \pi_{j,k}(X,t)| +  |D_X \pi_{j,k}^\perp \circ D_X \sigma_k(X,t) |  \lesssim \ve_k'(X,t).
        \end{align}
\end{lem}

The idea for the bi-Lipschitz estimates for $f$ is to track the change in length of curves as we push them forward by $\sigma_k$. The following lemma tells us that they increase by at most a factor $(1+\ve_k(X,t)^2)$ with an additional error term depending on the size on the time-derivative. This is why it is enough to assume square summability of the $\ve_k$ coefficients. 

\begin{lem}\label{l:ve-2-distortion}
	For $(X,t) \in V_k^8$ and $v = (v_X,v_t) \in T\Sigma_k(X,t)$, we have 
	\begin{align}
		\left| \| D\sigma_k(X,t) \cdot v \| - \|v\| \right| &\lesssim \ve_k(X,t)^2|v_X| + r_k^{-1} \ve_k(X,t) |v_t| \\
		&\leq \ve_k(X,t)^2\| v\| + r_k^{-1} \ve_k(X,t) |v_t|  . 
	\end{align}
\end{lem}

To prove Lemma \ref{l:ve-2-distortion}, we will need several preliminary lemmas which we now gather.

\begin{lem}\label{l:length-factor-projections}
	Let $(X,t) \in \R^{n+1}$ and let $i,j \in J_k$ such that $(X,t) \in 10B_{i,k} \cap 10B_{j,k} \neq \emptyset$. For any $v_X \in \R^n$, we have 
	\begin{align}
		\| D_X\pi_{i,k} \circ D_X\pi_{j,k} \circ D_X \pi_{i,k} \cdot v_X - D_X\pi_{i,k} \cdot v_X \| \lesssim \ve_k'(X,t)^2|v_X|.
	\end{align}
\end{lem}

\begin{proof}
	Fix some data as in the statement of the lemma. For any vertical plane $P$, we have that $D_X\pi_P \cdot v_X = (D_X \tilde{\pi} \cdot v_X, 0)$. In particular, it suffices to estimate 
    \begin{align}
        |  (D_X\tilde{\pi}_{i,k} \circ D_X\tilde{\pi}_{j,k} \circ D_X\tilde{\pi}_{i,k} - D_X\tilde{\pi}_{i,k}) \cdot v_X |.
    \end{align}
    David and Toro show in \cite[Equation (7.35)]{david2012reifenberg}, that any planes $V,W$ in $\R^n$ at an angle at most $\delta$ satisfy the estimate 
    \begin{align}\label{e:David-Toro-length-factor}
         | D_X V \circ D_X W \circ D_X V - D_X V | \lesssim \delta^2. 
    \end{align}
    By definition, the planes $V \coloneqq P'_{j,k}$ and $W \coloneqq P'_{i,k}$ make an angle at most $\ve_k'(X,t)$. Thus, the lemma follows directly from \eqref{e:David-Toro-length-factor} with this choice of planes.  
\end{proof}

\begin{lem}\label{l:length-factor-sigma}
	For $(X,t) \in \Sigma_k \cap V_k^8$ and $v = (v_X,v_t) \in \R^{n+1}$, we have 
	\begin{align}\label{e:length-factor-sigma}
			| D_X\tilde{\pi}_{i,k} \circ D_X\tilde{\sigma}_k(X,t) \circ D_X\tilde{\pi}_{i,k} \cdot v_X - D_X\tilde{\pi}_{i,k} \cdot v_X | \lesssim \ve_k(X,t)^2 |v_X| . 
	\end{align}
\end{lem}

\begin{proof}
	Fix $(X,t)$ and $v$ as in the statement. 
	We begin by writing the expression on the left-hand side of \eqref{e:length-factor-sigma} (inside the absolute values) as $u_1 + u_2$, where
	\begin{align}
		u_1 &\coloneqq \sum_{j \in J_k} \theta_{j,k}(X,t)\left[ D_X\tilde{\pi}_{i,k} \circ D_X\tilde{\pi}_{j,k} \circ D_X\tilde{\pi}_{i,k} - D_X\tilde{\pi}_{i,k} \right] \cdot v_X ; \\
		u_2 &\coloneqq \sum_{j\in J_k} \left( D_X\tilde{\pi}_{i,k} \cdot \tilde{\pi}_{j,k}(X) \right) \left[ D_X\theta_{j,k}(X,t) \circ D_X\tilde{\pi}_{i,k} \right] \cdot v_X. 
	\end{align}
	It follows directly from Lemma \ref{l:length-factor-projections} and the fact that the ball $\{10B_{j,k}\}$ have bounded overlap that 
    \[|u_1| \lesssim \ve_k'(X,t)^2 |v_X| \leq \ve_k(X,t)^2 |v_X|.\] 
    
    To estimate $|u_2|$, we first use the fact that $\sum_{j \in J_k} D\theta_{j,k}(X,t) = 0$ on $V_k^8$ to write 
	\begin{align}
		u_2 = \sum_{j \in J_k} \left( D_X\tilde{\pi}_{i,k} \cdot (\tilde{\pi}_{j,k}(X) - X) \right) \left[ D_X\theta_{j,k}(X,t) \circ D_X\tilde{\pi}_{i,k} \right] \cdot v_X .
	\end{align}
	Since $\tilde{\pi}_{j,k}(X) - X \in (P_{j,k}')^\perp$ and $P_{i,k}$ and $P_{j,k}$ make an angle at most $C\ve_k(X,t)$, we have 
	\begin{align}\label{e:D(pi-x)}
		|D\tilde{\pi}_{i,k} \cdot (\tilde{\pi}_{j,k}(X) - X)| \lesssim \ve'_k(X,t) |\tilde{\pi}_{j,k}(X) - X|. 
	\end{align}
    
    We claim that 
    \begin{align}\label{e:distance-to-plane}
        |\tilde{\pi}_{j,k}(X) - X| \lesssim \ve_k(X,t). 
    \end{align}
    Indeed, let $(Y,t) \in \Sigma_{k-1}$ such that $(X,t) = \sigma_{k-1}(Y,t)$ and let $l \in J_{k-1}$ such that $(X_{i,k},t_{i,k}) \in 2B_{l,k-1}$. By the triangle inequality and Theorem \ref{t:DT} \eqref{DT8}, it is easy to check that $(Y,t) \in 4B_{l.k-1}$. In particular, by applying Lemma \ref{l:parallel-to-plane} and \eqref{e:ve-prime-k-1-smaller-than-ve-k}, we have
    \begin{align}
        |(X,t) - \pi_{l,k-1}(Y,t)| =  | \sigma_{k-1}(Y,t) - \pi_{l,k-1}(Y,t) | \lesssim \ve_{k-1}'(Y,t) \lesssim \ve_k(X,t). 
    \end{align}
    Thus,  
    \begin{align}
        |(X,t) - \pi_{i,k}(X,t)| \leq |(X,t) - \pi_{l,k-1}(Y,t)| + C\ve_k(X,t) \lesssim \ve_k(X,t).  
    \end{align}
    The claim follows immediately.

	Now, using \eqref{e:D(pi-x)} and \eqref{e:distance-to-plane}, along with the derivative estimate for $\theta_{j,k}$, we have 
	\begin{align}
		|u_2| \lesssim \ve_k(X,t)^2 \, r_k \sum_{j \in J_k} r_k^{-1} | v_X|  \lesssim \ve_k(X,t)^2  |v_X| ,  
	\end{align}
	where in the final inequality we used the fact that $\{10B_{j,k}\}$ has bounded overlap. This estimate, combined with the corresponding estimate for $|u_1|$, finishes the proof. 
\end{proof}

\begin{lem}\label{l:small-normal}
    Let $(X,t) \in \Sigma_k \cap V_k^8$ and $v = (v_t,v_X) \in T\Sigma_k(X,t)$. For any $i \in J_k$ such that $(X,t) \in 10B_{i,k}$, we have 
    \begin{align}
        |D_X \tilde{\pi}_{i,k}^\perp \cdot v_X | = | v_X - D_X \tilde{\pi}_{i,k} \cdot v_X| \lesssim \ve_k(X,t). 
    \end{align}
\end{lem}

\begin{proof}
    Let $l \in J_{k-1}$ such that $(X_{i,k},t_{i,k}) \in 2B_{l,k-1}$, $(Y,t) \in \Sigma_{k-1}$ such that $(X,t) = \sigma_k(Y,t)$ and $u = (u_t,u_X) \in T\Sigma_{k-1}(X,t)$ such that $v_X = D_X\tilde{\sigma}_{k-1}(Y,t) \cdot u_X$. Using Theorem \ref{t:DT} \eqref{DT8} and the triangle inequality, we have that $(Y,t) \in 4B_{l,k-1}$. In particular, applying Lemma \ref{l:parallel-to-plane} and using that $\sigma_k$ is Lipschtz and $D_X \tilde{\pi}_{l,k-1}^\perp \circ D_X \tilde{\pi}_{l,k-1} = 0$, we get 
    \begin{align}
        | D_X \tilde{\pi}_{i,k}^\perp \cdot v_X | &\lesssim |D_X \tilde{\pi}_{l,k-1}^\perp \cdot v_X | = | D_X \tilde{\pi}_{l,k-1}^\perp \circ D_X \tilde{\sigma}_{k-1}(Y,t) \cdot u_X |  \\
        &= | D_X \tilde{\pi}_{l,k-1}^\perp \circ (D_X \tilde{\sigma}_{k-1}(Y,t) - D_X \tilde{\pi}_{l,k-1} )\cdot u_X | \\
        &\lesssim \ve_k'(Y,t)|u_X| \lesssim \ve_k(X,t)|v_X|,
    \end{align}
    as required. 
\end{proof}

\begin{proof}[\textbf{Proof of Lemma \ref{l:ve-2-distortion}}]
	Since the map $\sigma_k$ fixes the time component, we need only estimate its spatial component $\tilde{\sigma}_k$ (that is, $\sigma_k(X, t) = (\tilde{\sigma}_k(X, t), t)$). First, we can estimate 
	\begin{align}\label{e:|Dsigma.v|}
		\left| | D\tilde{\sigma}_k(X,t) \cdot v | - |v_X | \right| \leq |\partial_t \tilde{\sigma}_k(X,t) \, v_t | + \left| | D_X\tilde{\sigma}_k(X,t) \cdot v_X | - |v_X | \right|  .
	\end{align}
	By Lemma \ref{l:change-in-tangent-super}, we have 
	\begin{align}
		|\partial_t \tilde{\sigma}_k(X,t) | \lesssim r_k^{-1}\ve_k(X,t).
	\end{align}
	Thus, it only remains to estimate the final term in \eqref{e:|Dsigma.v|}. If $v_X = 0$, this term is trivially zero, so we may assume that $|v_X| \neq \emptyset$. 
    
    Defining $u_X \coloneqq v_X/|v_X|$, it suffices to show that 
	\begin{align}
		\left| |D_X \tilde{\sigma}_k(X,t) \cdot u_X | - 1 \right| \lesssim \ve_k(X,t)^2.  
	\end{align}
	Since $|u_X| = 1$, it further suffices to show that  
    \begin{align}\label{e:suffices-D-sigma} 
        | |D_X \tilde{\sigma}_k(X,t) \cdot u_X|^2 - 1 | \lesssim \ve_k(X,t)^2. 
    \end{align}
    
    To prove this final estimate, fix some index $i \in J_k$ such that $(X,t) \in 8B_{i,k}$ and begin by writing  
    \begin{align}
        D_X \tilde{\sigma}_k(X,t) \cdot u_X = w_1 + w_2 + w_3, 
    \end{align}
    where 
    \begin{align}
        w_1 &\coloneqq D_X\tilde{\pi}_{i,k}^{\perp} \circ D_X \tilde{\sigma}_k(X,t) \cdot u_X ; \\
        w_2 &\coloneqq D_X\tilde{\pi}_{i,k} \circ D_X \tilde{\sigma}_k(X,t) \circ D_X\tilde{\pi}_{i,k} \cdot u_X \\
        w_3 &\coloneqq D_X\tilde{\pi}_{i,k} \circ D_X \tilde{\sigma}_k(X,t) \circ D_X\tilde{\pi}_{i,k}^\perp \cdot u_X. 
    \end{align}
    It follows directly from Lemma \ref{l:parallel-to-plane} that 
    \[|w_1| \lesssim \ve_k(X,t)\] 
    and from Lemma \ref{l:length-factor-projections} that 
    \begin{align}
        |w_2 - D_X \tilde{\pi}_{i,k} \cdot u_X| = | (D_X\tilde{\pi}_{i,k} \circ D_X \tilde{\sigma}_k(X,t) \circ D_X\tilde{\pi}_{i,k}  -    D_X \tilde{\pi}_{i,k}) \cdot u_X | \lesssim \ve_k(X,t) . 
    \end{align}
    Furthermore, by using Lemma \ref{l:parallel-to-plane} and Lemma \ref{l:small-normal}, we have 
    \begin{align}
        |w_3| \leq |D_X\tilde{\pi}_{i,k} \circ ( D_X \tilde{\sigma}_k(X,t) - D_X\tilde{\pi}_{i,k}) \circ D_X\tilde{\pi}_{i,k}^\perp \cdot u_X| \lesssim \ve_k(X,t) |D_X \tilde{\pi}_{i,k}^\perp \cdot u_X| \lesssim \ve_k(X,t)^2. 
    \end{align}
    Putting these estimates together, using the fact that $w_1$ is orthogonal to $w_2 + w_3$ and applying the estimates $|a+b|^2 \lesssim |a|^2 + |b|^2$, we get 
    \begin{align}
         | |D_X \tilde{\sigma}_k(X,t) \cdot u_X|^2 - 1 | = | |w_1|^2 + |w_2 + w_3|^2 - 1 | &\leq | | D_X \tilde{\pi}_{i,k} \cdot u_X|^2 - 1 |  + C\ve_k(X,t)^2 \\
         &= |D_X\tilde{\pi}_{i,k}^\perp \cdot u_X|^2 + C\ve_k(X,t)^2. 
    \end{align}
    Recalling from above that $|D_X\tilde{\pi}_{i,k}^\perp \cdot u_X| \lesssim \ve_k(X,t)$, this finishes the proof of \eqref{e:suffices-D-sigma}. 
\end{proof}

\subsection{Bi-Lipschitz estimates for $f$} With the preliminary observation regarding $\sigma_k$ out of the way, we can start proving the bi-Lipschitz estimates from Theorem \ref{t:bi-Lip}. Assume for the remainder of the section that \eqref{e:ve-summable} holds for some constant $M \geq 0$. By translating and rotating, we can assume that \begin{align}
    P_0 = \R^{n-1} \times \{0\} \times \R. 
\end{align}
We start by showing that $f$ is Lipschitz.

\begin{lem}[$f$ is Lipschitz]  \label{l:f-upper-reg}
   For each $(X,t), (Y,s) \in P_0$, we have 
	\begin{align}
		\| f(X,t) - f(Y,s) \| \lesssim_M d_p((X,t),(Y,s)). 
	\end{align}
\end{lem}

\begin{proof}
	Let $(X,t),(Y,s) \in P_0$ and suppose to begin with that $d_p((X,t),(Y,s)) \geq 1$. In this case, recalling from Theorem \ref{t:DT} \eqref{DT3} that $d_p(f(Z,u) , (Z,u)) \lesssim \ve$ for all $(Z,u) \in P_0$, it is easy to check that $d_p((X,t),(Y,s)) \approx d_p(f(X,t) , f(Y,s))$. This is sufficient for the bi-Lipschitz estimate. \\

    Suppose now that $d_p((X,t),(Y,s)) < 1$. In the following, for each $k \geq 0$ we let $X_k ,Y_k \in \R^n$ denote the points such that 
    \begin{align}
        (X_k , t ) = f_k(X , t) \quad \mbox{ and } \quad (Y_k , s) = f_k(Y,s) . 
    \end{align}
    We start by showing the following. \\
	
    \noindent\textbf{Claim 1:} There is a largest integer, which we label $m$, such that $d_p((X_m,t),(Y_m,s)) \leq r_m$. Moreover, it satisfies the estimates
	\begin{align}\label{e:m-bounded}
		m \lesssim \log \left( \frac{1}{d_p((X,t),(Y,s))} \right)
	\end{align}
	and 
	\begin{align}\label{e:f_m-f}
		d_p(f_m(X,t) , f_m(Y,s)) \approx d_p(f(X,t) , f(Y,s)). 
	\end{align}
	
	Indeed, let $k \in \N$ be any index such that $d_p((X_k,t),(Y_k,s)) \leq r_k$. By applying Lemma \ref{l:sigma_k-bi-lip} ($k$ times), we obtain the estimate 
	\begin{align}\label{e:d-bounded}
		d_p((X,t),(Y,s)) \leq (1+C\ve)^k \, d_p((X_k,t),(Y_k,s)) \leq (1+C\ve)^k \, r_k . 
	\end{align}
	Since $(1+C\ve)^k \, r_k \to 0$ as $k \to \infty$ (recall that $r_k = 10^{-k}$), this implies the existence of the maximal value $m$. The estimate \eqref{e:m-bounded} follows by taking logarithms of \eqref{e:d-bounded} for the index $k = m$. 
    
    We now prove \eqref{e:f_m-f}. First, by the maximality of $m$, we have $d_p((X_{m+1} , t ) , (Y_{m+1} , t )) > r_{m+1}$. On the other hand, after iterating the estimate \eqref{Sigma_ksigma_k} from Theorem \ref{t:DT}, it is not difficult to show that $d_p((X_m,t) , f(X,t)) \lesssim \ve r_m$ and $d_p((Y_m,s) , f(Y,s)) \lesssim \ve r_m$. The inequality follows easily from these estimates and the triangle inequality, so long as $\ve$ is sufficiently small. \\

	Now having established the claim, let $\gamma \colon [0,1] \to P_0$ be a smooth constant speed parametrisation of the straight line from $(X,t)$ to $(Y,s)$. In particular, 
	\begin{align}
		\gamma'(u) = (Y-X, s-t) \quad \mbox{ for all } u \in [0,1]. 
	\end{align} 
	Furthermore, for $k \in \N$ and $u \in [0,1]$, let 
	\begin{align}
		v_k(u) \coloneqq Df_k(\gamma(u)) \cdot \gamma'(u). 
	\end{align}
	In this way, for each $k \geq 0$ and $u \in [0,1]$, one has 
	\begin{align}
		(v_k(u))_t = s-t, 
        \quad \text{and} \quad 
		v_{k+1}(u) = D\sigma_{k}(f_k(\gamma(u))) \cdot v_k(u).  
	\end{align}

    \noindent\textbf{Claim 2:} For each $u \in [0,1]$, we have  
    \begin{align}
		\| v_{m}(u) \| \lesssim \left[\prod_{k=0}^{m-1}(1+C\ve_k(X_k,t)^2)\right]\|\gamma'(u)\| + \sum_{k=0}^{m-1} r_k^{-1}\ve_k(X_k,t)|s-t|. 
    \end{align}

    Indeed, let $0 \leq l \leq m-1$ be the largest integer such that $f_l(\gamma(u)) \in V_l^8$. By Lemma \ref{l:ve-2-distortion}, for each $0 \leq k \leq l$, we have  
        \begin{align}
		\| v_{k+1}(u) \| \leq \big(1+C\ve_k(X_k,t)^2\big) \, \|v_k(u)\| + Cr_k^{-1} \, \ve_k(X_k,t) \, |s-t|. 
	\end{align}
    Hence, the iteration of this estimate up to the index $l$ gives
    \begin{align}
        \|v_{l+1}(u)\| \leq \left[\prod_{k=0}^{l}(1+C\ve_k(X_k,t)^2)\right]\|\gamma'(u)\| + C\sum_{k=0}^{l} r_k^{-1} \, \ve_k(X_k,t) \, |s-t|.
    \end{align}
    Since $f_{l+1}(\gamma(u)) \not \in V_{l+1}^{8}$ we cannot apply Lemma \ref{l:ve-2-distortion} any further. However, we can apply Lemma \ref{l:change-in-tangent-super} to get
    \begin{align}\label{e:l+2}
        \|v_{l+2}(u)\| \lesssim \left[\prod_{k=0}^{l}(1+C\ve_k(X_k,t)^2)\right]\|\gamma'(u)\| + \sum_{k=0}^{l} r_k^{-1} \, \ve_k(X_k,t) \, |s-t|.
    \end{align}
    By Definition \ref{d:CCBP} \eqref{CCBP1} and Theorem \ref{t:DT} \eqref{DT8}, we know that $f_k(\gamma(u)) \not\in V_k^{10}$ for all $k \geq l+2$. This, and Theorem \ref{t:DT} \eqref{DT4}, imply that $v_k(u) = v_{l+2}(u)$ for all $k \geq l+2$. In particular, $v_m(u) = v_{l+2}(u)$. This observation and \eqref{e:l+2} are enough to prove the claim. \\

	By the Fundamental Theorem of Calculus, Claim 2 and Cauchy-Schwarz, we observe that  
	\begin{align}
		|X_m - Y_m| &= \left| \int_0^1 D\tilde{f}_m(\gamma(u)) \cdot \gamma'(u) \, du  \right| \leq \int_0^1 \| v_m(u) \| \, du \\
		&\lesssim \left[\prod_{k=0}^{m-1} (1+C\ve_k(X_k,t)^2)\right] \int_0^1 \|\gamma'(u)\| \, du  + \sum_{k=0}^{m-1} r_k^{-1} \ve_k(X_k,t) |t-s| \\
		&\leq C(M) \, d_p((X,t),(Y,s)) + \left( \sum_{k=0}^{m-1} \ve_k(X_k,t)^2 \right)^{1/2} \left(  \sum_{k=0}^{m-1} r_k^{-2}\right)^{1/2} |t-s| \\
		&\lesssim_M d_p((X,t),(Y,s)) + r_m^{-1} |t-s| 
	\end{align} 
    (in the third and fourth lines, we used the assumption \eqref{e:ve-summable}). Recalling from Claim 1 that $d_p((X,t),(Y,s)) \leq r_m$, so that $|t-s|^{1/2} \leq r_m$, it follows direct from the above inequality that
	\begin{align} \label{eq:f_m_Lip}
		\hspace{1em} d_p(f_m(X,t) , f_m(Y,s)) = |X_m - Y_m| \lesssim_M d_p((X,t),(Y,s)) + |t-s|^\frac{1}{2} \lesssim d_p((X,t),(Y,s)). 
	\end{align}
    The proof is complete after recalling \eqref{e:f_m-f}.  
\end{proof}

To obtain the full lower Lipschitz estimates, we will first consider the lower Lipschitz estimates for $f$ on each time-slice. The proof will be similar to the proof of Lemma \ref{l:f-upper-reg}.

\begin{lem}[$f$ is bi-Lipschitz in every time slice] \label{lem:bilip_slice}
	For each $x, y \in \R^{n-1}$ and $t \in \R$, we have 
	\begin{align}
		\| f(x,0,t) - f(y,0,t) \| \gtrsim_M |x-y|. 
	\end{align}
\end{lem}

\begin{proof}
	Let $m \in \N$, $\gamma \colon [0,1] \to \R^n$ and $\{v_k(u)\}$ be as in the proof of Lemma \ref{l:f-upper-reg}, with $s = t$. Since $\gamma'(u)$ is purely spacial for each $u$, the same is true of the vectors $v_k(u)$. In particular, by the Fundamental Theorem of Calculus, we have 
	\begin{align}
			|X_m - Y_m| &= \left| \int_0^1 D\tilde{f}_m(\gamma(u)) \cdot \gamma'(u) \, du  \right| \leq \int_0^1 \| v_m(u) \| \, du.
	\end{align}
	Analogously to the proof of Claim 2 in Lemma \ref{l:f-upper-reg}, by applying Lemma \ref{l:ve-2-distortion} (for purely spacial vector) and Lemma \ref{l:change-in-tangent-super} (once), we obtain
	\begin{align}
		|X_m-Y_m| \gtrsim \left( \prod_{k=0}^{m-1} (1+C\ve_k(X_k,t)^2) \right)^{-1} \int_0^1 \|\gamma'(u)\| \, du \gtrsim_M |X-Y|. 
	\end{align} 
	The lemma follows directly from this and the estimate \eqref{e:f_m-f}.  
\end{proof}

We now show how to combine the global upper Lipschitz estimate with the lower Lipschitz estimate on each time slice to obtain the desired global bi-Lipschitz estimates. The following lemma, combined with Lemma \ref{l:f-upper-reg}, finishes the proof of Theorem \ref{t:bi-Lip}. 

\begin{lem}[$f$ is bi-Lipschitz]
	For each $(X,t) , (Y,s) \in P_0$, we have 
	\begin{align}
		\| f(X,t) - f(Y,s) \| \gtrsim_M d_p((X,t),(Y,s)). 
	\end{align}
\end{lem}

\begin{proof}
	Let $\alpha \in (0,1)$ (to be chosen sufficiently small depending on $M$) and fix some $(X,t), (Y,s)  \in P_0$. The proof split into two cases. First, if $|t-s|^{1/2} \geq \alpha |X-Y|$, then, since $f$ fixes the time component, 
	\begin{align}
		\| f(X,t) - f(Y,s) \| \geq |t-s|^{1/2} \gtrsim_\alpha d_p((X,t),(Y,s)). 
	\end{align}
	Since we will choose $\alpha$ depending only on $M$, this concludes the proof in this case. 
    
    Suppose instead that $|t-s|^{1/2} < \alpha |X-Y|$. In this case, by applying the lower Lipschitz estimate on each time-slice (Lemma~\ref{lem:bilip_slice}) and the general upper Lipschitz estimates (Lemma~\ref{l:f-upper-reg}), we get 
	\begin{align}
		\|f(X,t) - f(Y,s) \| &\geq \| f(X,t) - f(Y,t) \| - \| f(Y,t) - f(Y,s) \| \\
		&\geq c_1(M)|X-Y| - C_2(M)|t-s|^{1/2} \\
		&\geq c_1(M)|X-Y| - C_2(M) \,\alpha \, |X-Y| \\
		&\gtrsim_M |X-Y| \\
        &\gtrsim d_p((X,t),(Y,s)), 
	\end{align} 
	where in the second to last inequality we chose $\alpha = \alpha(M)$ sufficiently small. 
\end{proof}

\subsection{Additional regularity} In Section \ref{sec:regular} below, we will apply Theorem \ref{t:bi-Lip} in a specific context that will also allow us to upgrade the map $f$ to being regular parabolic bi-Lipschitz map. In order to achieve this, we will need the following regularity statement for the surfaces $\Sigma_k$, which is essentially a strengthening of Theorem \ref{t:DT} \eqref{DT5}. With Lemma \ref{l:Dorronsoro-for-C^2-maps} in mind, this will be extremely useful for upgrading the regularity of $f$.  

\begin{lem}[Second-derivative bounds]  \label{l:regular-graphs}
    Let $\ve_0$ be as in Theorem \ref{t:DT} and assume that 
	\[(P_0,\{B_{j,k}\},\{P_{j,k}\})_{k \geq k_0, \ j \in J_k} \quad \text{is an $\ve$-CCBP}\]
	for some $0 < \ve < \ve_0$. Then, for each $k \geq k_0$ and $j \in J_k$, the map $A_{j,k}$ from Theorem \ref{t:DT} \eqref{DT5} additionally satisfies the derivative estimates 
    \begin{align}
        |\partial_t A_{j,k}(X,t)| + |D_X^2 A_{j,k}(X,t)| \lesssim r_k^{-1}
    \end{align} 
    for all $(X,t) \in P_{j,k} \cap 49B_{j,k}$.
\end{lem}

Before discussing the proof of Lemma \ref{l:regular-graphs}, we show the $\sigma_k$ satisfies corresponding derivative estimates. These will be essential in the argument below. 

\begin{lem}[Derivative estimates for $\sigma_k$, restricted to spatial variables] \label{l:second-or-sigma-estimates}
    There exists $\Csigma \geq 1$ such that the following holds. For each $k \geq 0$, $(X,t) \in \Sigma_k$ and purely spatial tangent vectors $u,v \in T\Sigma_k(X,t)$ with $|u| , |v| \leq 1$, we have 
    \begin{align}
        | D\sigma_k(X,t) \cdot u - u| \leq \Csigma \, \ve, 
        \qquad \text{and} \qquad 
        |D^2 \sigma_k(X,t)[u,v]| \leq \Csigma \, \frac{\ve}{r_k}.
    \end{align}
\end{lem}

\begin{proof}
    Consider $(X, t) \in \Sigma_k$. The first estimate follows directly from Lemma \ref{l:change-in-tangent-super}. For the second estimate, we differentiate $\sigma_k$ twice to obtain 
    \begin{align*}
        D^2 \sigma_k(X, t)[u, v]
        & \leq 
        \sum_{j \in J_k} D^2 \theta_{j, k}(X, t) [u, v] \left(\pi_{j, k}(X, t) - (X, t) \right)
        \\ & \quad + \sum_{j \in J_k} (\nabla \theta_{j, k}(X, t) \cdot u) \left(D \pi_{j, k}(X, t)[v] - v\right)
        \\ & \quad + \sum_{j \in J_k} (\nabla \theta_{j, k}(X, t) \cdot v) \left(D \pi_{j, k}(X, t)[u] - u\right)
        \\ & \quad + \sum_{j \in J_k} \theta_{j, k}(X, t) \, D^2 \pi_{j, k} (X, t)[u, v]
        \\ & =:
        I + II + III + IV.
    \end{align*}
    Using \eqref{e:deriv-theta} we obtain 
    \begin{equation*}
        \abs{I} 
        \leq 
        \sum_{\substack{j \in J_k \\ (X, t) \in \supp \theta_{j, k}}} \frac{C}{r_k^2} \abs{u}\abs{v} \dist_p((X, t), P_{j, k})
        \leq 
        \frac{C\varepsilon}{r_k} \abs{u}\abs{v},
    \end{equation*}
    because we can estimate the distance from $(X, t)$ to $P_{j, k}$ using \eqref{Sigma_ksigma_k} because $(X, t) \in \Sigma_k$, and the contributing $j$ are bounded because of \eqref{e:supp_theta}. In turn, \eqref{e:deriv-theta} and Lemma~\ref{l:tangent-small-angle-to-plane} imply 
    \begin{equation*}
        \abs{II} 
        \leq 
        \frac{C}{r_k} \abs{u} \varepsilon \abs{v}
        \leq 
        \frac{C\varepsilon}{r_k} \abs{u}\abs{v},
    \end{equation*}
    and the same computation holds for $III$. Lastly, $IV = 0$ because $D^2 \pi_{j, k} = 0$. All these imply the required bound on the second derivatives. 
\end{proof}

We will prove Lemma \ref{l:regular-graphs} by induction on $k$. For this, we will need to recall some elements of the construction in \cite{hoffman2026parabolic}. First, the maps corresponding to the generation $k = 0$ are defined by 
\begin{align}\label{e:A-k-0}
    A_{j,0} := 0 \quad \mbox{ for all } j \in J_0. 
\end{align}
We then construct maps 
\begin{align}
    F_{j,k} \colon P_{j,k} \cap 40B_{j,k} \to P_{j,k}^\perp \quad  \mbox{ for each } k \geq 0 \mbox{ and } j \in J_k,
\end{align}
and 
\begin{align}
    A_{j,k} \colon P_{j,k} \cap 49B_{j,k} \to P_{j,k}^\perp \quad  \mbox{ for each } k \geq 1 \mbox{ and } j \in J_k,
\end{align}
inductively according to the following rules. \\

\noindent\textbf{Defining $F$.} Let $k \geq 0$ and suppose that $\{A_{j,k}\}_{j \in j_k}$ have been defined. For each $j \in J_k$ and $(X,t) \in P_{j,k} \cap 40B_{j,k}$, we define 
\begin{align}\label{e:def-F}
	 F_{j,k}(X,t) \coloneqq \pi^\perp_{j,k} \left( \sigma_k \left( h_{j,k}^{-1}(X,t) + A_{j,k}(h_{j,k}^{-1}(X,t))\right)\right),
\end{align}
where $\pi_{j,k}$ denotes the orthogonal projection onto $P_{j,k}$ and $h_{j,k}$ is a map from $P_{j,k} \cap 41B_{j,k}$ into $P_{j,k}$ whose image contains $P_{j,k} \cap 40B_{j,k}$ and is defined by 
\begin{align}\label{d:h}
	 h_{j,k}(X,t) := \pi_{j,k}\left[\sigma_k\left((X,t)+A_{j,k}(X,t)\right)\right]. 
\end{align}

It was shown in \cite{hoffman2026parabolic} that the map $h_{j,k}$ satisfies the derivative estimates 
\begin{align}\label{e:first-order-h}
    |Dh(X,t)\cdot v - v| \lesssim \ve 
\end{align}
for all $(X,t) \in 41B_{j,k}$ and all purely spacial vectors $v \in P_{j,k}'$ with $|v| \leq 1$, see \cite[(3.34)]{hoffman2026parabolic}. \\  

\noindent\textbf{Defining $A$.} Let $k \geq 0$ and suppose that $\{F_{i,k}\}_{i\in J_k}$ have already been defined. For each $j \in J_{k+1}$, choose an arbitrary $i \in J_k$ such that 
\begin{align}\label{e:inclusion}
    (X_{j,k+1},t_{j,k+1}) \in 2B_{i,k}. 
\end{align}
Then, for $(X,t) \in 49B_{j,k+1}$, let 
\begin{align}\label{e:def-A}
    A_{j,k+1}(X,t) 
    := 
    \pi_{j,k+1}^\perp \left[ \pi_{i,k}(g_{j,k+1}^{-1}(X,t)) + F\left( \pi_{i,k}(g_{j,k+1}^{-1}(X,t))\right) \right] - \pi_{j,k+1}^\perp(X_0,0).
\end{align}
where $(X_0,0)$ is an arbitrary point in $P_{j,k+1}$, and $g_{j,k+1}$ is a map from $P_{j,k+1} \cap 39B_{i,k}$ into $P_{j,k+1}$ whose image contains $P_{j,k+1} \cap 49B_{j,k+1}$, and is defined by 
\begin{align}\label{e:definitionf}
	 g_{j,k+1}(X,t) := \pi_{j,k+1} \left[ \pi_{i,k}(X,t) + F_{i,k}(\pi_{i,k}(X,t))\right]. 
\end{align}
It was shown in \cite[Lemma 3.13]{hoffman2026parabolic} that the map $g_{j,k+1}$ satisfies the estimates 
\begin{align}\label{e:first-order-g}
    | Dg(X,t) - I | \lesssim \ve \quad \mbox{ and } \quad  |g(X,t) - (X,t)| \lesssim \ve r_k 
\end{align}
for all $(X,t) \in P_{j,k+1} \cap 39B_{j,k+1}$ and all purely spatial vectors $v \in P_{j,k+1}'$. \\ 

Lemma \ref{l:regular-graphs} follows immediately from the following lemma which also provides, along with the estimate for the $A$ maps, the corresponding information for the $F$ maps.

\begin{lem}\label{l:regular-graphs-A-F}
    Let $\ve_0$ be as in Theorem \ref{t:DT} and assume that 
	\[(P_0,\{B_{j,k}\},\{P_{j,k}\})_{k \geq k_0, \ j \in J_k} \; \text{is an $\ve$-CCBP }\]
	for some 
    \begin{align}\label{e:epsilon}
    0 < \ve \ll \min\left\{\ve_0 , \frac{1}{100(1+\Csigma)}\right\}.
    \end{align}
    Letting
    \begin{align}
        \CA \coloneqq 4+8\Csigma, \quad \CF \coloneqq 12 + 30\Csigma \quad \mbox{ and } \quad  \CFprime \coloneqq 1+ 2\Csigma,
    \end{align}
    we have, for each $k \geq k_0$ and $j \in J_k$, that   
    \begin{align}\label{e:derivaties-A}
        |\partial_t A_{j,k}(X,t)| + |D_X^2 A_{j,k}(X,t)| \leq \CA r_k^{-1}, 
        \qquad 
        \text{for all $(X,t) \in P_{j,k} \cap 49B_{j,k}$},
    \end{align}
    \begin{align}\label{e:derivaties-F}
        |\partial_t F_{j,k}| + |D_X^2 F_{j,k}| \leq  \CF r_k^{-1}, 
        \qquad 
        \text{for all $(X,t) \in P_{j,k} \cap 40B_{j,k}$,}
    \end{align}
    and
    \begin{align}\label{e:derivatives-F-2}
         |\partial_t F_{j,k}| + |D_X^2 F_{j,k}| \leq  \CFprime r_k^{-1}, 
         \qquad 
         \text{for all $(X,t) \in P_{j,k} \cap 8B_{j,k}$.}
    \end{align}    
\end{lem}

\begin{rem}
    It is worth mentioning the role of \eqref{e:derivatives-F-2}. The proof of Lemma \ref{l:regular-graphs-A-F} is carried out by induction. We get derivative estimates for $F$ at generation $k$ depending on $A$ at generation $k$ and then we get derivative estimates for $A$ at generation $k+1$ depending on $F$ at generation $k$. The point is that,  given \eqref{e:derivaties-A} for some index $k$, we are able to obtain the estimates \eqref{e:derivatives-F-2} in this restricted domain for some constant $C_3$ completely independent of $C_1$. If we didn't have this, the constants in the construction would blow up.  
\end{rem}

It suffices to prove Lemma \ref{l:regular-graphs-A-F} in the case that $k_0 = 0$. We fix 
\[(P_0,\{B_{j,k}\},\{P_{j,k}\})_{k \geq 0, \ j \in J_k}\] 
as in the statement of Lemma \ref{l:regular-graphs} for the remainder of the sub-section. We will prove Lemma \ref{l:regular-graphs-A-F} by induction. Note that \eqref{e:derivaties-A} holds trivially for the base case $k = 0$ by \eqref{e:A-k-0}. Thus, it suffices to show that \eqref{e:derivaties-A} for index $k$ implies \eqref{e:derivaties-F} for the index $k$ and \eqref{e:derivaties-F} for the index $k$ implies \eqref{e:derivaties-A} for the index $k+1$. We split this into two steps.  \\

\noindent\textbf{Step 1.} The main goal of this step is to show the following. 

\begin{lem}\label{l:3.15-to-3.16}
    Let $k \geq 0$ and suppose that \eqref{e:derivaties-A} holds for $k$, then \eqref{e:derivaties-F} and \eqref{e:derivatives-F-2} hold for $k$. 
\end{lem}

Fix some $j \in J_k$. For the remainder of the proof of Lemma \ref{l:3.15-to-3.16} we will  drop the dependency on $j,k$ from our notation, i.e., we will write 
\begin{align}
    A \coloneqq A_{j,k} , \quad  F \coloneqq F_{j,k} , \quad P = P_{j,k} , \quad  \dots  
\end{align}

Our first goal will be to obtain second order derivative estimates on the map $h$ and $h^{-1}$ from \eqref{d:h}. As other times before, we denote $\tilde{h}$ the spatial part of $h$.

\begin{lem}\label{l:derivatives-h-edit}
    We have 
    \begin{align}\label{e:derivatives-h-edit}
         \max\left\{|\partial_t \tilde{h}(X,t)| ,  | D^2_Xh(X,t) | \right\} \leq 2(\Csigma + \CA) \, \ve \, r_k^{-1}
    \end{align}
    for each $(X,t) \in P$. 
    In particular, 
    \begin{align}\label{e:derivatives-h-1-edit}
         \max\left\{ |\partial_t [\tilde{h}^{-1}](X,t)| , | D^2_X[h^{-1}](X,t) | \right\} \leq 3(\Csigma + \CA) \, \ve \, r_k^{-1}
    \end{align}
    for each $(X,t) \in P$. 
\end{lem}

\begin{proof}
    The estimates in \eqref{e:derivatives-h-1-edit} follow directly from the estimates in \eqref{e:derivatives-h-edit}, \eqref{e:first-order-h} and the local inversion theorem. Thus, we will only be concerned with proving \eqref{e:derivatives-h-edit}. Fix $(X,t) \in P$ for the remainder of the proof and let us write ${A_*}$ for the graph map $(Y,s) \mapsto (Y,s) + A(Y,s)$. In particular, 
    \begin{align}\label{e:tilde-A}
        \partial_t \tilde{A}_* = \partial_t \tilde{A}, \quad D {A_*} = I + DA \quad \mbox{ and } \quad D^2{A_*} = D^2 A. 
    \end{align}
    
    We start with the time derivative estimate. Recalling the definition of $h$ from \eqref{d:h}, we can compute 
    \begin{align}
        \partial_t \tilde{h}(X,t) = S_1 + S_2 + S_3, 
    \end{align}
    where
    \begin{align}
        S_1 &\coloneqq \partial_t \tilde{\pi}\left(\sigma_k(A_*(X,t))\right) ; \\
        S_2 &\coloneqq D_X \tilde{\pi}\left(\sigma_k(A_*(X,t))\right) \partial_t \tilde{\sigma}_k(A_*(X,t)) ; \\
        S_3 &\coloneqq D_X \tilde{\pi}\left(\sigma_k(A_*(X,t))\right) D_X \tilde{\sigma}_k(A_*(X,t)) \, \partial_t \tilde{A}_*(X,t). 
    \end{align}
    Since the plane $P$ is vertical (so that $\tilde{\pi}$ is only a function of space), we have $S_1 =0$. By using the fact that $\tilde{\pi}$ is 1-Lipschitz and the estimates for $\sigma_k$ from Lemma \ref{l:second-or-sigma-estimates}, we can estimates $|S_2| \leq \Csigma \ve r_k^{-1}$. Using these same properties, the first equality in \eqref{e:tilde-A} and the induction hypothesis on $A$, we have $|S_3| \leq (1+C\ve)\CA \ve r_k^{-1} \leq 2\CA \ve r_k^{-1}$. Thus, we have the first estimate in \eqref{e:derivatives-h-edit}. 

    Now we deal with the spatial derivatives. Fix some purely spatial vectors $u,v \in P'$ (the linear space parallel to $P$) with $|u| , |v| \leq 1$. It suffices to show that  
    \begin{align}\label{e:derivatives-h-1}
        | D^2[h^{-1}](X,t)(u,v) | \leq 2(\Csigma+\CA) \ve r_k^{-1}. 
    \end{align}
    To prove this, 
    by applying the chain rule for second order differentials, we compute 
    \begin{align}
        D^2h(X,t)(u,v) = D\pi \cdot D^2[\sigma_k \circ {A_*}](X,t)(u,v) = D\pi \cdot ( S_1 + S_2), 
    \end{align}
    where 
    \begin{align}
        S_4 &\coloneqq D^2\sigma_k(X,t)\left(D{A_*}(X,t) \cdot u , D{A_*}(X,t) \cdot v\right) ; \\  
        S_5 &\coloneqq D\sigma_k({A_*}(X,t)) \cdot D^2{A_*}(X,t)(u,v). 
    \end{align}
    Since the map ${A_*}$ fixes the time component, the vectors $D{A_*}(X,t) \cdot u$  and $D{A_*}(X,t) \cdot v$ are purely spatial elements of $T\Sigma_{k}$. Moreover, by since $A$ is $C\ve$-Lipschitz, we have $|D{A_*}(X,t) \cdot u| \leq (1+C\ve)|u|$ and $|D{A_*}(X,t) \cdot v| \leq (1+C\ve)|v|$. Using these facts, along with Lemma \ref{l:second-or-sigma-estimates}, we get 
    \begin{align}
        |S_4| \leq (1+C\ve) \Csigma \ve r_k^{-1} \leq   2 \Csigma \ve r_k^{-1}
    \end{align}
    On the other hand, using Lemma \ref{l:second-or-sigma-estimates} and the induction hypothesis on $A$, we have 
    \begin{align}
        |S_5| \leq (1+C\ve) \CA \ve r_k^{-1} \leq 2  \CA \ve r_k^{-1}. 
    \end{align}
    Since $\pi$ is 1-Lipschitz, the above two estimate imply $|D^2h(X,t)(u,v)| \leq 2(\Csigma + \CA)\ve r_k^{-1}$, as required.
\end{proof}

We now have all the ingredients needed to estimates the derivatives of $F$. The following finishes the proof of Lemma \ref{l:3.15-to-3.16}.

\begin{lem}
    For each point $(X,t) \in P \cap 40B$, we have 
    \begin{align}\label{e:derivative-F-3-edit}
        \max\left\{ | \partial_t F(X,t)| \ ,  \ |D_XF^2(X,t)| \right\} \leq \CF \ve r_k^{-1}.
    \end{align}
    For each point $(X,t) \in P \cap 7B$, so long as $\ve$ is sufficiently small, we have 
    \begin{align}\label{e:derivative-F-4-edit}
        \max\left\{ | \partial_t F(X,t)| \ ,  \ |D_XF^2(X,t)| \right\}  \leq  \CFprime \ve r_k^{-1}.
    \end{align}
\end{lem}

\begin{proof}
    Fix $(X,t) \in P \cap 40B$ for the remainder of the proof. We start by estimating the time derivative of $F$. Recalling the definition of $F$ from \eqref{e:def-F}, we can compute 
    \begin{align}
        \partial_t F(X,t) = S_1 + S_2 + S_3, 
    \end{align}
    where 
    \begin{align}
        S_1 &\coloneqq D_X \pi^\perp \partial_t  \tilde{\sigma}_k \left(A_*(h^{-1}(X,t))\right) \\
        S_2 &\coloneqq D_X \pi^\perp D_X \tilde{\sigma}_k\left(A_*(h^{-1}(X,t))\right) \partial_t \tilde{A}_*(h^{-1}(X,t)) \\
        S_3 &\coloneqq D_X \pi^\perp D_X \tilde{\sigma}_k\left(A_*(h^{-1}(X,t))\right) D_X \tilde{A}_*(h^{-1}(X,t)) \, \partial_t [h^{-1}](X,t).
    \end{align}
    Using that $\pi^\perp$ is 1-Lipschitz, the derivative estimates for $\sigma_k$ from Lemma \ref{l:second-or-sigma-estimates} and for $h^{-1}$ from Lemma \ref{l:derivatives-h-edit}, and the induction hypotheses for $A$, we can estimate 
    \begin{align}
        |\partial_t F(X,t)| \leq \left[ \Csigma + (1+C\ve) \CA + 3(1+C\ve)(\CA + \Csigma)   \right] \ve r_k^{-1}. 
    \end{align}
    Recalling the definition of $\CF$, the above inequalities are enough to implies the time-derivative estimate in \eqref{e:derivative-F-3-edit}. 
    
    To obtain the improved time-derivative estimate in \eqref{e:derivative-F-4-edit}, in the case that $(X,t) \in P \cap 7B$, we first observe that this assumption implies $(Z,t) \coloneqq A_*(h^{-1}(X,t)) \in V_k^8$. In particular, by applying Lemma \ref{l:parallel-to-plane}, we have in this case that 
    \begin{align}\label{e:S_2+S+3}
        |S_2 + S_3| &\leq C(1+C\ve)\CA\ve^2 r_k^{-1} +  3C(1+C\ve)(\CA + \Csigma) \ve^2 r_k^{-1} \\
        &\leq C(1+C\ve)(4\CA + 3\Csigma) \ve^2 r_k^{-1} = C\ve(1+C\ve)(16+35\Csigma) \ve^2 r_k^{-1}.
    \end{align}
    Our choice of $\ve$ from \eqref{e:epsilon} implies that 
    \begin{align}\label{e:choice-of-epsilon}
        C\ve(1+C\ve)(16+35\Csigma) \leq 1. 
    \end{align} 
    which, together with \eqref{e:S_2+S+3}, implies the desired estimate in \eqref{e:derivative-F-3}.

    We now estimate the spatial derivatives. First, By applying the chain rule for second order derivatives, we obtain 
    \begin{align}\label{e:F}
        D^2F(X,t)(u,v) = D\pi^\perp \cdot \left[ S_4 + D\sigma_k(\tilde{A}(h^{-1}(X,t))) \cdot \left[ S_5 + S_6 \right] \right] ,
    \end{align}
    where 
    \begin{align}
        S_4 &\coloneqq D^2\sigma_k\left(\tilde{A}(h^{-1}(X,t))\right) \left( D[\tilde{A} \circ h^{-1}](X,t) \cdot u , D[\tilde{A} \circ h^{-1}](X,t) \cdot v \right); \\
        S_5 &\coloneqq D^2\tilde{A}(h^{-1}(X,t)) \left(D[h^{-1}](X,t) \cdot u , D[h^{-1}](X,t) \cdot v\right); \\
        S_6 &\coloneqq D\tilde{A}(h^{-1}(X,t)) \cdot D^2[h^{-1}](X,t)(u,v). 
    \end{align}
    Consider $S_4$ to begin with. Since the map $\tilde{A}$ and $h^{-1}$ both fix the time component, both vectors on which $D^2\sigma_k(\tilde{A}(h^{-1}(X,t)))$ is applied are purely spatial element of $T\Sigma_k(\tilde{A}(h^{-1}(X,t)))$. Thus, by applying Lemma \ref{l:second-or-sigma-estimates}, we have 
    \begin{align}\label{e:S_4}
        |S_4| \leq (1+C\ve)\Cseparation \ve r_k^{-1}. 
    \end{align}
    For $S_5$, we apply the induction hypothesis regarding the derivatives for $A$ and use the fact that $h^{-1}$ is $(1+C\ve)$-Lipschitz to get that 
    \begin{align}\label{e:S_5}
        |S_5|  \leq (1+C\ve)\CA \ve r_k^{-1}. 
    \end{align}
    For $S_6$, we use that fact that $\tilde{A}$ is $(1+C\ve)$-Lipschitz and Lemma \ref{l:derivatives-h-edit} to obtain 
    \begin{align}\label{e:S_6}
        |S_6| \leq 3(1+C\ve)(\Csigma + \CA) \ve r_k^{-1}. 
    \end{align}
    Recalling \eqref{e:F}, the estimate \eqref{e:derivative-F-3} follows from the estimates \eqref{e:S_4}, \eqref{e:S_5} and \eqref{e:S_6}, along with the fact that $\pi^\perp$ is 1-Lipschitz and $\sigma_k$ is $(1+C\ve)$-Lipschitz. 
    
    To obtain the improved estimate \eqref{e:derivative-F-4-edit} in the case that $(X,t) \in P \cap 7B$, we recall from above that $(Z,t) \coloneqq \tilde{A}(h^{-1}(X,t)) \in V_k^8$. In particular, by applying Lemma \ref{l:parallel-to-plane}, we have
    \begin{align}\label{e:improved-S_5-S_6}
        \left|D\pi^\perp \cdot D\sigma_k(Z,s) \cdot [ S_5 + S_6 ]\right| \leq C\ve |S_5 + S_6| &\leq C\ve (1+C\ve)(4\CA + 3\Csigma) \ve r_k^{-1} \\
        &= C\ve(1+C\ve)(16+35\Csigma) \ve r_k^{-1} \\
        &\leq \ve r_k^{-1},
    \end{align}
    where in the last inequality we used \eqref{e:choice-of-epsilon}. Thus, we have by \eqref{e:F}, \eqref{e:S_4} and \eqref{e:improved-S_5-S_6}, the desired estimate \eqref{e:derivative-F-4-edit}.
\end{proof}

\noindent\textbf{Step 2.} The main goal of this step is to show the following iteration.

\begin{lem}\label{l:3.16-to-3.15}
    Let $k \geq 0$ and suppose that \eqref{e:derivaties-F} and \eqref{e:derivatives-F-2} holds for $k$, then \eqref{e:derivaties-A} holds for $k+1$. 
\end{lem}

Fix some $j \in J_{k+1}$ and let $i \in J_k$ be the corresponding index from the construction of $A_{j,k+1}$ such that \eqref{e:inclusion} holds. For the remainder of the proof of Lemma \ref{l:3.16-to-3.15}, we will drop the dependency on $i,j$ and $k$ from the maps $g_{j,k+1}$ and $F_{i,k}$ defining $A_{j,k+1}$, i.e., we will write 
\begin{align}
    A \coloneqq A_{j,k+1}, \quad F \coloneqq F_{i,k} \quad \mbox{ and } \quad g \coloneqq g_{j,k+1}, 
\end{align}
but we will keep the dependency on the remaining objects. Analogous to the proof of Lemma \ref{l:3.15-to-3.16}, our first goal will be to obtain derivative estimates the map $g$. 

\begin{lem}\label{l:derivatives-g}
    For each point $(X,t) \in P_{j,k+1} \cap 7B_{i,k}$, We have 
    \begin{align}\label{e:derivatives-g-edit}
         \max\left\{|\partial_t \tilde{g}(X,t)| ,  | D^2_Xg(X,t) | \right\} \leq \CF \ve r_k^{-1}
    \end{align}
    for each $(X,t) \in P$. 
    In particular, 
    \begin{align}\label{e:derivatives-g-1-edit}
         \max\left\{ |\partial_t [\tilde{g}^{-1}](X,t)| , | D^2_X[g^{-1}](X,t) | \right\} \leq 2\CF \ve r_k^{-1}
    \end{align}
    for each $(X,t) \in P$. 
\end{lem}

\begin{proof}
If $(X,t) \in P_{j,k+1} \cap 49B_{j,k+1}$, then \eqref{e:inclusion} and \eqref{e:first-order-g} imply that $g^{-1}(X,t) \in 7B_{i,k}$. In particular, the estimates \eqref{e:derivatives-g-1-edit} on $P_{j,k+1} \cap 49B_{j,k+1}$ follows directly from the estimates \eqref{e:derivatives-g-edit}on $P_{j,k+1} \cap 7B_{i,k}$, \eqref{e:first-order-g} and the local inversion theorem. Thus, we will only be concerned with proving \eqref{e:derivatives-g-edit}. 

Fix $(X,t)$ as in the statement of the lemma. Let us write $F_*$ for the graph map $(Y,s) \mapsto (Y,s) + F(Y,s)$. In particular, 
    \begin{align}\label{e:tilde-F}
        |\partial_t F_* | = | \partial_t F| \quad \mbox{ and }  \quad D^2F_* = D^2 F. 
    \end{align} 
    Starting with the time-derivative, we compute 
\begin{align}
    \partial_t \tilde{g}(X,t) = D_X \tilde{\pi}_{j,k+1} \partial_t \tilde{F}_*(\pi_{i,k}(X,t)) .
\end{align}
Here, we have used the fact that $\partial_t \tilde{\pi}_{j,k+1} = \partial_t \tilde{\pi}_{i,k}=  0$. Since $(X,t) \in 7B_{i,k}$ and $P_{i,k}$ contains $(X_{i,k} , t_{i,k})$, we have $\pi_{i,k}(X,t) \in 7B_{i,k}$. Now, using that orthogonal projections are 1-Lipschitz, the identity $|\partial_t \tilde{F}_*| = | \partial_t F|$ and the induction hypothesis on the derivatives of $F$ inside $P_{i,k} \cap 7B_{i,k}$, it follows that 
\begin{align}
    |\partial_t \tilde{g}(X,t)| \leq \CFprime \ve r_k^{-1}. 
\end{align}

To prove the estimate on the spatial derivatives, it suffices to show 
\begin{align}
    |D^2 g(X,t)(u,v)| \leq \CFprime \ve r_k^{-1} 
\end{align}
for any points of purely spatial vectors $u,v$ with $|u| , |v| \leq 1$, which we now fix. By applying the chain rule for second order differentials, we compute 
    \begin{align}
        D^2g(X,t)(u,v) &= D\pi_{j,k+1} \cdot D^2 \tilde{F}(\pi_{i,k}(X,t))(D\pi_{i,k} \cdot u , D\pi_{i,k} \cdot v) \\
        &= D\pi_{j,k+1} \cdot D^2 F(\pi_{i,k}(X,t))(D\pi_{i,k} \cdot u , D\pi_{i,k} \cdot v)
    \end{align}
    By a similar reasoning as in the time-derivative estimate, it follows from the induction hypothesis on $F$ that 
    \begin{align}
        |D^2g(X,t)(u,v)| \leq \CFprime \ve r_{k}^{-1}, 
    \end{align}
    as required. 
\end{proof}

We now have all the ingredients needed to estimates the derivatives of $A$. The following finishes the proof of Lemma \ref{l:3.16-to-3.15}. 

\begin{lem}
    For each point $(X,t) \in P_{j,k+1} \cap 49B_{j,k+1}$, we have 
    \begin{align}\label{e:derivative-F-3}
        \max\left\{ |\partial_t A(X,t)| \ , \ |D_XA^2(X,t)| \right\} \leq \CA \ve r_k^{-1}.
    \end{align}
\end{lem}

\begin{proof}
    Fix $(X,t)$ as in the statement of the lemma. Notice to begin with that, since $(X,t) \in 49B_{j,k+1}$, it follows from \eqref{e:inclusion} and \eqref{e:first-order-g} that 
    \begin{align}
        d_p(g^{-1}(X,t),(X_{i,k} , t_{i,k})) \leq 49r_{k+1} + 2r_k + \ve r_k \leq 7r_k,
    \end{align}
    that is, $g^{-1}(X,t) \in 7B_{i,k}$. Since $P_{i,k}$ contains $(X_{i,k} , t_{i,k})$, it also follows that 
    \begin{align}\label{e:in-7-ball}
    \pi_{i,k}(g^{-1}(X,t)) \in 7B_{i,k}.
    \end{align} 
    
    We start by estimating the time-derivatives. Recalling the definition of $A$ in terms of $F$. By applying the chain rule and using the fact that $\partial_t \pi$ for any orthogonal projection $\pi$ onto a vertical plane, we have  
    \begin{align}
        \partial_t A(X,t) &= D_X \pi_{j,k+1}^\perp \partial_t \tilde{F}_*(\pi_{i,k}(g^{-1}(X,t))) \\ 
        &\hspace{2em} + D_X \pi_{j,k+1}^\perp D_X \tilde{F}_*(\pi_{i,k}(g^{-1}(X,t))) D_X \tilde{\pi}_{i,k} \partial_t[g^{-1}](X,t) .
    \end{align}
    By \eqref{e:in-7-ball} and the induction hypotheses on $F$, the first term above is at most $\CFprime \ve r_k^{-1}$ in absolute value. Since $\CA = 4\CFprime$, the time-derivative estimates now follows.  Using these same properties and the derivative estimates for $g^{-1}$ from Lemma \ref{l:derivatives-g}, the second term is at most $2(1+C\ve)\CFprime \ve r_k^{-1}$ in absolute value. 
    
    To estimate the spatial derivatives, it suffices to show that 
    \begin{align}
        |D^2A(X,t)(u,v)| \leq \CA \ve r_k^{-1} 
    \end{align}
    for any choice of purely spatial vectors $u,v \in P_{j,k+1}'$ with $|u|,|v| \leq 1$, which we now fix. We first apply the chain rule for second order derivatives to obtain  
    \begin{align}
        D^2A(X,t)(u,v) = D\pi_{j,k+1}^\perp \cdot \left[ S_1 + D\tilde{F}(\pi_{i,k}(g^{-1}(X,t)) \cdot S_2 \right], 
    \end{align}
    where 
    \begin{align}
        S_1 &\coloneqq D^2\tilde{F}(\pi_{i,k}(g^{-1}(X,t)))( D[\pi_{i,k} \circ g^{-1}](X,t) \cdot u , D[ \pi_{i,k} \circ g^{-1}](X,t) \cdot v) ; \\
        S_2 &\coloneqq D\pi_{i,k} \cdot D^2[g^{-1}](X,t)(u,v). 
    \end{align}
    To estimate the first term, note that it follows from \eqref{e:in-7-ball}, the induction hypothesis on $F$ and the fact that $g^{-1}$ is $(1+C\ve)$-Lipschitz, that 
    \begin{align}
        |D\pi_{j,k+1}^\perp \cdot S_1| \leq |S_1| \leq (1+C\ve) \CFprime \ve r_k^{-1}. 
    \end{align}
    On the other hand, it follows from Lemma \ref{l:derivatives-g} and the fact that $\tilde{F}$ is $(1+C\ve)$-Lipschitz that  
    \begin{align}
        |D\pi_{j,k}^\perp \cdot D\tilde{F}(\pi_{i,k}(g^{-1}(X,t)) \cdot S_2| \leq (1+C\ve)|S_2| \leq 2(1+C\ve)\CFprime \ve r_k^{-1}. 
    \end{align}
    Thus, putting these last estimate together, we have shown
    \begin{align}
        |D^2A(X,t)(u,v)| \leq 4\CFprime \ve r_k^{-1}  = \CA \ve r_k^{-1}, 
    \end{align}
    as required. 
\end{proof}

%%%%%%%%%%%%%%%%%%%%%%%%%%%%%%%%%%%%%%%%%%%%%%%%%%%%%%%%%%%%%%%%%
%%%%%%%%%%%%%%%%%%%%%%%%%%%%%%%%%%%%%%%%%%%%%%%%%%%%%%%%%%%%%%%%%
%%%%%%%%%%%%%%%%%%%%%%%%%%%%%%%%%%%%%%%%%%%%%%%%%%%%%%%%%%%%%%%%%
%%%%%%%%%%%%%%%%%%%%%%%%%%%%%%%%%%%%%%%%%%%%%%%%%%%%%%%%%%%%%%%%%

\newpage 
\section{Regular parabolic bi-Lipschitz maps for stopping-time regions} \label{sec:regular}

Given a measure $\mu$ and a stopping-time region $\calS \subset \cubes$ in which the density of $\mu$ is approximately constant, each cube is balanced and the Jones function is bounded (summed over cubes in the stopping-time region), we show in this section how to construct a corresponding CCBP which satisfies the hypotheses of Theorem \ref{t:bi-Lip}. In particular, this allows us to construct an approximating bi-Lipschitz map which, in fact, we can also show to be \textit{regular} parabolic Lipschitz. \\

The main result of this section is the following. 

\begin{prop}[Regular bi-Lipschitz approximation in a single stopping-time regime]  \label{p:construct-regular-images}
	For each $\alpha,\eta,\tau  \in (0,1)$ and $A > 1$ there is $\ve \in (0,1)$ such that the following holds. Let $\mu$ be a measure on $\R^{n+1}$ and $\calS \subseteq \cubes$ a stopping-time region satisfying the following conditions. 
	\begin{enumerate}
		\item \label{b:top-doubling} The maximal cube in $\calS$ is doubling, namely
		\begin{align}
			Q(\calS) \in \cubes^{db}.  
		\end{align}
        
		\item \label{b:doubling-are-balanced} The stopping-time regime is composed of balanced cubes, namely, for every $Q \in \calS \cap \cubes^{db}$,
		\begin{align}
			\tfrac{1}{28}B_Q \mbox{ is } \tau^2\mbox{-balanced}. 
		\end{align}
        
		\item \label{b:density-estimates} Densities are bounded at the scales of the stopping-time regime, namely, for every $Q \in \calS$, 
		\begin{align}
			\tau \Theta_\mu(2B_{Q(\calS)}) \leq \Theta_\mu(2B_Q) \leq A \Theta_\mu(2B_{Q(\calS)}).
		\end{align}
        
		\item \label{b:summability} The beta square function is small if restricted to the doubling scales of the stopping-time regime, namely, for each $Q \in \calS \cap \cubes^{db}$ there is a base point $(Z_Q,u_Q) \in Q$ and a radius $\rho_Q \in [59r(Q),60r(Q)]$ such that  
		\begin{align}
			\hspace{2em}
            \Theta_\mu(2B_{Q(\calS)})^{-1} \sum_{Q \in \calS \cap \cubes^{db}}  \beta_\mu((Z_Q,u_Q),\rho_Q)^2 \, \mathds{1}_Q(X,t) \leq \ve, \qquad (X,t) \in \R^{n+1}. 
		\end{align}
	\end{enumerate}
	Then, there exists a set $G \subseteq Q(S)$ with 
    \begin{align}
        \mu(Q(S) \setminus G) \leq \alpha \mu(Q(S))
    \end{align}
    and a \textbf{regular} parabolic bi-Lipschitz map $f \colon \R^n \to \R^{n+1}$ such that \begin{align}\label{e:contains-floor}
        z(\calS) \subseteq f(\R^n). 
    \end{align}
    and
	\begin{align}\label{e:f-is-close}
		\sup_{(X, t) \in G \cap B_{p}((Z_Q,u_Q) , \rho_Q)} \dist_p((X,t), f(\R^n)) \leq \eta \, \ell(Q), 
        \qquad 
        \text{ for all $Q \in \calS$.}
	\end{align}
    Moreover, in the case that $\mu$ is Ahlfors regular, we can obtain the stronger estimate that 
    \begin{align}\label{e:f-is-close-modified}
		\sup_{(X, t) \in \supp(\mu) \cap B_{p}((Z_Q,u_Q) , \rho_Q/2)} \dist_p((X,t), f(\R^n)) \leq \eta \, \ell(Q), 
        \qquad 
        \text{ for all $Q \in \calS$.}
	\end{align}
\end{prop}

The map $f$ from Proposition \ref{p:construct-regular-images} will be obtained via Theorem \ref{t:bi-Lip}, so we need to construct a CCBP. We will do so in Subsection~\ref{sec:construct_CCBP}, but first we need some preliminary results and extra definitions. In these two next subsections, we will follow the general scheme of \cite[Subsection 6.2]{azzam2015characterization}, strengthening it to work in the parabolic setting and actually give a bi-Lipschitz map. Later, in Subsection~\ref{subsec:regular}, new arguments will be used to show regularity of this maps (this is a purely parabolic feature related to the fact that parabolic rectifiability relies on this regularity, as opposed to its Euclidean analogue).

%%%%%%%%%%%%%%%%%%%%%%%%%%%%%%%%%%%%%%%%%%%%%%%%%%%%

\subsection{Preliminary arguments before the construction of the CCBP} 

As an elementary consequence of the properties of our stopping-time regime, we deduce that $\mu$ exhibits the desired polynomial scaling over doubling cubes. 

\begin{lem}[Scaling of $\mu$ over doubling cubes]  \label{l:volume-estimates-doubling}
    For each $Q \in \calS \cap \cubes^{db}$, 
    \begin{align}\label{e:volume-estimates-doubling}
        \mu(Q) \approx \mu(B_p((Z_Q,u_Q),\rho_Q)) \approx \Theta_\mu(2B_{Q(\calS)}) \, \rho_{Q}^{n+1}.
    \end{align}
\end{lem}
\begin{proof}
    It is easy to check, by the triangle inequality, that 
    \begin{align}
        B(Q) \;\; \subseteq \;\;  Q , \  2B_Q , \ B_p((Z_Q,u_Q),\rho_Q) \;\; \subseteq \;\; 100B(Q). 
    \end{align}
    Since $Q \in \cubes^{db}$, so that $\mu(B(Q)) \approx \mu(100B(Q))$, it then follows that 
    \begin{align}\label{e:volume-estimates-doubling-2}
        \mu(Q) \approx \mu(2B_Q) \approx \mu(B_p((Z_Q,u_Q),\rho_Q)). 
    \end{align}
    The lemma now follows after recalling that $\mu(2B_Q) = \Theta_\mu(2B_Q) \, \ell(Q)^{n+1}$, $\Theta_\mu(2B_Q) \approx \Theta_\mu(2B_{Q(\calS)})$ and $\ell(Q) \approx \rho_Q$. 
\end{proof}

Our goal now is to extend this estimate to (at least a portion of) the non-doubling cubes of the stopping-time regime. For that purpose, we will identify the cubes which have nearby doubling cubes, to be able to invoke the last result.

First, we show that any cube in $\calS$ has an ancestor, also in $\calS$, which is doubling and has comparable size. 

\begin{lem}[There exist not-too-big doubling ancestor] \label{l:nearby-doubling-cube}
There exists $\compdoub = \compdoub(\udensity,\ldensity) \geq 1$ such that the following holds. For each $Q \in \calS$ there exists $Q' \in \calS \cap \cubes^{db}$ such that $Q \subseteq Q'$ and $\ell(Q') \leq \compdoub \ell(Q)$. 
\end{lem}
\begin{proof}
    Since $Q(\calS) \in \cubes^{db}$, there always exists a doubling ancestor in $\calS$. And in fact, the least doubling ancestor cannot be too large because density decays along non-doubling chains (by Lemma~\ref{lemcad22}) and is at the same time bounded from below (by \eqref{b:density-estimates}); for more details, see \cite[Lemma 6.3]{azzam2015characterization}.
\end{proof}

To be sure that our cubes of interest also have doubling cubes below, we will restrict our attention to a subset of the stopping-time regime $\calS$.

\begin{notation}[Sub-stopping-time regime of cubes which have doubling cubes below]\label{n:T}
    Let 
\begin{align}
    \calT \coloneqq \{R \in \calS \colon R \supseteq P \mbox{ for some } P \in \calS \cap \cubes^{db} \setminus \min(\calS)\}. 
\end{align}
For each $Q \in \calS$, let $\widehat{Q}$ denote the smallest cube in $\calT \cap \cubes^{db}$ such that $Q \subseteq \widehat{Q}$.
\end{notation}

Now, since cubes in $\calT$ are designed to contain and be contained in doubling cubes in the original stopping-time regime $\calS$, they also enjoy similar estimates as the doubling ones, obtained in Lemma~\ref{l:volume-estimates-doubling}.

\begin{lem}[Scaling of $\mu$ over the subtree $\calT$]  \label{l:volume-estimates-non-doubling}
    For each $Q \in \calT$, we have 
    \begin{align}\label{e:volume-estimates-non-doubling}
        \mu(Q) \approx \mu(B_p((Z_{\widehat{Q}},u_{\widehat{Q}}),\rho_{\widehat{Q}})) \approx \Theta_\mu(2B_{Q(\calS)})\, \rho_{\widehat{Q}}^{n+1}. 
    \end{align}
\end{lem}

\begin{proof}
    Fix $Q \in \calT$. The second comparability estimate follows directly from Lemma \ref{l:volume-estimates-doubling}, so we need only deal with the first. If $Q \in \cubes^{db}$, then the first estimate follows again from Lemma \ref{l:volume-estimates-doubling}. Thus, we may suppose that $Q \not\in \cubes^{db}$.
    Let $Q'$ be the largest cube in $\calS \cap \cubes^{db}$ such that $Q' \subseteq Q$ (such a cube exist by the definition of $\calT$). By Lemma \ref{l:nearby-doubling-cube}, $\ell(Q) \approx \ell(Q') \approx \ell(\widehat{Q})$. Thus, by Lemma \ref{l:volume-estimates-doubling}, 
    \begin{align}
        \mu(Q') \approx \mu(B_p((Z_{\widehat{Q}},u_{\widehat{Q}}),\rho_{\widehat{Q}})). 
    \end{align}
    Since $Q' \subseteq Q \subseteq B_p((Z_{\widehat{Q}},u_{\widehat{Q}}) ,\rho_{\widehat{Q}})$, the above estimate implies \eqref{e:volume-estimates-non-doubling}.  
\end{proof}

These results allow us to make a choice of approximating plane for each $Q \in \calT$ that will be useful for our analysis later. Indeed, for any $Q \in \calT$, we consider its doubling ancestor $\widehat{Q}$: intuitively, it will be more beneficial to define the approximating plane at the scale of $Q$ with the direction of $L_{\widehat{Q}}$ (a plane minimizing $\beta$ at the scale of $\widehat{Q}$) because $\widehat{Q}$ carries more $\mu$-measure because of being doubling. In the next result, we show that, in fact, we can choose a plane $L^Q$ (notice the super-script) that passes through $Q$, is parallel to $L_{\widehat{Q}}$, and is still close to minimizing the corresponding $\beta$ numbers (despite having made a translation).

\begin{lem}[Translation and definition of the new approximating planes $L^Q$]  \label{l:new-base-point-and-plane}
    For each $Q \in \calT$, there is a point $(Y_Q,s_Q) \in Q$ such that 
    \begin{align}\label{e:new-base-point-and-plane-1}
        \dist((Y_Q,s_Q),L_{\widehat{Q}}) \lesssim \Theta_\mu(2B_{Q(\calS)})^{-\frac{1}{2}} \beta_\mu((Z_{\widehat{Q}},u_{\widehat{Q}}),\rho_{\widehat{Q}}) \, \rho_{\widehat{Q}}. 
    \end{align}
    Letting $L^Q$ denote the plane through $(Y_Q,s_Q)$ which is parallel to $L_{\widehat{Q}}$, we then have
    \begin{align}\label{e:new-base-point-and-plane-2}
        \rho_{\widehat{Q}}^{-n-1}\int_{B_p((Z_{\widehat{Q}},u_{\widehat{Q}}),\rho_{\widehat{Q}})} \left( \frac{\dist((X,t),L^Q)}{\rho_{\widehat{Q}}}\right)^2 \, d\mu(X, t)
        \lesssim \beta_\mu((Z_{\widehat{Q}},u_{\widehat{Q}}),\rho_{\widehat{Q}})^2. 
    \end{align}
\end{lem}

\begin{proof}
    Fix $Q \in \calT$ for the remainder of the proof. We start by proving \eqref{e:new-base-point-and-plane-1}. Let $M \geq 1$ (to be chosen sufficiently large) and suppose towards a contradiction that 
    \begin{align}
        Q \subseteq \left\{(X,t)  \colon \dist((X,t),L^Q) \geq M\Theta_\mu(2B_{Q(\calS)})^{-\frac{1}{2}}\beta_\mu((Z_{\widehat{Q}},u_{\widehat{Q}}),\rho_{\widehat{Q}}) \, \rho_{\widehat{Q}}\right\}. 
    \end{align}
By Chebyshev's inequality, the fact that $Q \subseteq B_p((Z_{\widehat{Q}},u_{\widehat{Q}}),\rho_{\widehat{Q}})$, and Lemma \ref{l:volume-estimates-non-doubling}, we obtain 
    \begin{align}
        \mu(Q) \leq \frac{\Theta_\mu(2B_{Q(\calS)})}{M^2\beta_\mu((Z_{\widehat{Q}},u_{\widehat{Q}}),\rho_{\widehat{Q}})} \beta_\mu((Z_{\widehat{Q}},u_{\widehat{Q}}),\rho_{\widehat{Q}}) \rho_{\widehat{Q}}^{n+1}  \approx \frac{\mu(Q)}{M^2}. 
    \end{align}
    However, this is a contradiction for $M$ sufficiently large. 
    
    To prove \eqref{e:new-base-point-and-plane-2}, we first note that, by \eqref{e:new-base-point-and-plane-1}, 
    \begin{align}
        \dist((X,t),L^Q) \leq \dist((X,t),L_{\widehat{Q}}) + C \Theta_\mu(2B_{Q(\calS)})^{-\frac{1}{2}} \beta_\mu((Z_{\widehat{Q}},u_{\widehat{Q}}),\rho_{\widehat{Q}}) \, \rho_{\widehat{Q}}
    \end{align}
    for each $(X,t) \in B(Z_{\widehat{Q}},u_{\widehat{Q}},\rho_{\widehat{Q}})$. Using this, and the inequality $(a+b)^2 \lesssim a^2 + b^2$ for $a,b \geq 0$, it follows that the left-hand side of \eqref{e:new-base-point-and-plane-2} is at most at constant multiple of 
    \begin{align}
        \beta_\mu ((Z_{\widehat{Q}},u_{\widehat{Q}}),\rho_{\widehat{Q}})^2 + \frac{\mu(B_p((Z_{\widehat{Q}},u_{\widehat{Q}}),\rho_{\widehat{Q}}))}{\Theta_\mu(2B_{Q(\calS)})  
        \, \rho_{\widehat{Q}}^{n+1}} \beta_\mu((Z_{\widehat{Q}},u_{\widehat{Q}}),\rho_{\widehat{Q}})^2. 
    \end{align}
    By Lemma \ref{l:volume-estimates-doubling}, this quantity is at most a constant multiple of $\beta_\mu((Z_{\widehat{Q}},u_{\widehat{Q}}),\rho_{\widehat{Q}})^2$, as required. 
\end{proof}

In the previous lemma, we have defined planes $L^Q$ at the scale of each $Q \in \calT$, but we have made translations with respect to the $\beta$-minimizing plane whenever $Q \notin \cubes^{db}$. In order to be able to later control the angles that get created between these planes along the stopping-time regime, we need to make sure that we have appropriate control of the angles at the doubling cubes (which provide the directions to the rest of the $L^Q$ planes). In the next result, which is a parabolic adaptation of \cite[Lemma 6.5]{azzam2015characterization}, we provide the control of the angles at doubling scales, in terms of beta numbers (for which we have control on how they accumulate, see \eqref{b:summability}).

\begin{lem}[Control of angles between planes associated to doubling cubes in the tree] \label{l:angle-control}
    Let $Q_1,Q_2 \in \calT \cap \cubes^{db}$ such that $\ell(Q_1) \approx \ell(Q_2)$ and let $P \in \calT \cap \cubes^{db}$ be the smallest cube such that 
    \begin{align}
        B_p((Z_{Q_1},u_{Q_1}),\rho_{Q_1}) \cup B_p((Z_{Q_2},u_{Q_2}),\rho_{Q_2}) \subseteq B_p((Z_P,u_P),\rho_{P}). 
    \end{align}
    Then, for each $i,j \in \{1,2\}$ and $(Y, s) \in L^{Q_i}$, we have 
    \begin{align}
        \dist_p((Y, s),L^{Q_j}) \lesssim \Theta_\mu(2B_{Q(\calS)})^{-\frac{1}{2}}\beta_\mu((Z_P,u_P),\rho_P) \big( \ell(Q_1) + \ell(Q_2) + \dist(y,Q_1) + \dist(y,Q_2) \big). 
    \end{align}
\end{lem}

\begin{proof}
    By symmetry, it suffices to estimate $\dist((Y, s), L^{Q_1})$ for $(Y, s) \in L^{Q_2}$. We will do so by using $L^P$ as an intermediary.
    
    \textbf{Step 1.} Let us start by comparing the beta numbers at $Q_1$ with those at $P$. Indeed, compute, using the definition of $\beta$ and the hypothesis of the lemma:
    \begin{align*}
        \beta_\mu ((Z_{Q_1}, u_{Q_1}), \rho_{Q_1})^2 
        & \leq 
        \rho_{Q_1}^{-(n+3)} \int_{B_{\rho_{Q_1}}(Z_{Q_1}, u_{Q_1})} \dist_p((X, t), L_p)^2 \, d\mu(X, t)
        \\ & \leq 
        \rho_{Q_1}^{-(n+3)} \int_{B_{\rho_P}(Z_P, u_P)} \dist_p((X, t), L_p)^2 \, d\mu(X, t)
        \\ & =
        \left( \frac{\rho_P}{\rho_{Q_1}} \right)^{n+3} \beta_\mu ((Z_P, u_P), \rho_P)^2
        \\ & \approx 
        \beta_\mu ((Z_P, u_P), \rho_P)^2,
    \end{align*}
    where in the last steps we have used that $L_p$ is the best approximating plane at the suggested scale (notice the subscript, so this plane is not the translated one, namely $L^P$), and that $\rho_P \approx \rho_{Q_1}$ (this can be argued using Lemma~\ref{l:nearby-doubling-cube}).

    \textbf{Step 2.} Using that $Q_1$ is balanced and doubling, we can estimate 
    \begin{equation*}
        \mu(B_{t_0r(Q_1)}(X_i, t_i))
        \gtrsim_\tau 
        \mu(B(Q_1))
        \approx 
        \mu(Q_1)
        \approx_{A, \tau}
        \Theta_\mu(2B_{Q(\calS)}) \rho_{Q_1}^{n+1},
    \end{equation*}
    where the last step is due to Lemma~\ref{l:volume-estimates-doubling}.

    On the other hand, we can estimate, using the inclusion in the statement and then \eqref{e:new-base-point-and-plane-2} (that is, our choice of the translated plane $L^p$)
    \begin{align*}
        \int_{B_{t_0r(Q_1)}(X_i, t_i) \cap B(Q_1)} \!\!\!\!\!\!\!\! \dist_p((X, t), L^p)^2 \, d\mu(X, t)
        & \leq 
        \int_{B_{\rho_P}(Z_p, u_P)} \!\!\!\!\! \dist_p((X, t), L^p)^2 \, d\mu(X, t)
        \\ & \lesssim 
        \rho_P^{n+3} \beta_\mu((Z_P, u_P), \rho_P)^2
        \\ & \approx 
        \rho_{Q_1}^{n+3} \beta_\mu((Z_P, u_P), \rho_P)^2,
    \end{align*}
    where in the last step we have again used that $\rho_P \approx \rho_{Q_1}$. Similarly, using an inclusion already shown in Lemma~\ref{l:volume-estimates-doubling} and \eqref{e:new-base-point-and-plane-2}, 
    \begin{align*}
        \int_{B_{t_0r(Q_1)}(X_i, t_i) \cap B(Q_1)} \!\!\!\!\!\!\!\!\!\!\!\!\!\! \dist_p((X, t), L^{Q_1})^2 \, d\mu(X, t)
        & \leq 
        \int_{B_{\rho_{Q_1}}(Z_{Q_1}, t_{Q_1})} \!\!\!\!\!\!\!\!\!\!\!\!\! \dist_p((X, t), L^{Q_1})^2 \, d\mu(X, t)
        \\ & \lesssim 
        \rho_{Q_1}^{n+3} \beta_\mu((Z_{Q_1}, u_{Q_1}), \rho_{Q_1})^2
        \\ & \lesssim 
        \rho_{Q_1}^{n+3} \beta_\mu((Z_P, u_P), \rho_P)^2,
    \end{align*}
    where the last estimate was proved in Step 1.

    The bounds from the last 3 displays together yield 
    \begin{multline*}
        \fint_{B_{t_0r(Q_1)}(X_i, t_i) \cap B(Q_1)} \left( \dist_p((X, t), L^P)^2 + \dist_p((X, t), L^{Q_1})^2 \right) \, d\mu(X, t)
        \\ \lesssim 
        \Theta_\mu(2B_{Q(\calS)})^{-1} \beta_\mu((Z_P, u_P), \rho_P)^2 \rho_{Q_1}^2,
    \end{multline*}
    so we can use Chebyshev's inequality to find $(Y_i, s_i) \in B_{t_0r(Q_1)}(X_i, t_i) \cap B(Q_1)$ such that 
    \begin{equation*}
        \max \left\{ \dist_p((X, t), L^P), \dist_p((X, t), L^{Q_1}) \right\}
        \lesssim 
        \Theta_\mu(2B_{Q(\calS)})^{-1/2} \beta_\mu((Z_P, u_P), \rho_P) \rho_{Q_1}.
    \end{equation*}

    \textbf{Step 3.} The fact that both $L^{Q_1}$ and $L^P$ are vertical planes imply that 
    \begin{equation*}
        \dist_p((Y, s), L^p) 
        =
        \dist_{\eucl}(Y, L^p \cap \{t=0\}), 
        \qquad 
        \forall \, (Y, s) \in L^{Q_1}.
    \end{equation*}
    We use this fact to rely on the elliptic result from \cite{azzam2015characterization} (namely Lemma~\ref{lem:angle_planes}) to estimate angles. 

    Indeed, let us apply Lemma~\ref{lem:angle_planes} to the spatial planes $L^{Q_1} \cap \{t=0\}$, $L^P \cap \{t=0\}$, choosing the points $\mathcal{X} := \{Y_0, Y_1, \ldots, Y_n\}$ constructed in the last Step. Since $(Y_i, s_i) \in B_{t_0r(Q_1)}(X_i, t_i) \cap B(Q_1)$ and $(X_i, t_i)$ satisfy the balancedness condition, it holds 
    \begin{equation*}
        \dist_{\eucl}(Y_i, \spann(\mathcal{X} \setminus \{Y_i\})) \gtrsim_\tau r(Q_1), 
        \qquad 
        \diam_{\eucl} (\mathcal{X}) \approx r(Q_1).
    \end{equation*}
    This means that we can apply Lemma~\ref{lem:angle_planes} with 
    \begin{equation*}
        \eta(\mathcal{X}) \approx r(Q_1), 
        \qquad 
        \varepsilon(\mathcal{X}) \approx \Theta_\mu(2B_{Q(\calS)})^{-1/2} \beta_\mu((Z_P, u_P), \rho_P)
    \end{equation*} 
    by the conclusion of Step 2 (this can be made as small as we desired by \eqref{b:summability}), which yields, for any $(Y, s) \in L^{Q_1}$, 
    \begin{multline*}
        \dist_p((Y, s), L^p)
        =
        \dist_{\eucl}(Y, L^p \cap \{t=0\})
        \\ \lesssim 
        \Theta_\mu(2B_{Q(\calS)})^{-1/2} \beta_\mu((Z_P, u_P), \rho_P) \left( \dist_{\eucl}(Y, \spann(\{Y_0, \ldots, Y_n\})) + r(Q_1) \right).
    \end{multline*}

    Now note that for any $(Y, s) \in L^{Q_1}$, we can estimate, since $(Y_i, s_i) \in B(Q_1) \subset Q_1$,
    \begin{equation*}
        \dist_{\eucl}(Y, \spann(\{Y_0, \ldots, Y_n\}))
        \lesssim
        \dist_p((Y, s), Q_1) + \ell(Q_1), 
    \end{equation*}
    which simplifies the previous display into 
    \begin{equation*}
        \dist_p((Y, s), L^p)
        \lesssim 
        \Theta_\mu(2B_{Q(\calS)})^{-1/2} \beta_\mu((Z_P, u_P), \rho_P) \left( \ell(Q_1) + \dist_p((Y, s), Q_1) \right), 
    \end{equation*}
    for every $(Y, s) \in L^{Q_1}$. 
    Interchanging the roles of $L^{Q_1}$ and $L^P$, we get 
    \begin{equation*}
        \dist_p((Y, s), L^{Q_1})
        \lesssim 
        \Theta_\mu(2B_{Q(\calS)})^{-1/2} \beta_\mu((Z_P, u_P), \rho_P) \left( \ell(Q_1) + \dist_p((Y, s), Q_1) \right), 
    \end{equation*}
    for every $(Y, s) \in L^P$.

    \textbf{Step 4.} Fix $(Y, s) \in L^{Q_2}$, and call $(Z, u)$ its projection into $L^P$. This way, we can compute 
    \begin{multline*}
        \dist_p((Y, s), L^{Q_1})
        \leq 
        \norm{(Y, s) - (Z, u)} + \dist((Z, u), L^{Q_1})
        \\ \lesssim 
        \dist((Y, s), L^P) + \Theta_\mu(2B_{Q(\calS)})^{-1/2} \beta_\mu((Z_P, u_P), \rho_P) \left( \ell(Q_1) + \dist_p((Z, u), Q_1) \right)
    \end{multline*}
    after using the result from Step 3. Noting that 
    \begin{multline*}
        \dist_p((Z, u), Q_1) 
        \leq 
        \norm{(Y, s) - (Z, u)} + \dist_p((Y, s), Q_1)
        \\ =
        \dist_p((Y, s), L^P) + \dist_p((Y, s), Q_1),
    \end{multline*}
    so the last estimate gets continued to 
    \begin{align*}
        \dist_p((Y, s), L^{Q_1})
        & \lesssim 
        \dist((Y, s), L^P)
        \\ & \quad + \Theta_\mu(2B_{Q(\calS)})^{-1/2} \beta_\mu((Z_P, u_P), \rho_P) \left( \ell(Q_1) + \dist_p((Y, s), Q_1) \right) 
        \\ & \quad + \Theta_\mu(2B_{Q(\calS)})^{-1/2} \beta_\mu((Z_P, u_P), \rho_P) \dist_p((Y, s), L^P)
        \\ & =:
        I + II + III.
    \end{align*} 

    At this point, we note that $II$ is already the term that we were seeking in the statement of the lemma. So it suffices to estimate $III$ and $I$. And indeed, $\Theta_\mu(2B_{Q(\calS)})^{-1/2} \beta_\mu((Z_P, u_P), \rho_P)$ can be made small as desired because of \eqref{b:summability}, which allows us to estimate $III \lesssim I$, and $I$ satisfies the desired bound by the results of Step 3. This finishes the proof.
\end{proof}

%%%%%%%%%%%%%%%%%%%%%%%%%%%%%%%%%%%%%%%%%%%%%%%%%%%%%%%%%%%%%%%%

\subsection{Construction of the CCBP associated to the stopping-time regime} \label{sec:construct_CCBP}

After all the information that we have deduced about what happens at the scales of our stopping-time regime $\calS$, we are now in a position to define a collection of balls and planes $(L_0,\{B_{j,k}\},\{L_{j,k}\})$ that provides nice approximations by planes at these relevant scales of the stopping-time regime. We will actually want it to define a CCBP in the sense of Definition~\ref{d:CCBP}, hence the need to properly control the angles that close-by planes form with each other: for that, the translations from Lemma~\ref{l:new-base-point-and-plane} to rely on doubling cubes for the construction, and the angle control from Lemma~\ref{l:angle-control}, will be very relevant.

\begin{notation}[The definition of the candidate to CCBP] \label{not:CCBP}
We define the collection of balls and planes $(L_0,\{B_{j,k}\},\{L_{j,k}\}\}$ step by step:
\begin{itemize}
    \item To start with, we define the approximating plane at the maximal scale $Q(\calS)$ by  
    \begin{align}
        L_0 \coloneqq L^{Q(\calS)},
    \end{align}
    where the latter plane is defined by the procedure in Lemma~\ref{l:new-base-point-and-plane}. 
    
    \item Then, to correctly combine the generations of the Reifenberg theorem (Theorems~\ref{t:DT} and \ref{t:bi-Lip}) and the David-Mattila cubes (Lemma~\ref{lem:cubes}), for each $k \geq 0$, we let $s(k) \in \Z$ be such that 
    \begin{align}
        56 \cubefactor \cuberatio^{-s(k)} \leq r_k < 56 \cubefactor \cuberatio^{-s(k)+1},
    \end{align}
    where we recall that $r_k = 10^{-k}$. We also pick $k_0$ such that $\ell(Q(\calS)) \leq r_{k_0} < A_0 \ell(Q(\calS))$, keeping the conventions from Notation~\ref{not:david-mattila}.

    \item At this point, for each $k \geq k_0$, we choose $\{(X_{j,k},t_{j,k})\}_{j \in J_k}$ as a maximal $r_k$-separated net within the selected points (found in Lemma~\ref{l:new-base-point-and-plane}) associated to the cubes of the $s(k)$-th generation of the stopping-time regime $\calT$
    \begin{align}
        \mathscr{C}_k \coloneqq \left\{(Y_Q,s_Q) \colon Q \in \calT_{s(k)} \coloneqq \calT \cap \cubes_{s(k)}\right\}
    \end{align}
    and let $\{B_{j,k}\}_{j \in J_k}$ be the corresponding collection of balls as in Definition \ref{d:CCBP}, namely $B_{j,k} := B_{r_k}(X_{j,k},t_{j,k})$.
    
    \item Finally, for each $k \geq k_0$ and $j \in J_k$, we propose as approximating plane at the scale of $B_{j,k}$ the following plane produced by the method of Lemma~\ref{l:new-base-point-and-plane}:
    \[L_{j,k} \coloneqq L^{Q_{j,k}},\] 
    where $Q_{j,k} \in \calT \cap \cubes_{s(k)}$ is the cube used in last paragraph, which satisfies $(X_{j,k},t_{j,k}) = (Y_{Q_{j,k}} , s_{Q_{j,k}})$. 
\end{itemize}
\end{notation}

We shall spend the remainder of this subsection proving the following result, that asserts that this construction actually gives rise to a CCBP in the sense of Definition~\ref{d:CCBP}, with additional control on the oscillations of close-by approximating planes.

\begin{prop} \label{p:CCBP_check}
    The triple $(L_0 , \{B_{j,k}\} , \{L_{j,k}\})$ just constructed defines a $C\ve$-CCBP which satisfies the summability condition \eqref{e:ve-summable} with constant $C\ve$. 
\end{prop}

The proof of this result will follow from the rest of results in this subsection.

We start by checking Definition \ref{d:CCBP} \eqref{CCBP1} and \eqref{CCBP2}, which follow just by unraveling the definitions.

\begin{lem}[Verification of \eqref{CCBP1}, nested generations]
For each $k \geq 1$ and $j \in J_k$, 
    \begin{align}\label{e:in-V-2}
        (X_{j,k},t_{j,k}) \in V_{k-1}^2. 
    \end{align}
\end{lem}
\begin{proof}
Fix some $k \geq 0$ and $j \in J_k$. Suppose to begin with that $s(k-1) = s(k)$. In this case, it follows that $Q_{j,k} \in \calT_{s(k-1)}$ and $(X_{j,k},t_{j,k}) \in \mathscr{C}_{k-1}$. By maximality, there exists $i \in J_{k-1}$ such that $(X_{j,k},t_{j,k}) \in B_{i,k-1}$, which finishes the proof of \eqref{e:in-V-2} in this case. Suppose instead that $s(k-1) < s(k)$. Let $R \in \calT_{s(k-1)}$ such that $Q_{j,k} \subseteq R$ and let $i \in J_{k-1}$ such that $(Y_R,s_R) \in B_{i,k-1}$. It then follows from the triangle inequality that 
\begin{align}
    d_p((X_{j,k},t_{j,k}),(X_{i,k-1},t_{i,k-1})) \leq \diam(R) + r_{k-1} \leq 56\cubefactor\cuberatio^{-s(k-1)} + r_{k-1} \leq 2r_{k-1},
\end{align}
as required.
\end{proof}

\begin{lem}[Verification of \eqref{CCBP2}, initial angle]
If $j \in J_{0}$, then 
\begin{align}
    d_{100B_{j,0}}(L_{j,0} , L_0) = 0. 
\end{align}
\end{lem}

\begin{proof}
    By construction, there is only a single index $j_0 \in J_0$ and this index satisfies $(X_{j_0,0},t_{j_0,0}) = (Y_{Q(\calS)},s_{Q(\calS)})$. In particular, $Q_{j_0,0} = Q(\calS)$ and $L_{j_0,0} = L^{Q(\calS)} = L_0$. 
\end{proof}

Before proving directly conditions \eqref{CCBP3} and \eqref{CCBP4} in Definition \ref{d:CCBP} about the changes of angles of close-by scales, we prove the following useful preliminary estimate, which gives us a way to control $\ve_k$ (recall its definition at the beginning of Subsection~\ref{sec:bilip}) at points nearby the CCBP by the $\beta$ numbers of cubes in the stopping-time.

\begin{lem}[Control of the $\varepsilon$'s from the CCBP construction by $\beta$'s from the stopping-time]  \label{l:ve-control}
    Let 
    \begin{align}
        \lambda_1 \coloneqq 100\cubefactor\cuberatio. 
    \end{align}
    Then, for each $k \geq 0$ and $(X,t) \in V_k^{100}$ there exists $Q \in \calS \cap \cubes^{db}$ satisfying $\ell(Q) \approx r_k$ such that 
    \begin{align}\label{e:ve-control}
    (X,t) \in B_p((Z_Q,u_Q),\lambda_1  \rho_Q), 
    \quad \text{and} \quad 
    \ve_k(X,t) \lesssim \Theta_\mu(2B_{Q(\calS)})^{-1/2}\beta_\mu((Z_Q,u_Q),\rho_Q). 
    \end{align}
\end{lem}

\begin{proof}
   Let $k,l \geq 0$ such that $|l -k| \leq 2$ and let $j \in J_k$ and $i \in J_l$ such that $(X,t) \in 100B_{j,k} \cap 100B_{i,l}$ (such indices exist by the assumption $(X,t) \in V_k^{100}$). Let $Q$ be the smallest cube in $\calS \cap \cubes^{db}$ such that 
   \begin{align}
       B_p((Z_{\widehat{Q}_{j,k}},u_{\widehat{Q}_{j,k}}), \rho_{\widehat{Q}_{j,k}}) 
       \cup 
       B_p((Z_{\widehat{Q}_{i,l}},u_{\widehat{Q}_{i,l}}), \rho_{\widehat{Q}_{i,l}}) 
       \subseteq 
       B_p((Z_{Q},u_{Q}),\rho_{Q}) .
   \end{align}
    Since $(X_{j,k},t_{j,k}) = (Y_{Q_{j,k}},s_{Q_{j,k}}) \in Q \subseteq B_p((Z_{\widehat{Q}_{j,k}},u_{\widehat{Q}_{j,k}}),\rho_{\widehat{Q}_{j,k}})$ and $(X,t) \in 100B_{j,k}$, it follows that  
   \begin{align}
       d_p(X,t,Z_Q,u_Q) \leq 100r_k + \rho_Q \leq 100 \cdot 56 \cubefactor \cuberatio r(Q) + \rho_Q \leq \lambda_1 \rho_Q.
   \end{align}
   The assertions $\ell(Q) \approx r_k$ and the estimate in \eqref{e:ve-control} follow directly from Lemma \ref{l:angle-control}. 
\end{proof}

Definition \ref{d:CCBP} \eqref{CCBP3} and \eqref{CCBP4} now follow at once, as the following lemma shows. 

\begin{lem}[Verification of \eqref{CCBP3} and \eqref{CCBP4}]  \label{l:angle-control-CCBP}
    Let $k \geq 0$, $l \in \{k,k+1\}$ and consider indices $j \in J_k$ and $i\in J_l$ such that $d_p((X_{j,k},t_{j,k}),(X_{i,l},t_{i,l})) \leq 100r_l$. Then,  
    \begin{align}\label{e:angle-control-CCBP}
        d_{20B_{j,k}}(L_{j,k} , L_{i,l}) \lesssim d_{100B_{j,k}}(L_{j,k} , L_{i,l}) \lesssim \ve_k(X_{j,k},t_{j,k}) \lesssim \ve.
    \end{align}
\end{lem}

\begin{proof}
    The first inequality in \eqref{e:angle-control-CCBP} is clear from the definition of the bilateral distance $d_B$. The second inequality follows from the definition of $\ve_k$ and the fact that, by assumption, $(X_{j,k},t_{j,k}) \in 100B_{j,k} \cap 100B_{i,l}$. The final inequality follows directly from Lemma \ref{l:ve-control} since, for each $Q \in \calS \cap \cubes^{db}$, it holds $\Theta_\mu(2B_{Q(\calS)})^{-1/2} \beta_\mu(Z_Q,u_Q,\rho_Q) \lesssim \ve$. 
\end{proof}
    
To finish the proof of Proposition~\ref{p:CCBP_check}, it only remains to check the summability condition \eqref{e:ve-summable}. We start by showing that the summability condition Proposition \ref{p:construct-regular-images} \eqref{b:summability} still holds for a related but larger function.

\begin{lem}[Summability of $\beta$'s associated to our CCBP]  \label{l:modified-function-summable}
    We have 
    \begin{align}
        \Theta_\mu(2B_{Q(\calS)})^{-1} \sum_{Q \in \calS \cap \cubes^{db}} \beta((Z_Q,u_Q),\rho_Q)^2 \mathds{1}_{B_p((Z_Q,u_Q),2\lambda_1\rho_Q)}(X,t) \lesssim \ve 
    \end{align}
    for all $(X,t) \in \R^{n+1}$. 
\end{lem}

\begin{proof}
    Fix $(X,t) \in \R^{n+1}$ and $k_0 \in \N$. Let $k_1$ be the largest integer $0 \leq k \leq k_0$ such that $(X,t) \in B_p((Z_Q,u_Q),2\lambda_1\rho_Q)$ for some $Q \in \calS \cap \cubes^{db}_k$. Then, choose an arbitrary cube in $\calS \cap \cubes^{db}_{k_1}$, which we denote by $Q_{k_1}$, such that 
    \begin{align}
        (X,t) \in B_p((Z_{Q_{k_1}}, u_{Q_{k_1}}), 2\lambda_1 \rho_{Q_{k_1}}). 
    \end{align}
    We will prove that 
    \begin{align}\label{e:up-to-k_1}
        \sum_{k=0}^{k_1}\sum_{Q \in \calS \cap \cubes^{db}_k} \beta_\mu((Z_Q,u_Q),\rho_Q)^2 \mathds{1}_{B_p((Z_Q,u_Q),2\lambda_1\rho_Q)}(X,t) \\
        \lesssim  \sum_{Q \in \calS \cap \cubes^{db}} \beta_\mu((Z_Q,u_Q),\rho_Q)^2 \mathds{1}_{Q}(Z_{Q_{k_1}},u_{Q_{k_1}}) .
    \end{align}
     Recalling from Proposition \ref{p:construct-regular-images} \eqref{b:summability} that the right-hand side of \eqref{e:up-to-k_1} is at most $C\ve^2$, this is sufficient to prove the lemma after taking $k_0 \to \infty$.

    To prove \eqref{e:up-to-k_1}, fix some index $0 \leq k \leq k_1$ and some cube $Q \in \calS \cap \cubes_k^{db}$ such that $(X,t) \in B_p((Z_{Q} , u_Q) , 2\lambda_1\rho_Q)$.
    We construct below a corresponding cube $R = R(Q) \in \calS \cap \cubes^{db}$ such that 
    \begin{align}\label{e:R-from-Q}
        \hspace{2em}(Z_{Q_{k_1}},u_{Q_{k_1}}) \in R, \quad \ell(R) \approx \ell(Q) \quad \mbox{ and } \quad B_p((Z_Q,u_Q),\rho_Q) \subseteq B_p((Z_R,u_R),\rho_R).
    \end{align}
     Indeed, if $\ell(Q) \geq 10^4\cubefactor\lambda_1\ell(Q(\calS))$, then we can set $R = Q(\calS)$ and the properties in \eqref{e:R-from-Q} are immediate. Thus, we will assume that 
     \begin{align}\label{e:Q-small}
        \ell(Q) \leq 10^4\cubefactor\lambda_1\ell(Q(\calS)).
     \end{align}
     In this case, we let $R$ be the smallest cube in $\calS \cap \cubes^{db}$ which contains $Q_{k_1}$ and satisfies $\ell(R) \geq 10^4\cubefactor\lambda_1\ell(Q)$ (which exists by \eqref{e:Q-small}). Note, in particular, that 
     \begin{align}
         120 \lambda_1 r(Q) \leq 2 r(R). 
     \end{align}
     The first property in \eqref{e:R-from-Q} is immediate from the definition of $R$ and the second property follows from Lemma \ref{l:nearby-doubling-cube}. The final properties follows by the triangle inequality. Indeed, if $(Z,u) \in B_p((Z_Q,u_Q),\rho_Q)$, then  
     \begin{align}
        d_p((Z,u),(Z_R,u_R)) &\leq d_p((Z,u),(X,t)) + d_p((X,t),(Z_{Q_{k_1}},u_{Q_{k_1}})) + d_p((Z_{Q_{k_1}},u_{Q_{k_1}}) , (Z_R,u_R) ) \\
        &\leq 4\lambda_1 \rho_Q + \lambda_1\rho_{Q_{k_1}} + 28r(R) \\
        &\leq 240\lambda_1r(Q) + 60\lambda_1r(Q_{k_1}) + 28r(R) \\
        &\leq 4r(R) + r(R) + 28r(R) \leq \rho_R. 
     \end{align}

     Combining the three properties in \eqref{e:R-from-Q}, one obtains
     \begin{align}
         \beta_\mu((Z_Q,u_Q),\rho_Q)^2 \mathds{1}_{B_p((Z_Q,u_Q),2\lambda_1\rho_Q)}(X,t) \lesssim \beta_\mu((Z_{R(Q)},u_{R(Q)}),\rho_{R(Q)})^2 \mathds{1}_{R(Q)}(Z_{Q_{k_1}},u_{Q_{k_1}}). 
     \end{align}
     The second and third properties in \eqref{e:R-from-Q} imply, by a simple volume argument, that 
     \begin{align}
         \left|\left\{ Q \in \calS \cap \cubes^{db} \colon R(Q) = R \right\} \right| \lesssim 1 \quad \mbox{ for each } R \in \calS \cap \cubes^{db}. 
     \end{align}
     The above two properties imply \eqref{e:up-to-k_1}.
\end{proof}

We can finally prove \eqref{e:ve-summable}, which finishes the proof of Proposition~\ref{p:CCBP_check}.

\begin{lem}[Verification of \eqref{e:ve-summable}, summability of angles]
    For each $(X,t) \in L_0$, 
    \begin{align}
        \sum_{k=0}^\infty \ve_k(f_k(X,t))^2 \lesssim \ve^2. 
    \end{align}
\end{lem}

\begin{proof}
    Fix $(X,t) \in L_0$ and $k_0 \in \N$. We will show that 
    \begin{align}\label{e:ve-sum-truncated}
        \sum_{k=0}^{k_0} \ve_k(f_k(X,t))^2 \lesssim \ve^2
    \end{align}
    for some implicit constant independent of $k_0$. This is sufficient for the proof after taking $k_0 \to \infty$. 
    
    If there is no $k \in \{0,\dots,k_0\}$ such that $f_k(X,t) \in V_k^{100}$, then the left-hand side of \eqref{e:ve-sum-truncated} is zero and there is nothing to show. Suppose otherwise, and let $k_1$ denote the largest integer $0 \leq k \leq k_0$ such that $f_k(X,t) \in V_k^{100}$. By Definition \ref{d:CCBP} \eqref{CCBP1} and maximality, we have 
    \begin{align}
    \begin{split}\label{e:in-out-of-100}
        f_k(X,t) \in V_k^{100} \quad &\mbox{ for all } \quad 0 \leq k \leq k_1 \\
        f_k(X,t) \not\in V_k^{100} \quad &\mbox{ for all } \quad k_1 < k \leq k_0. 
        \end{split}
    \end{align}
    By Lemma \ref{l:ve-control}, for each $0 \leq k \leq k_1$ there exists $Q_k \in \calS \cap \cubes^{db}$ satisfying $\ell(Q_k) \approx r_k$ and 
    \begin{align}
        f_k(X,t) \in  B_p((Z_{Q_k},u_{Q_k}),\lambda_1\rho_{Q_k}), 
        \quad \text{and} \quad 
        \ve_k(f_k(X,t))^2 \lesssim \Theta_\mu(2B_{Q(\calS)})^{-1} \beta_\mu((Z_{Q_k},u_{Q_k}),\rho_{Q_k}). 
    \end{align}
    Since $d_p(f_k(X,t) , f_{k_1}(X,t)) \lesssim \ve r_k \approx \ve \rho_{Q_k}$, it follows that $f_{k_1}(X,t) \in B_p((Z_{Q_k},u_{Q_k}),2\lambda_1\rho_{Q_k})$. Thus, we actually have 
    \begin{align}\label{e:f_{k_1}}
        \hspace{2em} \ve_k(X,t)^2 \lesssim \Theta_\mu(2B_{Q(\calS)})^{-1} \beta_\mu((Z_{Q_k},u_{Q_k}),\rho_{Q_k})^2 \, \mathds{1}_{B_p((Z_{Q_k},u_{Q_k}),2\lambda_1\rho_{Q_k})}(f_{k_1}(X,t)). 
    \end{align}
    Since $\ell(Q_k) \approx r_k$, there is a uniformly bounded number of indices $0 \leq k \leq k_1$ associated to a single cube in $T \in \calS \cap \cubes^{db}$ (in the sense that $Q_k = T$). Using this fact, along with \eqref{e:f_{k_1}} and Lemma \ref{l:modified-function-summable}, we now get  
    \begin{align}
        \sum_{k=0}^{k_1} \ve_k(X,t)^2 
        \lesssim \sum_{Q \in \calS \cap \cubes^{db}} \Theta_\mu(2B_{Q(\calS)})^{-1} \beta_\mu((Z_Q,u_Q),\rho_Q)^2 \, \mathds{1}_{B_p((Z_Q,u_Q),2\lambda_1\rho_Q)}(f_{k_1}(X,t))) \lesssim \ve. 
    \end{align}
    The sum over the indices $k_1 < k \leq k_0$ is trivially zero by \eqref{e:in-out-of-100}. 
\end{proof}

%%%%%%%%%%%%%%%%%%%%%%%%%%%%%%%%%%%

\subsection{Regularity of the constructed map} \label{subsec:regular}

In the last subsection, we have defined a CCBP in Notation~\ref{not:CCBP}, and later we have checked in Proposition~\ref{p:CCBP_check} that it satisfies the hypotheses of Theorem~\ref{t:bi-Lip}. Therefore, this Theorem informs us that the associated $f \colon P_0 \to \R^{n+1}$ as in Theorem~\ref{t:DT} is actually bi-Lipschitz. Our goal now, in order to finish the proof of Proposition \ref{p:construct-regular-images}, is to show that $f$ is \textit{regular} parabolic Lipschitz, and that it satisfies \eqref{e:f-is-close}. We start with the question of regularity, which is the key novelty. 

\begin{prop}[Regularity of $f$]  \label{p:f-regular-Lipschitz}
    The map $f \colon P_0 \to \R^{n+1}$ is regular parabolic Lipschitz. 
\end{prop}

First, let us identify $P_0$ with $\R^n$ by a rotation only in the spatial variables. Next, note that bi-Lipschitzianity of $f$ easily implies that 
\begin{align}
    \Sigma \coloneqq f(\R^d) 
\end{align}
is parabolic Ahlfors-David regular (see Definition~\ref{def:UR}). Therefore, putting together Definition~\ref{def:regular} and Remark~\ref{rem:UR_with_betas}, our goal in Proposition~\ref{p:f-regular-Lipschitz} is to show that $\Sigma$ is parabolic uniformly rectifiable, which boils down to showing the strong geometric lemma for the $\beta$-coefficients of $\Sigma$, that is, that the estimate \eqref{eq:UR_via_betas} holds for 
\begin{equation} \label{def:beta_Sigma}
    \beta_\Sigma := \beta_{\mathcal{H}^{n+1}_p \lfloor \Sigma},
\end{equation} 
where the latter are defined as in Definition~\ref{def:beta}.

After an elementary discretization argument, Proposition~\ref{p:f-regular-Lipschitz} gets reduced to showing
\begin{align} \label{claim:UR_discrete}
	\sum_{\substack{I \in \Delta \\ I \subset I_0}} \beta_\Sigma(3I)^2|I| \lesssim |I_0| \quad \mbox{ for all } I_0 \in \Delta.  
\end{align}
where $\Delta$ is the standard grid of dyadic parabolic cubes in $\R^n = \R^{n-1} \times \R$ (recall that we identified $P_0$ with $\R^n$). Notice that these are not the David-Mattila cubes $\cubes$, which are adapted to $\mu$ instead: $\Delta$ is composed of usual parabolic cylinders $I = B_{2^{-k}}(X_k, t_k)$, centered at dyadic points $(X_k, t_k)$, and with dyadic side-lengths (in this case, $\ell(I) := 2^{-k}$). Note that, because of the parabolic metric, $|I| \approx \ell(I)^{n+1}$. To ease the notation, we abbreviate 
\begin{equation*}
    \Delta(I_0) := \big\{ I \in \Delta : I \subseteq I_0 \big\}.
\end{equation*}

\begin{rem}
    It is easy to note that 
    \begin{equation} \label{eq:beta_less_Omega}
        \beta_\Sigma(3I) \lesssim \Omega_f(3I)
        \qquad \text{for any $I \in \Delta$,}
    \end{equation} 
    for the $\beta$ numbers only care about the geometric shape of the image $\Sigma = \text{im}(f)$, but $\Omega$ numbers also care about its parametrization. Indeed, if $A \cdot X + b$ is a competitor in the definition of $\Omega_f(3I)$, we can consider the vertical plane $P := \{(A\cdot X + b, t) : X \in \R^{n-1}, t \in \R\}$ and it will hold 
    \begin{equation*}
        \dist_p(f(X, t), P) \leq \abs{f(X, t) - (A\cdot X + b, t)},
    \end{equation*}
    which jointly with the fact that $|3I| \approx \ell(I)^{n+1}$ and bi-Lipschitzianity of $f$ give \eqref{eq:beta_less_Omega}. Thus, if we are able to prove that, for some $I_0$, it holds 
    \begin{align} 
    	\sum_{\substack{I \in \Delta \\ I \subset I_0}} \Omega_f(3I)^2|I| \lesssim |I_0|,
    \end{align}
    the corresponding estimate for the beta numbers as in \eqref{claim:UR_discrete} will follow. In fact, we will be able to prove the stronger estimate for the Omega numbers in several \textit{easier} cases, depending on the position of $I_0$.
\end{rem}

The idea to prove the estimate \eqref{claim:UR_discrete} relies on dividing into different cases, which are treated separately:
\begin{itemize}
    \item Controlling the flatness of graph balls which are much larger than the top ball, or live very far from it, will be easy since the map $f$ is the identity far from the top ball.
    \item For graph balls which are near the top ball and which are close (in scale and location) to some ball in the CCBP, we will be able to control the $\beta$-numbers using the summability condition \eqref{e:ve-summable}.  
    \item For balls which are near the top ball but are very far from the CCBP, we will show that the graph can be written as a the graph of $C^2$ map with good enough estimates of the time and spacial derivatives. We will then use Taylor's Theorem to control the $\beta$-numbers.
\end{itemize}

Let us put this outline into action and prove \eqref{claim:UR_discrete} along the next pages. For simplicity of the notation in what follows, assume, without loss of generality, that 
\begin{equation*}
    r_{k_0} = 1, \quad \text{that is} \quad k_0 = 0,
\end{equation*}
where we recall that $k_0$ was chosen in Notation~\ref{not:CCBP}.

Fix some $I_0 \in \Delta$. We begin by showing that the sum over large cubes in controlled. Indeed, let 
\begin{align}
	\calL \coloneqq \{I \in \Delta(I_0) \colon \ell(I) > 1\}.
\end{align}

\begin{lem}[Large cubes] \label{l:L}
	We have 
	\begin{align}
		\sum_{I \in \calL}  \Omega_f(3I)^2|I| \lesssim |I_0|.
	\end{align}
\end{lem}

\begin{proof}
    Given any $I \in \calL$, we can estimate $\Omega_f(3I)$ by approximating with the identity (i.e., choosing $A = I$, $b=1$ in \eqref{def:Omega_number}) and applying the inequality from Theorem~\ref{t:DT} \eqref{DT3}, which yields 
    \begin{equation*}
        \Omega_f(3I)^2 
        \lesssim 
        \varepsilon^2 \ell(I)^{-2}.
    \end{equation*}
    Therefore, arranging the cubes in $\calL$ by generations according to their side-length, we can estimate 
    \begin{equation*}
        \sum_{I \in \calL}  \Omega_f(3I)^2 |I|
        =
        \sum_{k \geq 0} \sum_{\substack{I \subset I_0 \\ \ell(I) = 2^k}}  \Omega_f(3I)^2 |I|
        \lesssim 
        \varepsilon^2 \sum_{k \geq 0} (2^{-k})^2 \sum_{\substack{I \subset I_0 \\ \ell(I) = 2^k}} |I|
        \leq 
        \varepsilon^2 |I_0| \sum_{k \geq 0} (2^{-k})^2
        \approx 
        \varepsilon^2 |I_0|,
    \end{equation*}
    where we have used that cubes from a fixed generation of $\Delta$ are disjoint.
\end{proof}

To deal with the small cubes, $\Delta(I_0) \setminus \calL$, we will need to identify further sub-families. First, we partition $\Delta(I_0) \setminus \calL$ into cubes which are near $B_{j_0,0}$ and those which are far from it. Indeed, let
\begin{align}
	\mathcal{N} \coloneqq \{I \in \Delta(I_0) \setminus \calL \colon 10B_{j_0,0} \cap 3I \neq\empty \emptyset\}, 
    \qquad 
	\calF \coloneqq \{I \in \Delta(I_0) \setminus \calL \colon 10B_{j_0,0} \cap 3I = \emptyset\}.
\end{align}

\begin{lem}[Far away cubes]  \label{l:F}
	We have 
	\begin{align}
		\sum_{I \in \calF} \Omega_f(3I)^2|I| = 0. 
	\end{align}
\end{lem}

\begin{proof}
    The key observation is that $f$ is just the identity in the regions meaningful for $\calF$. To show that, first note that, given any $(X_{j,k}, t_{j, k})$ for $k \geq 0$ and $j \in J_k$, there exists, by Definition~\ref{d:CCBP}.\eqref{CCBP1}, some $(X_{i_1, k-1}, t_{i_1, k-1})$ with $i_1 \in J_{k-1}$ such that $|(X_{j,k}, t_{j, k}) - (X_{i_1, k-1}, t_{i_1, k-1})| < 2 r_{k-1}$. In turn, for that new point $(X_{i_1, k-1}, t_{i_1, k-1})$, there exists a further point $(X_{i_2, k-2}, t_{i_2, k-2})$ with $i_2 \in J_{k-2}$ such that $|(X_{i_1, k-1}, t_{i_1, k-1}) - (X_{i_2, k-2}, t_{i_2, k-2})| < 2r_{k-2}$. Iterating this we obtain a chain of points in successive generations, until we reach $(X_{j_0, 0}, t_{j_0, 0})$, and it holds 
    \begin{align*}
        |(X_{j,k}, t_{j, k}) - (X_{j_0, 0}, t_{j_0, 0})|
        & \leq 
        |(X_{j,k}, t_{j, k}) - (X_{i_1, k-1}, t_{i_1, k-1})|
        \\ & \quad + |(X_{i_1, k-1}, t_{i_1, k-1}) - (X_{i_2, k-2}, t_{i_2, k-2})|
        \\ & \quad + \ldots 
        \\ & \quad + |(X_{i_{k-1}, 1}, t_{i_{k-1}, 1}) - (X_{j_0, 0}, t_{j_0, 0})|
        \\ & < 
        2 r_{k-1} + 2 r_{k-2} + \ldots + 2r_{0}
        \\ & \leq 
        4,
    \end{align*}
    recalling that $r_j = 10^{-j}$ for any $j$. This, along with \eqref{e:supp_theta}, imply that
    \begin{equation*}
        \supp \theta_{j, k}
        \subset 
        10 B_{j, k}
        =
        B_{p} (X_{j,k}, t_{j, k} , 10r_k)
        \subset 
        B_{p} (X_{j_0, 0}, t_{j_0, 0}, 10r_k + 4)
        \subset 
        10 B_{j_0, 0}
    \end{equation*}
    for $k \geq 1$. Also, $\supp \theta_{j_0, 0} \subset 10 B_{j_0, 0}$ directly by \eqref{e:supp_theta}.

    Now pick $(X, t) \in 3I$ for $I \in \calF$. By definition of $\calF$, this implies that $(X, t) \notin 10B_{j_0, 0}$, so by the conclusions in last paragraph, it holds $\theta_{j, k}(X, t) = 0$ for any $k \geq 0$ and $j \in J_k$. Thus, $\sigma_k(X, t) = (X, t)$ (see \eqref{e:sigma_k}), whence $f_k(X, t) = (X, t)$ (see \eqref{e:deff}), for any $k \geq 0$. By the convergence $f_k \to f$ asserted in Theorem~\ref{t:DT}, it holds $f(X, t) = (X, t)$. Inserting this conclusion in the definition \eqref{def:Omega_number} with $A = I, b=1$ yields $\Omega_f(3I) = 0$, from which the result follows.
\end{proof}

It now only remains to deal with the sum over $\mathcal{N}$. For this, we will need to take a further subdivision. This subdivision will be based on whether the image of the cube is close (in scale and location) to some ball in the CCBP, or whether it is instead far from any ball in the CCBP. To be more precise, we introduce the following  stopping-time distance associated to the balls from the CCBP by
\begin{align}
	\rho(X,t) \coloneqq \inf_{\substack{k \geq 0 \\ j \in J_k}} \big[ r_k + \dist((X,t),B_{j,k}) \big],
    \qquad 
    (X,t) \in \R^{n+1},
\end{align}
and we also use the notation 
\begin{align}
	\rho(I) \coloneqq \sup_{(X,t) \in I} \rho(f(X,t)), 
    \qquad 
    I \subseteq \Delta.
\end{align}

Let $\tau > 0$ be small, to be determined later, and let $\calM$ be the maximal (with respect to inclusion) set of cubes in $\Delta$ such that 
\begin{align}
	\ell(I) < \tau \rho(I), 
\end{align}
so that by maximality, it follows easily that 
\begin{align}
	\ell(I) \approx \tau \rho(I), 
    \qquad 
    I \in \mathcal{M}.
\end{align}
Cubes below $\calM$ will be considered to be far from the relevant scales of the CCBP: the most relevant behavior will happen above $\calM$. To define this more precisely, we let $\mathcal{N}_0$ be the set of maximal cubes in $\mathcal{N}$. In particular, by our definition of $\calL$ above, it holds 
\begin{align}
    \mathcal{N}_0 = \big\{J_0 \in \mathcal{N}_0 \colon \; \ell(J_0) = 1\big\}
\end{align}
Now, for $J_0 \in \mathcal{N}_0$, we define the set of cubes close to the scales of the CCBP by
\begin{align}
	\Tree(J_0) := \big\{ J \in \Delta(J_0) \colon J \mbox{ is not strictly contained in any cube from } \calM \big\}.
\end{align}

Given all these definitions, we now show that, for each $J_0 \in \mathcal{N}_0$ and $J \in \Tree(J_0)$, there is a ball in the CCBP which is close to the scale of $J$, in the sense that it has comparable size and contains $f(3J)$. 

\begin{lem}[Cubes in the stopping-time correspond to scales of the CCBP]  \label{l:correspondence-CCBP}
	For each $J_0 \in \mathcal{N}_0$ and $I \in \Tree(J_0)$, there exists $k =k(I) \geq 0$ and $j = j(I) \in J_k$ such that 
	\begin{align}\label{e:correspondence-CCBP}
		\ell(I) \approx r_k \quad \mbox{ and } f(3I) \subseteq 10B_{j,k}. 
	\end{align} 
	Moreover, there is bounded correspondence of the mapping $I \mapsto (k,j)$, i.e., 
	\begin{align}\label{e:bounded-correspondence}
		|\{I \in \Tree(J_0) \colon k(I) = k \mbox{ and } j(I) = j\}| \lesssim 1 
	\end{align}
	for each $k \geq 0$ and $j \in J_k$, with an implicit constant which is independent of $k, j$.
\end{lem}

\begin{proof}
    Fix $I$ as in the statement. By definition of $\Tree$, it holds $\rho(I) \leq \tau^{-1} \ell(I)$. This means that there exists $k^{(0)} \geq 0$ and $j^{(0)} \in J_{k^{(0)}}$ such that 
    \begin{equation*}
        \dist(f(X_I, t_I), B_{j^{(0)}, k^{(0)}}) + r_{k^{(0)}}
        \leq
        \tau^{-1} \ell(I),
    \end{equation*}
    where $(X_I, t_I)$ is the center of $I$ (although any point works). Using \eqref{CCBP1} iteratively, we can find $k^{(1)}, k^{(2)}, \ldots, k$ and associated $j^{(1)}, j^{(2)}, \ldots, j$ such that 
    \begin{equation*}
        \norm{(X_{j^{(i)}, k^{(i)}}, t_{j^{(i)}, k^{(i)}})
        - (X_{j^{(i+1)}, k^{(i+1)}}, t_{j^{(i+1)}, k^{(i+1)}})}
        \leq 
        2 r_{k^{(i+1)}}.
    \end{equation*}
    We stop the first time that it holds $r_k \geq 4\ell(I) / \tau$ (the cubes for which $\ell(I) > \tau / 4$ can be handled more easily associating them to $J_0$, so let us concrentrate on the more interesting case when $\ell(I) \leq \tau / 4$). Thus, the by minimality of this choice of $k$, it holds 
    \begin{equation*}
        \frac{4\ell(I)}{\tau} 
        \leq 
        r_k
        =
        10 r_{k+1}
        < 
        10 \frac{4\ell(I)}{\tau},
    \end{equation*}
    which already gives the first item in the statement.

    To prove the second item, fix $(X, t) \in 3I$. Then, we may estimate, using that $f$ is Lipschitz (Theorem~\ref{t:bi-Lip}, Lemma~\ref{l:f-upper-reg}), 
    \begin{align*}
        \norm{f(X, t) - (X_{j, k}, t_{j, k})}
        & \leq 
        \norm{f(X, t) - f(X_I, t_I)}
        + \norm{f(X_I, t_I) - (X_{j^{(0)}, k^{(0)}}, t_{j^{(0)}, k^{(0)}})}
        \\ & \qquad + \norm{(X_{j^{(0)}, k^{(0)}}, t_{j^{(0)}, k^{(0)}}) - (X_{j, k}, t_{j, k})}
        \\ & \leq 
        \norm{f(X, t) - f(X_I, t_I)}
        + \dist_p(f(X_I, t_I), B_{j^{(0)}, k^{(0)}}) + r_{k^{(0)}}
        \\ & \qquad + \sum_i \norm{(X_{j^{(i)}, k^{(i)}}, t_{j^{(i)}, k^{(i)}}) - (X_{j^{(i+1)}, k^{(i+1)}}, t_{j^{(i+1)}, k^{(i+1)}})}
        \\ & \leq
        C(M) \, \diam_p(3I) + \tau^{-1} \ell(I) + 2 (r_{k^{(0)}} + r_{k^{(1)}} + \ldots + r_k)
        \\ & \leq 
        C(n, M)\, \ell(I) + \frac{r_k}{4} + 4 r_k
        \\ & \leq 
        \frac{C(n, M)}{4} \tau r_k + \frac{r_k}{4} + 4 r_k
        \\ & \leq 
        6 r_k
    \end{align*}
    just by choosing $\tau > 0$ small enough. This implies $f(3I) \subset 10B_{j, k}$.

    Lastly, given $I \in \Tree(J_0)$ and its associated $k, j$, let us estimate the amount of $I' \in \Tree(J_0)$ that can share the same $k, j$. First note that, in such case, it must hold $\ell(I') \approx r_k \approx \ell(I)$; which means that $I'$ is in a dyadic generation which is boundedly away from that of $I$. Moreover, since both $f(3I)$ and $f(3I')$ are contained in $10B_{j, k}$, and $f$ is bi-Lipschitz (Theorem~\ref{t:bi-Lip}), we can estimate 
    \begin{equation*}
        \norm{(X_I, t_I) - (X_{I'}, t_{I'})}
        \lesssim
        \norm{f(X_I, t_I) - f(X_{I'}, t_{I'})}
        \leq 
        \diam_p(10B_{j, k})
        \lesssim 
        r_k.
    \end{equation*}
    This geometric closeness, together with the closeness in dyadic generation between $I$ and $I'$ established above, easily imply (by elementary properties of dyadic grids) that the amount of $I'$ than can share $k, j$ with $I$ is bounded uniformly, which finishes the proof of the lemma.
\end{proof}

We now show that our assignment $I \mapsto (k(I), j(I))$ is, somehow, geometrically stable. Concretely, we show that neighboring cubes are assigned neighboring balls, in the concrete sense below.

\begin{lem}[Nearby cubes get assigned nearby balls]  \label{l:M-cubes-small}
    There exists $C_* \geq 1$ such the following holds. For each $J \in \mathcal{N}_0$ and $I \in \Tree(J_0)$, for each $J \in \calM$ it holds
    \begin{align}
        3I \cap J \neq \emptyset
        \quad \implies \quad 
        \ell(J) \lesssim \ell(I)
        \quad \text{and} \quad 
        10B_{j(J),k(J)} \subseteq C_*B_{j(I),k(I)}. 
    \end{align}
\end{lem}
\begin{proof}
    Choose $(X,t)\in3I\cap J$. By Lemma~\ref{l:correspondence-CCBP},
    $f(X,t)\in10B_{j(I),k(I)}$, so by definition of $\rho$ we get
    \begin{equation*}
        \rho(f(X,t))
        \leq r_{k(I)}
           +\dist_p\bigl(f(X,t),B_{j(I),k(I)}\bigr)
        \lesssim r_{k(I)}.
    \end{equation*}
    Thus, using that $\rho$ is Lip(1, 1/2) with constant 1, and $f$ is also Lip(1, 1/2) (see Theorem~\ref{t:bi-Lip} or Lemma~\ref{l:f-upper-reg}), we get 
    \begin{multline*}
        \rho(J)
        \leq 
        \rho(f(X, t)) + \sup_{(Y, s) \in J} \abs{\rho(f(Y, s)) - \rho(f(X, t))}
        \leq 
        C r_{k(I)} + C(M) \,d_p\bigl((Y,s),(X,t)\bigr)
        \\ \leq 
        C r_{k(I)} + C(n, M) \, \ell(J)
        < 
        C r_{k(I)} + C(n, M) \, \tau \rho(J),
    \end{multline*}
    where in the last estimate we have used the fact that $J \in \calM$. Thus, choosing $\tau$ small enough, we can hide the last term to infer (again using $J \in \calM$)
    \begin{equation*}
        \ell(J) < \tau \rho(J) \lesssim r_{k(I)} \approx \ell(I),
    \end{equation*}
    where the last estimate comes from Lemma~\ref{l:correspondence-CCBP}. This gives the first assertion of the lemma.

    To prove the second assertion, consider $J_0' \in \mathcal{N}_0$ the dyadic ancestor of $J$ of side
    length $1$. By definition of $\Tree$ and $\calM$, it holds $J \in \Tree(J_0')$. Hence, applying Lemma~\ref{l:correspondence-CCBP}, we obtain
    \begin{equation*}
        f(3J) \subset 10B_{j(J),k(J)},
        \quad \text{and} \quad 
        r_{k(J)}
        \lesssim
        \ell(J)
        \lesssim
        \ell(I)
        \lesssim 
        r_{k(I)}.
    \end{equation*}
    Thus, for any $(Y,s)\in10B_{j(J),k(J)}$, it holds 
    \begin{multline*}
        d_p((Y,s), (X_{j(I),k(I)}, t_{j(I),k(I)}))
        \leq
        d_p((Y,s), f(X,t))
        + d_p(f(X,t), (X_{j(I),k(I)}, t_{j(I),k(I)}))
        \\ 
        \lesssim 
        r_{k(J)} + r_{k(I)}
        \lesssim 
        r_{k(I)},
    \end{multline*}
    where we have used that $f(X, t) \in f(J) \subset 10B_{j(J), k(J)}$, and $f(X, t) \in 10B_{j(I), k(I)}$. From this, the second assertion in the lemma follows.
\end{proof}

On the other hand, let us state and quickly prove some elementary properties of the balls obtained after the assignment from Lemma~\ref{l:correspondence-CCBP}.

\begin{lem}[Elementary properties of $B_{j(I), k(I)}$]  \label{l:B-j-k-properties}
    Let $I \in \Tree(J_0)$ and set $k = k(I)$ and $j = j(I)$. First, we have 
    \begin{align}\label{e:B-j-k-contains-a-cube}
        2B_{Q_{j,k}} \subseteq 2B_{j,k}. 
    \end{align}
    Second, if $C_*r_{k} \leq 1$, there exists $Q \in \calS \cap \cubes^{db}$ with $\ell(Q) \approx r_{k}$ such that 
    \begin{align}\label{e:B-j-k-contained-in-a-cube}
        C_*B_{j,k} \subseteq B_p((Z_Q,u_Q),\rho_Q)
    \end{align}
    and also 
    \begin{align}\label{e:volume-B-j-k}
        \mu(\lambda B_{j,k}) \approx \Theta_\mu(2B_{Q(\calS)}) \, r_k^{n+1}, 
        \qquad 
        \text{for any $\lambda \in [2,C_*]$.}
    \end{align}
    Moreover, there is an absolute constant $C_1 \geq 1$ (i.e., independent of $I$ and $J_0$) such that 
    \begin{align}
        Q \subseteq B^* \coloneqq B_p(f(X_{I_0},t_{I_0}) , C_1\ell(I_0)).
    \end{align}
    
\end{lem}

\begin{proof}
    Note that $Q_{j,k} \in \calT$ by Notation~\ref{not:CCBP}. The inclusion \eqref{e:B-j-k-contains-a-cube} follows the triangle inequality after unraveling the appropriate definitions. For \eqref{e:B-j-k-contained-in-a-cube}, we let $Q$ be the smallest ancestor of $Q_{j,k}$ in $\calS \cap \cubes^{db}$ such that $C_* r_k \leq r(Q)$ (such cube always exists because of the assumption that $C_*r_k \leq 1$, recalling that $1 \leq r(Q(\calS))$). Since 
    \begin{align}
        (X_{j,k},t_{j,k}) \in Q_{j,k} \subseteq Q \quad \mbox{ and } \quad (Z_Q,u_Q) \in Q,
    \end{align}
    for any $(X,t) \in C_*B_{j,k}$ we have 
    \begin{align}
        d_p((X,t),(Z_Q,u_Q)) \leq C_*r_k + 56r(Q) \leq 57r(Q) \leq \rho_Q,
    \end{align}
    which is \eqref{e:B-j-k-contained-in-a-cube}. In turn, noting that $Q \in \calS \cap \cubes^{db}$ and $\ell(Q) \approx \ell(Q_{j,k}) \approx r_k$, the estimate \eqref{e:volume-B-j-k} is an immediate corollary of Lemma \ref{l:volume-estimates-doubling}, Lemma \ref{l:volume-estimates-non-doubling}, \eqref{e:B-j-k-contains-a-cube} and \eqref{e:B-j-k-contained-in-a-cube}. The final implication holds by using the fact that $I \in \mathcal{N}$ (it is close to $I_0$) and the comparability estimates $\ell(Q) \sim r_k \sim \ell(I) \leq \ell(I_0)$. 
\end{proof}

The fact that we have established, at least for the small cubes in these $\Tree$ (those for which $C_*r_k \leq 1$) are related to the cylinders $B_p((Z_Q, u_Q), \rho_Q)$, for which our assumption \eqref{b:summability} from Proposition~\ref{p:construct-regular-images} gives summability of their associated beta numbers, implies that we can sum also sum the beta numbers associated to a subtree.

\begin{lem}[Summability of betas for small cubes in the trees]  \label{l:Tree'}
    For $J_0 \in \mathcal{N}_0$, let 
    \begin{align}
        \Tree'(J_0) := \big\{I \in \Tree(J_0) \colon C_*r_{k(I)} \leq 1 \big\}.
    \end{align}
    Then, 
    \begin{align}
        \Theta_\mu(2B_{Q(\calS)})^{-1} \sum_{J_0 \in \mathcal{N}_0}\sum_{I \in \Tree'(J_0)} \beta_\mu(C_*B_{j(I),k(I)})^2 \, r_{k(I)}^{n+1} \lesssim \ell(I_0)^{n+1}.  
    \end{align}
\end{lem}

\begin{proof}
    Since there are only finitely many cubes in $\mathcal{N}_0$, it suffices to bound the sum over a single family $\Tree'(J_0)$ for some $J_0$ which we now fix. By Lemma \ref{l:B-j-k-properties}, for each $I \in \Tree'(J_0)$ there exists $Q_I \in \calS \cap \cubes^{db}$ such that 
    \begin{align}\label{e:Q_I}
        \ell(Q_I) \approx r_{k(I)}, \quad Q_I \subseteq B^* \quad \mbox{and} \quad C_*B_{j(I),k(I)} \subseteq B_p((Z_{Q_I},u_{Q_I}),\rho_{Q_I}).
    \end{align}
    The first and third properties, together with Lemma \ref{l:volume-estimates-doubling}, imply that \begin{align}
        \beta_\mu(C_*B_{j(I),k(I)}) \, r_{k(I)}^{n+1} 
        \lesssim 
        \beta_\mu((Z_{Q_I},u_{Q_I}),\rho_{Q_I})\, \Theta_\mu(2B_{Q(\calS)})^{-1} \, \mu(Q_I).
    \end{align}
    Moreover, they also imply that the mapping $I \mapsto Q_I$ has uniformly bounded correspondence. The two facts, together with the second property in \eqref{e:Q_I} and Proposition~\ref{p:construct-regular-images}.\eqref{b:summability}, yield
    \begin{align}
    \sum_{I \in \Tree'(J_0)} & \Theta_\mu(2B_{Q(\calS)})^{-1} \beta_\mu(C_*B_{j(I),k(I)})^2 \, r_{k(I)}^{n+1} 
    \\ & \lesssim 
    \sum_{\substack{Q \in \calS \cap \cubes^{db} \\ Q \subseteq B^*}} \Theta_\mu(2B_{Q(\calS)})^{-2}\beta_\mu((Z_Q,u_Q),\rho_Q)^2 \mu(Q)
    \\ & \leq 
    \int_{Q(\calS) \cap B^*} \sum_{Q \in \calS \cap \cubes^{db}} \Theta_\mu(2B_{Q(\calS)})^{-2}\beta_\mu((Z_Q,u_Q),\rho_Q)^2 \mathds{1}_Q(X,t) \, d\mu(X,t)
    \\ & \leq 
    \ve \, \Theta_\mu(2B_{Q(\calS)})^{-1} \mu(Q(\calS) \cap B^*) 
    \\ & \lesssim 
    \ve \, \ell(I_0)^{n+1}, 
    \end{align}
    where the last estimate is due to our hypothesis on density estimates in Proposition~\ref{p:construct-regular-images}.\eqref{b:density-estimates}, recalling that $r(B^*) = C_1\ell(I_0)$ from Lemma~\ref{l:B-j-k-properties}.
\end{proof}

Now that we have obtained many properties about the assignment $I \mapsto B_{j(I), k(I)}$, let us see how can estimate the associated $\beta_\Sigma$ sum (from \eqref{claim:UR_discrete}) for the cubes in $\mathcal{N}$. We start by the following estimate, which follows upon a covering argument by cubes in $\mathcal{M}$ and taking into account only the definition of $B_{j(I), k(I)}$, and not the last properties that we have inferred for these balls.

\begin{lem}[Elementary estimate for $\beta_\Sigma$ in the regime $\mathcal{N}_0$]  \label{l:near}
	For each $J_0 \in \mathcal{N}_0$ and $I \in \Tree(J_0)$, if $k = k(I)$ and $j = j(I)$ are as in Lemma \ref{l:correspondence-CCBP}, then
	\begin{align}\label{e:omega-estimate-local}
		\beta_\Sigma(f(3I))^2|I| 
        \lesssim 
        \Theta_\mu(2B_{Q(\calS)})^{-1} \beta_\mu(C_*B_{j,k})^2r_k^{n+1} 
        + \sum_{\substack{J \in \calM \\ J \cap 3I \neq \emptyset}} \left( \frac{\ell(J)}{\ell(I)} \right)^2|J|.
	\end{align}
\end{lem}

\begin{proof}
	Fix $J_0 \in \mathcal{N}_0$, $I \in \Tree(J_0)$ and $k=k(I)$, $j = j(I)$ as above. Let $\ve > 0$ and choose a $d$-dimensional plane $P$ such that 
	\begin{align}
		\beta_\mu(C_*B_{j,k},P)^2 \leq \beta_\mu(C_*B_{j,k})^2 + \ve. 
	\end{align} 
    
	To get an estimate for $\beta_\Sigma(f(3I))$, let us cover $3I$ with cubes from $\calM$, namely define 
	\begin{align}
		\calM_I \coloneqq \{J \in \calM \colon J \cap 3I \neq \emptyset\},
	\end{align}
    and, in fact, we partition $\calM_I$ into two sets, those whose image is (relatively) near to $P$ and those whose image is far from $P$. Indeed, let $\lambda \geq 1$ (to be chosen sufficiently large) and set 
	\begin{align}
			\calM_I^n \coloneqq \big\{J \in \calM \colon \dist(f(J) , P) < \lambda\ell(J)\big\}, 
            \qquad 
			\calM_I^f \coloneqq \big\{J \in \calM \colon \dist(f(J) , P) \geq \lambda\ell(J)\big\},
	\end{align}
    so that, letting $\nu\coloneqq \calH_{\rm p}^{n+1}\lfloor{\Sigma}$ and recalling the definition of $\beta_\Sigma$ in \eqref{def:beta_Sigma}, we can estimate 
	\begin{align}
        \beta_\Sigma(f(3I))^2 \, |I| 
        \lesssim 
        \sum_{J \in \calM_I^n} \int_{f(J)} \left( \frac{\dist((X,t),P)}{\ell(I)} \right)^2 \, d\nu + \sum_{J \in \calM_I^f} \int_{f(J)} \left( \frac{\dist((X,t),P)}{\ell(I)} \right)^2 \, d\nu
        =:
		I_1 + I_2.
	\end{align}
    
	Using the definition of $\calM_I^n$ and the properties $\diam(f(J)) \approx \ell(J)$ and $\nu(f(J)) \approx |J|$, we can estimate 
	\begin{align}
		I_1 \lesssim_\lambda \sum_{J \in \calM_I^n} \int_{f(J)} \left( \frac{\ell(J)}{\ell(I)} \right)^2 \, d\nu \leq \sum_{\substack{J \in \calM \\ J \cap 3I \neq \emptyset}} \left( \frac{\ell(J)}{\ell(I)} \right)^2|J|.
	\end{align}
    
	In turn, to estimate $I_2$ we proceed as follows. First, for each $J \in \calM_I^f$, consider the index pair $(k(J),j(J))$ obtained from Lemma \ref{l:correspondence-CCBP}. If $\lambda$ is sufficiently large, depending on the Lipschitz constant of $f$ and the implicit constant in \eqref{e:correspondence-CCBP}, then 
    \begin{equation*}
        d_p((X,t),f(X_J,t_J)) \leq \lambda \ell(J)/2 \qquad \text{for all $(X,t) \in f(J) \cup 10B_{j(J),k(J)}$.}
    \end{equation*} 
    Therefore, since $J \in \calM_I^f$, we have $\dist(f(X_J,t_J),P) \geq \lambda \ell(J)$ and we obtain
    \begin{align}
        \dist_p((X,t),P) \approx \dist(f(X_J,t_J) , P) \qquad \mbox{ for all } (X,t) \in f(J) \cup 10B_{j(J),k(J)}. 
    \end{align}
    On the other hand, by \eqref{e:volume-B-j-k} and Lemma \ref{l:correspondence-CCBP}, we have 
	\begin{align}
		\mu(10B_{j(J),k(J)}) \approx \Theta_\mu(2B_{Q(\calS)}) \, r_{k(J)}^{n+1} \approx \Theta_\mu(2B_{Q(\calS)}) \, |J|.
	\end{align}
	Combing the above two estimates, we obtain 
	\begin{align}
		 \int_{f(J)} \left( \frac{\dist((X,t),P)}{\ell(I)} \right)^2 \, d\nu \approx \Theta_\mu(2B_{Q(\calS)})^{-1} \int_{10B_{j(J),k(J)}} \left( \frac{\dist((X,t),P)}{\ell(I)} \right)^2 d\mu, 
	\end{align}
	so summing over $J \in \calM_I^f$, and later using \eqref{e:bounded-correspondence} and the fact that $10B_{j(J),k(J)} \subseteq C_*B_{j,k}$ for each $J \in \calM_I^f$ (recall Lemma \ref{l:M-cubes-small}), we have 
	\begin{align}
		I_2 &\lesssim_A \Theta_\mu(2B_{Q(\calS)})^{-1} \sum_{J \in \calM_I^f}   \int_{10B_{j(J),k(J)}} \left( \frac{\dist((X,t),P)}{\ell(I)} \right)^2 d\mu \\
		&\lesssim \Theta_\mu(2B_{Q(\calS)})^{-1} \int_{20B_{j,k}}  \left( \frac{\dist((X,t),P)}{\ell(I)} \right)^2 d\mu\\
		&\lesssim \Theta_\mu(2B_{Q(\calS)})^{-1} \left( \beta_\mu(C_*B_{j,k})^2 + \ve \right) r_k^{n+1}.
	\end{align}
	Taking $\ve \to 0$ finishes the proof of \eqref{e:omega-estimate-local}.
\end{proof}

Putting together the last elementary estimate for individual terms $\beta_\Sigma(f(3I))$, with the rest of more refined estimates for the balls $B_{j, k}$ that we obtained, for instance, in Lemma~\ref{l:Tree'} using the key hypotheses of Proposition~\ref{p:construct-regular-images}, we can finally obtain the summability for the beta numbers associated all the cubes inside $\Tree$ structures rooted in $\mathcal{N}$.

\begin{lem}[Summability of betas across the $\Tree$ structures rooted at $\mathcal{N}$]
    We have  
	\begin{align}\label{e:omega-estimate-global}
		\sum_{J_0 \in \mathcal{N}_0} \sum_{I \in \Tree(J_0)} \beta_\Sigma(f(3I))^2 \, |I| \lesssim |I_0|. 
	\end{align}
\end{lem}

\begin{proof}
    Since there are only finitely many cubes in $\mathcal{N}_0$, we need only estimate the sum over a single family $\Tree(J_0)$ for some $J_0$ which we now fix. Recalling from Lemma \ref{l:correspondence-CCBP} that $\ell(I) \approx r_{k(I)}$, any $I \in \Tree(J_0) \setminus \Tree'(J_0)$ must satisfy $\ell(I) \approx 1 = \ell(I_0)$, i.e., the cubes from $\Tree(J_0) \setminus \Tree(J_0)$ only belong to finitely many generations. Since $\Sigma$ is $C\ve$-Reifenberg flat, it must be that $\beta_\Sigma(f(3I)) \lesssim \ve$ for all $I \in \Tree(J_0)$. These two facts together imply 
    \begin{align}
        \sum_{I \in \Tree(J_0) \setminus \Tree'(J_0)} \beta_\Sigma(f(3I))^2|I| \lesssim |I_0|. 
    \end{align}

    To deal with the sum over cubes in $\Tree(J_0)$, we first apply Lemmas \ref{l:Tree'} and \ref{l:near} to get 
    \begin{align}
        \sum_{I \in \Tree'(J_0)} \beta_\Sigma(f(3I))^2|I| &\lesssim \sum_{I \in \Tree'(J_0)} \Theta_\mu(2B_{Q(\calS)})^{-1} \beta_\mu(C_*B_{j(I),k(I)})^2 r_{k(I)}^{n+1} \\
        &\hspace{4em}+ \sum_{I \in \Tree'(J_0)} \sum_{\substack{J \in \calM \\ J \cap 3I \neq \emptyset}} \left( \frac{\ell(J)}{\ell(I)}\right)^2|J| \\
        &\lesssim |I_0| + \sum_{I \in \Tree'(J_0)} \sum_{\substack{J \in \calM \\ J \cap 3I \neq \emptyset}} \left( \frac{\ell(J)}{\ell(I)}\right)^2|J|. 
    \end{align}
    The last term is estimated by swapping the order of summation and summing over a geometric series. More precisely, consider cubes $J \in \calM$ and $I \in \Tree(J_0)$ such that $J \cap 3I \neq \emptyset$. Since $I \subseteq J_0$, we must have $J \cap 3J_0 \neq \emptyset$. Furthermore, by Lemma \ref{l:M-cubes-small} there exists $k_0 \in \N$ such that $\ell(I) \geq 2^{-k_0} \ell(J)$. It follows from this last estimate and doubling that, for each $J \in \calM$ and $k \in \N$, there are a uniformly bounded number of cubes $I \in \Tree(J_0) \cap \Delta_k$ such that $J \cap 3I \neq \emptyset$. Putting all of this together, we get 
    \begin{align}
    \sum_{I \in \Tree'(J_0)} \sum_{\substack{J \in \calM \\ J \cap 3I \neq \emptyset}} \left( \frac{\ell(J)}{\ell(I)}\right)^2|J| &\leq \sum_{\substack{J \in \calM \\ J \cap 3J_0 \neq \emptyset }} \sum_{\substack{I \in \Tree'(J_0) \\ J \cap 3I \neq \emptyset\\ \ell(I) \geq 2^{-k_0}\ell(J)}} \left( \frac{\ell(J)}{\ell(I)}\right)^2|J|\\
    &\lesssim \sum_{\substack{J \in \calM \\ J \cap 3J_0 \neq \emptyset }} \sum_{k=-k_0}^\infty 2^{-2k}|J| \\
    &\lesssim \sum_{\substack{J \in \calM \\ J \cap 3J_0 \neq \emptyset}} |J| \lesssim |J_0| = |I_0|.  
    \end{align}
\end{proof}

Now that we have obtained an estimate for the $\Tree$'s rooted at $\mathcal{N}$, we are left to obtain estimates for the scales which are below all these $\Tree$'s. Recalling their definition, this means that we are left to obtain estimates below the family $\calM$. In this case, by definition of $\calM$, the relevant scales are much below those of the CCBP, and this means that the construction of the function $f$ stopped at a finite generation $f_k$, and we will be able to use the $C^2$ smoothness estimates of these finite iterates $f_k$ (see Theorem~\ref{t:DT}) to deduce the summability estimates for the associated Omega numbers.

\begin{lem}[Cubes much below the scale of the CCBP, below $\calM$] \label{l:far}
	For each $I \in \calM$ such that $I \subseteq I_0$, we have 
	\begin{equation*}
		\sum_{J \subseteq I} \Omega_f(3J)^2|J| \lesssim |I|. 
	\end{equation*}
\end{lem}

\begin{proof}
    \textbf{Step 1:} the construction of $f$ stabilizes over $3I$. Let $(X, t) \in 3I$. We will prove that, if $k \geq K$ for $K$ large enough (we can take $K$ such that it holds $\frac{\tau^{-1}\ell(I)}{1000} \leq r_K < \frac{\tau^{-1}\ell(I)}{100}$), then $f_k(X, t) \notin V_k^{10}$. For the sake of a contradiction, assume that $f_k(X, t) \in 10B_{j, k}$ for $j \in J_k$.

    First, note that by definition of $\rho$ and $\calM$, there exists $(X_I, t_I) \in I$ such that $\rho(f(X_I, t_I)) > \tau^{-1}\ell(I)$. Moreover, using that $f$ is Lipschitz (with constant $C(M)$) by Theorem~\ref{t:bi-Lip}, we obtain 
    \begin{equation*}
        \abs{\rho(f(X_I, t_I)) - \rho(f(X, t))}
        \leq 
        \abs{f(X_I, t_I) - f(X, t)}
        \leq 
        C(M) \norm{(X_I, t_I) - (X, t)}
        \leq 
        C(M, n) \ell(I),
    \end{equation*}
    so that if $\tau$ is chosen small enough depending on $M$ and $n$, we have 
    \begin{equation} \label{eq:rho_below_3I}
        \rho(f(X, t)) \geq \frac{\tau^{-1}}{2} \ell(I).
    \end{equation}

    On the other hand, since $f_k$ converges to $f$ uniformly (see Theorem~\ref{t:DT}), we can estimate 
    \begin{align*}
        \rho(f(X, t))
        & \leq 
        \rho(f_k(X, t)) + \abs{\rho(f(X, t)) - \rho(f_k(X, t))}
        \\ & \leq 
        \left( \dist_p(f_k(X, t), B_{j, k}) + r_k \right) + \abs{f(X, t) - f_k(X, t)}
        \\ & \leq 
        \left( 10r_k + r_k \right) + \sum_{m=K}^\infty \abs{f_{m+1}(X, t) - f_m(X, t)}
        \\ & \leq 
        11 r_k + C\varepsilon r_K
        \leq 
        20 r_K,
    \end{align*}
    where we have used \eqref{Sigma_ksigma_k} to estimate the difference between two consecutive stages of the construction of $f$, and chosen $\varepsilon$ small. If $K$ is chosen large enough so that $r_K \ll \ell(I)$, this is clearly incompatible with \eqref{eq:rho_below_3I}, which proves the claim. \\  

    \noindent\textbf{Step 2:} the derivative estimates. Step 1 implies, by Theorem~\ref{t:DT}.\eqref{DT4} implies that $f = f_K$ over $3I$. Therefore, $\Omega_f(3J) = \Omega_{f_K}(3J)$ for any $J \subset I$. By Lemma~\ref{l:Dorronsoro-for-C^2-maps}, we are left to show that 
    \begin{equation} \label{claim:D2FK}
        \abs{\partial_t F_K} \lesssim \varepsilon r_K^{-1}, 
        \qquad 
        |D^2_X F_K| \lesssim \varepsilon r_K^{-1}, 
    \end{equation}
    where $F_K$ is the spatial component of $f_K$, namely $f_k(X, t) = (F_K(X, t), t)$. 
    
    To prove \eqref{claim:D2FK}, recalling that $f_{k+1} = \sigma_k \circ f_k$, we need estimates for derivatives of $\sigma_k$ in $\Sigma_k$ (the image of $f_k$, according to Theorem~\ref{t:DT}). In Lemma~\ref{l:ve-2-distortion} we already saw that, over $\Sigma_k$, 
    \begin{equation} \label{eq:dt_sigma}
        \abs{\partial_t \tilde{\sigma}_k} \leq \frac{C\varepsilon}{r_k},
    \end{equation}
    where $\tilde{\sigma}_k$ is the spatial component of $\sigma_k$. Moreover, in Lemmas~\ref{l:change-in-tangent-super} and \ref{l:second-or-sigma-estimates} we proved that
    \begin{equation} \label{eq:dx_sigma}
        \abs{D_x \tilde{\sigma}_k} \leq 1+C\varepsilon,
        \quad 
        |D^2_x \tilde{\sigma}_k| \leq \frac{C \varepsilon}{r_k}, 
        \quad 
        \text{over }\Sigma_k.
    \end{equation}

    \textbf{Step 3}: iteration. Let us now iterate the estimates on $k$. First, we note that, as a consequence of \eqref{eq:f_m_Lip}, one can show
    \begin{equation} \label{eq:dx_Fk} 
        \abs{D_x F_k} \leq L, 
        \qquad \text{over $\Sigma_k$},
    \end{equation}
    where $L$ is uniform on $k$ (depends on $M$). We now prove by induction that 
    \begin{equation*}
        \abs{\partial_t F_k} \leq \frac{\widetilde{C}}{r_k}, 
        \qquad 
        |D^2_x F_k| \leq \frac{\widetilde{C}}{r_k}, 
        \qquad 
        \text{over $\Sigma_k$},
    \end{equation*}
    where $\widetilde{C}$ is independent of $k$, which will prove the claimed \eqref{claim:D2FK}, hence finishing the proof. 
    
    By differentiating $f_{k+1} = \sigma_k \circ f_k$, we obtain 
    \begin{equation*}
        \partial_t F_{k+1} (X, t)
        =
        \nabla \tilde{\sigma}_k(f_k(X, t)) \, \partial_t F_k(X, t)
        + \partial_t \tilde{\sigma}_k (f_k(X, t), t),
    \end{equation*}
    so that by using the induction hypothesis, \eqref{eq:dx_sigma} and \eqref{eq:dt_sigma} 
    \begin{equation*}
        \abs{\partial_t F_{k+1}}
        \leq 
        \left(1+C\varepsilon\right) \frac{\widetilde{C}}{r_k} 
        + \frac{C\varepsilon}{r_k}
        =
        \frac{1}{r_{k+1}} \left( \frac{(1+C\varepsilon)\widetilde{C} + C\varepsilon}{10} \right)
        \leq 
        \frac{\widetilde{C}}{r_{k+1}}
    \end{equation*}
    simply by making $\varepsilon$ small enough so that $C \varepsilon \leq 1$. A similar computation works for the second spatial derivatives, using \eqref{eq:dx_Fk} this time, too:
    \begin{equation*}
        |D^2_x F_{k+1}|
        \leq 
        |D^2 \tilde{\sigma}_k| \abs{D_x F_k}^2 + \abs{D \tilde{\sigma}_k} |D^2_x F_k| 
        \leq 
        \frac{C}{r_k} L^2 + (1+C\varepsilon) \frac{\widetilde{C}}{r_k}
        \leq 
        \frac{\widetilde{C}}{r_{k+1}},
    \end{equation*}
    by choosing $\widetilde{C}$ large enough depending on $L$ (which ultimately depends on $M$)\footnote{Note that in the fact that $L$ is uniform is relevant to make the choice of $\widetilde{C}$, hence the importance of \eqref{eq:dx_Fk}. If we had estimated $\abs{D_x F_k}$ iterating \eqref{eq:dx_sigma}, the uniform estimate would not have followed.}. This finishes the proof by induction, and therefore the lemma.
\end{proof}

Gathering all the estimates for the beta (or sometimes Omega) numbers that we have developed over the last pages, we have covered all the possible cases, and we can simply wrap up to obtain the estimate \eqref{claim:UR_discrete}, and hence finish the proof of Proposition~\ref{p:f-regular-Lipschitz}.

\begin{proof}[\textbf{Proof of Proposition \ref{p:f-regular-Lipschitz}}]
    Recall that we already explained that Proposition \ref{p:f-regular-Lipschitz} reduces to the estimate \eqref{claim:UR_discrete}. To prove this, fix an arbitrary $I_0 \in \Delta$. First, by Lemmas \ref{l:L} and \ref{l:F} (and \eqref{eq:beta_less_Omega}), we have 
    \begin{align}
        \sum_{I \in \Delta(I_0)} \beta_\Sigma(3I)^2|I| \lesssim |I_0| + \sum_{I \in \mathcal{N}} \beta_\Sigma(3I)^2|I|. 
    \end{align}
    By construction, we have 
    \begin{align}
        \mathcal{N} \subseteq \Biggl(\; \bigcup_{J_0 \in \mathcal{N}_0} \Tree(J_0)\Biggr) 
        \cup 
        \Biggl(\; \bigcup_{\substack{I \in \calM \\ J \subseteq I_0}} \Delta(J) \Biggr)  .
    \end{align}
    Thus, by applying Lemmas \ref{l:near} and \ref{l:far}, we have 
    \begin{align}
        \sum_{I \in \mathcal{N}} \beta_\Sigma(3I)^2|I| &\leq \sum_{J_0 \in \mathcal{N}_0} \sum_{I \in \Tree(J_0)} \beta_\Sigma(3I)^2|I| + \sum_{\substack{J \in \calM \\ J \subseteq I_0}} \sum_{I \in \Delta(J)} \beta_\Sigma(3I)^2|I| \\
        &\lesssim \sum_{J_0 \in \mathcal{N}_0} |J_0| + \sum_{\substack{J \in \calM \\ J \subseteq I_0}} |J| \lesssim |I_0|,
    \end{align}
    where in the final inequality we used that $\mathcal{N}_0$ and $\{J \in \mathcal{M} \colon J \subseteq I_0\}$ are both disjoint collections of cubes contained in $I_0$. Combining the above estimates, \eqref{claim:UR_discrete} follows.
\end{proof}

Now that we have verified that the map $f$ is regular parabolic bi-Lipschitz, it only remains to check the items \eqref{e:contains-floor} -- \eqref{e:f-is-close-modified}. We start by showing \eqref{e:contains-floor}, i.e., that the image of $f$ contains every point where we never stop. 

\begin{lem}
    The inclusion \eqref{e:contains-floor} holds. 
\end{lem}

\begin{proof}
    Fix some $(X,t) \in z(\calS)$. By definition, $(X,t)$ is contained in infinitely many cubes from $\calS$. Recalling the definition of $\calT$, it follows that $(X,t)$ is also contained in infinitely many cubes from $\calT$. Let $\{Q_i\}_{i \in \N} \subseteq \calT$ be a sequences of such cubes such that $\ell(Q_i) \to 0$. By maximality in the definition of $\{(X_{j,k},t_{j,k})\}$, for each $i$ there exists $k(i) \in \N$ and $j(i) \in J_{k(i)}$ such that 
    \begin{align}
        d_p((X_{Q_i},t_{Q_i}),(X_{j(i),k(i))},t_{j(i),k(i)})) &\leq \ell(Q_i) + d_p((Y_{Q_i},s_{Q_i}), (X_{j(i),k(i)},t_{j(i),k(i)})) \\
        &\leq \ell(Q_i) + r_{k(i)} \lesssim \ell(Q_i). 
    \end{align}
Since $(X,t) \in Q_i$ for each $i$, we have that $(X,t)$ is the limit of this sequence of points. In particular, recalling Theorem \ref{t:DT} \eqref{DT1}, we have  $(X,t) \in E_\infty \subseteq f(\R^n)$. 
\end{proof}

We finish this section by defining the good set $G$ and showing \eqref{e:f-is-close}. Moreover, we will show how to obtain \eqref{e:f-is-close-modified} in the case the $\mu$ is Ahlfors regular. Indeed, let $\gamma \in (0,1)$ (to be chosen below sufficiently small depending on $\eta$). Then, define
\begin{align}
    B &\coloneqq \big\{ (X,t) \in Q(S) \colon \dist((X,t),L^Q)\, \mathds{1}_{B_{p}((Z_Q,u_Q) , \rho_Q)}(X,t) \geq \gamma  \ell(Q) \mbox{ for some } Q \in \calT \cap \cubes^{db} \big\};  \\
    G &\coloneqq Q(\calS) \setminus B. 
\end{align}

We first show that the closeness estimate \eqref{e:f-is-close} holds with this definition of $G$. 

\begin{lem}
    The inequality \eqref{e:f-is-close} holds. 
\end{lem}

\begin{proof}
    Fix $Q \in S$ and $(X,t) \in G \cap B_{p}((Z_Q,u_Q) , \rho_Q)$. So long as $\gamma$ and $\ve$ are sufficiently small depending on $\eta$, it suffices to find $k \in \Z$ and $j \in J_k$ such that 
    \begin{align}\label{e:j-k}
        r_k \approx \ell(Q), \quad (X,t) \in 10B_{j,k} \quad \mbox{ and } \quad \dist((X,t),L_{j,k}) \lesssim (\gamma + \ve) \ell(Q). 
    \end{align}
    Indeed, supposing such indices exist, let $(X',t') \in L_{j,k} \cap 11B_{j,k}$ be the point in $L_{j,k}$ closest to $(X,t)$, so that $d_p((X,t),(X',t')) \lesssim (\gamma + \ve) \ell(Q)$. Recalling from Theorem \ref{t:DT} that $\Sigma_k$ is locally a graph over $L_{j,k}$ of the map $A_{j,k}$, whose $L^\infty$ norm is at most $C\ve r_k \approx C \ve \ell(Q)$, there exists a point $(X'',t'') \in \Sigma_k$ such that $d_p((X',t'),(X'',t'')) \lesssim \ve \ell(Q)$. By Theorem \ref{t:DT} \eqref{DT8}, we know that this last point satisfies $\dist_p((X'',t''),\Sigma) \lesssim \ve \ell(Q)$. The inequality \eqref{e:f-is-close} now follows by the triangle inequality after we take $\ve$ and $\gamma$ small enough.

    Our objective now is to find indices satisfying \eqref{e:j-k}. Let $\widehat{Q} \in \calT \cap \cubes^{db}$ be the corresponding cube from Notation \ref{n:T}. Recall that this means 
    \begin{align}\label{e:Q-hat-proof}
        Q \subseteq \widehat{Q} \quad \mbox{ and } \quad \ell(Q) \approx \ell(\widehat{Q}). 
    \end{align}
    It is then easy to check that $(X,t) \in B_{\rho_{\widehat{Q}}}(Z_{\widehat{Q}} , u_{\widehat{Q}})$, which, by the definition of $G$, implies 
    \begin{align}\label{e:dist-L-hat}
        \dist_p((X,t),L^{\widehat{Q}}) \leq \gamma \ell(\widehat{Q}) \approx \gamma \ell(Q). 
    \end{align}
    Let $k \in \Z$ such that $(Y_{\widehat{Q}},s_{\widehat{Q}}) \in \mathscr{C}_k$ (recall Notation \ref{not:CCBP}). It clear from \eqref{e:Q-hat-proof} and the definition of $\mathscr{C}_k$ that 
    \begin{align}
        \ell(Q) \approx \ell(\widehat{Q}) \approx r_k.
    \end{align}
    Thus, we have the first requirement of \eqref{e:j-k}. By maximality, there exists $j \in J_k$ such that 
    \begin{align}
        d_p( (Y_{\widehat{Q}},s_{\widehat{Q}}) , (X_{j,k} , t_{j,k})) \leq r_k. 
    \end{align}
    Since $(X,t) \in B_{\rho_{\widehat{Q}}}(Z_{\widehat{Q}},u_{\widehat{Q}})$ and $(Z_{\widehat{Q}},u_{\widehat{Q}}) , (Y_{\widehat{Q}},s_{\widehat{Q}}) \in \widehat{Q}$, follow from the triangle inequality that comparability between $r_k$ and the size of $\widehat{Q}$ from Notation \ref{not:CCBP} that  
    \begin{align}\label{e:X,t,Xj-tj}
        d_p((X,t),(X_{j,k},t_{j,k})) &\leq d_p((X,t),(Z_{\widehat{Q}},u_{\widehat{Q}})) + d_p((Z_{\widehat{Q}},u_{\widehat{Q}}) , (Y_{\widehat{Q}},u_{\widehat{Q}})) + d_p((Y_{\widehat{Q}},u_{\widehat{Q}}) , (X_{j,k},t_{j,k}))  
        \nonumber
        \\
        &\leq \rho_{\widehat{Q}} + 56r(\widehat{Q}) + r_k \leq 5r_k,
    \end{align}
    i.e., the indices satisfy the second requirement of \eqref{e:j-k}. It now only remains to check that $(X,t)$ is close to the plane $L_{j,k}$.  
    Recall that $Q_{j,k}$ is the cube in $\calT$ such that $(X_{j,k},t_{j,k}) = (Y_{Q_{j,k}} , s_{Q_{j,k}})$ and that 
    \begin{align}\label{e:identification-of-planes}
        L_{j,k} = L^{Q_{j,k}} = L^{\widehat{Q}_{j,k}}.
    \end{align}
    Since $(X_{j,k},y_{j,k}) \in Q_{j,k} \subseteq \widehat{Q}_{j,k}$, it follows from \eqref{e:X,t,Xj-tj} and the triangle inequality that  
    \begin{align}
        d_p((X,t),(X_{\widehat{Q}_{j,k}} , t_{\widehat{Q}_{j,k}})) \leq 5r_k + \ell(\widehat{Q}_{j,k}) \lesssim \ell(Q). 
    \end{align}
    Since both $\widehat{Q}$ and $\widehat{Q}_{j,k}$ are doubling cubes in $\calT$, we have by Lemma \ref{l:angle-control}, \eqref{e:dist-L-hat} and the identity \eqref{e:identification-of-planes} that 
    \begin{align}
        \dist((X,t),L_{j,k}) \lesssim (\gamma + \ve) \ell(Q),
    \end{align}
    which completes the proof of \eqref{e:j-k}. 
\end{proof}

Now we show that the complement of $G$ in $Q(S)$ has small measure. 

\begin{lem}
    For $\ve$ sufficient small, depending on $\alpha$ and $\gamma$, the set $G$ satisfies 
    \begin{align}
        \mu(Q(S) \setminus G) = \mu(B) \leq \alpha \mu(Q(S)). 
    \end{align}
\end{lem}

\begin{proof}
    By the Chebyshev and Cauchy-Schwarz inequalities, we can first estimate 
    \begin{align}
         \mu(B) &\leq \gamma^{-1}\int_{Q(S)}\left( \sum_{Q \in S \cap \cubes^{db}} \left( \frac{\dist((X,t),L^Q)}{\ell(Q)} \right)^2 \mathds{1}_{B_{p}((Z_Q,u_Q) , \rho_Q)}(X,t) \right)^{1/2}\, d\mu \\
        &\leq \gamma^{-1}\left(\int_{Q(S)} \sum_{Q \in S \cap \cubes^{db}} \left( \frac{\dist((X,t),L^Q)}{\ell(Q)} \right)^2 \mathds{1}_{B_{p}((Z_Q,u_Q) , \rho_Q))}(X,t) \, d\mu \right)^{1/2 }\mu(Q(S))^{1/2} \\
        &\leq  \gamma^{-1}\left(\sum_{Q \in S \cap \cubes^{db}} \int_{B_{p}((Z_Q,u_Q) , \rho_Q)} \left( \frac{\dist((X,t),L^Q)}{\ell(Q)} \right)^2 \, d\mu \right)^{1/2} \mu(Q(S))^{1/2}. 
    \end{align}
    Recall from \eqref{e:volume-estimates-doubling} that $\mu(Q) \approx \Theta_\mu(2B_{Q(S)}) \rho_Q^{n+1}$. Combining this estimate with \eqref{e:new-base-point-and-plane-2}, by carrying on from the above we obtain 
    \begin{align}
        \mu(B) &\lesssim  \gamma^{-1}\left(\sum_{Q \in S \cap \cubes^{db}}\Theta_\mu(2B_{Q(S)})^{-1} \beta_\mu((Z_Q,u_Q),\rho_Q)^2 \mu(Q) \right)^{1/2} \mu(Q(S))^{1/2} \\
        &= \gamma^{-1} \left(\sum_{Q \in S \cap \cubes^{db}} \int_{Q(S)} \Theta_\mu(2B_{Q(S)})^{-1} \beta_\mu((Z_Q,u_Q),\rho_Q)^2 \mathds{1}_{Q}(X,t) \, d\mu  \right)^{1/2} \mu(Q(S))^{1/2}.
    \end{align}
    By the summability hypothesis \eqref{b:summability} of Proposition \ref{p:construct-regular-images}, this last expression is at most $\gamma^{-1}\ve \mu(Q(S))$. Thus, the lemma follows so long as $\ve = \ve(\alpha , \gamma)$ is sufficiently small. 
\end{proof}

\begin{lem}
    If $\mu$ is Ahlfors regular and $\ve$ is sufficiently small, depending on $\eta$ and the regularity constant, we have \eqref{e:f-is-close-modified}.  
\end{lem}

\begin{proof}
    So long as $\ve$ is chosen small enough with the dependencies indicated in the statement, this follows directly from the inequality 
    \begin{align}
        r(B)^{-1}\inf_L\sup_{(X,t) \in \supp(\mu) \cap B} \dist(X,t,L) \lesssim \beta_\mu(2B)^\frac{1}{n+2}
    \end{align}
    for any ball $B$ centred on $\supp(\mu)$. Here, the infimum is taken over all vertical planes $L$. For a proof of this inequality we direct the reader to \cite[Lemma 2.12]{azzam2018analyst} where it is proved in the Euclidean metric. The proof in this setting is exactly the same up to notation (choosing $d=n+1$ and restricting to vertical planes).    
\end{proof}

%%%%%%%%%%%%%%%%%%%%%%%%%%%%%%%%%%%%%%%%%%%%%%%%%%%%%%%%%%%%
%%%%%%%%%%%%%%%%%%%%%%%%%%%%%%%%%%%%%%%%%%%%%%%%%%%%%%%%%%%%
%%%%%%%%%%%%%%%%%%%%%%%%%%%%%%%%%%%%%%%%%%%%%%%%%%%%%%%%%%%%
%%%%%%%%%%%%%%%%%%%%%%%%%%%%%%%%%%%%%%%%%%%%%%%%%%%%%%%%%%%%
%%%%%%%%%%%%%%%%%%%%%%%%%%%%%%%%%%%%%%%%%%%%%%%%%%%%%%%%%%%%
%%%%%%%%%%%%%%%%%%%%%%%%%%%%%%%%%%%%%%%%%%%%%%%%%%%%%%%%%%%%

\newpage

\section{The proof that finite Jones function implies rectifiability}\label{sec:Jones-implies-rectifiability}

With the bi-Lipschitz Reifenberg theorem and the construction of regular parabolic Lipschitz for stopping-times with small Jones function out of the way, we can now more directly tackle the proof of the main result, Theorem \ref{thm:main}. In this section we deal with the implication \textit{finite square function of Jones' beta numbers implies rectifiability}. We follow the outline of \cite[Section 5]{azzam2015characterization} using, as the key tool, the construction of regular parabolic bi-Lipschitz maps from Proposition \ref{p:construct-regular-images}. \\

By standard arguments, it suffices to show that any compact subset $F$ of $E$ with positive $\mu$-measure which satisfies 
\begin{align}\label{e:beta-bounds}
    \int_F \int_0^1 \beta_\mu((X,t),r)^2 \, \frac{dr}{r} d\mu  < \infty
\end{align}
and (as a consequence of our assumption that densities are finite $\mu$-a.e.)
\begin{align}\label{e:density-bounds}
    \sup_{r > 0} \frac{\mu(B_p((X,t),r)}{r^{n+1}} < \infty, \quad \mbox{ for all } \quad (X,t) \in F, 
\end{align}
has a rectifiable subset of positive $\mu$-measure. We fix an $F$ satisfying the above conditions for the remainder of the section. \\

Let $\kappa \in (0,1)$, to be chosen sufficiently small below. Recalling our assumption on the densities, let $(X_0,t_0)$ be a density point of $F$, and then choose a radius $r_0 > 0$ such that the ball $B_0 = B_{p}((X_0,t_0),r_0)$ is almost completely contained in $F$ and, at the same time, satisfies a doubling estimate in the form
\begin{align}\label{e:B_0}
    \mu(B_0 \setminus F) \leq \kappa^2 \mu(B_0) \quad \mbox{ and } \quad \mu\left(\frac{1}{2} B_0 \right) \gtrsim \mu(B_0).  
\end{align}
Note that, for $\kappa$ small enough, \eqref{e:B_0} implies 
\begin{align}\label{e:F-large-in-half-ball}
    \mu\left(F \cap \frac{1}{2} B_0 \right) \gtrsim \mu(B_0). 
\end{align}
This means that, since we are only interested in obtaining a rectifiable subset of $F$ with positive $\mu$-measure, it suffices to find this subset within $B_0$. \\

Let $\cubes(B_0) \subseteq \cubes$ denote the set of cubes contained in $B_0$ and let  $\mathcal{B}$ denote the set of cubes $R \in \cubes(B_0)$ which contain a much smaller portion of $F$ than $B_0$ did, in the sense that 
\begin{align}
    \mu(R \setminus F) \geq \kappa \mu(R). 
\end{align}
By the same argument as in \cite[Lemma 5.1]{azzam2015characterization}, it quickly follows from \eqref{e:density-bounds} that these bad cubes $\mathcal{B}$ occupy only a small fraction of $B_0$: 
\begin{align}\label{e:small-mass-bad}
    \mu\left( \bigcup_{R \in \calB} R \right) \lesssim \kappa \mu(B_0). 
\end{align}

Our goal next will be to construct Reifenberg maps (as in Proposition~\ref{p:construct-regular-images}) that cover a portion of the support of $\mu$ within $B_0$. This will quickly imply rectifiability of that part of $\supp \mu$ because of all the properties developed in the previous sections. To be able to apply Proposition~\ref{p:construct-regular-images}, we first need to select an initial cube which is doubling, which we shall do next. It will also be useful to choose this cube small enough, so we also fix $\delta \in (0,1)$, to be chosen sufficiently small below. 

\begin{lem}
    There exists $R_0 \in \cubes^{db} \cap \cubes(B_0)$ such that 
\begin{align}\label{e:R_0}
    \ell(R_0) \leq \delta, \qquad \delta^{-1}B_{R_0} \subseteq \frac{9}{10}B_0 \qquad \mbox{and} \qquad \mu\left((F \cap R_0) \setminus \bigcup_{R \in \calB} R \right) > 0 . 
\end{align}
\end{lem}

\begin{proof}
    Let $\calM$ be a maximal collection of cubes in $\cubes^{db} \cap \cubes(B_0)$ which have non-empty intersection with $\tfrac{1}{2}B_0$ and which have side-length at most $\delta/10$. By Lemma \ref{lem:cover_doubling}, we have 
    \begin{align}\label{e:M-covers}
        \mu\left( \frac{1}{2} B_0 \setminus \bigcup_{T \in \calM} T \right) = 0.
    \end{align}
    Moreover, it is clear from the definition and the triangle inequality that each cube in $\calM$ satisfies the first two properties in \eqref{e:R_0}. We now only need to find some cube in $R_0 \in \calM$ which satisfies the right-most property. Suppose, towards a contradiction, that no such cube exists. Using this assumption, along with \eqref{e:small-mass-bad} and \eqref{e:M-covers}, one observes that 
    \begin{align}
        \mu\left(F \cap \frac{1}{2}B_0\right) 
        & = 
        \sum_{T \in \calM} \mu \left( F \cap T \right)
        =
        \sum_{T \in \calM}\left[\mu\left( F \cap T \cap \bigcup_{R \in \calB} R \right) +  \mu\left( (F \cap T)\setminus \bigcup_{R \in \calB} R \right)  \right] 
        \\
        & = \sum_{T \in \calM} \mu\left( F \cap T \cap \bigcup_{R \in \calB} R \right)
        \leq \mu\left( \bigcup_{R \in \calB} R \right) \lesssim \kappa \mu(B_0).  
    \end{align}
    For $\kappa$ small enough, this directly contradicts \eqref{e:F-large-in-half-ball}. 
\end{proof}

The idea now is to apply Proposition \ref{p:construct-regular-images} iteratively, constructing a sequence $\{\calS_i\}_{i \in \N}$ of stopping-time regions and a set $F' \subseteq F \cap R_0$ satisfying
\begin{align}\label{e:F'}
    \mu(F') > 0 \quad \mbox{ and } \quad  F' \subseteq \bigcup_{i \in \N} z(\calS_i),
\end{align}
where $z(\calS_i)$ is the set of points in $Q(\calS_i)$ where we never stop, namely
\begin{align}
    z(\calS_i) \coloneqq Q(\calS_i) \setminus \bigcup_{Q \in \min(\calS_i)} Q.
\end{align} 
That this is sufficient to conclude the required implication of Theorem \ref{thm:main} follows from the next lemma. 

\begin{lem}
    Let $\calS$ be a stopping-time region which satisfies the hypotheses of Proposition \ref{p:construct-regular-images} and let $\Sigma = f(\R^n)$ denote the image of the associated regular bi-Lipschitz mapping. Then, 
    \begin{align}\label{e:ac}
         \mu\lfloor{z(\calS)} \ll \calH_p^{n+1}\lfloor{\Sigma}
    \end{align}
    In particular, $\mu\lfloor{z(\calS)}$ is parabolic rectifiable. 
\end{lem}

\begin{proof}
    The conclusion that $\mu\lfloor{z(\calS)}{z(\calS)}$ is parabolic rectifiable follows directly from condition \eqref{e:ac} and the fact that $\Sigma$ is parabolic uniformly rectifiable, as shown in Proposition~\ref{p:construct-regular-images} (see also Remark~\ref{rem:UR_with_betas}). Therefore, we are left to prove \eqref{e:ac}. For that purpose, it clearly suffices to show that 
    \begin{align}\label{e:density}
        \limsup_{r \to 0} \frac{\mu(B_p(X,t,r))}{\calH_p^{n+1}(\Sigma \cap B_p(X,t,r))} \lesssim 1, 
        \qquad 
        \text{uniformly for all $(X,t) \in z(\calS)$.}
    \end{align}
    
    Fix some $(X,t) \in z(\calS)$. Since, by definition of the contact set $z(\calS)$, $(X,t)$ is contained in infinitely many cubes from $\calS$ (with radii tending to 0), for each $r > 0$ small enough, there exists $Q(r) \in \calS$ such that 
    \begin{align}
        B_p(X,t,r) \subseteq 2B_{Q(r)} \quad \mbox{ and } \quad r \approx \ell(Q(r)).
    \end{align} 
    Since
    $(X,t) \in \Sigma$ (by \eqref{e:contains-floor} applied along the sequence of cubes in $\calS$ containing $(X, t)$ and radii converging to 0) and $\Sigma$ is Ahlfors regular (because $f$ from Proposition~\ref{p:construct-regular-images} is bi-Lipschitz), we have
    \begin{align}
        \calH_p^{n+1}(\Sigma \cap B_p(X,t,r)) \approx r^{n+1} \approx \ell(Q(r))^{n+1}
    \end{align}
    uniformly for each such $r$. Using the above properties, along with the density conditions in the hypotheses of Proposition \ref{p:construct-regular-images}, we conclude that 
    \begin{align}
        \frac{\mu(B_p(X,t,r))}{\calH_p^{n+1}(\Sigma \cap B_p(X,t,r))} \lesssim \frac{\mu(2B_{Q(r)})}{\ell(Q(r))^{n+1}} \lesssim_A \frac{\mu(2B_{Q(\calS)})}{\ell(Q(\calS))^{n+1}}
    \end{align}
    with all implicit constant independent of both $(X,t)$ and $r$. This implies \eqref{e:density} immediately. 
\end{proof}

In light of the previous discussion, our goal now is to be able to apply iteratively Proposition~\ref{p:construct-regular-images} below the chosen $R_0$, until we find appropriate contact regions $z(\calS_i)$. For that, first of all, we need some work to find the appropriate stopping-time regions which satisfy the hypotheses of Proposition~\ref{p:construct-regular-images}.

Following the construction of \cite[Sections 6.1 \& 6.3]{azzam2015characterization} verbatim, for every small enough doubling cube $R$ which contains a large proportion of $F$, it is possible to construct a stopping-time region $\mathcal{S}_R'$, with explicit control over the volume over a sub-class of minimal cubes and such that, after removing the set of minimal cubes, the resulting stopping-time region $\mathcal{S}_R$ satisfies the hypotheses of Proposition \ref{p:construct-regular-images}. More concretely, the following result holds.

\begin{lem}\label{l:stopping-time-with-set}
    Let $A > 1$ and $\tau \in (0,1)$ such that $\tau \ll A^{-1}$. Let $R \in \cubes^{db}$ such that $\ell(R) \leq \delta$ and $\mu(R \setminus F) \leq \kappa$. Then, for $\delta,\kappa$ sufficiently small (depending only on $\tau$ and $A$), there exists a constant $C$, depending only on $A$ and $\tau$, and stopping-time regions, which we denote by $\calS_R \subseteq \mathcal{S}_R'$, such that $\calS_R$ satisfies the hypotheses of Proposition \ref{p:construct-regular-images} and such that
    \begin{align}\label{e:stop-control}
        \sum_{Q \in \min(\calS_R') \setminus (\mathrm{HD}(R) \cup \mathrm{UB}(R))} \mu(Q) 
        & \leq \tau^{1/2} \mu(R) \\
        & \quad + \frac{C}{\Theta_\mu(2B_R)} \sum_{Q \in \calS_R'} \int_{F \cap \delta^{-1}B_Q}  \int_{\delta \ell(Q)}^{\delta^{-1}\ell(Q)} \beta_\mu((X,t),r) \, \frac{dr}{r} d\mu .
    \end{align} 
    Here, $\mathrm{HD}(R)$ are certain high density minimal cubes, and $\mathrm{UB}(R)$ are certain non-balanced minimal cubes. The stopping-time region $\mathcal{S}_R$ is obtained from $\mathcal{S}_R'$ by removing the minimal cubes.   
\end{lem}

The precise definitions of $\mathrm{HD}(R)$ and $\mathrm{UB}(R)$ are not needed in the current exposition. The important point is that, by slightly modifying the stopping-time procedure, one can use the estimate \eqref{e:stop-control} to construct a sequence of trees whose top cubes satisfy a Carleson condition on the corresponding densities. Indeed, the following lemma is obtained by directly applying the abstract stopping-time procedure of \cite[Sections 5.2 - 5.5]{azzam2015characterization}. One starts at the top cube $R_0$ and stops if they encounter a cube in $\mathrm{HD}^*(R_0)$, $\mathrm{UB}^*(R_0)$ or $\mathrm{O}^*(R_0)$, the first two sets being a modification of the sets $\mathrm{HD}(R_0)$, $\mathrm{UB}(R_0)$ from above, and the latter set (where $\mathrm{O}$ stand for other), being modification of the set of cubes in \eqref{e:stop-control}). Below these stopped cubes, one looks for a maximal family of cubes $R$ satisfying the hypotheses of Lemma \ref{l:stopping-time-with-set} and then restarts. The family of stopping-times is obtained by iterating this procedure. Stopping because one encounters an $\mathrm{O}^*(R)$ cubes is controlled by \eqref{e:stop-control}. Stopping because one encounters $\mathrm{HD}(R_0)$ or $\mathrm{UB}(R_0)$ cubes is controlled by arguing that the density increases.

\begin{lem}\label{l:sequence-of-stopping-times}
    There is a sequence of stopping-time regions $\{\calT_i\}_{i \in \N}$ such that the following conditions hold. 
    \begin{enumerate}
        \item \label{c:sequence-of-stopping-times-1} For each $i \in \N$, we have $Q(\calT_i) \in \cubes^{db}(R_0)$.
        
        \item \label{c:sequence-of-stopping-times-2} For each $i \in \N$, if $\calS_i \coloneqq \calS_{Q(\calT_i)}$ is the stopping-time region constructed in Lemma \ref{l:stopping-time-with-set}, then
        \begin{align}
            z(\calT_i) \subseteq z(\calS_i). 
        \end{align}
        
        \item \label{c:sequence-of-stopping-times-2.5} The stopping-time regions $\{\calT_i\}_{i \in \N}$ can be partitioned into layers as 
        \[\{\calT_i\}_{i \in \N} = \bigcup_{N \in \N} \mathrm{Layer}(N)\] 
        such that the following covering property holds: For each $N_0 \in \N$, up to a set of $\mu$-measure zero, we have the inclusion 
        \begin{align}
            R_0 \setminus \bigcup_{R \in \calB} R 
            \; \subseteq \;
            \bigcup_{\calT_i \in \mathrm{Layer}(N_0)} Q(\calT_i) \, \cup \,  \bigcup_{N=1}^{N_0-1} z(\calT_i) . 
        \end{align}
        
        \item \label{c:sequence-of-stopping-times-3} We have control over the top cubes by the estimate 
        \begin{align}
            \sum_{i \in \N} \Theta_\mu(2B_{Q(\calT_i)}) \, \mu(Q(\calT_i)) 
            \; \lesssim \; 
            \mu(R_0) + \int_F \int_0^1 \beta_\mu((X,t)S,r) \, \frac{dr}{r} d\mu . 
        \end{align}
    \end{enumerate}
\end{lem}

We now have all the ingredients we need to define a set $F'$ which satisfies the conditions \eqref{e:F'}. Recall that we already discussed constructing such a set is sufficient to prove Theorem \ref{thm:main}.

First, by Lemma \ref{l:sequence-of-stopping-times} \eqref{c:sequence-of-stopping-times-2.5}, let $A_1$ be a full $\mu$-measure subset of $R_0 \setminus \bigcup_{R \in \calB} R$ such that 
\begin{align}
    A_1 \subseteq \bigcup_{\calT_i \in \mathrm{Layer}(N_0)} Q(\calT_i) \, \cup \, \bigcup_{N=1}^{N_0-1} z(\calT_i), 
    \qquad 
    \text{for all $N_0 \in \N$.}
\end{align}
Then, let 
\begin{align}
    A_2 \coloneqq \left\{(X,t) \in A_1 \colon \; \Theta^{n+1,*}(X,t,\mu) > 0 \right\},
\end{align}
so that, recalling the hypotheses of Theorem \ref{thm:main}, $A_2$ is a set of full $\mu$-measure in $A_1$. 

Next, by our initial assumption \eqref{e:beta-bounds}, the condition Lemma \ref{l:sequence-of-stopping-times} \eqref{c:sequence-of-stopping-times-3} implies that 
\begin{align}
    \sum_{i \in \N} \Theta_\mu(2B_{Q(\calT_i)}) \mathds{1}_{Q(\calT_i)}(X,t) < \infty \quad \mbox{ for } \mu\mbox{-a.e. } (X,t) \in R_0. 
\end{align}
As a consequence, by following directly the short are argument at the start of \cite[Section 5.6]{azzam2015characterization}, we deduce that $\mu$-a.e. point $(X,t) \in A_2$ is contained in only finitely many cubes from $\{Q(\calT_i)\}_{i \in \N}$. Namely, if we denote by $A_3$ this set of points
\begin{align}
    A_3 \coloneqq \left\{(X,t) \in A_2 \colon \sum_{i \in \N} \mathds{1}_{Q(\calT_i)}(X,t) < \infty \right\}, 
\end{align}
it has full $\mu$-measure in $A_2$. 

The fact that any point in $A_3$ belongs to only finitely many stopping-time regimes, intuitively means that it is a contact set for the smallest one of this regimes. More concretely, it follows from Lemma \ref{l:sequence-of-stopping-times} \eqref{c:sequence-of-stopping-times-2} and \eqref{c:sequence-of-stopping-times-2.5} that 
\begin{align}\label{e:A_2-a-subset}
    A_3 \; \subseteq \; \bigcup_{N=1}^\infty z(\calT_i) \; \subseteq \; \bigcup_{N=1}^\infty z(\calS_i). 
\end{align}
Indeed, right-most inclusion is immediate from Lemma \ref{l:sequence-of-stopping-times} \eqref{c:sequence-of-stopping-times-2}, so we need only prove the left-most inclusion. Consider some $(X,t) \in A_3$, and let $Q(\calT_i)$ be the smallest cube such that $(X,t) \in Q(\calT_i)$ (which exists by definition of $A_3$). Let also $N_0 \in \N$ be such that $\calT_i \in \text{Layer}(N_0)$. Since $(X,t) \in A_1$,
\begin{align}
    (X,t) \in \bigcup_{\calT_j \in \mathrm{Layer}(N_0+1)} Q(\calT_j) \cup \bigcup_{N=1}^{N_0} z(\calT_j). 
\end{align}
However, by minimality, we know that $(X,t) \not\in Q(\calT_j)$ for any $\calT_j \in \mathrm{Layer}(N_0+1)$. This finishes the proof of \eqref{e:A_2-a-subset}. 

We finally define 
\begin{align}
    F' \coloneqq F \cap A_3. 
\end{align}
It is clear from the construction that $F'$ is a set of full $\mu$-measure in 
\begin{align}
    (F \cap R_0) \setminus \bigcup_{R \in \calB} R, 
\end{align}
so, recalling \eqref{e:R_0}, it holds $\mu(F') > 0$, which is the first part of \eqref{e:F'}. The second part of \eqref{e:F'}, namely that $F' \subseteq \bigcup_{i \in \N} z(\calS_i)$, follows directly from \eqref{e:A_2-a-subset}. This completes the verification of \eqref{e:F'}, and hence the proof of this implication of Theorem~\ref{thm:main} by the previous discussion.

%%%%%%%%%%%%%%%%%%%%%%%%%%%%%%%%%%%%%%%%%%%%%%%%%%%%%%%%%%%%
%%%%%%%%%%%%%%%%%%%%%%%%%%%%%%%%%%%%%%%%%%%%%%%%%%%%%%%%%%%%
%%%%%%%%%%%%%%%%%%%%%%%%%%%%%%%%%%%%%%%%%%%%%%%%%%%%%%%%%%%%
%%%%%%%%%%%%%%%%%%%%%%%%%%%%%%%%%%%%%%%%%%%%%%%%%%%%%%%%%%%%
%%%%%%%%%%%%%%%%%%%%%%%%%%%%%%%%%%%%%%%%%%%%%%%%%%%%%%%%%%%%
%%%%%%%%%%%%%%%%%%%%%%%%%%%%%%%%%%%%%%%%%%%%%%%%%%%%%%%%%%%%
%%%%%%%%%%%%%%%%%%%%%%%%%%%%%%%%%%%%%%%%%%%%%%%%%%%%%%%%%%%%

\newpage

\section{Preliminaries for rectifiability implies finite Jones function} \label{sec:preliminaries_2}

We now turn our attention to the other implication of Theorem \ref{thm:main}, namely, \textit{rectifiability implies finiteness of the Jones square function}. Let us first give some further notation before diving into the proofs.

\begin{notation}
    In what follows, we will sometimes use notation that somehow differs from that from the previous sections. We summarize the differences, and also recall some relevant facts that will be useful later.

    \begin{itemize}
        \item (Boldface space-time points) We will often write points in space-time $\R^n \times \R$ as boldface Roman letters, e.g., $\mathbf{X}$. For a cube $Q$, we will write the centre point as 
        \begin{align}
            \mathbf{X}_Q \coloneqq (X_Q,t_Q). 
        \end{align}

        \item (Parabolic balls) In this section, Section~\ref{proof-of-the-main-theorem} and the appendix, we call the cylinders $B_{p}({\X},r)$ of
        Section~\ref{sec:setting} \textit{balls}. Note that the balls $B_{p}({\X},r)$ are not the balls associated to the distance $d_{p}$. We will use the following two facts many times. First,
        \begin{equation}\label{eq:ball-distance}
            d_{p}({\X},{\Y})\le r\quad\Longrightarrow\quad {\Y}\in B_{p}({\X},r)
            \quad\Longrightarrow\quad d_{p}({\X},{\Y})\le2r.
        \end{equation}
        Second, if ${\Y}\in B_{p}({\X},\rho)$, then
        \begin{equation}\label{eq:ball-inclusion}
            B_{p}({\Y},r)\subset B_{p}({\X},r+\rho).
        \end{equation}
        Indeed, the spatial distances add up, and $r^2+\rho^2\le(r+\rho)^2$.

        \item (Vertical planes) Vertical planes are the sets
        \[
            V_{\xi,c}\coloneqq \{(X,t)\in\mathbb R^{n+1}:\xi\cdot X=c\},
            \qquad\xi\in S^{n-1},\quad c\in\mathbb R.
        \]
        For a point ${\X}=(X,t)$, the closest point of $V_{\xi,c}$ has the same
        time coordinate, and so
        \[
            d_{p}({\X},V_{\xi,c})=|\xi\cdot X-c|.
        \]
        In particular, if a vertical plane $V$ has non-empty intersection with $B_{p}({\X},r)$, then
        $d_{p}({\X},V)\le r$.
        \\
        For a vertical plane $V$, we will write 
        \begin{align}
            \lambda_V := \mathcal H_{p}^{n+1}\lfloor V.
        \end{align}

        \item (Lipschitz graphs) After a rotation of the spatial variables, a Lip(1,1/2) graph is a set of the form
        \[
        \Gamma=\{(X,A(X,t),t):(X,t)\in\mathbb R^{n-1}\times\mathbb R\},
        \]
        where $A:\mathbb R^{n-1}\times\mathbb R\to\mathbb R$ is parabolic
        Lipschitz (Definition~\ref{regfngraph.def}). We call $\Gamma$ regular
        when $A$ is regular (Definition~\ref{def:regular}).

        \item For a finite measure $\nu$, we write $\|\nu\|$ for its total variation.  Also, $\supp f$ is the closed support of a function $f$.
    \end{itemize}
\end{notation}

%%%%%%%%%%%%%%%%%%%%%%%%%%%%%%%%%%%%%%%%%%%%%%%%%%%%%%%%%%%%%%%%
%%%%%%%%%%%%%%%%%%%%%%%%%%%%%%%%%%%%%%%%%%%%%%%%%%%%%%%%%%%%%%%%
%%%%%%%%%%%%%%%%%%%%%%%%%%%%%%%%%%%%%%%%%%%%%%%%%%%%%%%%%%%%%%%%
%%%%%%%%%%%%%%%%%%%%%%%%%%%%%%%%%%%%%%%%%%%%%%%%%%%%%%%%%%%%%%%%
%%%%%%%%%%%%%%%%%%%%%%%%%%%%%%%%%%%%%%%%%%%%%%%%%%%%%%%%%%%%%%%%
%%%%%%%%%%%%%%%%%%%%%%%%%%%%%%%%%%%%%%%%%%%%%%%%%%%%%%%%%%%%%%%%

\newpage

\section{The main lemma for rectifiability implies finite Jones function}
\label{sec:graph}

In this section we establish the main lemma for this direction of Theorem~\ref{thm:main}: when restricted to the family of regular parabolic Lipschitz graphs, the conclusion of Theorem~\ref{thm:main} is true with the alpha numbers, namely, the Jones function for any finite measure $\mu$, defined in terms of $\alpha$-numbers (Definition \ref{d:alpha-numbers}), is finite almost everywhere (with respect to surface measure).

\begin{thm}[Main Lemma: summability of alpha numbers over regular Lip(1,1/2) graphs]\label{thm:graph-alpha}
Let $\Gamma$ be a regular Lip(1,1/2) graph,
and define its surface measure by $\sigma := \mathcal H_{p}^{n+1}\lfloor\Gamma$. Let $\mu$ be a finite Borel measure on $\mathbb R^{n+1}$. Then
\[
\int_0^1\alpha_\mu(\mathbf{{\X}},r)^2\,\frac{dr}{r}<\infty
\qquad\text{for }\sigma\text{-a.e. }{\X}\in\Gamma.
\]
\end{thm}

We will introduce alpha numbers soon in the next subsection. 

Before proving this main lemma, let us draw the following simple corollary which will be used in Section~\ref{proof-of-the-main-theorem}, which just ammounts to using correctly the definition of parabolic rectifiability.

\begin{cor}[Extension to parabolic rectifiable]  \label{cor:graph-alpha}
Let $\mu$ be a finite Borel measure on $\mathbb R^{n+1}$ which is parabolic rectifiable. Then
\begin{align}\label{e:dini-alpha-corollary}
\int_0^1\alpha_\mu({\X},r)^2\,\frac{dr}{r}<\infty
\qquad\text{for }\mu\text{-a.e. }{\X}\in\mathbb R^{n+1}.
\end{align}
\end{cor}

\begin{proof}
Since $\mu$ is parabolic rectifiable, its support can be covered up to a set of $\mu$-measure zero by a countable collection of regular Lip$(1,1/2)$ graph $\{\Gamma_j\}_j$. Fix $j$ and put $\sigma_j=\mathcal H_{p}^{n+1}\lfloor\Gamma_j$. By
Theorem~\ref{thm:graph-alpha}, the integral in \eqref{e:dini-alpha-corollary} is finite for every $\X \in \Gamma_j \setminus N_j$, for some
$N_j\subset\Gamma_j$ with $\sigma_j(N_j)=0$. Since $\mu \ll \mathcal H_{p}^{n+1}$ and 
$\mathcal H_{p}^{n+1}(N_j)=0$, it follows that $\mu(N_j)=0$. Since the graphs
$\Gamma_j$ cover $\mu$-almost all of $\mathbb R^{n+1}$, the corollary
follows.
\end{proof}

We now begin the proof of Theorem \ref{thm:graph-alpha} in earnest. In the rest of this section, we assume that $\mu$ and $\Gamma$ are as
in Theorem~\ref{thm:graph-alpha}. We emphasize that the measure $\mu$ does not need to be
supported on $\Gamma$.

\subsection{The graph measure}
\label{the-graph-measure}
After a spatial rotation, we can assume that $\Gamma$ is the graph of a regular Lipschitz function $A$ over the plane
\[
\mathbb P := \mathbb R^{n-1}\times\mathbb R, 
\]
which we will call the \textit{base}. This amounts to saying that the height of the graph gives the $n$-th spatial coordinate.

\begin{notation}[Restricted metric in $\Pee$]
    We endow $\Pee$ with the natural metric induced by $d_p$, namely $\rho$ defined by 
    \begin{align}
        \rho((u,t),(v,\tau))
        :=|u-v|+|t-\tau|^{1/2}.
    \end{align}
    We also define cylinders, for $z=(u,t)\in\mathbb P$ and $r>0$, by
    \[
        B_{\mathbb P}(z,r)
        :=
        \big\{(v,\tau)\in\mathbb P:|u-v|\le r,\ |t-\tau|\le r^2\big\},
    \]
    and naturally \eqref{eq:ball-distance} and \eqref{eq:ball-inclusion} hold restricted to $\Pee$.
\end{notation}

\begin{notation}[Graph over $\Pee$]
    We define a map $\Phi$ on $\Pee$ by 
    \[
    \Phi(u,t)\coloneqq (u,A(u,t),t),
    \]
    so by definition, 
    \[ \Gamma=\Phi(\mathbb P).\]
    Moreover, if $L$ is the parabolic Lipschitz
    constant of $A$, then
    \begin{align}\label{e:Phi-bi-Lipschitz}
    \rho(z,z')
    \le d_{p}(\Phi(z),\Phi(z'))
    \le C_L\rho(z,z').
    \end{align}
\end{notation}

\begin{notation}[Measures over $\Pee$]
    We also let 
    \begin{equation*}
        \text{$m:=du\,dt$ \quad be Lebesgue measure on $\mathbb P$.}
    \end{equation*} 
    This way, it holds
    $m(B_{\mathbb P}(z,r))\approx r^{n-1}r^2=r^{n+1}$,
    and $(\mathbb P,\rho)$ has homogeneous dimension $n+1$. By the bi-Lipschitzianity from \eqref{e:Phi-bi-Lipschitz}, it holds 
    \begin{equation}\label{eq:graph-measure-comparison}
    C_L^{-1}m(H)\le \sigma(\Phi(H))\le C_Lm(H)
    \qquad\text{for Borel }H\subset\mathbb P.
    \end{equation}
    so by the Radon--Nikodym theorem,
    \begin{equation}
        \sigma=\Phi_\#(w\,m)
    \label{eq:beta_alpha:2-3}
    \end{equation}
    for a measurable density $w$ satisfying
    $0<c_L\le w\le C_L<\infty$ $m$-a.e. Note, in particular, that this implies that $\Gamma$ is parabolic Ahlfors--David regular
    (Definition~\ref{def:UR}), or AD-regular for short.
\end{notation} 

\begin{figure}[h!]
    \centering
    \begin{tikzpicture}[scale=0.78,>=Stealth]
      \def\A(#1){4.1+0.35*sin(1.1*(#1) r)+0.06*(#1)}
      \draw[->,black] (-0.3,0) -- (9.4,0) node[below] {$u$};
      \draw[->,black] (-0.3,0) -- (-0.3,7.9) node[left] {$x_n$};
      \pgfmathsetmacro{\yQ}{\A(4.1)}
      \draw[dashed] (4.1,\yQ) circle (3.1);
      \node[font=\scriptsize] at (6.6,6.85) {$B_Q$};
      \fill[blue!20] (4.1,\yQ) circle (0.25);
      \draw (4.1,\yQ) circle (0.25);
      \draw[thick] plot coordinates {
    (0.00000,4.10000) (0.22500,4.19924) (0.45000,4.29326) (0.67500,4.37715) (0.90000,4.44661) (1.12500,4.49824) (1.35000,4.52971) (1.57500,4.53993) (1.80000,4.52910) (2.02500,4.49870) (2.25000,4.45141) (2.47500,4.39093) (2.70000,4.32176) (2.92500,4.24896) (3.15000,4.17777) (3.37500,4.11336) (3.60000,4.06048) (3.82500,4.02317) (4.05000,4.00453) (4.27500,4.00652) (4.50000,4.02983) (4.72500,4.07388) (4.95000,4.13680) (5.17500,4.21557) (5.40000,4.30623) (5.62500,4.40406) (5.85000,4.50393) (6.07500,4.60058) (6.30000,4.68893) (6.52500,4.76442) (6.75000,4.82329) (6.97500,4.86275) (7.20000,4.88124) (7.42500,4.87844) (7.65000,4.85535) (7.87500,4.81420) (8.10000,4.75832) (8.32500,4.69194) (8.55000,4.61992) (8.77500,4.54748) (9.00000,4.47986)
    };
      \node[font=\scriptsize] at (8.6,{\A(9)+0.35}) {$\Gamma$};
      \draw[line width=2.2pt,blue!70!black] (3.6,0) -- (4.6,0);
      \node[below,font=\scriptsize] at (4.1,0) {$Q_0$};
      \draw[dotted] (3.6,0) -- (3.6,{\A(3.6)});
      \draw[dotted] (4.6,0) -- (4.6,{\A(4.6)});
      \draw[line width=2.2pt,blue!70!black] plot coordinates {
    (3.60000,4.06048) (3.68333,4.04466) (3.76667,4.03116) (3.85000,4.02013) (3.93333,4.01172) (4.01667,4.00602) (4.10000,4.00314) (4.18333,4.00314) (4.26667,4.00606) (4.35000,4.01192) (4.43333,4.02071) (4.51667,4.03240) (4.60000,4.04693)
    };
      \fill (4.1,\yQ) circle (1.4pt);
      \node[font=\scriptsize,anchor=south] at (4.1,{\yQ+0.55}) {${\X}_Q=\Phi(z_Q)$};
      \node[font=\scriptsize,blue!70!black,anchor=north west] at (4.6,{\A(4.6)-0.05}) {$Q=\Phi(Q_0)$};
      \fill (4.1,0) circle (1.4pt) node[above left=-1pt,font=\scriptsize] {$z_Q$};
      \draw[->,blue!70!black] (2.4,0.5) to[bend left=20] node[left,font=\scriptsize] {$\Phi$} (3.45,{\A(3.45)-0.2});
    \end{tikzpicture}
    \caption{The slice $\{t=t_Q\}$ of a $\Gamma$-cube $Q=\Phi(Q_0)$, for
    $n=2$. Since $\Phi$ does not change $t$, it maps the slice
    of $Q_0$ onto the slice of $Q$. The shaded disc and the dashed disc are
    the slices of $B_{p}({\X}_Q,\ell(Q)/4)$ and of
    $B_Q=B_{p}({\X}_Q,3\diam_{p}(Q))$.}
    \label{fig:gamma-cube}
\end{figure}
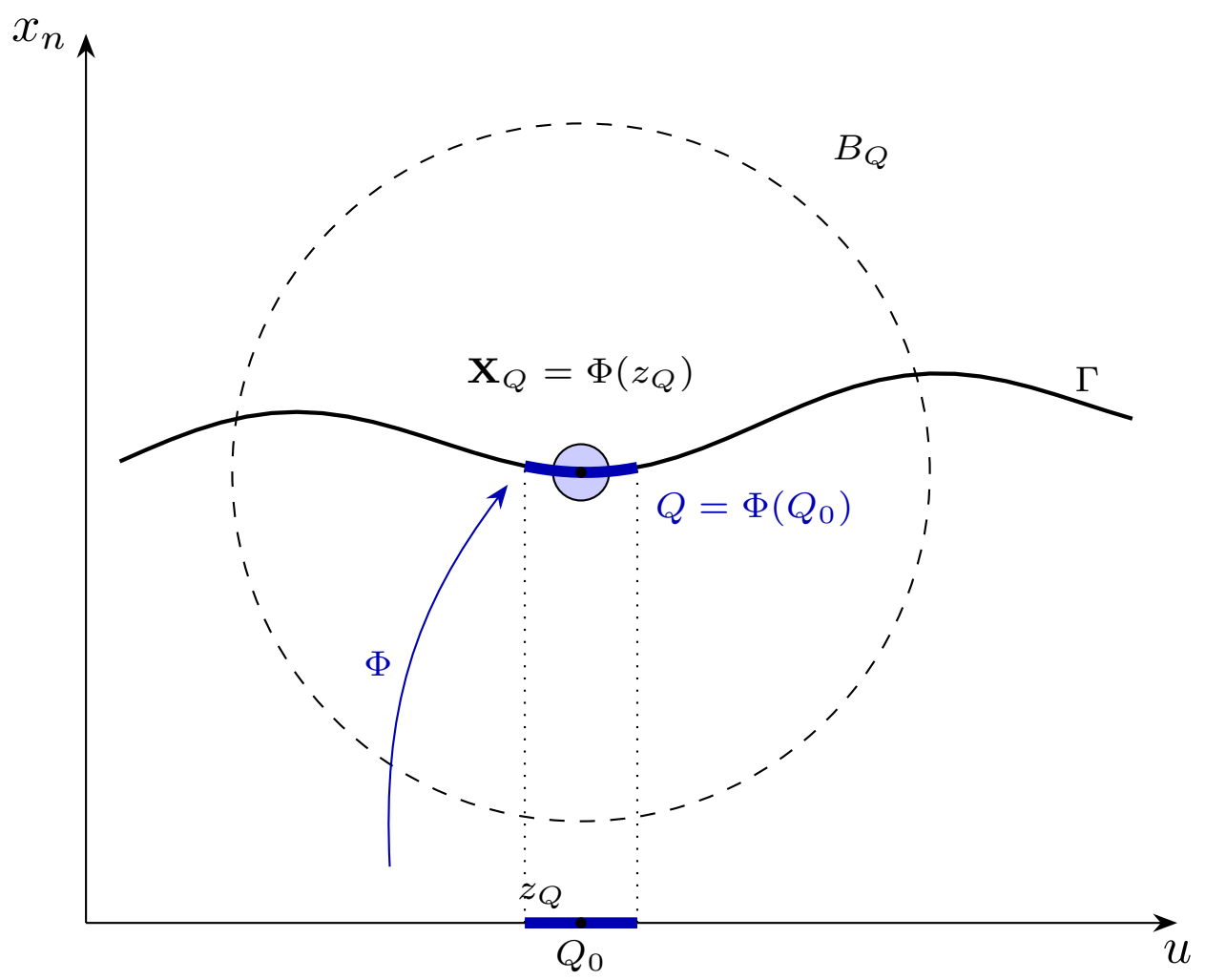

\begin{notation}[Dyadic cubes over $\Gamma$]
    It will be convenient to introduce a collection cubes on $\Gamma$. As in \cite{tolsa2015characterization}, these will be the cubes induced by the cubes on the base.
    \begin{itemize}
        \item The parabolic dyadic cubes of $\mathbb P$ are the half-open cubes with spatial side $\ell=2^{-j}$ and time side $\ell^2$, $j\in\mathbb Z$.

        \item A \emph{$\Gamma$-cube} is a set $Q=\Phi(Q_0)$, where $Q_0$ is a parabolic dyadic cube of $\mathbb P$.

        \item We let $\ell(Q)=\ell(Q_0)$ and $z_Q=(u_Q,t_Q)$ be the center of $Q_0$, and we call ${\X}_Q:=\Phi(z_Q)$ the center of $Q$.

        \item We write $\mathcal D_\Gamma$ for the family of all $\Gamma$-cubes and $\mathcal D_\Gamma^0$ for those cubes with $\ell(Q)\le1$. These are not the cubes $\mathcal D$ of Section~\ref{sec:cubes}.

        \item The $\Gamma$-cubes of a fixed side length are disjoint and cover $\Gamma$, and two $\Gamma$-cubes are either nested or disjoint.

        \item Since $\Phi$ maps the center of the base cube to the center of the graph cube, one has that $\Phi(z)\in B_{p}({\X}_Q,r)$ implies
        $z\in B_{\mathbb P}(z_Q,r)$. Together with the distance bounds for
        $\Phi$ and the size of the base cubes, this gives
        \[
        \Gamma\cap B_{p}({\X}_Q,\ell(Q)/4)\subset Q,
        \qquad
        \ell(Q)\le\diam_{p}(Q)\le C\ell(Q),
        \qquad
        \sigma(Q)\approx\ell(Q)^{n+1}.
        \]

        \item We define the ball $B_Q$ associated to $Q$ by 
        \[B_Q:=B_{p}({\X}_Q,3\diam_{p}(Q))\] 
        which contains $Q$ and is centered on $\Gamma$ (see
        Figure~\ref{fig:gamma-cube}).
    \end{itemize}       
\end{notation}

Now that we have introduced a collection of cubes on $\Gamma$, it will also be convenient to define the corresponding $\alpha$-numbers. Moreover, it will be useful to quantity how flat the the map $A$ (which defines $\Gamma$) is near the corresponding base cube. 

\begin{notation}[Alpha numbers]
    Recalling Definition~\ref{d:alpha-numbers}, for each $Q\in\mathcal D_\Gamma$ we set
    \begin{align}
        \alpha_\mu(Q)\coloneqq\alpha_\mu({\X}_Q,r(B_Q)).
    \end{align}
    For a fixed vertical plane $V$, we
    also write
    \[
    \alpha_{\mu,V}(Q)\coloneqq \frac1{r(B_Q)^{n+2}}\inf_{a\ge0}
    \dist_{B_Q}(\mu,a\lambda_V),
    \]
    so that $\alpha_\mu(Q)=\inf_V\alpha_{\mu,V}(Q)$, and in sums $Q\subset R$
    includes $Q=R$.
\end{notation}

\begin{notation}[Gamma numbers]
For $Q\in\mathcal D_\Gamma$, we set
\[D_Q\coloneqq B_{\mathbb P}(z_Q,2r(B_Q)),\] 
so that 
\begin{align}
    \Phi^{-1}(B_Q) \subseteq D_Q.  
\end{align}
Then, set
\[
\gamma_Q^2 :=  \frac{1}{\ell(Q)^{n+3}}
\inf_{\substack{p\in\mathbb R^{n-1}, \\ q\in\mathbb R}}
\int_{D_Q}|A(u,t)-p\cdot(u-u_Q)-q|^2\,dm(u,t).
\]
The infimum is over affine functions of $u$, which do not depend on $t$. This is the analogue for the ball $D_Q$ of the $\gamma$-number of Definition~\ref{def:gamma}. 
\end{notation}

After all the above notation and definitions, let us go back to our goal in this section, which is to provide a proof of the Main Lemma~\ref{thm:graph-alpha}. An elementary reduction is to note that, to prove the Main Lemma, it suffices to show the following discretized estimate:
\begin{equation}\label{eq:main-cubes}
\sum_{Q\in\mathcal D_\Gamma^0:\,{\X}\in Q}\alpha_\mu(Q)^2<\infty
\qquad\text{for }\sigma\text{-a.e. }{\X}\in\Gamma.
\end{equation}
Indeed, let ${\X}\in Q$ and $\ell(Q)/2<r\le\ell(Q)$.  Since
$d_{p}({\X},{\X}_Q)\le\diam_{p}(Q)$ and
$r\le\diam_{p}(Q)$, we have $B_{p}({\X},r)\subset B_Q$.
Restricting the test functions to this smaller ball, we get
$\alpha_\mu({\X},r)\le(r(B_Q)/r)^{n+2}\alpha_\mu(Q)\lesssim\alpha_\mu(Q)$.
Integrating over $\ell(Q)/2<r\le\ell(Q)$ and summing over the cubes
of $\mathcal D_\Gamma^0$ which contain ${\X}$ gives the Main Lemma at ${\X}$. \\

A first step towards proving \eqref{eq:main-cubes} is to show that any measure which is absolute continuous with respect to surface measure on a regular parabolic graph, with uniformly bounded density, satisfies a Carleson condition for the $\alpha$-numbers. Moreover, we want the best plane in the definition of the $\alpha$-numbers to be chosen depending only on the underlying geometry of $\Gamma$. We prove this in Proposition \ref{prop:beta-alpha} below for graphs which satisfy a Carleson condition on the $\gamma_Q$ coefficients from above. It turns out that regular parabolic graphs naturally satisfy this assumption, see Proposition \ref{prop:graph-measure}. More precisely, for a bounded Borel function $g\ge0$ on $\Gamma$, we compare $g\sigma$ with the class of flat measures. In Proposition~\ref{prop:beta-alpha}, we prove a Carleson condition for the quantities $\alpha_{g\sigma,V_Q}(Q)$, where the vertical planes $V_Q$ depend only on $\Gamma$.

\begin{prop}[Control on gamma numbers implies control on alpha numbers]  \label{prop:beta-alpha}
Let $\Gamma$ be a Lip(1,1/2) graph and
$\sigma=\mathcal H_{p}^{n+1}\lfloor\Gamma$. Suppose that
\begin{align}\label{e:gamma-carleson}
\sum_{Q\subset R}\gamma_Q^2\sigma(Q)\le C_\gamma\sigma(R)
\qquad\text{for every }R\in\mathcal D_\Gamma.
\end{align}
Then there are vertical planes $V_Q$ and numbers $c_Q\in[0,C]$,
$Q\in\mathcal D_\Gamma$, such that the following holds for
every $R\in\mathcal D_\Gamma$.
\begin{enumerate}
\item The numbers
$\kappa_Q=r(B_Q)^{-(n+2)}\dist_{B_Q}(\sigma,c_Q\lambda_{V_Q})$
satisfy
\[
\sum_{Q\subset R}\kappa_Q^2\sigma(Q)\lesssim\sigma(R).
\]
\item For every bounded Borel function $g\ge0$ on $\Gamma$,
\[
\sum_{Q\subset R}\alpha_{g\sigma,V_Q}(Q)^2\sigma(Q)
\lesssim\|g\|_{L^\infty(\sigma)}^2\,\sigma(R).
\]
\end{enumerate}
The constant $C$ depends only on $n$ and $L$. The implicit constants
depend on $n$, $L$, $C_\gamma$ and the fixed dyadic parameters.
\end{prop}

\begin{proof}
 The proof follows that of \cite[Theorem~1.1]{TolsaUR}. The idea is as follows: we take $V_Q$ as the graph of the best-approximating affine function for $A$ over $D_Q$ in the definition of $\gamma_Q$. Since this function does not depend on $t$, $V_Q$ is a vertical plane. Then, we pass from $g\sigma$ to a flat measure in two steps. First, we move each point $\Phi(z)$ of $\Gamma$ to the point of $V_Q$ over the same base point $z$. This moves it by the difference between the heights of $\Gamma$ and $V_Q$ over $z$, and we estimate the total cost using $\gamma_Q$. Second, we replace the density by its average and we estimate the resulting error by martingale differences. Finally, we use the assumed square Carleson estimate for $\gamma_Q$ and prove the corresponding estimate for the density term using martingale differences. See Figure~\ref{fig:move-flatten}. 

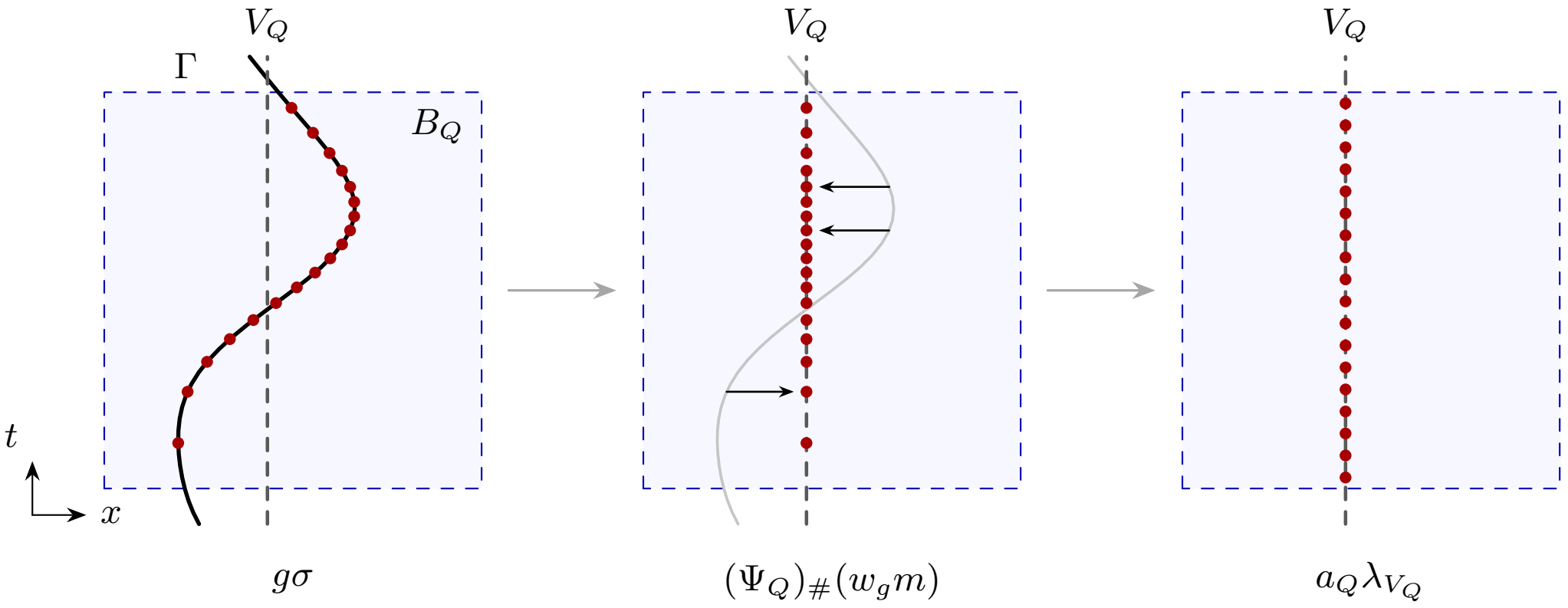
\begin{figure}[H]
\centering
\begin{tikzpicture}[x=1.55cm,y=1.55cm,font=\small,>=Stealth,
  line cap=round,line join=round]
  \def\q{-0.1420126989}
  \def\offset{0.1947091712}
  \def\transportdots{-0.4438/-0.8494,-0.3916/-0.5640,-0.2833/-0.3974,-0.1563/-0.2712,-0.0256/-0.1650,0.1005/-0.0705,0.2164/0.0166,0.3180/0.0989,0.4026/0.1783,0.4681/0.2561,0.5127/0.3335,0.5356/0.4118,0.5358/0.4922,0.5131/0.5761,0.4674/0.6658,0.3985/0.7642,0.3061/0.8770,0.1869/1.0155}
  \begin{scope}
    \fill[blue!3] (-1.05,-1.1025) rectangle (1.05,1.1025);
    \draw[blue!70!black,dashed] (-1.05,-1.1025) rectangle (1.05,1.1025);
    \node[anchor=north east,inner sep=4pt] at (1.05,1.1025) {$B_Q$};
    \draw[black,line width=1pt]
      plot coordinates {
(-0.52263,-1.30000) (-0.55537,-1.23500) (-0.58128,-1.17000) (-0.60143,-1.10500) (-0.61676,-1.04000) (-0.62793,-0.97500) (-0.63526,-0.91000) (-0.63857,-0.84500) (-0.63730,-0.78000) (-0.63047,-0.71500) (-0.61687,-0.65000) (-0.59519,-0.58500) (-0.56426,-0.52000) (-0.52317,-0.45500) (-0.47152,-0.39000) (-0.40953,-0.32500) (-0.33810,-0.26000) (-0.25887,-0.19500) (-0.17414,-0.13000) (-0.08678,-0.06500) (0.00000,0.00000) (0.08280,0.06500) (0.15833,0.13000) (0.22358,0.19500) (0.27611,0.26000) (0.31416,0.32500) (0.33677,0.39000) (0.34383,0.45500) (0.33604,0.52000) (0.31481,0.58500) (0.28210,0.65000) (0.24019,0.71500) (0.19154,0.78000) (0.13847,0.84500) (0.08307,0.91000) (0.02700,0.97500) (-0.02862,1.04000) (-0.08315,1.10500) (-0.13649,1.17000) (-0.18889,1.23500) (-0.24087,1.30000)
};
    \draw[gray!70!black,dashed,line width=.85pt] (\q,-1.3)--(\q,1.3);
    \node[above] at (\q,1.3) {$V_Q$};
    \node[anchor=east,black] at (-.45,1.25) {$\Gamma$};
    \foreach \px/\pt in \transportdots{
      \fill[red!65!black] ({\px-\offset},\pt) circle (1.45pt);
    }
    \node at (0,-1.62) {$g\sigma$};
  \end{scope}
  \begin{scope}[shift={(3,0)}]
    \fill[blue!3] (-1.05,-1.1025) rectangle (1.05,1.1025);
    \draw[blue!70!black,dashed] (-1.05,-1.1025) rectangle (1.05,1.1025);
    \draw[gray!45,line width=.7pt]
      plot coordinates {
(-0.52263,-1.30000) (-0.55537,-1.23500) (-0.58128,-1.17000) (-0.60143,-1.10500) (-0.61676,-1.04000) (-0.62793,-0.97500) (-0.63526,-0.91000) (-0.63857,-0.84500) (-0.63730,-0.78000) (-0.63047,-0.71500) (-0.61687,-0.65000) (-0.59519,-0.58500) (-0.56426,-0.52000) (-0.52317,-0.45500) (-0.47152,-0.39000) (-0.40953,-0.32500) (-0.33810,-0.26000) (-0.25887,-0.19500) (-0.17414,-0.13000) (-0.08678,-0.06500) (0.00000,0.00000) (0.08280,0.06500) (0.15833,0.13000) (0.22358,0.19500) (0.27611,0.26000) (0.31416,0.32500) (0.33677,0.39000) (0.34383,0.45500) (0.33604,0.52000) (0.31481,0.58500) (0.28210,0.65000) (0.24019,0.71500) (0.19154,0.78000) (0.13847,0.84500) (0.08307,0.91000) (0.02700,0.97500) (-0.02862,1.04000) (-0.08315,1.10500) (-0.13649,1.17000) (-0.18889,1.23500) (-0.24087,1.30000)
};
    \draw[gray!70!black,dashed,line width=.85pt] (\q,-1.3)--(\q,1.3);
    \node[above] at (\q,1.3) {$V_Q$};
    \draw[black,line width=.5pt,-{Stealth[length=1.3mm]}]
      ({-.3916-\offset},-.5640)--({\q-.07},-.5640);
    \draw[black,line width=.5pt,-{Stealth[length=1.3mm]}]
      ({.5127-\offset},.3335)--({\q+.07},.3335);
    \draw[black,line width=.5pt,-{Stealth[length=1.3mm]}]
      ({.5131-\offset},.5761)--({\q+.07},.5761);
    \foreach \px/\pt in \transportdots{
      \fill[red!65!black] (\q,\pt) circle (1.45pt);
    }
    \node at (0,-1.62) {$(\Psi_Q)_\#(w_gm)$};
  \end{scope}
  \begin{scope}[shift={(6,0)}]
    \fill[blue!3] (-1.05,-1.1025) rectangle (1.05,1.1025);
    \draw[blue!70!black,dashed] (-1.05,-1.1025) rectangle (1.05,1.1025);
    \draw[gray!70!black,dashed,line width=.85pt] (\q,-1.3)--(\q,1.3);
    \node[above] at (\q,1.3) {$V_Q$};
    \foreach \k in {1,...,18}{
      \fill[red!65!black]
        (\q,{-1.1025+(\k-.5)*2.205/18}) circle (1.45pt);
    }
    \node at (0,-1.62) {$a_Q\lambda_{V_Q}$};
  \end{scope}
  \draw[->,gray!70,line width=.6pt] (1.2,0)--(1.8,0);
  \draw[->,gray!70,line width=.6pt] (4.2,0)--(4.8,0);
  \draw[black,->] (-1.45,-1.25)--(-1.15,-1.25) node[right] {$x$};
  \draw[black,->] (-1.45,-1.25)--(-1.45,-.95) node[above left] {$t$};
\end{tikzpicture}
\caption{The two terms in \eqref{eq:beta_alpha:2-22} for $n=1$.
Red dots show density with respect to $dt$: the middle panel preserves their time coordinates,
and the last has constant density $(w_g)_{J_Q}$.}
\label{fig:move-flatten}
\end{figure} 

Let us now get more concrete and provide a rigorous proof. We start by showing the second implication. At the end of the proof we comment on how to obtain the first implication. 

For each $Q$, let \[L_Q(u) :=  p_Q\cdot(u-u_Q)+q_Q\] 
be the least-squares minimizer in $\gamma_Q$. Since $A$ is continuous, $A\in L^2(D_Q,m)$. The affine functions of $u$
form a finite dimensional subspace of $L^2(D_Q,m)$, so the minimizer
exists and is unique. An affine function vanishing a.e. on $D_Q$ is zero, so $p_Q$ and $q_Q$ are unique as well. We set
\[
V_Q:=\{(u,L_Q(u),t):(u,t)\in\mathbb P\},
\qquad
\Psi_Q(u,t):=(u,L_Q(u),t).
\]
Note that these choices depend only on the geometry of $\Gamma$, and not on $g$.

We want show now that $\Psi_Q$ is bi-Lipschitz with constants depending only on $n$ and $L$. For $n=1$, $L_Q$ is constant so there is nothing to show. For $n\ge2$, let $r_Q:=2r(B_Q)$ be the
radius of $D_Q$. Then, $D_Q=\{|u-u_Q|\le r_Q\}\times\{|t-t_Q|\le r_Q^2\}$ and 
\[
\int_{|u-u_Q|\le r_Q}(u-u_Q)\,du=0,\qquad
\frac1{m(D_Q)}\int_{D_Q}(u_i-u_{Q,i})(u_j-u_{Q,j})\,dm=\frac{\delta_{ij}r_Q^2}{n+1}.
\]
By the first identity and Fubini's theorem, every function of $t$ alone, such as
$q_Q$ or $A(u_Q,t)$, is orthogonal in $L^2(D_Q,m)$ to each $u_i-u_{Q,i}$. Thus the
normal equations give
\[
p_Q=\frac{n+1}{r_Q^2m(D_Q)}\int_{D_Q}(u-u_Q)\bigl(A(u,t)-A(u_Q,t)\bigr)\,dm(u,t).
\]
Since $|u-u_Q|\le r_Q$ and $|A(u,t)-A(u_Q,t)|\le L|u-u_Q|$ on $D_Q$, we get
$|p_Q|\le(n+1)L$. Since $\Psi_Q$ fixes the $(u,t)$ coordinates and $|p_Q|\le(n+1)L$, we have
\[
\rho(z,z')\le d_{p}(\Psi_Q(z),\Psi_Q(z'))\le(1+|p_Q|)\rho(z,z'),
\]
as required. 

To compare $g\sigma$ with a flat measure on $V_Q$, we write
$\lambda_{V_Q}$ in the coordinates of the base. The measure $H\mapsto\lambda_{V_Q}(\Psi_Q(H))$ on
$\mathbb P$ is translation invariant, since translations in the base
give translations along $V_Q$. By the same argument as in \eqref{eq:graph-measure-comparison}, applied to
$\Psi_Q$ instead of $\Phi$, it is comparable to $m$. Thus, it is finite and
positive on balls. Translation invariant Radon measures on
$\mathbb P=\mathbb R^n$ are unique up to a constant factor, so it equals $\Lambda_Q m$ for a constant
$\Lambda_Q$, and
\begin{equation}\label{eq:beta_alpha:2-20}
\lambda_{V_Q}=\Lambda_Q(\Psi_Q)_\#m,
\qquad C_L^{-1}\le\Lambda_Q\le C_L.
\end{equation}

Now fix a bounded Borel function $g\ge0$ on $\Gamma$, and write
$\|g\|_\infty=\|g\|_{L^\infty(\sigma)}$. Set 
\[w_g=(g\circ\Phi)w,\]
with $w$ as in \eqref{eq:beta_alpha:2-3}. 
Then, 
$g\sigma=\Phi_\#(w_gm)$.
Since $\sigma$ and $\Phi_\#m$ have the same null sets,
\[
\|w_g\|_{L^\infty(m)}\le C_L\|g\|_\infty.
\]

As in \cite[Section~4]{TolsaUR}, we expand this density in martingale
differences. To do this, we must first place each $D_Q$ in
a dyadic cube of comparable size.
For this, we use finitely many adjacent dyadic systems on $(\mathbb P,\rho)$
from \cite[Theorem~4.1]{HK}, with a fixed scale parameter $0<\delta\le\frac{1}{96}$.
For a cube $I$ of generation $k$, we set $\ell(I)=\delta^k$.
The cubes in each system form nested partitions and, for each $I$, 
there is a point $z_I\in I$ such that
\[
B_{\mathbb P}(z_I,C^{-1}\ell(I))\subset I\subset B_{\mathbb P}(z_I,C\ell(I)),
\]
where $C>0$ is fixed. In particular,
$\diam_\rho I\lesssim\ell(I)$ and $m(I)\approx\ell(I)^{n+1}$.
By the covering property in \cite[Theorem~4.1]{HK}, for each $Q$ we choose a cube $J_Q$ in one of these systems such that 
\[
D_Q\subset J_Q,\qquad\ell(J_Q)\approx\ell(Q).
\]

For cubes $I \subset J$ in 
one of these systems and
$\varphi\in L^2(J,m)$, define 
\[
\Delta_I\varphi
=\sum_{I'}(\varphi_{I'}-\varphi_I)\mathbf 1_{I'},
\]
where $I'$ runs over the children of $I$ in the same system and
$\varphi_I$ is the $m$-average of $\varphi$ on $I$.
These martingale differences have mean zero on $I$ and are pairwise
orthogonal within that system. By standard properties of martingale differences, one has
\[
\varphi-\varphi_J=\sum_{I\subset J}\Delta_I\varphi
\quad\text{in }L^2(J,m),
\qquad
\sum_{I\subset J}\|\Delta_I\varphi\|_{L^2(m)}^2\le\int_J|\varphi|^2\,dm.
\]

We apply these martingale differences with $\varphi=w_g\in L^2(J_Q,m)$. For each $Q \in \mathcal{D}_\Gamma$, we put
\[
\omega_Q\coloneqq \sum_{I\subset J_Q}
\frac{\ell(I)^{1+\frac{n+1}{2}}}{\ell(Q)^{n+2}}\,\|\Delta_Iw_g\|_{L^2(m)},
\]
where the cubes $I$ belong to the system of $J_Q$. Since
$r(B_Q)\approx\ell(Q)$, the denominator is, up to constants, the
normalization in the $\alpha$-coefficients. The factors of $\ell(I)$ naturally appear in the computions below. The factor
$\ell(I)$ comes from cancellation against a Lipschitz function
in \eqref{eq:beta_alpha:2-24} below. The factor $\ell(I)^{(n+1)/2}$ comes from
the Cauchy--Schwarz inequality and $m(I)\approx\ell(I)^{n+1}$.

We claim that
\begin{equation}\label{eq:beta_alpha:2-11}
\sum_{Q\subset R}\omega_Q^2\sigma(Q)
\lesssim_L\|g\|_\infty^2\sigma(R).
\end{equation}

Fix $Q$ for the moment and let $k$ be the generation of $J_Q$, so that
$\delta^k=\ell(J_Q)\approx\ell(Q)$. The cubes $I\subset J_Q$ of
generation $k+j$ are disjoint and each have measure comparable to
$\delta^{(k+j)(n+1)}$. Since $m(J_Q)\approx\delta^{k(n+1)}$, there are at most
$C\delta^{-j(n+1)}$ of them.
Hence
\[
\sum_{I\subset J_Q}\left(\frac{\ell(I)}{\ell(Q)}\right)^{n+2}
\lesssim\sum_{j\ge0}\delta^{-j(n+1)}\delta^{j(n+2)}
=\sum_{j\ge0}\delta^j\lesssim1.
\]
Write
\[
\left(\frac{\ell(I)}{\ell(Q)}\right)^{1+\frac{n+1}{2}}
=\left(\frac{\ell(I)}{\ell(Q)}\right)^{\frac{n+2}{2}}
\left(\frac{\ell(I)}{\ell(Q)}\right)^{\frac{1}{2}}.
\]
Since $\sigma(Q)\approx\ell(Q)^{n+1}$, the Cauchy--Schwarz inequality and the preceding
estimate give
\[
\begin{aligned}
\omega_Q^2\sigma(Q)
&\approx\left(\sum_{I\subset J_Q}
\left(\frac{\ell(I)}{\ell(Q)}\right)^{1+\frac{n+1}{2}}
\|\Delta_Iw_g\|_{L^2(m)}\right)^2\\
&\le\sum_{I\subset J_Q}\left(\frac{\ell(I)}{\ell(Q)}\right)^{n+2}
\sum_{I\subset J_Q}\frac{\ell(I)}{\ell(Q)}\|\Delta_Iw_g\|_{L^2(m)}^2
\lesssim\sum_{I\subset J_Q}\frac{\ell(I)}{\ell(Q)}
\|\Delta_Iw_g\|_{L^2(m)}^2.
\end{aligned}
\]

Below, we will want to take a sum of the above inequality over cubes $Q \subset R$. To estimate the resulting right-hand side we will need two facts: each cube $J$ equals $J_Q$ for boundedly
many $Q$ 
and, for each $I$ in the same system as $J_Q$,  
\begin{align}\label{e:geometric-series}
\sum_{Q:\,J_Q\supset I}\frac{\ell(I)}{\ell(Q)}
\lesssim 1.
\end{align}
By changing the
order of summation, this reduces the estimate to orthogonality. To obtain the above facts, fix a system of cubes and a cube $J$ in it. If $J=J_Q$, then $\ell(Q)\approx\ell(J)$, so
$\ell(Q)$ takes boundedly many values.
For each value, the base cubes $Q_0\subset D_Q\subset J$ are disjoint.
Since $m(Q_0)=\ell(Q)^{n+1}\approx m(J)$, there are boundedly many of them. Now fix a cube $I$, and
consider the cubes $Q$ such that $J_Q$ lies in the system of $I$ and
$J_Q\supset I$. Since the cubes of a system are nested, $J_Q$ is $I$ or
an ancestor of $I$, and there is one such cube in each generation. Thus,
by the bounded multiplicity and $\ell(Q)\approx\ell(J_Q)$, we have
\[
\sum_{Q:\,J_Q\supset I}\frac{\ell(I)}{\ell(Q)}
\lesssim\sum_{j\ge0}\delta^j\lesssim1.
\]

Now we are ready to take the sum. Fix $R$ and let $R_0$ denote the corresponding base cube. To apply orthogonality, we place the
cubes $J_Q$, $Q\subset R$, in boundedly many cubes of side length comparable
to $\ell(R)$. Let $k_R$ be the largest
integer such that $\delta^{k_R}\ge\ell(J_Q)$ for every $Q\subset R$.
Then $\delta^{k_R}\approx\ell(R)$. For $Q\subset R$, the cube $J_Q$ lies
in a cube $T$ of generation $k_R$ of its own system. Since $T$ contains
$z_Q\in R_0$ and $\diam_\rho T\lesssim\ell(R)$, it lies in
$B_{\mathbb P}(z_R,C\ell(R))$. In each system, these cubes $T$ are disjoint
and have measure $\approx\ell(R)^{n+1}$, so there are boundedly many of
them. Each cube $I$ with $I\subset J_Q$ for some $Q\subset R$ lies in
exactly one of these cubes $T$ of its own system. In the following
display, $I$ and $J_Q$ belong to the system of $T$. By Tonelli's
theorem, we may change the order of summation. Since
$\|w_g\|_{L^\infty(m)}\le C_L\|g\|_\infty$ and
$m(T)\approx\ell(R)^{n+1}\approx\sigma(R)$, we get
\[
\begin{aligned}
\sum_{Q\subset R}\omega_Q^2\sigma(Q)
\lesssim\sum_T\sum_{I\subset T}\|\Delta_Iw_g\|_{L^2(m)}^2
\sum_{\substack{Q\subset R\\ I\subset J_Q\subset T}}\frac{\ell(I)}{\ell(Q)}
\lesssim\sum_T\int_T|w_g|^2\,dm
\lesssim_L\|g\|_\infty^2\sigma(R),
\end{aligned}
\]
where the second step uses \eqref{e:geometric-series} and orthogonality on each $T$.
This proves \eqref{eq:beta_alpha:2-11}.

For each $Q$, we choose the coefficient $a_Q$ such that $a_Q\lambda_{V_Q}$ has
the constant density $(w_g)_{J_Q}$ in the coordinates of the base. That is, we set
\[
a_Q=\frac{(w_g)_{J_Q}}{\Lambda_Q}.
\]
Although not explicit in our notation, we note that $a_Q$ depends on the choice of $g$. Since $w_g\ge0$, we have $a_Q\ge0$.
Fix a $1$-Lipschitz function $f$
supported in $B_Q$. Put $h=f\circ\Psi_Q$. By
$g\sigma=\Phi_\#(w_gm)$ and \eqref{eq:beta_alpha:2-20}, we see that
\[
\begin{aligned}
\int f\,d(g\sigma)&=\int_{\mathbb P}(f\circ\Phi)\,w_g\,dm,\\
a_Q\int f\,d\lambda_{V_Q}&=a_Q\,\Lambda_Q\int_{\mathbb P}h\,dm
=(w_g)_{J_Q}\int_{\mathbb P}h\,dm.
\end{aligned}
\]
By subtracting the two identities and adding and subtracting
$\int_{\mathbb P}h\,w_g\,dm$, we obtain 
\begin{equation}\label{eq:beta_alpha:2-22}
\int f\,d(g\sigma)-a_Q\int f\,d\lambda_{V_Q}
= \int_{\mathbb P}(f\circ\Phi-f\circ\Psi_Q)w_g\,dm
+\int_{\mathbb P}h\bigl(w_g-(w_g)_{J_Q}\bigr)\,dm.
\end{equation}
The maps $\Phi$ and $\Psi_Q$ keep the base point $z_Q$. So, if
$\Phi(z)\in B_Q$ or $\Psi_Q(z)\in B_Q$, then
\[
z\in B_{\mathbb P}(z_Q,r(B_Q)).
\]
Thus $f\circ\Phi$ and $h$ vanish outside $B_{\mathbb P}(z_Q,r(B_Q))$, and
$B_{\mathbb P}(z_Q,r(B_Q))\subset D_Q\subset J_Q$.

For the first term in \eqref{eq:beta_alpha:2-22}, the points
$\Phi(z)$ and $\Psi_Q(z)$, with $z=(u,t)$, differ only in the $x_n$-coordinate.
Since $f$ is $1$-Lipschitz, we have
\[
|f(\Phi(z))-f(\Psi_Q(z))|\le |A(u,t)-L_Q(u)|.
\]
Thus, since $f\circ\Phi$ and $h$ vanish outside $D_Q$, the first term is bounded in absolute value by
$\int_{D_Q}|A-L_Q|w_g\,dm$. Since $w_g\le C_L\|g\|_\infty$, the Cauchy--Schwarz
inequality, $m(D_Q)\approx\ell(Q)^{n+1}$ and the definition of $\gamma_Q$,
with $L_Q$ as its minimizer, give
\begin{equation}\label{eq:beta_alpha:2-23}
\int_{D_Q}|A-L_Q|w_g\,dm
\lesssim_L\|g\|_\infty m(D_Q)^{\frac{1}{2}}
\left(\int_{D_Q}|A-L_Q|^2\,dm\right)^{\frac{1}{2}}
\lesssim_L\|g\|_\infty\gamma_Q\ell(Q)^{n+2}.
\end{equation}
Now consider the second term. Since $f$ is $1$-Lipschitz and
$\Psi_Q$ is bi-Lipschitz, $\operatorname{Lip}_\rho h\lesssim_L1$.
Also, $h$ is bounded and vanishes outside $J_Q$, so $h\in L^2(J_Q,m)$.
The expansion of $w_g-(w_g)_{J_Q}$ converges in $L^2(J_Q,m)$, so we
can integrate it against $h$ term by term.
Since each $\Delta_Iw_g$ vanishes outside $I$ and has mean zero on $I$,
with $z_I$ as above, we get
\[
\int_{\mathbb P}h\bigl(w_g-(w_g)_{J_Q}\bigr)\,dm
=\sum_{I\subset J_Q}\int_Ih\,\Delta_Iw_g\,dm
=\sum_{I\subset J_Q}\int_I\bigl(h-h(z_I)\bigr)\Delta_Iw_g\,dm.
\]
Now for $z\in I$, we have $\rho(z,z_I)\le C\ell(I)$, and so
$|h(z)-h(z_I)|\lesssim_L\ell(I)$. Thus, by the Cauchy--Schwarz inequality
on each $I$ and $m(I)\approx\ell(I)^{n+1}$, we get
\begin{equation}\label{eq:beta_alpha:2-24}
\left|\int_{\mathbb P}h\bigl(w_g-(w_g)_{J_Q}\bigr)\,dm\right|
\lesssim_L\sum_{I\subset J_Q}\ell(I)\,m(I)^{\frac{1}{2}}
\|\Delta_Iw_g\|_{L^2(m)}
\approx\omega_Q\ell(Q)^{n+2}.
\end{equation}

Taking the supremum over $f$ and using
\eqref{eq:beta_alpha:2-22}--\eqref{eq:beta_alpha:2-24} and
$r(B_Q)\approx\ell(Q)$, we get
\begin{equation}\label{eq:beta_alpha:2-25}
\frac{\dist_{B_Q}(g\sigma, a_Q\lambda_{V_Q})}{r(B_Q)^{n+2}}
\lesssim_L\|g\|_\infty\gamma_Q+\omega_Q.
\end{equation}
Since $a_Q\ge0$, the same bound holds for $\alpha_{g\sigma,V_Q}(Q)$.
Squaring \eqref{eq:beta_alpha:2-25}, summing over $Q\subset R$ and
using the assumed Carleson estimate for $\gamma_Q$ and
\eqref{eq:beta_alpha:2-11}, we get the second conclusion. 

For the
first conclusion, we take $g=1$ and let $c_Q$ be the $a_Q$ from above, corresponding to this choice.
Then $0\le c_Q\le C$ by the bounds for $w$ and $\Lambda_Q$.
For $g=1$, $\kappa_Q$ is the left-hand side of \eqref{eq:beta_alpha:2-25}, so
the same sum gives the bound for $\kappa_Q$.
\end{proof}

We next verify the hypothesis of Proposition~\ref{prop:beta-alpha}
for regular parabolic Lipschitz graphs.

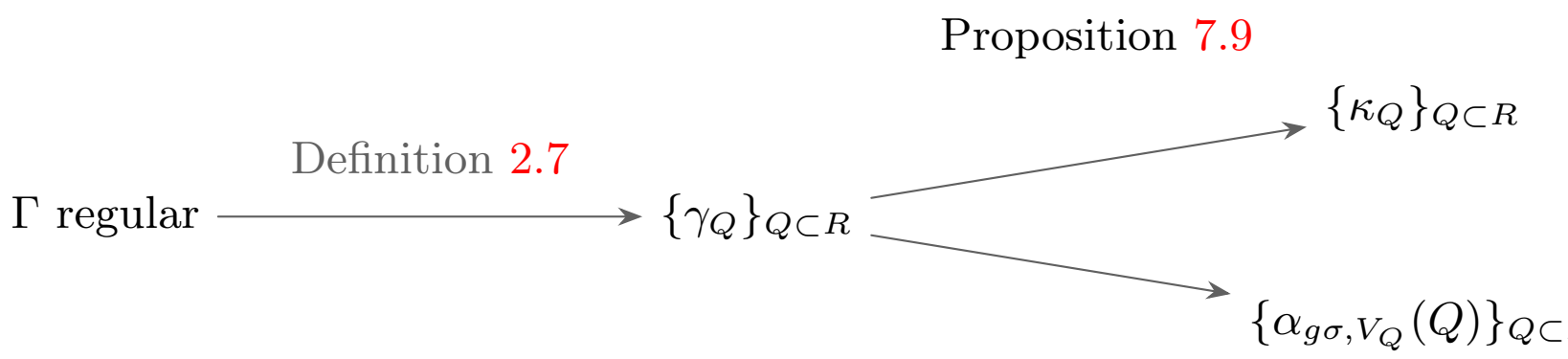
\begin{figure}[H]
\centering
\begin{tikzpicture}[font=\small,>=Stealth,line cap=round]
  \node (regular) at (0,0) {$\Gamma\ \text{regular}$};
  \node (beta) at (4.5,0) {$\{\gamma_Q\}_{Q\subset R}$};
  \node (surface) at (9.1,.75) {$\{\kappa_Q\}_{Q\subset R}$};
  \node (density) at (9.1,-.75) {$\{\alpha_{g\sigma,V_Q}(Q)\}_{Q\subset R}$};
  \draw[->,gray!75!black] (regular)--node[above=5pt]
    {Definition~\ref{def:regular}} (beta);
  \draw[->,gray!75!black] (beta)--(surface);
  \draw[->,gray!75!black] (beta)--(density);
  \node at (6.85,1.23) {Proposition~\ref{prop:beta-alpha}};
\end{tikzpicture}
\caption{Regularity gives the $\gamma_Q$-Carleson estimate, by
Definition~\ref{def:regular} and the comparison in the proof below;
Proposition~\ref{prop:beta-alpha} gives the estimates for
$\kappa_Q$ and $\alpha_{g\sigma,V_Q}$ on the same planes.}
\label{fig:graph-carleson}
\end{figure}

\begin{prop}[{Regular Lip(1,1/2) graphs satisfy the gamma control in Proposition~\ref{prop:beta-alpha}}] \label{prop:graph-measure}
Let $\Gamma$ be a regular Lip(1,1/2) graph,
with $\sigma=\mathcal H_{p}^{n+1}\lfloor\Gamma$. Then the
conclusions of Proposition~\ref{prop:beta-alpha} hold, with constants
depending on $n,L~\& ~M$ where  $M$ is the constant in Definition~ \ref{def:regular}. In particular,
\[
\sum_{Q\subset R}\alpha_\sigma(Q)^2\sigma(Q)
\lesssim\sigma(R)
\qquad\text{for every }R\in\mathcal D_\Gamma.
\]
\end{prop}

\begin{proof}
By Proposition~\ref{prop:beta-alpha} (since $\alpha_\sigma(Q)\le\kappa_Q$), it suffices to prove that
$\Gamma$ satisfies the Carleson condition \eqref{e:gamma-carleson}. 
For $z=(u,t)\in\mathbb P$ and $r>0$, the $\gamma$-numbers of $A$ from
Definition~\ref{def:gamma} are
\[
\gamma_A(z,r)^2=\frac{1}{r^2m(B_{\mathbb P}(z,r))}
\inf_{p\in\mathbb R^{n-1},\,q\in\mathbb R}\int_{B_{\mathbb P}(z,r)}
|A(v,\tau)-p\cdot(v-u)-q|^2\,dm(v,\tau).
\]
Since $A$ is regular, Definition~\ref{def:regular} gives the Carleson
estimate
\begin{equation}\label{eq:gamma-carleson-continuous}
\int_{B_{\mathbb P}(z_0,r_0)}\int_0^{r_0}
\gamma_A(z,r)^2\,\frac{dr}{r}\,dm(z)\lesssim r_0^{n+1}
\end{equation}
for all $z_0\in\mathbb P$ and $r_0>0$.

To pass from \eqref{eq:gamma-carleson-continuous} to the sum over the
$\Gamma$-cubes, we compare $\gamma_Q$ with $\gamma_A(z,r)$ when $z$ is in the
base cube of $Q$ and $r\approx\ell(Q)$. For $z\in Q_0=\Phi^{-1}(Q)$ and $C_D\ell(Q)\le r<2C_D\ell(Q)$, we
have $D_Q\subset B_{\mathbb P}(z,r)$ if $C_D$ is large enough, depending
only on $n$ and $L$. We fix such a $C_D$. The infima in $\gamma_Q$ and
$\gamma_A(z,r)$ are over the same affine functions. Thus, restricting the
integral defining $\gamma_A(z,r)$ to $D_Q$ and using $r\approx\ell(Q)$ gives
$\gamma_Q^2\lesssim_L\gamma_A(z,r)^2$. As $\sigma(Q)\approx m(Q_0)$,
we obtain
\[
\sum_{Q\subset R}\gamma_Q^2\sigma(Q)
\lesssim_L
\sum_{Q\subset R}\int_{Q_0}\int_{C_D\ell(Q)}^{2C_D\ell(Q)}
\gamma_A(z,r)^2\,\frac{dr}{r}\,dm(z)
\lesssim\sigma(R).
\]
For the last inequality, the base cubes of each generation are
disjoint, as are the radius intervals for different generations.
All the integrals are over base points in $\Phi^{-1}(R)$ and radii
at most $2C_D\ell(R)$, so we apply
\eqref{eq:gamma-carleson-continuous} on the ball
$B_{\mathbb P}(z_R,2C_D\ell(R))$, which contains $\Phi^{-1}(R)$. This
proves the proposition.
\end{proof}

\subsection{The case when \texorpdfstring{$\mu$}{mu}
is supported on \texorpdfstring{$\Gamma$}{Gamma}}
\label{the-supported-case}

Before proving the Main Lemma~\ref{thm:graph-alpha} (and more concretely, \eqref{eq:main-cubes}) in full generality (that is, for general measures $\mu$, not necessarily supported on $\Gamma$), we first prove a stronger form of
\eqref{eq:main-cubes} for measures which \textit{are} supported on $\Gamma$, where we compute the $\alpha$-numbers with respect to the fixed planes
$V_Q$ from Proposition~\ref{prop:beta-alpha}, instead of taking the infimum over all vertical planes. This corresponds to
\cite[Lemma~2.3]{tolsa2015characterization}.

In the rest of this section, 
\begin{equation*}
    \text{$V_Q$, $c_Q$ and $\kappa_Q$ are as in Proposition~\ref{prop:beta-alpha},}
\end{equation*} 
which applies to $\Gamma$ by
Proposition~\ref{prop:graph-measure}. 

\begin{lem}[Control of alphas for measures supported on $\Gamma$]  \label{lem:supported}
Let $\nu$ be a finite Borel measure supported on $\Gamma$.
Then
\begin{align}\label{e:alpha-estimate-supported-on-Gamma}
\sum_{Q\in\mathcal D_\Gamma^0:\,{\X}\in Q}
\alpha_{\nu,V_Q}(Q)^2<\infty
\qquad\text{for }\sigma\text{-a.e. }{\X}\in\Gamma.
\end{align}
\end{lem}

If $\nu=g\sigma$ with $g$ bounded, this follows from
Proposition~\ref{prop:graph-measure}. In general, we use a
Calder\'on--Zygmund decomposition of $\nu$ with respect to $\sigma$. Its
bad part has total mass zero on each Calder\'on--Zygmund cube $P$, so its
integral against a test function depends only on the oscillation of that
function on $P$.

\begin{figure}[ht]
        \centering
        \begin{tikzpicture}[x=1.4cm,y=1.4cm,>=Stealth,font=\small]
          \fill[blue!7] (4.5,0) circle (2.4);
          \draw[blue!60!black,dashed] (4.5,0) circle (2.4);
          \node[blue!60!black] at (3.3,1.5) {$B_Q$};
          \draw[orange!80!black,densely dotted] (6.6,0) circle (0.6);
          \draw[orange!80!black,dotted,thick] (6.6,0) circle (1.8);
          \node[orange!80!black] at (8.75,1.2) {$3B_P$};
          \node[orange!80!black] at (6.6,0.85) {$B_P$};
          \draw[thick] (1.0,0) -- (9.1,0) node[right] {$\Gamma$};
          \draw[line width=2.6pt,blue!60!black] (4.1,0) -- (4.9,0);
          \node[blue!60!black,below=3pt] at (4.5,0) {$Q$};
          \draw[line width=2.6pt,orange!80!black] (6.5,0) -- (6.7,0);
          \node[orange!80!black,below=3pt] at (6.6,0) {$P$};
          \foreach \x in {1.15,1.65,2.15,2.65,3.15,3.65,4.15,4.65,5.15,5.65,6.15,6.52,6.60,6.68,7.15,7.65,8.15,8.65}{\fill[red!65!black] (\x,0) circle (1pt);}
          \fill (4.3,0) circle (1.6pt);
          \node[above=2pt] at (4.3,0) {${\X}$};
        \end{tikzpicture}
        \caption{A bad cube $P$ and a cube $Q$, in a fixed-time slice
        of a flat graph ($n=2$). Here $\nu(P)>\lambda\sigma(P)$,
        $P\cap B_Q\ne\varnothing$, $r(B_P)\le r(B_Q)$, and
        ${\X}\in Q\setminus3B_P$. The red dots illustrate $\nu$.}
    \label{fig:cz}
\end{figure}
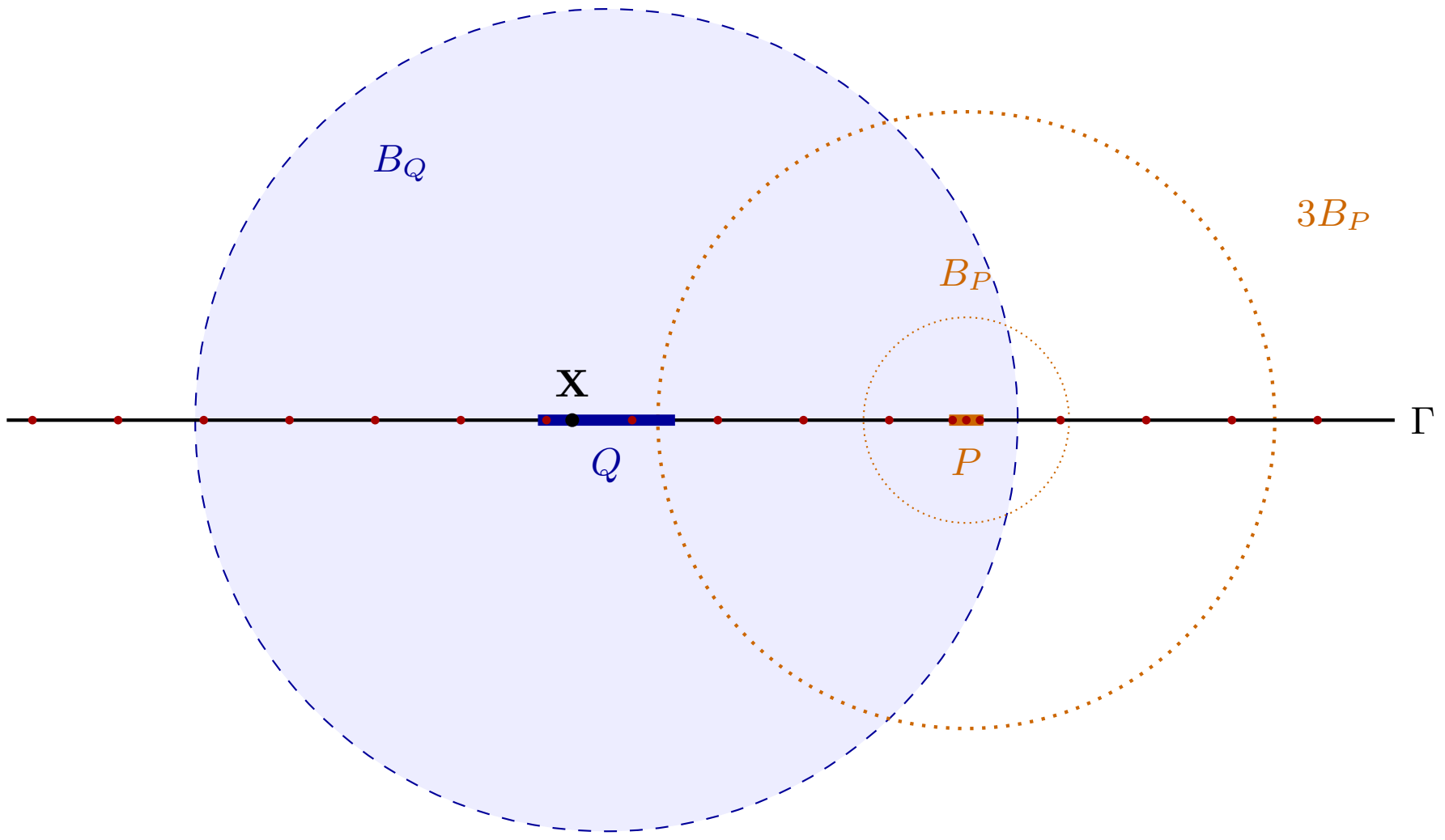

\begin{proof}
Suppose first that $\nu=g\sigma$ with $g\ge0$ bounded. Let
$R\in\mathcal D_\Gamma$ with $\ell(R)=1$. By
Proposition~\ref{prop:graph-measure} and Tonelli's theorem,
\[
\int_R\sum_{Q\subset R:\,{\X}\in Q}\alpha_{g\sigma,V_Q}(Q)^2\,d\sigma({\X})
=\sum_{Q\subset R}\alpha_{g\sigma,V_Q}(Q)^2\sigma(Q)
\lesssim\|g\|_\infty^2\sigma(R)<\infty.
\]
Every $Q\in\mathcal D_\Gamma^0$ which contains a point ${\X}\in R$ lies
in $R$, and these cubes $R$ cover $\Gamma$. So the lemma holds in this
case.

Now let $\nu$ be arbitrary, and fix $\lambda>0$.  Consider the collection of $\Gamma$-cubes $P$ satisfying $\nu(P)>\lambda\sigma(P)$. Since $\nu$ is a finite measure and $\Gamma$ is AD-regular, the diameters of these cubes are uniformly bounded, so each of them lies in a maximal cube with respect to inclusion. Let $\mathcal B$ be the family of these maximal cubes. We call these the \textit{bad cubes}.
By maximality, the cubes in $\mathcal B$ are disjoint. In particular, 
\begin{align}\label{e:sum-of-P-bounded}
    \sum_{P \in \mathcal B} \nu(P) \leq \|\nu\|. 
\end{align}
Furthermore, if
$\widehat P$ is the parent of $P\in\mathcal B$, then by maximality and AD-regularity of $\Gamma$,  
\begin{align}\label{e:nu-less-than-sigma}
\nu(P)\le\lambda\sigma(\widehat P)\lesssim\lambda\sigma(P).\end{align}
Every
$\Gamma$-cube $Q$ which meets $G=\Gamma\setminus\bigcup_{P\in\mathcal B}P$
satisfies $\nu(Q)\le\lambda\sigma(Q)$, otherwise $Q$ would lie in a
bad cube. Hence,
\[\nu\lfloor G\le\lambda\,\sigma\lfloor G.\] 
Indeed, let $E\subset G$ be Borel and $U\supset E$ be open.
Consider the maximal $\Gamma$-cubes with side at most $1$ which lie in $U$ and meet $E$.
They are disjoint and cover $E$, so $\nu(E)\le\lambda\sigma(U)$.
Letting $\sigma(U)\to\sigma(E)$ proves the claim. Thus
\[
\nu=g\sigma+\sum_{P\in\mathcal B}\nu_P,
\qquad
\nu_P=\nu\lfloor P-\frac{\nu(P)}{\sigma(P)}\,\sigma\lfloor P,
\]
where $g\sigma=\nu\lfloor G+\sum_{P}\frac{\nu(P)}{\sigma(P)}\sigma\lfloor P$
and $0\le g\lesssim\lambda$. Each $\nu_P$ is concentrated on $P$ and satisfies
$\nu_P(\Gamma)=0$ and $\|\nu_P\|\le2\nu(P)$.

Consider the set 
\[\widetilde G=\Gamma\setminus\bigcup_{P\in\mathcal B}3B_P.\]
By
AD-regularity and \eqref{e:sum-of-P-bounded}, 
\[
\sigma(\Gamma\setminus\widetilde G)
\le\sum_{P\in\mathcal B}\sigma(3B_P)
\lesssim\sum_{P\in\mathcal B}\sigma(P)
\le\frac1\lambda\sum_{P\in\mathcal B}\nu(P)
\le\frac{\|\nu\|}{\lambda}.
\]
Since $\lambda$ can be taken arbitrarily large, it suffices to prove
the estimate \eqref{e:alpha-estimate-supported-on-Gamma} for $\sigma$-a.e. ${\X}\in\widetilde G$. Let $\mathcal G$ denote the family of $\Gamma$-cubes $Q$ such that  
\begin{align}\label{e:definition-calG}
    Q \not\subseteq \bigcup_{P \in \mathcal{B}} 3B_P.
\end{align}
In this way, any $\Gamma$-cube which contains a point ${\X}\in\widetilde G$ must be an element of $\mathcal{G}$.

Let $Q\in\mathcal G$, and let $f$ be $1$-Lipschitz with
$\supp f\subset B_Q$. For each $P \in \mathcal{B}$, since $\nu_P(\Gamma)=0$ and ${\X}_P\in P$, we have 
\[
\Bigl|\int f\,d\nu_P\Bigr|
=\Bigl|\int_P\bigl(f-f({\X}_P)\bigr)\,d\nu_P\Bigr|
\le\diam_{p}(P)\,\|\nu_P\|
\lesssim\ell(P)\,\nu(P).
\] 
Since $\int f\,d\nu_P=0$ whenever $P\cap B_Q=\varnothing$,
the triangle inequality for $\dist_{B_Q}$ and
$r(B_Q)\ge3\ell(Q)$ give 
\[
\alpha_{\nu,V_Q}(Q)\le\alpha_{g\sigma,V_Q}(Q)
+C\sum_{P\in\mathcal B:\,P\cap B_Q\ne\varnothing}
\frac{\ell(P)\,\nu(P)}{\ell(Q)^{n+2}}.
\]
By the first part of the proof, applied to $g\sigma$, it is enough to
show that
\begin{equation}\label{eq:supported-goal}
\sum_{Q\in\mathcal G\cap\mathcal D_\Gamma^0:\,{\X}\in Q}
\biggl(\sum_{P\in\mathcal B:\,P\cap B_Q\ne\varnothing}
\frac{\ell(P)\,\nu(P)}{\ell(Q)^{n+2}}\biggr)^2<\infty
\qquad\text{for }\sigma\text{-a.e. }{\X}\in\widetilde G.
\end{equation}

Consider for the moment some cube $Q\in\mathcal G$, and let $P\in\mathcal B$ be such that $P \cap B_Q \neq \varnothing$
(see Figure~\ref{fig:cz}). We claim that 
\begin{equation}\label{eq:supported-radii}
r(B_P)\le r(B_Q).
\end{equation}
Indeed, suppose towards a contradiction that $r(B_P)>r(B_Q)$. Let ${\Y}_0\in P\cap B_Q$. Then
\[
d_{p}({\Y}_0,{\X}_P)\le\diam_{p}(P)=\frac{r(B_P)}{3},
\]
so ${\Y}_0\in\frac13B_P$. Also, ${\X}_Q\in B_{p}({\Y}_0,r(B_Q))$. Thus, by \eqref{eq:ball-inclusion},
\[
B_Q\subset B_{p}({\Y}_0,2r(B_Q))
\subset B_{p}\Bigl({\X}_P,2r(B_Q)+\frac{r(B_P)}{3}\Bigr)\subset3B_P.
\]
It follows that $Q\subset B_Q\subset3B_P$, which contradicts the fact that $Q\in\mathcal G$, see \eqref{e:definition-calG}. We see also that 
\begin{align}\label{e:P-B_Q}
    P \subseteq 2B_Q \quad \mbox{ and } \quad \ell(P) \lesssim \ell(Q). 
\end{align}
Indeed, let ${\Y}_0\in P\cap B_Q$. Every ${\Y}\in P$ satisfies
$d_{p}({\Y},{\Y}_0)\le\diam_{p}(P)$. So, by \eqref{eq:ball-inclusion} and
\eqref{eq:supported-radii},
\[
P\subset B_{p}({\Y}_0,\diam_{p}(P))
\subset B_{p}\bigl({\X}_Q,r(B_Q)+\diam_{p}(P)\bigr)
\subset B_{p}\Bigl({\X}_Q,\frac{4}{3}r(B_Q)\Bigr).
\]
Thus $P\subset2B_Q$. Furthermore, 
\[
\ell(P)\le\diam_{p}(P)\le\diam_{p}(Q)\lesssim\ell(Q).
\]

Since the cubes in $\mathcal B$ are disjoint and $\nu(P)\lesssim\lambda\sigma(P)$ (see \eqref{e:nu-less-than-sigma}),
AD-regularity and \eqref{e:P-B_Q} give
\[
\sum_{P\in\mathcal B:\,P\cap B_Q\ne\varnothing}
\frac{\ell(P)\,\nu(P)}{\ell(Q)^{n+2}}
\lesssim\frac{\lambda}{\ell(Q)^{n+1}}
\sum_{P\in\mathcal B:\,P\subset2B_Q}\sigma(P)
\le\frac{\lambda\,\sigma(2B_Q)}{\ell(Q)^{n+1}}
\lesssim\lambda.
\]

Thus, the inner sum in \eqref{eq:supported-goal} is at most $C\lambda$  and hence its square is at most $C\lambda$ mulitplied by the sum itself. This makes it so that it is suffices to show
\begin{equation}\label{eq:supported-goal2}
\sum_{Q\in\mathcal G\cap\mathcal D_\Gamma^0:\,{\X}\in Q}
\sum_{P\in\mathcal B:\,P\cap B_Q\ne\varnothing}
\frac{\ell(P)\,\nu(P)}{\ell(Q)^{n+2}}<\infty
\qquad\text{for }\sigma\text{-a.e. }{\X}\in\Gamma.
\end{equation}

To prove \eqref{eq:supported-goal2}, we integrate over $\Gamma$. Since
$\sigma(Q)\approx\ell(Q)^{n+1}$, Tonelli's theorem gives
\[
\int_\Gamma\sum_{Q\in\mathcal G\cap\mathcal D_\Gamma^0:\,{\X}\in Q}
\sum_{P\in\mathcal B:\,P\cap B_Q\ne\varnothing}
\frac{\ell(P)\,\nu(P)}{\ell(Q)^{n+2}}\,d\sigma({\X})
\approx\sum_{P\in\mathcal B}\nu(P)
\sum_{Q\in\mathcal G\cap\mathcal D_\Gamma^0:\,P\cap B_Q\ne\varnothing}
\frac{\ell(P)}{\ell(Q)}.
\]
Fix $P\in\mathcal B$. By \eqref{e:P-B_Q}, the cubes $Q$ in
the last sum satisfy $\ell(Q)\gtrsim\ell(P)$. For each side length,
there are only boundedly many of them. Indeed, their balls $B_Q$ meet $P$, so
\[
d_{p}({\X}_Q,{\X}_P)\le2r(B_Q)+\diam_{p}(P)\lesssim\ell(Q).
\]
By AD-regularity, only boundedly many disjoint cubes of side $\ell$
have centers in a ball of radius $C\ell$.
The inner sum is therefore bounded by a geometric series, which gives 
\[
\sum_Q\frac{\ell(P)}{\ell(Q)}\lesssim1.
\]
Thus the integral is at most $C\sum_{P\in\mathcal B}\nu(P)\le C\|\nu\|<\infty$.
This proves
\eqref{eq:supported-goal2}, and the lemma follows.
\end{proof}

\subsection{The approximating measure for the general case}
\label{the-approximating-measure}

Let us now proceed and tackle the general case of Main Lemma~\ref{thm:graph-alpha}, that is, the case that the measure $\mu$ is not supported on $\Gamma$. In this case, we will need a way of approximating the parts of the measure which lie away from $\Gamma$. This approximation, which will be denoted by $\eta$, is defined in this section. \\ 

In order to define the approximation$\eta$ we will, as in \cite[Section~2.3]{tolsa2015characterization}, make use of a Whitney decomposition of
$\Gamma^c=\mathbb R^{n+1}\setminus\Gamma$. 

\begin{notation}[Whitney decomposition of $\Gamma^c$]
    We consider a family
    $\mathcal W$ of disjoint parabolic dyadic cubes\footnote{
        The parabolic dyadic cubes of $\mathbb R^{n+1}$ are the
        half-open cubes with spatial side $\ell=2^{-j}$ and time side $\ell^2$,
        $j\in\mathbb Z$. For $\theta>0$, $\theta P$ is the cube with the same center
        ${\X}_P$ as $P$, spatial side $\theta\ell(P)$ and time side
        $\theta^2\ell(P)^2$.
    } which cover $\Gamma^c$, and
    a constant $R_{\mathcal W}>20$ such that, for every $P\in\mathcal W$,
    \begin{enumerate}
    \item\label{it:whitney-1} $10P\subset\Gamma^c$,
    \item\label{it:whitney-2} $R_{\mathcal W}P\cap\Gamma\ne\varnothing$,
    \end{enumerate}
    and
    \begin{equation}\label{eq:whitney-small}
    \diam_{p}(P)<\dist_{p}(P,\Gamma).
    \end{equation}
    These are the \textit{Whitney cubes}. Note that, if $P\in\mathcal W$ meets $B_{p}({\Y},r)$ at a point ${\Y}_1$, with
    ${\Y}\in\Gamma$, then $\diam_{p}(P)<\dist_{p}(P,\Gamma)\le d_{p}({\Y}_1,{\Y})\le2r$. Thus
    $P\subset B_{p}({\Y}_1,2r)$, and so
    \begin{equation}\label{eq:whitney-inclusion}
    P\subset B_{p}({\Y},3r).
    \end{equation}
\end{notation}

The construction of this family of Whitney cubes is well-known in the Euclidean setting, and adapts without further problems to the current parabolic scenario.\footnote{
    The cubes can be constructed, for instance, as follows. Take the maximal parabolic dyadic cubes $P$ with
    $10P\cap\Gamma=\varnothing$. They cover $\Gamma^c$ since $\Gamma$ is
    closed and nonempty. For the parent $\widehat P$ of $P$, $10\widehat P$
    meets $\Gamma$ and $10\widehat P\subset22P$, so (2) holds with
    $R_{\mathcal W}=22$. By (1), $\dist_{p}(P,\Gamma)\ge\frac92\ell(P)$.
    So, replacing each cube, if necessary, by its descendants of a fixed
    generation, depending only on $n$, we get \eqref{eq:whitney-small}, with a larger
    $R_{\mathcal W}$ depending only on $n$. 
} \\

We denote the projection onto $\Gamma$ by
\[\Pi_\Gamma(x',x_n,t) := (x',A(x',t),t).\]

\begin{lem}
    The projection defined right above satisfies 
    \begin{align}\label{eq:projection}
        d_{p}({\Y},\Pi_\Gamma({\Y})) \leq (1+2L) \, \dist_{p}({\Y},\Gamma). 
    \end{align}
\end{lem}
\begin{proof}
    For ${\Y}=(y',y_n,t)$ and ${\bZ}=(u,A(u,\tau),\tau)\in\Gamma$,
    \[
    |y_n-A(y',t)|\le|y_n-A(u,\tau)|+L\bigl(|u-y'|+|\tau-t|^{1/2}\bigr)
    \le(1+2L)\,d_{p}({\Y},{\bZ}).
    \]
    Since $d_{p}({\Y},\Pi_\Gamma({\Y}))=|y_n-A(y',t)|$, taking the infimum over ${\bZ}$
    gives \eqref{eq:projection}.
\end{proof}

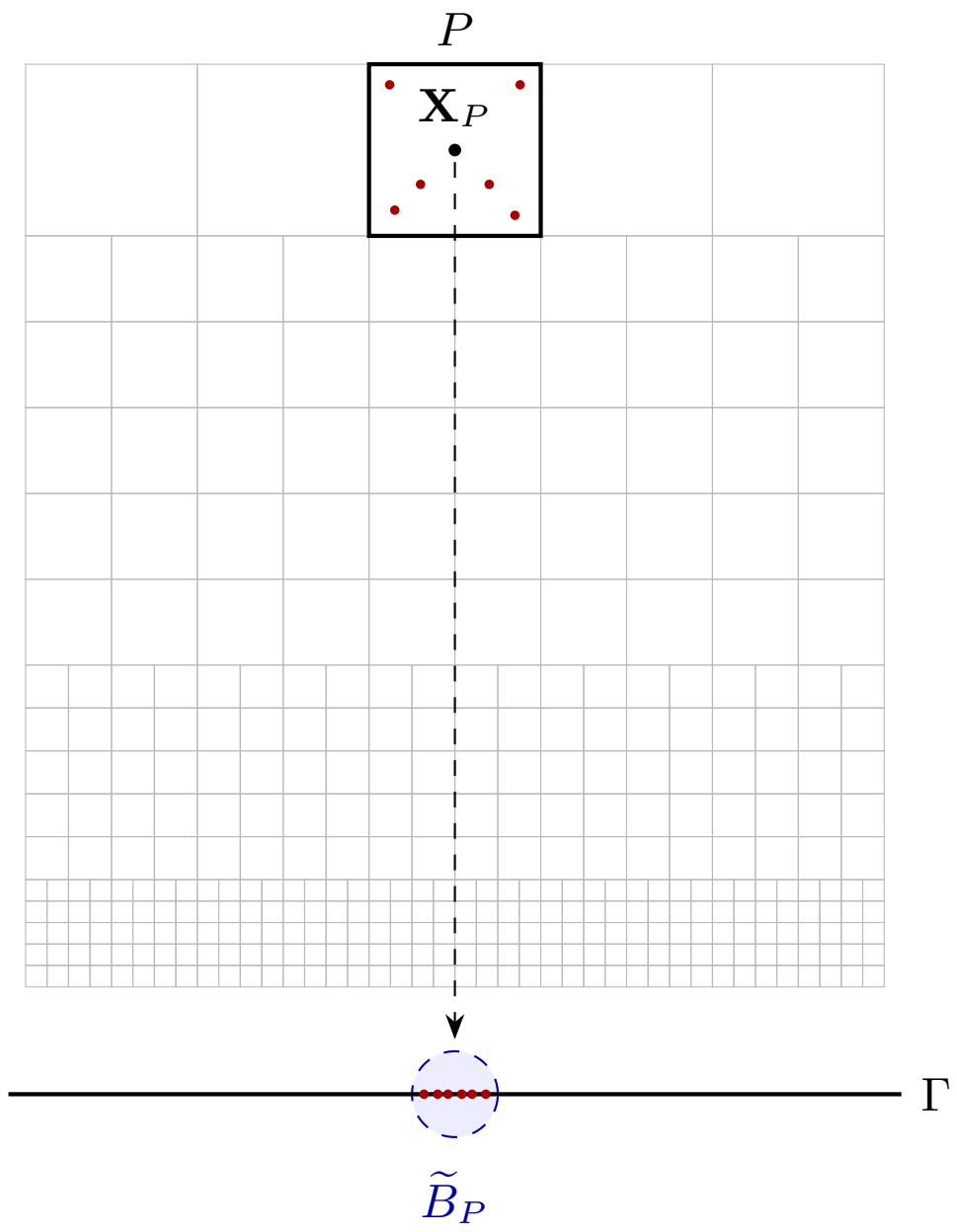
\begin{figure}[h!]
\centering
\begin{tikzpicture}[x=1.15cm,y=1.15cm,>=Stealth,font=\small]
  \foreach \h/\lo/\nr in {0.125/0.625/4,0.25/1.25/4,0.5/2.5/4,1/5/0}{
    \pgfmathtruncatemacro{\nc}{5/\h-1}
    \foreach \j in {0,...,\nr}{
      \foreach \i in {0,...,\nc}{
        \draw[gray!55,very thin] ({2+\i*\h},{\lo+\j*\h}) rectangle ({2+(\i+1)*\h},{\lo+(\j+1)*\h});
      }
    }
  }
  \fill[blue!7] (4.5,0) circle (0.25);
  \draw[blue!60!black,dashed] (4.5,0) circle (0.25);
  \draw[thick] (1.9,0) -- (7.1,0) node[right] {$\Gamma$};
  \draw[thick] (4,5) rectangle (5,6);
  \node[above] at (4.5,6) {$P$};
  \fill (4.5,5.5) circle (1.2pt);
  \node[above=1pt] at (4.5,5.5) {${\X}_P$};
  \draw[->,dashed] (4.5,5.43) -- (4.5,0.32);
  \node[blue!60!black,below=3pt] at (4.5,-0.25) {$\widetilde B_P$};
  \foreach \p in {(4.15,5.15),(4.3,5.3),(4.12,5.88),(4.85,5.12),(4.7,5.3),(4.88,5.88)}{\fill[red!65!black] \p circle (0.9pt);}
  \foreach \x in {4.32,4.4,4.46,4.54,4.6,4.68}{\fill[red!65!black] (\x,0) circle (0.9pt);}
\end{tikzpicture}
\caption{Whitney cubes for $\Gamma=\{x_2=0\}$, in the fixed-time slice
through ${\X}_P$ ($n=2$). The cubes of side $h$ fill the strip
$5h\le x_2<10h$. The red dots illustrate $\mu\lfloor P$ and $g_P\sigma$,
which have the same mass.}
\label{fig:whitney}
\end{figure}

\begin{notation}[The approximating measure $\eta$]
    To define $\eta$, we use the Whitney decomposition and several additional localization tools:
\begin{itemize}
    \item First, we move the mass of each cube $P\in\mathcal W$ to a ball centered on
    $\Gamma$, near $P$ and of comparable size, with a Lipschitz density with
    respect to $\sigma$. Indeed, to each $P\in\mathcal W$, we associate a ball 
    \[\widetilde B_P := B_{p}(\widetilde {\X}_P,\frac{\ell(P)}{4}), \quad \text{ where } \quad 
    \widetilde {\X}_P := \Pi_\Gamma({\X}_P).\] 
    
    \item By \eqref{it:whitney-2}, every point of $R_{\mathcal W}P$ is within $C\ell(P)$. So $\dist_{p}({\X}_P,\Gamma)\lesssim\ell(P)$. Thus, by \eqref{eq:projection},
    \[
    \widetilde B_P\subset B_{p}({\X}_P,C\ell(P)).
    \]
    
    \item Let us also fix a Lipschitz function $\psi:[0,\infty)\to[0,1]$ with $\psi=1$ on
    $[0,\frac12]$ and $\psi=0$ on $[1,\infty)$. For $P\in\mathcal W$ and
    ${\Y}\in\mathbb R^{n+1}$, we defined the following rescaled and normalized cut-offs:
    \[
    \varphi_P({\Y}) := \psi\Bigl(\frac{d_{p}({\Y},\widetilde {\X}_P)}
    {r(\widetilde B_P)}\Bigr),
    \qquad
    g_P := \frac{\mu(P)}{\|\varphi_P\|_{L^1(\sigma)}}\,\varphi_P.
    \]

    \item Since $\sigma(B_{p}(\widetilde {\X}_P,\frac{\ell(P)}{16}))\le
    \|\varphi_P\|_{L^1(\sigma)}\le\sigma(\widetilde B_P)$, AD-regularity gives
    $\|\varphi_P\|_{L^1(\sigma)}\approx\ell(P)^{n+1}$, and so
    \begin{equation}\label{eq:bump-bounds}
    \int g_P\,d\sigma=\mu(P),
    \qquad
    \|g_P\|_\infty\lesssim\frac{\mu(P)}{\ell(P)^{n+1}},
    \qquad
    \Lip_{d_{p}}g_P\lesssim\frac{\mu(P)}{\ell(P)^{n+2}}.
    \end{equation}
    Moreover, each $g_P$ vanishes outside $\widetilde B_P$. In Figure~\ref{fig:whitney}, we see the Whitney
    cubes for a flat $\Gamma$ and, for a particular choice of Whiteny cube $P$, the corresponding bal $\widetilde B_P$. 
    
    \item Lastly, we define the measure $\eta$, supported on $\Gamma$, which we will later use to approximate $\mu\lfloor\Gamma^c$ in the proof of the Main Lemma, by
    \begin{align}\label{e:def-eta}
    \eta\coloneqq\sum_{P\in\mathcal W}g_P\,\sigma,
    \end{align}
    noting that $\|\eta\|=\mu(\Gamma^c)$.
\end{itemize}
\end{notation}

A Whitney cube $P$ may miss $2B_Q$, for a $\Gamma$-cube $Q$, while its bump
$g_P$ is nonzero in $B_Q$. To control the contribution of these cubes, we
introduce the following function, similarly as in \cite[Section~2.3]{tolsa2015characterization}.

\begin{notation}[The function $G_{\theta_0}$]
    For a constant $\theta_0>10$, which we choose in
    Lemma~\ref{lem:five-terms} below, we consider
    \[
    G_{\theta_0}({\X}) := \sum_{P\in\mathcal W}\frac{\mu(P)}{\ell(P)^{n+1}}\,
    \mathbf 1_{\theta_0\widetilde B_P}({\X}).
    \]
    Note that, by Tonelli's theorem and AD-regularity,
    \begin{align}\label{e:G-finite}
    \|G_{\theta_0}\|_{L^1(\sigma)}\lesssim \theta_0^{n+1}\sum_{P\in\mathcal W}\mu(P)\le
    \theta_0^{n+1}\|\mu\|,
    \end{align}
    and thus $G_{\theta_0}({\X})<\infty$ for $\sigma$-a.e. ${\X}\in\Gamma$.
\end{notation}

%%%%%%%%%%%%%%%%%%%%%%%%%%%%%%%%%%%%%%%%%%%%%%%%%%%%%%%%%%%%%%%%%%%%%%%%%%
%%%%%%%%%%%%%%%%%%%%%%%%%%%%%%%%%%%%%%%%%%%%%%%%%%%%%%%%%%%%%%%%%%%%%%%%%%
%%%%%%%%%%%%%%%%%%%%%%%%%%%%%%%%%%%%%%%%%%%%%%%%%%%%%%%%%%%%%%%%%%%%%%%%%%
%%%%%%%%%%%%%%%%%%%%%%%%%%%%%%%%%%%%%%%%%%%%%%%%%%%%%%%%%%%%%%%%%%%%%%%%%%
%%%%%%%%%%%%%%%%%%%%%%%%%%%%%%%%%%%%%%%%%%%%%%%%%%%%%%%%%%%%%%%%%%%%%%%%%%

\subsection{The \texorpdfstring{$\alpha$}{alpha} coefficients
of the \texorpdfstring{$\Gamma$}{Gamma}-cubes}
\label{the-good-cubes}

The next lemma is the parabolic version of \cite[Lemma~2.4]{tolsa2015characterization}. It gives us a way to estimate the $\alpha$-numbers for $\mu$ in terms of the quantities $\eta$ and $G_{\theta_0}$, from the previous section, with some error terms which we will later be able to control.  

\begin{lem}[Control of alphas in terms of quantities involving $\mu \llcorner \Gamma, \eta$ and $G_{\theta_0}$]  \label{lem:five-terms}
    If $\theta_0$ is big enough, depending only on $n$ and $L$, every $Q\in\mathcal D_\Gamma$ satisfies 
    \begin{align}\label{e:five-terms}
    \begin{split}
    \alpha_\mu(Q)\lesssim{}&
    \alpha_{\mu\lfloor\Gamma,V_Q}(Q)
    +\alpha_{\eta,V_Q}(Q)
    +\int_{6B_Q}
    \frac{\dist_{p}({\Y},\Gamma)}{\ell(Q)^{n+2}}\,d\mu({\Y})\\
    &+\inf_{{\X}\in Q}G_{\theta_0}({\X})\,\kappa_Q
    +\sum_{\substack{P\in\mathcal W:\,P\cap2B_Q=\varnothing,\\
    \widetilde B_P\cap B_Q\ne\varnothing}}
    \frac{\mu(P)\,\ell(Q)}{\ell(P)^{n+2}},
    \end{split}
    \end{align}
    with $\kappa_Q$ as in Proposition \ref{prop:beta-alpha}. 
\end{lem}

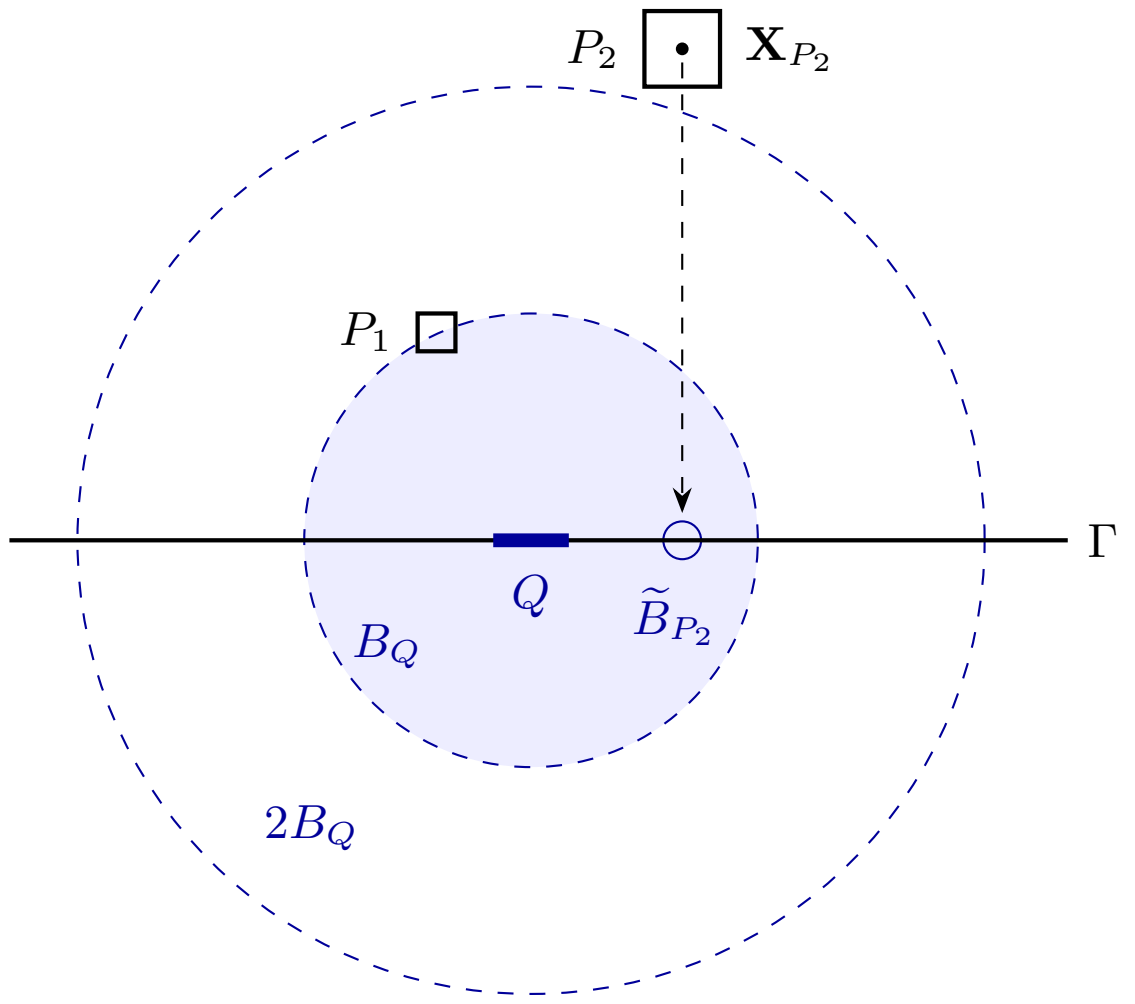
\begin{figure}[ht]
\centering
\begin{tikzpicture}[x=2cm,y=2cm,>=Stealth,font=\small]
  \fill[blue!7] (4.125,0) circle (0.75);
  \draw[blue!60!black,dashed] (4.125,0) circle (0.75);
  \draw[blue!60!black,dashed] (4.125,0) circle (1.5);
  \node[blue!60!black] at (3.65,-0.35) {$B_Q$};
  \node[blue!60!black] at (3.4,-0.95) {$2B_Q$};
  \draw[blue!60!black] (4.625,0) circle (0.0625);
  \draw[thick] (2.4,0) -- (5.9,0) node[right] {$\Gamma$};
  \draw[line width=2.6pt,blue!60!black] (4,0) -- (4.25,0);
  \node[blue!60!black,below=3pt] at (4.125,0) {$Q$};
  \node[blue!60!black,below=3pt] at (4.6,-0.04) {$\widetilde B_{P_2}$};
  \draw[thick] (3.75,0.625) rectangle (3.875,0.75);
  \node[left] at (3.73,0.69) {$P_1$};
  \draw[thick] (4.5,1.5) rectangle (4.75,1.75);
  \node[left] at (4.48,1.625) {$P_2$};
  \fill (4.625,1.625) circle (1.2pt);
  \node[right] at (4.77,1.625) {${\X}_{P_2}$};
  \draw[->,dashed] (4.625,1.58) -- (4.625,0.09);
\end{tikzpicture}
\caption{The two families in the proof of Lemma~\ref{lem:five-terms},
for $\Gamma=\{x_2=0\}$, in the fixed-time slice through ${\X}_{P_2}$
($n=2$). The cube $P_1$ meets $2B_Q$. The cube $P_2$ misses $2B_Q$, but
$\widetilde B_{P_2}$ meets $B_Q$.}
\label{fig:near-far}
\end{figure}

Before starting the proof, we give an indication of where the five terms come from. The first terms simply comes from the part of $\mu$ which is supported on $\Gamma$. The remaining four terms come from the part of $\mu$ supported away from $\Gamma$. Writing
$\mu\lfloor\Gamma^c=\eta+\sum_{P\in\mathcal W}(\mu\lfloor P-g_P\sigma)$, it is clear that the second term comes from the $\eta$ component. If
$P$ has non-empty intersection with $2B_Q$, then moving the mass of $P$ to $g_P\sigma$ costs at most
a constant times $\ell(P)\mu(P)$, and this gives the third term. If $P$
is disjoint from $2B_Q$, then $\int f\,d\mu\lfloor P=0$ whenever $\supp f\subset B_Q$,
but $\int fg_P\,d\sigma$ may be nonzero. Comparing $\sigma$ with $c_Q\lambda_{V_Q}$ gives the fourth term, and the
variation of these bumps on $B_Q$ gives the fifth. In
Figure~\ref{fig:near-far}, we see one cube of each kind. We now start the proof in earnest. 

\begin{proof}[Proof of Lemma \ref{lem:five-terms}]
Since $\mu=\mu\lfloor\Gamma+\mu\lfloor\Gamma^c$, the triangle inequality for
$\dist_{B_Q}$ gives
\[
\alpha_\mu(Q)\le\alpha_{\mu\lfloor\Gamma,V_Q}(Q)
+\alpha_{\mu\lfloor\Gamma^c,V_Q}(Q),
\]
so it remains to estimate the last coefficient. 

To this end, we consider a
$1$-Lipschitz function $f$ with $\supp f\subset B_Q$. Recalling the definition of $\eta$ from \eqref{e:def-eta} and using the fact that the Whitney
cubes $\mathcal{W}$ are disjoint and cover $\Gamma^c$, we can write 
\begin{align}
\int f\,d\mu\lfloor\Gamma^c
= \sum_{P \in \mathcal{W}} \int f\,d\mu\lfloor P =\sum_{P\in\mathcal W}\int f\,d(\mu\lfloor P-g_P\sigma)
+\int f\,d\eta.
\end{align}
Let $I_1$ be the family of the cubes of $\mathcal W$ which meet $2B_Q$,
and $I_2=\mathcal W\setminus I_1$. Fix $a, \bar a \ge0$.  Then
\begin{equation}\label{eq:off-graph-integral-bound}
\Bigl|\int f\,d\mu\lfloor\Gamma^c
-(a-\bar a)\int f\,d\lambda_{V_Q}\Bigr|
\le\Bigl|\int f\,d\eta-a\int f\,d\lambda_{V_Q}\Bigr|
+S_1+S_2.
\end{equation}
where
\begin{align}\label{e:eta-S_1-S_2}
S_1:=\sum_{P\in I_1}\Bigl|\int f\,d(\mu\lfloor P-g_P\sigma)\Bigr|,
\qquad
S_2:=\Bigl|\sum_{P\in I_2}\int f\,d(\mu\lfloor P-g_P\sigma)
+\bar a\int f\,d\lambda_{V_Q}\Bigr|.
\end{align}

Note that the number $a-\bar a$ may be negative. However, if $\nu$ is a positive
measure and $a'<0$, then
\[
\dist_{B_Q}(\nu,a'\lambda_{V_Q})\ge\dist_{B_Q}(\nu,0),
\]
as we see by testing with $|f|$ in place of $f$. Here $|f|$ is admissible
whenever $f$ is, and we use the same formula for $\dist_{B_Q}$ when the
coefficient is negative. All of this is to say that 
\begin{align}
    \alpha_{\mu\lfloor\Gamma^c,V_Q}(Q) &\leq r(B_Q)^{-n-2} \dist_{B_Q}(\mu\lfloor\Gamma^c , (a-\bar a) \lambda_{V_Q} ) \\
    & \leq r(B_Q)^{-n-2} \left[ \dist_{B_Q}(\eta , a \lambda_{V_Q}) + \sup_f S_1 + \sup_f S_2  \right], 
\end{align}
where the second inequality follows from \eqref{eq:off-graph-integral-bound} after taking the supremum over all admissible functions. 
Taking the
infimum over $a\ge0$, we get
\[
\alpha_{\mu\lfloor\Gamma^c,V_Q}(Q)\le\alpha_{\eta,V_Q}(Q)
+\frac{1}{r(B_Q)^{n+2}}\Bigl(\sup_fS_1+\sup_fS_2\Bigr).
\]
Note that the final two terms in the above inequality depend on $\bar a$. We will choose it later. 

First we deal with $S_1$. Let $P\in I_1$. The measures $\mu\lfloor P$ and
$g_P\sigma$ have the same mass $\mu(P)$, and both vanish outside
$B_{p}({\X}_P,C\ell(P))$. So, subtracting $f({\X}_P)$, we deduce that
\[
\Bigl|\int f\,d(\mu\lfloor P-g_P\sigma)\Bigr|
=\Bigl|\int\bigl(f-f({\X}_P)\bigr)\,d(\mu\lfloor P-g_P\sigma)\Bigr|
\lesssim\ell(P)\,\mu(P)
\le\int_P\dist_{p}({\Y},\Gamma)\,d\mu({\Y}),
\]
since $\dist_{p}({\Y},\Gamma)>\diam_{p}(P)\ge\ell(P)$ for
${\Y}\in P$ by \eqref{eq:whitney-small}. Since $P$ is assumed to meet $2B_Q$, \eqref{eq:whitney-inclusion} gives
$P\subset6B_Q$. Thus, since the cubes of $I_1$ are disjoint,
\[
S_1\lesssim\int_{6B_Q}\dist_{p}({\Y},\Gamma)\,d\mu({\Y}).
\]

Next we consider $S_2$. Since $\supp f \subset B_Q$, we have by the definition of $I_2$ that $\int f\,d\mu\lfloor P=0$ for all $P \in I_2$. Moreover, since each $g_P$ vanishes outside $\widetilde B_P$, we have $fg_P=0$ unless $\widetilde B_P$ meets $B_Q$. In particular, writing (as in \cite{tolsa2015characterization}) $P\sim Q$ if $P\in I_2$ and $\widetilde B_P\cap B_Q\ne\varnothing$, we have 
\begin{align}
    S_2 = \Bigl|\sum_{P\sim Q}\int fg_P\,d\sigma-\bar a\int f\,d\lambda_{V_Q}\Bigr|.
\end{align}
We first show that the cubes $P\sim Q$ are not much smaller than $Q$. Fix $P\sim Q$ for the moment. Since $P$ is disjoint from $2B_Q$, we have ${\X}_P\notin2B_Q$. On the other hand,
$\widetilde B_P\subset B_{p}({\X}_P,C\ell(P))$ has nonempty intersection with $B_Q$, and so, by \eqref{eq:ball-inclusion},
\[
{\X}_P\in B_{p}\bigl({\X}_Q,C\ell(P)+r(B_Q)\bigr).
\]
Hence $2r(B_Q)<C\ell(P)+r(B_Q)$. Thus,
\begin{align}\label{e:r(B_Q)-leq-ell(P)}
\ell(Q) \approx r(B_Q) \lesssim\ell(P).
\end{align}
Since the Whitney cubes are disjoint, this implies
\[
\sum_{P\sim Q}\frac{\mu(P)}{\ell(P)^{n+1}}\lesssim\frac{\|\mu\|}{\ell(Q)^{n+1}}.
\]
By \eqref{eq:bump-bounds}, \eqref{e:r(B_Q)-leq-ell(P)} and the
preceding bound, we have
$\sum_{P\sim Q}\|g_P\|_\infty<\infty$ and
$\sum_{P\sim Q}\Lip_{d_{p}}g_P<\infty$. Thus the series
\[G_Q\coloneqq \sum_{P\sim Q}g_P\] 
converges uniformly on $\mathbb R^{n+1}$, and $G_Q$
is Lipschitz. Let $c_Q$ be as in Proposition~\ref{prop:beta-alpha}. We can now estimate 
\[
\begin{aligned}
S_2&=\Bigl|\int fG_Q\,d\sigma-\bar a\int f\,d\lambda_{V_Q}\Bigr|\\
&\le\Bigl|\int fG_Q\,d(\sigma-c_Q\lambda_{V_Q})\Bigr|
+\int|f|\,\bigl| c_QG_Q-\bar a\bigr|\,d\lambda_{V_Q}
\eqqcolon S_{2,1}+S_{2,2}.
\end{aligned}
\]

We will estimate $S_{2,1}$ in terms of $\kappa_Q$ (recall Proposition \ref{prop:beta-alpha}). First, we need to bound the Lipschitz constant
of $F=fG_Q$. Notice that $\supp F\subset B_Q$. Since $f$ is $1$-Lipschitz and vanishes outside $B_Q$, we have $|f|\le r(B_Q)$, and $r(B_Q)\lesssim\ell(P)$ for $P\sim Q$ by \eqref{e:r(B_Q)-leq-ell(P)}. Then, 
\eqref{eq:bump-bounds} gives
\[
\Lip_{d_{p}}F
\le\|G_Q\|_\infty+r(B_Q)\Lip_{d_{p}}G_Q
\le\sum_{P\sim Q}\bigl(\|g_P\|_\infty+ r(B_Q)\Lip_{d_{p}}g_P\bigr)
\lesssim\sum_{P\sim Q}\frac{\mu(P)}{\ell(P)^{n+1}}.
\]
If $P\sim Q$ and ${\X}\in Q$, then
\[
d_{p}({\X},\widetilde {\X}_P)\le4r(B_Q)+2r(\widetilde B_P)\lesssim\ell(P),
\]
so ${\X}\in \theta_0\widetilde B_P$ if $\theta_0$ is big enough, depending only on $n$
and $L$. Indeed, the constant in the last
inequality comes from
\eqref{e:r(B_Q)-leq-ell(P)}, the Whitney property (2) and \eqref{eq:projection}. Therefore,
\begin{equation}\label{eq:far-bumps}
\sum_{P\sim Q}\frac{\mu(P)}{\ell(P)^{n+1}}\le\inf_{{\X}\in Q}G_{\theta_0}({\X}).
\end{equation}
If $F\ne0$, then $\frac{F}{\Lip_{d_{p}}F}$ is admissible in
the definition of $\dist_{B_Q}$, and if $F=0$, then $ S_{2,1}=0$. So, recalling the definition of $\kappa_Q$ from Proposition \ref{prop:beta-alpha}, we have 
\[
S_{2,1}\le\Lip_{d_{p}}F\,\dist_{B_Q}(\sigma,c_Q\lambda_{V_Q})
\lesssim\inf_{{\X}\in Q}G_{\theta_0}({\X})\,\kappa_Q\,r(B_Q)^{n+2}.
\]

Finally we turn our attention to $S_{2,2}$. On $B_Q$, we compare $G_Q$ with
its value at ${\X}_Q$. That is, we choose
$\bar a= c_QG_Q({\X}_Q)\ge0$. 
For
${\Y}\in B_Q$, \eqref{eq:bump-bounds} gives
\[
|G_Q({\Y})-G_Q({\X}_Q)|
\le2r(B_Q)\sum_{P\sim Q}\Lip_{d_{p}}g_P
\lesssim\sum_{P\sim Q}\frac{\mu(P)\,\ell(Q)}{\ell(P)^{n+2}}.
\]
Since $c_Q\le C$, $|f|\le r(B_Q)$ and $\lambda_{V_Q}(B_Q)\lesssim\ell(Q)^{n+1}$,
we get
\[
S_{2,2}\lesssim\sum_{P\sim Q}\frac{\mu(P)\,\ell(Q)^{n+3}}{\ell(P)^{n+2}}.
\]
Dividing by
$r(B_Q)^{n+2}\approx\ell(Q)^{n+2}$ gives the lemma.
\end{proof}

\subsection{Proof of the Main Lemma}
\label{the-main-lemma-proof}

After all the work in the previous subsections, we are ready to prove \eqref{eq:main-cubes}, and thus the Main Lemma \ref{thm:graph-alpha}. In the proof, we will need the maximal function
\[
M_{n+1}\mu({\X})\coloneqq \sup_{r>0}\frac{\mu(B_{p}({\X},r))}{r^{n+1}}.
\]
By the $5r$-covering lemma and the upper AD-regularity of $\sigma$,
\[
\sigma \big(\{{\X}\in\Gamma:M_{n+1}\mu({\X})>\lambda\}\big)\lesssim\frac{\|\mu\|}{\lambda}
\qquad \forall \, \lambda>0,
\]
so $M_{n+1}\mu({\X})<\infty$ for $\sigma$-a.e. ${\X}\in\Gamma$.

\begin{proof}[\textbf{Proof of the Main Lemma~\ref{thm:graph-alpha}}]
As already explained, we only need to verify \eqref{eq:main-cubes}.
For fixed ${\X}\in\Gamma$, we let $\mathcal D_\Gamma^0(\X)$ denote the set of cubes in $\cubes_\Gamma^0$ such that $\X \in Q$. By Lemma~\ref{lem:five-terms}, it
suffices to show that, for each of the five terms on the right-hand side of \eqref{e:five-terms}, the sum over $\cubes_\Gamma^0(\X)$ of this term squared
is finite for $\sigma$-a.e. ${\X}\in\Gamma$. 

For the first two terms, the $\sigma$-a.e. finiteness
follows from Lemma~\ref{lem:supported} applied to $\mu\lfloor\Gamma$ and
$\eta$. For the fourth term, we note that 
\[
\sum_{Q \in \cubes_\Gamma^0(\X)}\inf_{{\Y}\in Q}G_{\theta_0}({\Y})^2\kappa_Q^2
\le G_{\theta_0}({\X})^2\sum_{Q \in \cubes_\Gamma^0(\X)}\kappa_Q^2
\]
for every $\X \in \Gamma$. 
Recalling from \eqref{e:G-finite} and Proposition~\ref{prop:graph-measure} (as
in the first part of the proof of Lemma~\ref{lem:supported}), respectively, that $G_{\theta_0}({\X})$ and $\sum_{Q \in \cubes_\Gamma^0(\X)} \kappa_Q^2$ are finite for
$\sigma$-a.e. ${\X}$, we have control over the fourth term. 

Next we consider the fifth term. The cubes in the sum there are the cubes $P\sim Q$ as in
the proof of Lemma~\ref{lem:five-terms}. Thus, by \eqref{e:r(B_Q)-leq-ell(P)}, $\ell(Q)\lesssim\ell(P)$. It now follows from \eqref{eq:far-bumps}
\[
\sum_{\substack{P\in\mathcal W:\,P\cap2B_Q=\varnothing,\\
\widetilde B_P\cap B_Q\ne\varnothing}}
\frac{\mu(P)\,\ell(Q)}{\ell(P)^{n+2}}
\lesssim\sum_{\substack{P\in\mathcal W:\,P\cap2B_Q=\varnothing,\\
\widetilde B_P\cap B_Q\ne\varnothing}}
\frac{\mu(P)}{\ell(P)^{n+1}}
\le G_{\theta_0}({\X}).
\]
Thus, square sum of this term is at most $CG_{\theta_0}({\X})$ times the term itself, and it
suffices to show that
\[
\sum_{Q \in \cubes_\Gamma^0(\X)}  \sum_{\substack{P\in\mathcal W:\,P\cap2B_Q=\varnothing,\\
\widetilde B_P\cap B_Q\ne\varnothing}}
\frac{\mu(P)\,\ell(Q)}{\ell(P)^{n+2}}<\infty
\qquad\text{for }\sigma\text{-a.e. }{\X}\in\Gamma.
\]
By Tonelli's theorem and $\sigma(Q)\approx\ell(Q)^{n+1}$, the integral of the
left-hand side over $\Gamma$ is comparable to
\[
\sum_{P\in\mathcal W}\mu(P)
\sum_{\substack{Q\in\mathcal D_\Gamma^0:\,P\cap2B_Q=\varnothing,\\
\widetilde B_P\cap B_Q\ne\varnothing}}
\frac{\ell(Q)^{n+2}}{\ell(P)^{n+2}}.
\]
For fixed $P$, the cubes $Q$ in the inner sum satisfy
$\ell(Q)\lesssim\ell(P)$ as these are the cubes with $P\sim Q$ so from \eqref{e:r(B_Q)-leq-ell(P)} we get $\ell(Q)\approx r(B_Q)\lesssim \ell(P)$.
Since $B_Q$ meets $\widetilde B_P$ and
$r(B_Q)\lesssim\ell(P)$, we also have $Q\subset C\widetilde B_P$. Those of one side
length $\ell$ are disjoint, so they contribute at most
\[
\frac{C\sigma(C\widetilde B_P)\,\ell}{\ell(P)^{n+2}}\lesssim\frac{\ell}{\ell(P)}.
\]
Summing over the dyadic $\ell\lesssim\ell(P)$, the integral is at most
$C\|\mu\|$.

Finally we deal with the third term. By the Cauchy--Schwarz inequality,
\[
\biggl(\int_{6B_Q}\frac{\dist_{p}({\Y},\Gamma)}{\ell(Q)^{n+2}}
\,d\mu({\Y})\biggr)^2
\le\frac{\mu(6B_Q)}{\ell(Q)^{n+1}}
\int_{6B_Q}\frac{\dist_{p}({\Y},\Gamma)^2}{\ell(Q)^{n+3}}
\,d\mu({\Y}).
\]
For ${\X}\in Q$, we have $6B_Q\subset B_{p}({\X},7r(B_Q))$, and so
$\frac{\mu(6B_Q)}{\ell(Q)^{n+1}}\lesssim M_{n+1}\mu({\X})$. Since $M_{n+1}\mu({\X})<\infty$
for $\sigma$-a.e. ${\X}$, it suffices to show that
\[
\sum_{Q \in \cubes_\Gamma^0(\X)}\int_{6B_Q}
\frac{\dist_{p}({\Y},\Gamma)^2}{\ell(Q)^{n+3}}\,d\mu({\Y})<\infty
\qquad\text{for }\sigma\text{-a.e. }{\X}\in\Gamma.
\]
Since points of $\Gamma$ do not contribute to the above integral, Tonelli's theorem and
$\sigma(Q)\lesssim\ell(Q)^{n+1}$ bound the integral of the left-hand side over
$\Gamma$ by
\[
C\sum_{Q\in\mathcal D_\Gamma^0}\int_{6B_Q}
\frac{\dist_{p}({\Y},\Gamma)^2}{\ell(Q)^2}\,d\mu({\Y})
=C\int_{\Gamma^c}\dist_{p}({\Y},\Gamma)^2
\sum_{Q\in\mathcal D_\Gamma^0:\,{\Y}\in6B_Q}\frac{1}{\ell(Q)^2}\,d\mu({\Y}).
\]
If ${\Y}\in\Gamma^c$ and ${\Y}\in6B_Q$, then
\[
\dist_{p}({\Y},\Gamma)\le d_{p}({\Y},{\X}_Q)\le12r(B_Q)\lesssim\ell(Q).
\]
The centers of these cubes $Q$ with $\ell(Q)=\ell$ lie in
$B_{p}({\Y},C\ell)$. The sets
$\Gamma\cap B_{p}({\X}_Q,\frac{\ell}{4})\subset Q$ are disjoint with
$\sigma$-measure comparable to $\ell^{n+1}$, so there are boundedly many such
$Q$. Thus, summing over the dyadic $\ell\gtrsim\dist_{p}({\Y},\Gamma)$, we get
\[
\sum_{Q\in\mathcal D_\Gamma^0:\,{\Y}\in6B_Q}\frac{1}{\ell(Q)^2}
\lesssim\frac{1}{\dist_{p}({\Y},\Gamma)^2},
\]
and the last integral is at most $C\|\mu\|$.
Hence the five sums of squares are finite for $\sigma$-a.e. ${\X}\in\Gamma$, which
proves \eqref{eq:main-cubes}.
\end{proof}

%%%%%%%%%%%%%%%%%%%%%%%%%%%%%%%%%%%%%%%%%%%%%%%%%%%%%%%%%%%%%%%%%%%%%%
%%%%%%%%%%%%%%%%%%%%%%%%%%%%%%%%%%%%%%%%%%%%%%%%%%%%%%%%%%%%%%%%%%%%%%
%%%%%%%%%%%%%%%%%%%%%%%%%%%%%%%%%%%%%%%%%%%%%%%%%%%%%%%%%%%%%%%%%%%%%%
%%%%%%%%%%%%%%%%%%%%%%%%%%%%%%%%%%%%%%%%%%%%%%%%%%%%%%%%%%%%%%%%%%%%%%
%%%%%%%%%%%%%%%%%%%%%%%%%%%%%%%%%%%%%%%%%%%%%%%%%%%%%%%%%%%%%%%%%%%%%%

\newpage

\section{The proof that rectifiability implies finite Jones function}
\label{proof-of-the-main-theorem}

In this section, we prove the implication of Theorem~\ref{thm:main} that parabolic rectifiability of $\mu$ implies finiteness $\mu$-a.e. of the Jones' square function \eqref{e:finite-jones-beta-main-thm}. 

Although the bound demanded by \eqref{e:finite-jones-beta-main-thm} (and Theorem~\ref{thm:main}) is about the square function of the beta numbers, we will actually be able to control a larger square function, which is the one associated to the alpha numbers. In the following result, we show that this is indeed the case: controlling the alpha square function is harder than the beta one. For this relationship to be true, we need to be careful with the densities. Concretely, the lower and upper $(n+1)$-densities are
\[
\Theta^{n+1}_*({\X},\mu)
:=\liminf_{r\downarrow0}\frac{\mu(B_{p}({\X},r))}{r^{n+1}},
\qquad
\Theta^{n+1,*}({\X},\mu)
:=\limsup_{r\downarrow0}\frac{\mu(B_{p}({\X},r))}{r^{n+1}}.
\]
Note that the upper density is the one we already defined in Theorem~\ref{thm:main}. Moreover, along the section, $\alpha_\mu(Q)$ and $V_Q$ refer to ambient cubes
and the balls $B_Q^+$ defined below, not to the $\Gamma$-cubes
of Section~\ref{sec:graph}. 

\begin{prop}[Control of alpha numbers implies control of beta numbers]  \label{prop:alpha-beta}
Let $\mu$ be a finite Borel measure on $\mathbb R^{n+1}$. Suppose that, on a
measurable set $F$, the lower and upper $(n+1)$-densities of $\mu$ are
positive and finite and
\begin{align}\label{e:finite-jones-alpha}
\int_0^1\alpha_\mu({\X},r)^2\,\frac{dr}{r}<\infty
\quad\text{for }\mu\text{-a.e. }{\X}\in F.
\end{align}
Then
\begin{align}\label{e:finite-jones-beta}
\int_0^1\beta_\mu({\X},r)^2\,\frac{dr}{r}<\infty
\quad\text{for }\mu\text{-a.e. }{\X}\in F.
\end{align}
\end{prop}

Before proving Proposition \ref{prop:alpha-beta}, let us see why it implies the corresponding implication in Theorem \ref{thm:main}.

\begin{proof} [\textbf{Proof of the implication in Theorem~\ref{thm:graph-alpha} assuming Proposition \ref{prop:alpha-beta}}]
Our goal is to show \eqref{e:finite-jones-beta-main-thm}. Therefore, we fix, for the remainder of the proof, a finite Borel measure $\mu$ which is assume to be parabolic rectifiable. 

By Corollary~\ref{cor:graph-alpha}, the
$\alpha$-square function of $\mu$, as in \eqref{e:finite-jones-alpha}, is finite $\mu$-a.e. Thus, if we can show that the densities of $\mu$ are positive and finite $\mu$-a.e., an application of Proposition \ref{prop:alpha-beta} finishes the proof of \eqref{e:finite-jones-beta-main-thm}. 
Consider a covering graph $\Gamma_j$. Set
$\sigma_j=\mathcal H_{p}^{n+1}\lfloor\Gamma_j$. Since
$\mu\ll\mathcal H_{p}^{n+1}$ (recall that this is part of the definition of rectifiability), we can write
$\mu\lfloor\Gamma_j=f_j\sigma_j$. Put
$\nu_j=\mu\lfloor(\mathbb R^{n+1}\setminus\Gamma_j)$.
Since $\sigma_j$ is doubling, the Lebesgue differentiation theorem
applies to $f_j$.
We claim that
\[
\lim_{r\downarrow0}
\frac{\nu_j(B_{p}({\X},r))}{\sigma_j(B_{p}({\X},r))}=0
\qquad\text{for }\sigma_j\text{-a.e. }{\X}\in\Gamma_j.
\]
For $\lambda>0$, let $E_\lambda\subset\Gamma_j$ be the set where the upper limit exceeds
$\lambda$. For $\delta>0$, its points admit arbitrarily small balls
$B=B_{p}({\X},r)$ with $r<\delta$ and
$\nu_j(B)>\lambda\sigma_j(B)$. By the fivefold covering lemma,
we choose disjoint such balls $B_i$ whose fivefold dilates
cover $E_\lambda$. By AD-regularity,
\begin{equation}\label{eq:main-singular-vitali}
\begin{aligned}
\sigma_j(E_\lambda)
\le C\sum_i\sigma_j(B_i)
\le\frac{C}{\lambda}\sum_i\nu_j(B_i) \le\frac{C}{\lambda}
\nu_j\bigl(\{{\Y}:d_{p}({\Y},\Gamma_j)\le2\delta\}\bigr).
\end{aligned}
\end{equation}
This tends to zero as $\delta\downarrow0$, since $\Gamma_j$ is closed,
$\nu_j$ is finite and $\nu_j(\Gamma_j)=0$. Taking rational
$\lambda>0$ proves the claim.
Thus
\[
\lim_{r\downarrow0}
\frac{\mu(B_{p}({\X},r))}{\sigma_j(B_{p}({\X},r))}
=f_j({\X})
\qquad\text{for }\sigma_j\text{-a.e. }{\X}\in\Gamma_j.
\]
Since $\mu\lfloor\Gamma_j\ll\sigma_j$, this also holds for $\mu$-a.e. ${\X}\in\Gamma_j$.
Now $0<f_j({\X})<\infty$ at $\mu$-a.e. point of $\Gamma_j$, so
AD-regularity of $\sigma_j$ gives the required density bounds.

Let $F$ be the subset of $\bigcup_j\Gamma_j$ where the lower and upper
$(n+1)$-densities of $\mu$ are positive and finite. It is measurable
by Appendix~\ref{app:measurability}. The graphs cover
$\mu$-almost all of $\mathbb R^{n+1}$, so $F$ has full $\mu$-measure.
Thus Proposition~\ref{prop:alpha-beta} proves Theorem~\ref{thm:main}.
\end{proof}

We now prove Proposition~\ref{prop:alpha-beta}, following the
stopping argument of \cite[Section~3]{tolsa2015characterization}. We first fix a set $K$
with uniform density and $\alpha$-square-function bounds and stop at
the maximal cubes whose triples miss $K$. The remaining good cubes
have uniform density bounds and square-summable $\alpha$-coefficients
(Lemmas~\ref{lem:stopping-cubes} and~\ref{lem:alpha-sum}).
Lemmas~\ref{lem:plane-comparison}--\ref{lem:alpha-beta-key} compare
the approximating planes across scales and estimate the distance to
$V_Q$. Lemma~\ref{lem:beta-sum} then gives the weighted dyadic
$\beta_2$ estimate, from which we deduce the continuous
square-function bound.

\subsection{Preliminaries} \label{subsec:preliminaries_partI}

\begin{notation}[Parabolic dyadic cubes]
    Let $\mathcal D_{n+1}$ be the family of parabolic dyadic cubes of $\mathbb R^{n+1}$
    from Section~\ref{the-approximating-measure}. These are not the cubes
    $\mathcal D$ of Section~\ref{sec:cubes}, which are adapted to $\mu$. Write ${\bZ}_Q$ for the center of
    $Q\in\mathcal D_{n+1}$, and put 
    \[D_n:=2+\sqrt n.\] Then
    \[
    \sup_{{\X}\in Q}d_{p}({\X},{\bZ}_Q)\le D_n\ell(Q),
    \qquad
    d_{p}({\bZ}_Q,{\bZ}_P)\le D_n\ell(P)
    \]
    when $Q$ is a child of $P$. As in
    Section~\ref{the-approximating-measure}, $3Q$ consists of $3^n\cdot9$
    cubes of the same generation, and
    \[
    \sup_{{\X}\in3Q}d_{p}({\X},{\bZ}_Q)\le3D_n\ell(Q),
    \qquad 3Q\subset3P, \qquad\text{if }Q\subset P.
    \]
    In particular, $3Q\subset B_{p}({\bZ}_Q,3D_n\ell(Q))$.
\end{notation}

\begin{notation}[Choice of constants $L_- , L_+, C_*, C_1$ and $C_2$] \label{not:constants}
    For the construction below, we need to introduce constants $L_- , L_+, C_*, C_1$ and $C_2$. We declare now that these constants are chosen such that the following inequalities hold:
    \[
    L_->3D_n+2,\qquad C_*>L_-+2,\qquad C_2>2C_*+4D_n.
    \]
    and
    \begin{align}
        L_+>\max\{2D_n,C_2+2,3D_n,C_*+8D_n+3\} = \max\{C_2+2 , C_*+8D_n+3  \}.
    \end{align}
    and 
    \begin{align}\label{e:C_1}
    C_1>L_++3D_n.
    \end{align}
\end{notation}

In order to get lower and upper bounds for $\mu$ and $\alpha_\mu(Q)$, respectively, we introduce
\[
B_Q^- := B_{p}({\bZ}_Q,L_-\ell(Q)),\qquad
B_Q^+ := B_{p}({\bZ}_Q,L_+\ell(Q)).
\]

\begin{rem}
    If $Q$ is a child of $P$, then as $L_+ > 2D_n$, we get
    \[
    d_{p}({\bZ}_Q,{\bZ}_P)+L_+\ell(Q)
    \le\left(D_n+\frac{L_+}{2}\right)\ell(P)<L_+\ell(P).
    \]
    Thus
    \begin{equation}\label{eq:ambient-balls-nested}
    B_Q^+\subset B_P^+.
    \end{equation}
\end{rem}

For these ambient cubes, we define
\begin{equation}\label{eq:ambient-alpha}
\alpha_\mu(Q) := \frac{1}{\ell(Q)^{n+2}}
\inf_{a\ge0,V}\dist_{B_Q^+}(\mu,a\lambda_V).
\end{equation}

\subsection{The stopping cubes}\label{s:stopping-time}

By the assumptions of Proposition~\ref{prop:alpha-beta}, we can cover
$F$, up to a $\mu$-null set, by countably many measurable sets $K$ with
constants $c_K,C_K,r_0,N>0$ such that
\begin{equation}\label{eq:uniform-density-bounds}
c_Kr^{n+1}\le\mu(B_{p}({\X},r))\le C_Kr^{n+1}
\qquad \forall \, {\X}\in K,\ 0<r\le r_0,
\end{equation}
and
\begin{align}\label{eq:uniform-alpha-bounds}
\int_0^{r_0}\alpha_\mu({\X},r)^2\,\frac{dr}{r}\le N
\qquad \forall \, {\X}\in K.
\end{align}
One can in fact take $c_K$ and $r_0$ to be reciprocals of positive integers, and
$C_K$ and $N$ to be positive integers. Measurability follows from
Appendix~\ref{app:measurability}. We will fix such a $K$ for the remainder of the section. \\ 

\begin{figure}[H]
\centering
\begin{tikzpicture}[font=\small,>=Stealth,line cap=round,line join=round]
  \begin{scope}
    \node at (0,2.5) {(a) $\mathrm{HD}_0$};
    \filldraw[fill=blue!7,draw=blue!60!black,dashed,thick]
      (-1.4,-1.5) rectangle (1.4,1.5);
    \node[blue!60!black,above] at (-.9,1.5) {$B_Q^+$};
    \foreach \p in {(-1.1,-1.2),(-.65,-1.15),(-.2,-1.2),(.3,-1.1),
      (.85,-1.2),(1.15,-.9),(-1.1,-.65),(-.7,-.6),(-.25,-.75),
      (.4,-.7),(.95,-.55),(-1.15,-.05),(-.7,.05),(.6,-.05),
      (1.1,.1),(-1.05,.55),(-.6,.65),(-.15,.45),(.35,.6),
      (.9,.65),(1.15,.4),(-1.1,1.15),(-.65,1.05),(-.1,1.2),
      (.45,1.15),(1.05,1.1)}
      {\fill[red!65!black] \p circle (1.2pt);}
    \fill[black] (0,0) circle (1.2pt);
    \node[below=2pt] at (0,0) {${\bZ}_Q$};
    \draw[->,black] (-1.65,-1.7)--(1.75,-1.7) node[right] {$x$};
    \draw[->,black] (-1.65,-1.7)--(-1.65,1.95) node[above] {$t$};
    \node at (0,-2.25) {$\mu(B_Q^+)\ge M\ell(Q)^{n+1}$};
  \end{scope}
  \begin{scope}[shift={(4.7cm,0)}]
    \node at (0,2.5) {(b) $\mathrm{LD}_0$};
    \filldraw[fill=blue!7,draw=blue!60!black,dashed,thick]
      (-1.4,-1.5) rectangle (1.4,1.5);
    \node[blue!60!black,above] at (-.9,1.5) {$B_Q^-$};
    \foreach \p in {(-.7,.7),(.8,-.6)}
      {\fill[red!65!black] \p circle (1.2pt);}
    \fill[black] (0,0) circle (1.2pt);
    \node[below=2pt] at (0,0) {${\bZ}_Q$};
    \draw[->,black] (-1.65,-1.7)--(1.75,-1.7) node[right] {$x$};
    \draw[->,black] (-1.65,-1.7)--(-1.65,1.95) node[above] {$t$};
    \node at (0,-2.25) {$\mu(B_Q^-)\le\frac{1}{M}\ell(Q)^{n+1}$};
  \end{scope}
  \begin{scope}[shift={(9.4cm,0)}]
    \node at (0,2.5) {(c) $\mathrm{BA}_0$};
    \filldraw[fill=blue!7,draw=blue!60!black,dashed,thick]
      (-1.4,-1.5) rectangle (1.4,1.5);
    \node[blue!60!black,above] at (-.9,1.5) {$3Q$};
    \filldraw[fill=gray!14,draw=black,thick]
      (-.466667,-.166667) rectangle (.466667,.166667);
    \node at (0,0) {$Q$};
    \draw[black,line width=1.5pt] (1.7,-.65)--(1.7,1.2);
    \node[above] at (1.7,1.2) {$K$};
    \draw[->,black] (-1.65,-1.7)--(1.75,-1.7) node[right] {$x$};
    \draw[->,black] (-1.65,-1.7)--(-1.65,1.95) node[above] {$t$};
    \node at (0,-2.25) {$Q\notin\mathrm{HD}_0\cup\mathrm{LD}_0$};
    \node at (0,-2.85) {$3Q\cap K=\varnothing$};
  \end{scope}
\end{tikzpicture}
\caption{The three stopping conditions for $n=1$, with each panel rescaled separately. As usual red dots illustrate $\mu$.}
\label{fig:stopping-tests}
\end{figure}
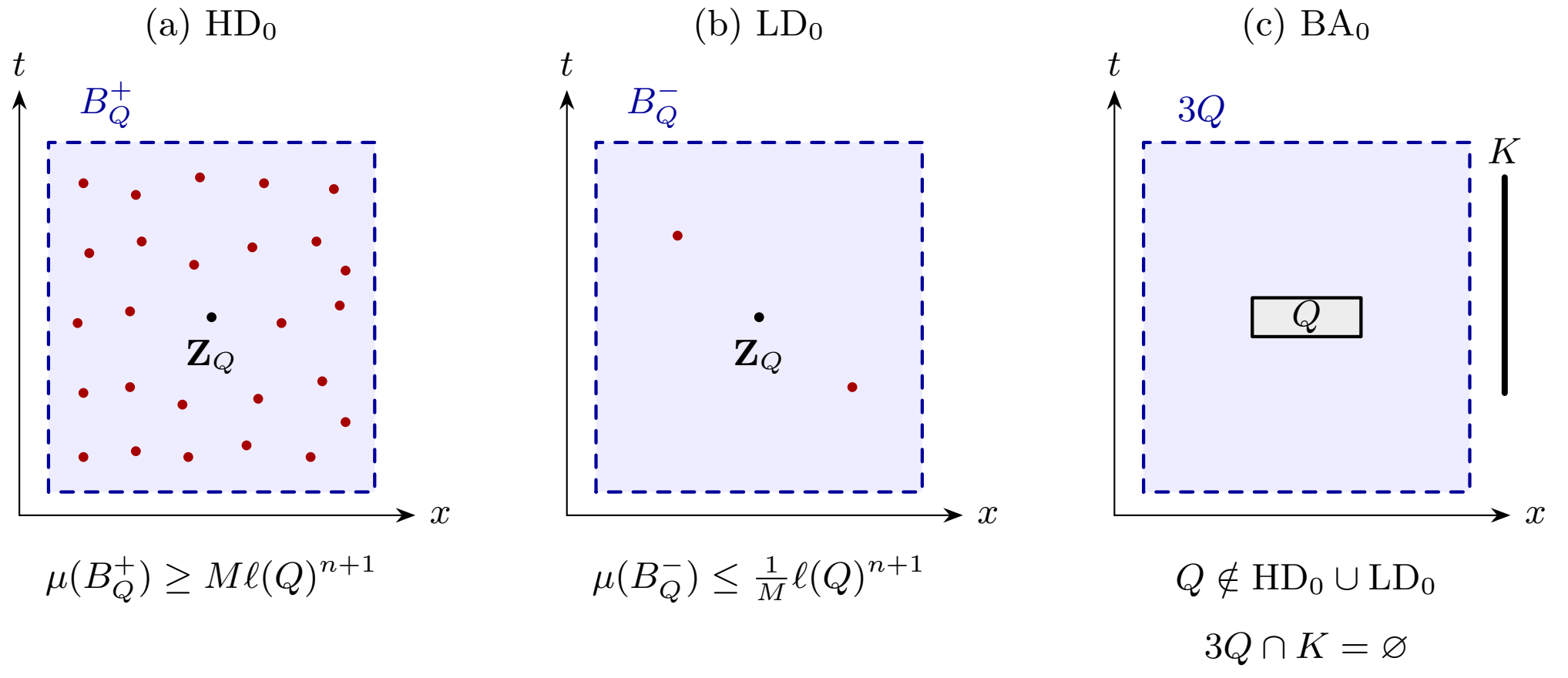

\begin{notation}[Stopping conditions]
    Fix a top cube $R$ meeting $K$ with
    $2C_1\ell(R)<r_0$. Set
    \[
    \mathcal D_{n+1}(3R) := \big\{Q\in\mathcal D_{n+1}:Q\subset3R,
    \ \ell(Q)\le\ell(R)\big\},
    \]
    that is, the union of the trees rooted at the $3^n\cdot9$ cubes of side
    $\ell(R)$ making up $3R$. Now, for $M\ge1$, we define in $\mathcal D_{n+1}(3R)$ the following stopping conditions (see also Figure~\ref{fig:stopping-tests}):
    \begin{itemize}
    \item $Q\in\mathrm{HD}_0$ (high density) if $\mu(B_Q^+)\ge M\ell(Q)^{n+1}$;
    \item $Q\in\mathrm{LD}_0$ (low density)\footnote{
        Using $B_Q^-$ in this $\mathrm{LD}_0$ condition will let us apply the lower density bound at a point of $3Q\cap K$.
    } if $\mu(B_Q^-)\le M^{-1}\ell(Q)^{n+1}$;
    \item $Q\in\mathrm{BA}_0$ (big alphas) if
    $Q\notin\mathrm{HD}_0\cup\mathrm{LD}_0$ and $3Q\cap K=\varnothing$.
    \end{itemize}
    Let $\operatorname{Stop}$ be the family of maximal (and thus disjoint) cubes in
    $\mathrm{HD}_0\cup\mathrm{LD}_0\cup\mathrm{BA}_0$.
\end{notation}

Although we stop according to the high density, lower density and big $\alpha$ condition, because of the assumptions on $K$, it turns out that this is equivalent to stopping whenever a cube is disjoint from $K$. This is the content of the following lemma, which corresponds to \cite[Lemma~3.1]{tolsa2015characterization}.

\begin{lem}[Stopping conditions are equivalent to disjointness with $K$]  \label{lem:stopping-cubes}
If $M$ is big enough, depending only on $c_K$, $C_K$, $n$ and
$L_\pm$, then
\[\text{$Q\in\mathcal D_{n+1}(3R)$ belongs to
$\mathrm{HD}_0\cup\mathrm{LD}_0\cup\mathrm{BA}_0$} 
\quad \iff \quad 
\text{$3Q\cap K=\varnothing$.}\] 
Thus, the stopping cubes are the maximal cubes
$Q\in\mathcal D_{n+1}(3R)$ with $3Q\cap K=\varnothing$. In particular,
\begin{equation}\label{eq:alpha_beta:2-18}
K\cap R\cap\bigcup_{S\in\operatorname{Stop}}S=\varnothing.
\end{equation}
\end{lem}

\begin{proof}
The implication $3Q \cap K = \emptyset$ implies $Q \in \mathrm{HD}_0\cup\mathrm{LD}_0\cup\mathrm{BA}_0$ is easy. Indeed, we may assume that $Q \not\in \mathrm{HD}_0\cup\mathrm{LD}_0$ otherwise there is nothing to show. However, $3Q \cap K = \emptyset$ and $Q \not\in \mathrm{HD}_0\cup\mathrm{LD}_0$ is precisely the definition of $\mathrm{BA}_0$, and we are done. 

We now deal with the implication $Q \in \mathrm{HD}_0\cup\mathrm{LD}_0\cup\mathrm{BA}_0$ implies $3Q \cap K = \emptyset$. Fix some cube in this collection and suppose towards a contradiction that there exists $\bZ \in 3Q \cap K$. Recalling the definition of $\mathrm{BA}_0$, we must have $Q \in \mathrm{HD}_0 \cup \mathrm{LD}_0$.
\begin{itemize}
    \item Consider first the case that $Q\in\mathrm{HD}_0$.
    Since $d_{p}({\bZ},{\bZ}_Q)\le3D_n\ell(Q)$,
    \[
    B_Q^+\subset B_{p}({\bZ},(L_++3D_n)\ell(Q)).
    \]
    The radius is less than $r_0$, so \eqref{eq:uniform-density-bounds} gives
    \[
    M\ell(Q)^{n+1}\le\mu(B_Q^+)
    \le C_K(L_++3D_n)^{n+1}\ell(Q)^{n+1}.
    \]
    However, this is a contradiction if $M>C_K(L_++3D_n)^{n+1}$.

    \item Consider now the case that $Q\in\mathrm{LD}_0$. Since
    $L_->3D_n+2$, we have
    \[
    B_{p}({\bZ},\ell(Q))\subset B_Q^-.
    \]
    Thus, since $\ell(Q)<r_0$, \eqref{eq:uniform-density-bounds} gives
    \[
    c_K\ell(Q)^{n+1}\le\mu(B_{p}({\bZ},\ell(Q)))
    \le\mu(B_Q^-)\le M^{-1}\ell(Q)^{n+1}.
    \]
    This is a contradiction if $M>c_K^{-1}$.
\end{itemize}
Thus, we have finished the proof of the first implication, which proves the description of $\operatorname{Stop}$.
For $S\in\operatorname{Stop}$, we have $S\subset3S$ and
$3S\cap K=\varnothing$. Thus we get
\eqref{eq:alpha_beta:2-18}.
\end{proof}

From now on, take, accordingly with the proof of the last result,
\[
M := 2\max\{1,c_K^{-1},C_K(L_++3D_n)^{n+1}\}
\]
so that the conclusions of Lemma~\ref{lem:stopping-cubes} hold. 

Let $\mathcal G$ be the cubes in
$\mathcal D_{n+1}(3R)$ not contained in a stopping cube, and put
\[
\mathcal G(R) := \big\{Q\in\mathcal G:Q\subset R\big\}.
\]
The cubes of $\mathcal G$ are the good cubes. By the stopping conditions
$\mathrm{HD}_0$ and $\mathrm{LD}_0$ we know that
\begin{equation}\label{eq:alpha_beta:2-9}
M^{-1}\ell(Q)^{n+1}\le\mu(B_Q^-),
\qquad
\mu(B_Q^+)\le M\ell(Q)^{n+1}, 
\qquad 
\text{for all $Q \in \mathcal G$.}
\end{equation}
Moreover, by Lemma~\ref{lem:stopping-cubes},
\begin{align}\label{e:good-intersect-K}
Q\in\mathcal G\quad\Longleftrightarrow\quad3Q\cap K\ne\varnothing.
\end{align}

Finally, we also note that 
\begin{equation}\label{eq:alpha-beta-stop-scale}
S\in\operatorname{Stop},\quad Q\in\mathcal G(R),\quad
S\cap3Q\ne\varnothing
\quad\Longrightarrow\quad\ell(S)\le2\ell(Q).
\end{equation}
Indeed, suppose towards a contradiction that $\ell(S)\ge4\ell(Q)$. The spatial and time sides of $3Q$
are at most $\frac{3}{4}\ell(S)$ and $\frac{9}{16}\ell(S)^2$.
These are less than the respective margins $\ell(S)$ and $4\ell(S)^2$
between $S$ and $3S$. Since $S\cap3Q\ne\varnothing$, we get $3Q\subset3S$.
This contradicts $3Q\cap K\ne\varnothing$ and
$3S\cap K=\varnothing$. \\

As well as the density conditions in \eqref{eq:alpha_beta:2-9}, we know by the stopping condition $\mathrm{BA}_0$ that the square-sum of $\alpha$-numbers of cubes containing a common point is bounded. This is the content of the next lemma, which corresponds to \cite[Lemma~3.2]{tolsa2015characterization}.

\begin{lem}[Alpha square function of the good cubes is controlled]  \label{lem:alpha-sum}
For every ${\X}\in3R$, we have
\begin{align}\label{eq:good-alpha-Jones}
\sum_{\substack{Q\in\mathcal G\\{\X}\in Q}}
\alpha_\mu(Q)^2\le CN,
\end{align}
where $N$ is the uniform bound in \eqref{eq:uniform-alpha-bounds} and $C$ depends only on $n$. In particular, 
\begin{align}\label{eq:good-alpha-carleson}
    \sum_{Q \in \mathcal G} \alpha_\mu(Q)^2\mu(Q) \leq CN \mu(3R). 
\end{align}
\end{lem}

\begin{proof}
The estimate \eqref{eq:good-alpha-carleson} follows directly from \eqref{eq:good-alpha-Jones} by integrating over $3R$, so we need only prove \eqref{eq:good-alpha-Jones}. 
We first estimate the coefficients along the good ancestors of a cube $P \in \mathcal G$. Fix some point ${\bZ}'_P\in3P\cap K$.
For every good cube $Q \in \mathcal G$ which contains $P$, we have ${\bZ}'_P\in3Q$.
Since $3Q\subset B_{p}({\bZ}_Q,3D_n\ell(Q))$ and $C_1>L_++3D_n$, for each 
$r\in[C_1\ell(Q),2C_1\ell(Q)]$ we get
\[
B_Q^+\subset B_{p}({\bZ}'_P,r).
\]
Using this inclusion and $r \approx \ell(Q)$, it also follows that 
\[\alpha_\mu(Q)\le C\alpha_\mu({\bZ}'_P,r).\] 
Each of the intervals $[C_1\ell(Q),2C_1\ell(Q)]$ are contained in $[0,r_0]$ and have disjoint
interiors. Thus, squaring and
integrating over these intervals with respect to $\frac{dr}{r}$, we get
\begin{equation}\label{eq:good-ancestor-sum}
\sum_{\substack{Q\in\mathcal G\\P\subset Q}}
\alpha_\mu(Q)^2
\le C\int_0^{r_0}\alpha_\mu({\bZ}'_P,r)^2\,\frac{dr}{r}\le CN.
\end{equation}

To deal with the sum over all good cubes containing a point ${\X}\in3R$, we proceed as follows. By construction, the good cubes containing ${\X}$ are nested.
If there are finitely many of them and at least one, we apply
\eqref{eq:good-ancestor-sum} to the smallest one.
If there are infinitely many, we write $Q_1\supset Q_2\supset\cdots$.
For each $k$, \eqref{eq:good-ancestor-sum} implies that the sum over
$Q_1,\dots,Q_k$ is bounded by $CN$. Letting $k\to\infty$ gives the lemma.
\end{proof}

%%%%%%%%%%%%%%%%%%%%%%%%%%%%%%%%%%%%%%%%%%%%%%%%%%%%%%%%%%%%%%%%%%
%%%%%%%%%%%%%%%%%%%%%%%%%%%%%%%%%%%%%%%%%%%%%%%%%%%%%%%%%%%%%%%%%%
%%%%%%%%%%%%%%%%%%%%%%%%%%%%%%%%%%%%%%%%%%%%%%%%%%%%%%%%%%%%%%%%%%
%%%%%%%%%%%%%%%%%%%%%%%%%%%%%%%%%%%%%%%%%%%%%%%%%%%%%%%%%%%%%%%%%%
%%%%%%%%%%%%%%%%%%%%%%%%%%%%%%%%%%%%%%%%%%%%%%%%%%%%%%%%%%%%%%%%%%

\subsection{A key estimate bounding angles of nearby planes by alpha numbers} \label{s:angles}

The main result of this subsection is Lemma \ref{lem:consecutive-planes}, in which we control angles between consecutive generation can be controlled by the $\alpha$-numbers. This is relevant because, to estimate $\beta$-numbers (as is our goal, for we strive to prove \eqref{e:finite-jones-beta-main-thm}), we will need to compute distance from points to planes. More concretely, following the strategy of Tolsa in \cite{tolsa2015characterization}, for a point contained in a good cube as above, we will compute the distance to the plane corresponding to some ancestor cube $Q$ by iteratively projecting onto the parent planes. In order to obtain meaningful estimates, we need to show that angles between consecutive generations can be controlled, and this is what we will accomplish in Lemma \ref{lem:consecutive-planes} below. 

Before proving the main result, we need several preliminary lemmas. In the first one we explain how, if we find two vertical planes whose surface measures are very close to each other in the sense of the transport-like distance between measures, then the planes are geometrically very close to each other, too.

\begin{lem}[Comparison of distances between planes in terms of their surface measures]  \label{lem:plane-comparison}
Let $r>0$, $0<R_0<R_1$ and $R_2>R_1+1$. Let $B:=B_{p}({\bZ},R_2r)$, and let
$V_i:=V_{\xi_i,c_i}$ be vertical planes which meet
$B_{p}({\bZ},R_0r)$. Suppose
\[
a_i\ge a_*>0,\qquad i=1,2,
\]
and
\begin{equation}
\dist_B(a_1\lambda_{V_1},a_2\lambda_{V_2})
\le\varepsilon r^{n+2}.
\label{eq:alpha_beta:2-2}
\end{equation}
Then
\begin{equation}
\begin{aligned}
\sup_{{\Y}\in V_1\cap B_{p}({\bZ},R_1r)}
d_{p}({\Y},V_2)
+
\sup_{{\Y}\in V_2\cap B_{p}({\bZ},R_1r)}
d_{p}({\Y},V_1)
\le \frac{C}{a_*}\varepsilon r.
\end{aligned}
\label{eq:alpha_beta:2-3}
\end{equation}
Here $C$ depends only on $n$, $R_0$ and $R_1$.
\end{lem}

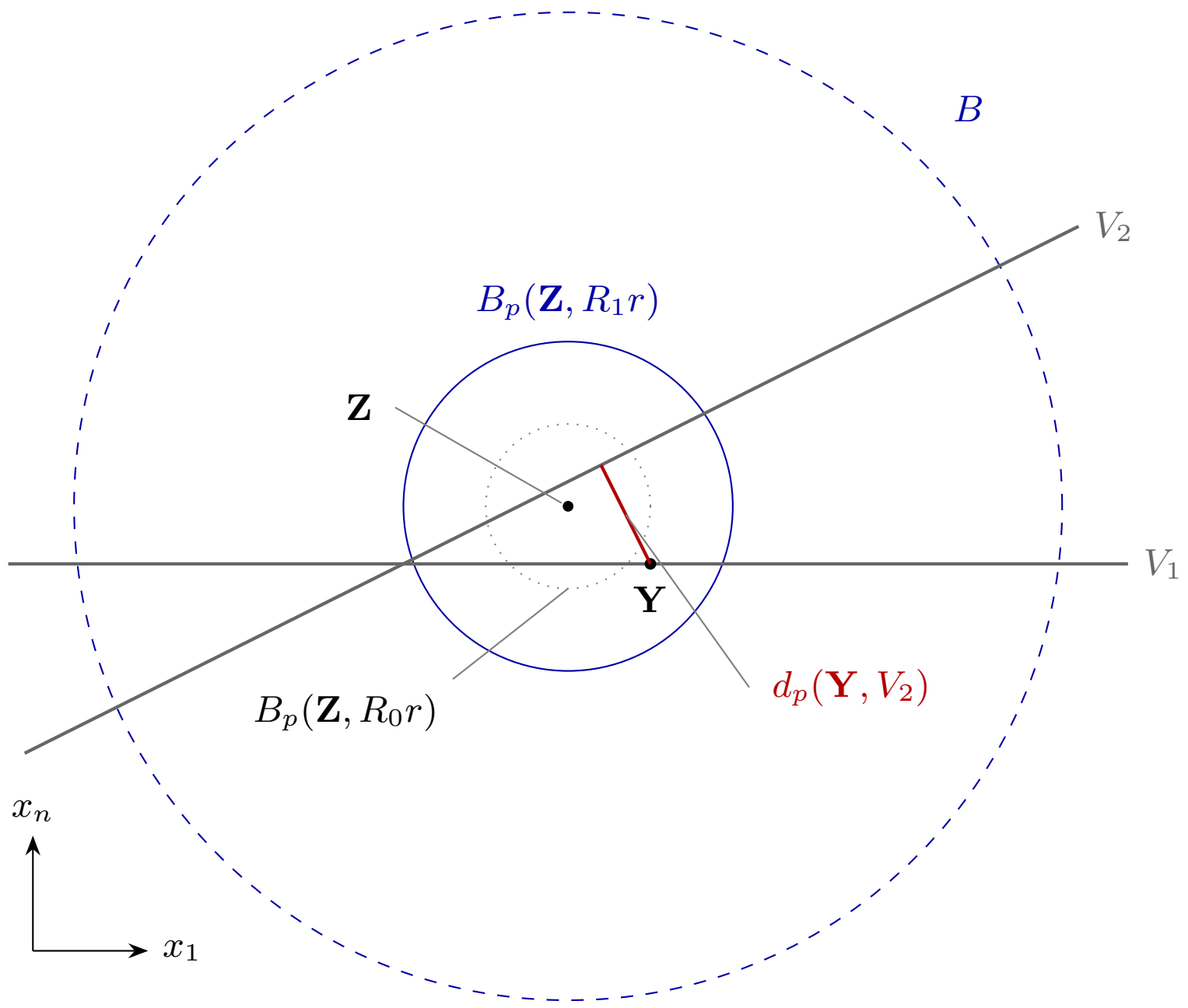
\begin{figure}[H]
\centering
\begin{tikzpicture}[x=1.4cm,y=1.4cm,>=Stealth,font=\small]
\draw[blue!65!black,dashed] (0,0) circle (3);
\draw[blue!65!black] (0,0) circle (1);
\draw[gray,dotted] (0,0) circle (.5);
\node[blue!65!black,anchor=south west] at (2.25,2.25) {$B$};
\node[blue!65!black,anchor=south] at (0,1.04)
{$B_{p}({\bZ},R_1r)$};
\draw[gray] (0,-.5)--(-.7,-1.05);
\node[anchor=north east] at (-.7,-1.05) {$B_{p}({\bZ},R_0r)$};
\draw[gray!80!black,thick] (-3.4,-.35)--(3.4,-.35)
node[right] {$V_1$};
\draw[gray!80!black,thick] (-3.3,-1.5)--(3.1,1.7)
node[right] {$V_2$};
\fill (0,0) circle (1.3pt);
\draw[gray] (-.04,.02)--(-1.05,.6);
\node[anchor=east] at (-1.1,.6) {${\bZ}$};
\fill (.5,-.35) circle (1.4pt) node[below=2pt] {${\Y}$};
\draw[red!70!black,thick] (.5,-.35)--(.2,.25);
\draw[gray] (.35,-.05)--(1.1,-1.1);
\node[red!70!black,anchor=west] at (1.15,-1.1)
{$d_{p}({\Y},V_2)$};
\draw[->,black] (-3.25,-2.7)--(-2.55,-2.7) node[right] {$x_1$};
\draw[->,black] (-3.25,-2.7)--(-3.25,-2.0) node[above] {$x_n$};
\end{tikzpicture}
\caption{The distance test in Lemma~\ref{lem:plane-comparison},
in a fixed-time slice ($n=2$). Here $r=1$, $R_0=\frac12$,
$R_1=1$ and $R_2=3$.}
\label{fig:plane-comparison-test}
\end{figure}

We test the two measures against $f=\chi\,d_{p}(\cdot,V_2)$, where $\chi=1$ on
$B_{p}({\bZ},R_1r)$ and $\chi=0$ outside $B$. In
Figure~\ref{fig:plane-comparison-test}, we see that $f$ vanishes on $V_2$,
while on $V_1$ it is the distance to $V_2$ where $\chi=1$.

\begin{proof}
Choose a smooth cutoff $\chi$ such that 
\[0\leq \chi \leq 1, \quad \text{$\chi=1$ on $B_p(\mathbf Z, R_1r)$}, \quad \supp\chi\subset B_p(\mathbf Z, (R_1+1)r)\subset B \quad \text{and} \quad \Lip_{d_p}\chi\le \frac{C}{r},\] 
and let $f:=\chi\,d_p(\cdot,V_2)$. Since $V_2 \cap B_p(\mathbf Z,R_0r)\neq \phi$, we deduce that $d_p(\cdot,V_2) \le Cr$ on $\supp\chi$, which implies that $\Lip_{d_p}f\le C$. Since $\supp f \subset B$, we conclude that $f/C$ is an admissible function in the definition of $\dist_B$. Testing against $f/C$, and using that fact that $f$ vanishes on $V_2$, we obtain 
\begin{equation}\label{eq:plane-distance-test}
a_1\int_{V_1} \chi(\mathbf Y)d_p(\mathbf Y,V_2)\,d\lambda_{V_1}(\mathbf Y) \le C\varepsilon r^{n+2} . 
\end{equation}

Now write $\mathbf Z=(Z,t_0)$ and $V_1=H_1\times\mathbb R$. Define also
$h(Y)=\xi_2\cdot Y-c_2$ for $Y\in H_1$, so clearly $h$ is affine and \begin{align}\label{e:h=distance}
d_p((Y,t),V_2)=|h(Y)| \quad \mbox{ on } \quad V_1. 
\end{align}
Let us also denote $Z_1$ the orthogonal projection of $Z$ onto $H_1$, so that it holds $|Z-Z_1|=d_p(\mathbf Z,V_1)\le R_0r$ and $|Y-Z|^2 =|Y-Z_1|^2 +|Z_1-Z|^2$ for any $Y \in H_1.$
Moreover, it is also true that
\[
V_1 \cap B_p(\mathbf Z,R_1r)=U\times I,
\]
where 
\begin{align}\label{e:U-times-I}
\hspace{2em} U:= \left\{Y\in H_1: |Y-Z_1| \le \rho := \sqrt{R_1^2r^2-|Z-Z_1|^2} \right\}, 
\qquad 
I:= \left[t_0-R_1^2r^2,t_0+R_1^2r^2\right],
\end{align}
and it holds $\sqrt{R_1^2-R_0^2} \, r \le \rho\le R_1r$ because $|Z-Z_1|\le R_0r$ and $R_0<R_1$.

Now, in \eqref{eq:beta_alpha:2-20} we have $d\lambda_{V_1} =\Lambda\,dY\,dt$ where $dY$ is the Lebsgue measure on $H_1$ and $\Lambda>0$ depends only on $n$. As $\chi=1$ on $U\times I$ and $|U||I| \approx r^{n+1}$, from \eqref{eq:plane-distance-test} we get 
\begin{equation}\label{eq:plane-average}
\frac{1}{|U|}\int_U|h(Y)| dY \leq \frac{C\varepsilon r}{a_1}.   
\end{equation}
Using the inequality $\|A\|_{L^\infty(B)} \lesssim_n |B|^{-1} \|A\|_{L^1(B)}$ whenenever $A$ is an affine function on $\R^n$ and $B$ is a ball in $\R^n$, we obtain from \eqref{eq:plane-average} that
\begin{align}
    \sup_U |h| \leq \frac{C\varepsilon r}{a_1}.
\end{align}
Recalling \eqref{e:h=distance} and \eqref{e:U-times-I}, and using that $a_1 > a_*$, we deduce that 
\[
\sup_{\mathbf Y\in V_1\cap B_p(\mathbf Z,R_1r)}d_p(\mathbf Y,V_2)
=\sup_U|h|\le Ca_*^{-1}\varepsilon r, 
\]
which gives the first estimate in \eqref{eq:alpha_beta:2-3}. Reversing the roles of $V_1$ and $V_2$ gives the other estimate in \eqref{eq:alpha_beta:2-3} and finishes the proof.
\end{proof}

Now we show that, whenever the densities are controlled from above and below, and the alpha numbers are small, as is the case for the good cubes from before, we can always choose a nice approximating vertical plane at the corresponding scale.

\begin{lem}[Choice of approximating planes to $\mu$ at the scales of good cubes]  \label{lem:flat-coefficients}
There exist constants $\tau_M>0$ and $0<c_M\le C_M$, depending only on $M$
and $n$, such that the following holds. Let $Q\in\mathcal D_{n+1}$ satisfy
\eqref{eq:alpha_beta:2-9} and $\alpha_\mu(Q)<\tau_M$. Then there are
$a_Q>0$ and a vertical plane $V_Q$ such that
\begin{equation}
\dist_{B_Q^+}(\mu,a_Q\lambda_{V_Q})
\le2\alpha_\mu(Q)\ell(Q)^{n+2},
\label{eq:alpha_beta:2-10}
\end{equation}
and
\begin{equation}
c_M\le a_Q\le C_M,
\qquad
V_Q\cap B_{p}({\bZ}_Q,C_*\ell(Q))\ne\varnothing.
\label{eq:alpha_beta:2-11}
\end{equation}
\end{lem}

\begin{proof}
Write $\ell=\ell(Q)$. Recall that $\lambda_V$ is a dimensional constant
times Lebesgue measure on $V$. Thus
\[
\lambda_V(B_{p}({\W},r))\approx r^{n+1}\quad\text{for }{\W}\in V,
\qquad
\lambda_V(B_{p}({\bZ},r))\lesssim r^{n+1}\quad\text{for every }{\bZ}.
\]
For the second estimate, if $V$ meets $B_{p}({\bZ},r)$ at a point ${\W}$,
then $B_{p}({\bZ},r)\subset B_{p}({\W},2r)$. Thus, this estimate follows from the upper regularity at points centered on $V$ from the first estimate. \\ 

Let $a\ge0$ and let $V$ be a vertical plane such that
\begin{equation}\label{eq:flat-near-minimizer}
\dist_{B_Q^+}(\mu,a\lambda_V)\le\frac{1}{2M}\,\ell^{n+2}.
\end{equation}
We claim that
\begin{equation} \label{claim:V}
    V\cap B_{p}({\bZ}_Q,(L_-+1)\ell)\ne\varnothing
    \qquad\text{and}\qquad
    c_M\le a\le C_M,
\end{equation}
where $c_M$ and $C_M$ depend only on $M$ and $n$.

To prove the claim, first we write ${\bZ}_Q := (z_Q,t_Q)$ and, for ${\Y}=(y,u)$, we test with
\[
f_1({\Y}) := \min\left\{\ell,\max\left\{0,(L_-+1)\ell-\max\left\{|y-z_Q|,|u-t_Q|^{1/2}\right\}\right\}\right\}.
\]
This function is $1$-Lipschitz, since $|y-z_Q|$ and $|u-t_Q|^{1/2}$ are $1$-Lipschitz
functions of ${\Y}$. It equals $\ell$ on $B_Q^-$ and vanishes outside
$B_{p}({\bZ}_Q,(L_-+1)\ell)$. Thus, it is supported in $B_Q^+$, since $L_-+1<L_+$. Let us verify the different points in \eqref{claim:V}:
\begin{itemize}
    \item By \eqref{eq:alpha_beta:2-9} and
    \eqref{eq:flat-near-minimizer},
    \[
    a\int f_1\,d\lambda_V\ge\ell\,\mu(B_Q^-)-\frac{1}{2M}\,\ell^{n+2}
    \ge\frac{1}{2M}\,\ell^{n+2}.
    \]
    Since the integral is positive and $f_1$ vanishes outside of $B_{p}({\bZ}_Q,(L_-+1)\ell)$, this implies the intersection property of \eqref{claim:V}.

    \item Next, since
    $f_1\le\ell$ and vanishes outside of $B_{p}({\bZ}_Q,(L_-+1)\ell)$, the upper bound for $\lambda_V$
    gives
    \[
    \int f_1\,d\lambda_V
    \le\ell\,\lambda_V\bigl(B_{p}({\bZ}_Q,(L_-+1)\ell)\bigr)
    \lesssim\ell^{n+2}.
    \]
    Thus $a\ge c_M$.

    \item Lastly, to see that $a \leq C_M$, choose a point $W \in V \cap B_{p}({\bZ}_Q,(L_-+1)\ell)$ and consider the test function 
    \[
    f_2 := \min\left\{\ell,\max\left\{0,2\ell-d_{p}(\cdot,{\W})\right\}\right\}.
    \]
    This function is $1$-Lipschitz, equals $\ell$ on $B_{p}({\W},\frac{\ell}{2})$ and vanishes
    outside of $B_{p}({\W},2\ell)$. By \eqref{eq:ball-inclusion} and the choice of constant from Notation \ref{not:constants},
    \[
    B_{p}({\W},2\ell)\subset B_{p}({\bZ}_Q,(L_-+3)\ell)\subset B_Q^+.
    \]
    In particular, $f_2$ is supported in $B_Q^+$. It now follows that 
    \[
    a\,\ell^{n+2}\lesssim a\int f_2\,d\lambda_V
    \le\ell\,\mu(B_Q^+)+\frac{1}{2M}\,\ell^{n+2}\le(M+1)\,\ell^{n+2},
    \]
    which gives $a\le C_M$ and finishes the proof of the claim \eqref{claim:V}. 
\end{itemize} 

Now that we have proved \eqref{claim:V}, take $\tau_M := \frac{1}{4M}$ and fix some $Q \in \cubes_{n+1}$ which satisfies \eqref{eq:alpha_beta:2-9} and $\alpha_\mu(Q)<\tau_M$. If
$\alpha_\mu(Q)>0$, the definition of $\alpha_\mu(Q)$ gives $a_Q\ge0$ and
a vertical plane $V_Q$ satisfying \eqref{eq:alpha_beta:2-10}. Since
\[
2\alpha_\mu(Q)<2\tau_M=\frac{1}{2M},
\]
the estimate \eqref{eq:flat-near-minimizer} holds. Since $L_-+1<C_*$, the claim \eqref{claim:V} gives
\eqref{eq:alpha_beta:2-11}. 

Suppose instead that $\alpha_\mu(Q)=0$. In such case, we choose $a_k\ge0$ and planes
$V_k=V_{\xi_k,c_k}$ with
\[
\dist_{B_Q^+}(\mu,a_k\lambda_{V_k})\to0.
\]
By the claim \eqref{claim:V}, for large $k$ we have 
\[
V_k\cap B_{p}({\bZ}_Q,(L_-+1)\ell)\ne\varnothing
\qquad\text{and}\qquad
c_M\le a_k \le C_M,
\]
Thus, by passing to a subsequence,
we may assume
\[
a_k\to a_Q,\qquad\xi_k\to\xi,\qquad c_k\to c.
\]
Defining $V_Q = V_{\xi,c}$, it is clear that  
\[
V_k\cap B_{p}({\bZ}_Q,(L_-+1)\ell)\ne\varnothing
\qquad\text{and}\qquad
c_M\le a_k \le C_M.
\]
Since $L_-+1<C_*$, \eqref{eq:alpha_beta:2-11} holds also in this case when $\alpha_\mu(Q) = 0$.

It remains to prove \eqref{eq:alpha_beta:2-10}. Since $\xi_k\to\xi$ and
$c_k\to c$, there are isometries $\rho_k$ of $\mathbb R^n$, converging to
the identity uniformly on compact sets, such that the maps
\[
T_k(x,t)=(\rho_k(x),t)
\]
send $V_Q$ onto $V_k$. These maps preserve $d_{p}$, so
$\lambda_{V_k}$ is the push-forward of $\lambda_{V_Q}$ under $T_k$. Let $f$ be
continuous with compact support in $\R^{n+1}$. The functions $f\circ T_k$ are uniformly
bounded, vanish outside a fixed compact set and converge to $f$ pointwise.
By dominated convergence,
\[
\int f\,d\lambda_{V_k}=\int f\circ T_k\,d\lambda_{V_Q}
\longrightarrow\int f\,d\lambda_{V_Q}.
\]
For every $1$-Lipschitz function $f$ with $\supp f\subset B_Q^+$, we let
$k\to\infty$ in
\[
\Bigl|\int f\,d\mu-a_k\int f\,d\lambda_{V_k}\Bigr|
\le\dist_{B_Q^+}(\mu,a_k\lambda_{V_k})
\]
to get $\dist_{B_Q^+}(\mu,a_Q\lambda_{V_Q})=0$, which is
\eqref{eq:alpha_beta:2-10}.
\end{proof}

\begin{notation}[Approximating vertical plane] \label{not:V_Q}
For each cube $Q$ satisfying \eqref{eq:alpha_beta:2-9}, fix
$(a_Q,V_Q)$ as in Lemma~\ref{lem:flat-coefficients} when
$\alpha_\mu(Q)<\tau_M$. If $\alpha_\mu(Q)\ge\tau_M$, let $V_Q$
be the vertical plane through ${\bZ}_Q$ with spatial normal $e_1$. In this case
any vertical plane through ${\bZ}_Q$ would do. In either case, (by \eqref{eq:alpha_beta:2-11} or by definition), we have
\[
V_Q\cap B_{p}({\bZ}_Q,C_*\ell(Q))\ne\varnothing.
\]
\end{notation}

We now demonstrate that the angles between consecutive planes from those chosen above can be controlled by the corresponding $\alpha$-numbers. 

\begin{lem}[Control of angles between consecutive approximating planes by alpha numbers] \label{lem:consecutive-planes}
Let $Q$ be a child of $P$, and suppose both cubes satisfy
\eqref{eq:alpha_beta:2-9}. Then
\begin{equation}
\sup_{{\Y}\in V_Q\cap B_{p}({\bZ}_Q,C_2\ell(Q))}
d_{p}({\Y},V_P)
\le
C_M\bigl(\alpha_\mu(Q)+\alpha_\mu(P)\bigr) \, \ell(P),
\label{eq:alpha_beta:2-13}
\end{equation}
where $C_2$ is the fixed projection radius chosen above.
Now let $Q$ and $Q'$ be cubes satisfying \eqref{eq:alpha_beta:2-9}, with
$Q'\subset3Q$ and $\ell(Q')=\ell(Q)$. Then
\begin{equation}\label{eq:alpha-beta-neighbour-planes}
\sup_{{\Y}\in V_{Q'}\cap B_{p}({\bZ}_{Q'},(C_*+4D_n+2)\ell(Q))}
d_{p}({\Y},V_Q)
\le C_M\bigl(\alpha_\mu(Q)+\alpha_\mu(Q')\bigr)\ell(Q).
\end{equation}
Here $C_M$ depends only on $M$ and $n$.
\end{lem}

\begin{proof}
We begin with  \eqref{eq:alpha_beta:2-13}. Suppose to begin with that either
\[\mbox{$\alpha_\mu(Q)\ge\tau_M$ or $\alpha_\mu(P)\ge\tau_M$}.\] Fix
\[
{\Y}\in V_Q\cap B_{p}({\bZ}_Q,C_2\ell(Q)),\qquad
{\W}\in V_P\cap B_{p}({\bZ}_P,C_*\ell(P))
\]
(where the second point exists by Notation \ref{not:V_Q}). By the triangle inequality,
\[
d_{p}({\Y},V_P)
\le d_{p}({\Y},{\bZ}_Q)+d_{p}({\bZ}_Q,{\bZ}_P)+d_{p}({\bZ}_P,{\W})
\le(C_2+D_n+2C_*)\ell(P).
\]
Then, since $\alpha_\mu(Q)+\alpha_\mu(P)\ge\tau_M$, this gives
\eqref{eq:alpha_beta:2-13} so long as  $C_M\ge(C_2+D_n+2C_*)\tau_M^{-1}$. \\

Suppose instead that \[\mbox{$\alpha_\mu(Q)<\tau_M$ and $\alpha_\mu(P)<\tau_M$.}\]
Put
${\bZ}={\bZ}_Q$, $r=\ell(Q)$ and choose parameters $R_0,R_1$ and $R_2$ (as in Lemma \ref{lem:plane-comparison}) according to 
\[
R_0:=2(C_*+D_n),\qquad R_1:=C_2,\qquad R_2=L_+.
\]
By the choice of the constants in Notation \ref{not:constants}, it is easy to check that $R_0<R_1$ and $R_2>R_1+1$. Also,
$B_{p}({\bZ},R_2r)=B_Q^+$. Since
\[
d_{p}({\bZ}_Q,{\bZ}_P)\le D_n\ell(P)=2D_nr,
\]
both planes $V_Q$ and $V_P$ meet $B_{p}({\bZ},R_0r)$. In summary, the radii $R_0,R_1,R_2$, the ball $B_Q^+$ and the planes $V_Q,V_P$ satisfy the hypotheses of Lemma \ref{lem:plane-comparison}. 

Since $B_Q^+\subset B_P^+$, the
triangle inequality and \eqref{eq:alpha_beta:2-10} give
\[
\dist_{B_Q^+}(a_Q\lambda_{V_Q},a_P\lambda_{V_P})
\le2\alpha_\mu(Q)\ell(Q)^{n+2}+2\alpha_\mu(P)\ell(P)^{n+2}
\le2^{n+3}\bigl(\alpha_\mu(Q)+\alpha_\mu(P)\bigr)r^{n+2}.
\]
Also, $a_Q\ge c_M$ and $a_P\ge c_M$ by \eqref{eq:alpha_beta:2-11}. Thus,
applying Lemma~\ref{lem:plane-comparison} with $a_*=c_M$,
$\varepsilon=2^{n+3}(\alpha_\mu(Q)+\alpha_\mu(P))$ and the above radii, we
get \eqref{eq:alpha_beta:2-13}. \\ 

Now we prove \eqref{eq:alpha-beta-neighbour-planes}. Write
$\ell :=\ell(Q)=\ell(Q')$. Since ${\bZ}_{Q'}\in3Q$, we have
\[
d_{p}({\bZ}_Q,{\bZ}_{Q'})<3D_n\ell.
\]
Suppose that $\alpha_\mu(Q)\ge\tau_M$ or $\alpha_\mu(Q')\ge\tau_M$. Fix
\[
{\Y}\in V_{Q'}\cap B_{p}({\bZ}_{Q'},(C_*+4D_n+2)\ell),\qquad
W\in V_Q\cap B_{p}({\bZ}_Q,C_*\ell).
\]
Then
\[
d_{p}({\Y},V_Q)
\le d_{p}({\Y},{\bZ}_{Q'})+d_{p}({\bZ}_{Q'},{\bZ}_Q)+d_{p}({\bZ}_Q,{\W})
\le(4C_*+11D_n+4)\ell.
\]
Since $\alpha_\mu(Q)+\alpha_\mu(Q')\ge\tau_M$, this gives
\eqref{eq:alpha-beta-neighbour-planes} if
$C_M\ge(4C_*+11D_n+4)\tau_M^{-1}$.

Now suppose that $\alpha_\mu(Q)<\tau_M$ and $\alpha_\mu(Q')<\tau_M$. The
balls $B_Q^+$ and $B_{Q'}^+$ need not be nested, so we use a smaller ball
which lies in both. Put ${\bZ}={\bZ}_{Q'}$, $r=\ell$ and
\[
R_0=C_*+4D_n,\qquad R_1=C_*+4D_n+2,\qquad R_2=L_+-4D_n,
\]
and let $B=B_{p}({\bZ},R_2\ell)$. Again $R_0<R_1$ and $R_2>R_1+1$. By
the bound for $d_{p}({\bZ}_Q,{\bZ}_{Q'})$, both planes meet
$B_{p}({\bZ},R_0\ell)$. Also,
\[
B\subset B_{Q'}^+\cap B_Q^+,
\]
since $(L_+-4D_n)\ell+3D_n\ell<L_+\ell$. Thus, by the triangle inequality
and \eqref{eq:alpha_beta:2-10},
\[
\dist_B(a_{Q'}\lambda_{V_{Q'}},a_Q\lambda_{V_Q})
\le2\bigl(\alpha_\mu(Q)+\alpha_\mu(Q')\bigr)\ell^{n+2}.
\]
Since $a_Q\ge c_M$ and $a_{Q'}\ge c_M$, applying
Lemma~\ref{lem:plane-comparison} with $a_*=c_M$,
$\varepsilon=2(\alpha_\mu(Q)+\alpha_\mu(Q'))$ and the above radii gives
\eqref{eq:alpha-beta-neighbour-planes}.
\end{proof}

We are now at the point where we can estimate the distance of a point $\X$ (which is contained in some good cube) to one of the planes associated to a good ancestor (Lemma \ref{lem:alpha-beta-key}). It will convenient to introduce that notation $\ell_\X$ which quantifies the size of the small good cube containing $\X$ (if it exists). \\  

For ${\X}\in3R$, set
\[
\ell_{\X}\coloneqq\begin{cases}
\ell(S),&{\X}\in S\text{ for some }S\in\operatorname{Stop},\\
0,&{\X}\notin\bigcup_{S\in\operatorname{Stop}}S.
\end{cases}
\]
The following lemma corresponds to \cite[Lemma 3.3]{tolsa2015characterization}. 

\begin{lem}\label{lem:alpha-beta-key}
For $Q\in\mathcal G(R)$ and ${\X}\in3Q$ we have 
\begin{equation}
\begin{aligned}
d_{p}({\X},V_Q)
\le C_M\alpha_\mu(Q)\ell(Q)
+C_M\sum_{\substack{P\in\mathcal G,\ {\X}\in P\subset3Q\\
\ell(P)\le\ell(Q)}}
\alpha_\mu(P)\ell(P)+C\ell_{\X}.
\end{aligned}
\label{eq:alpha_beta:2-20}
\end{equation}
Here $C_M$ depends only on $M$ and $n$, and $C$ depends only on $n$.
\end{lem}

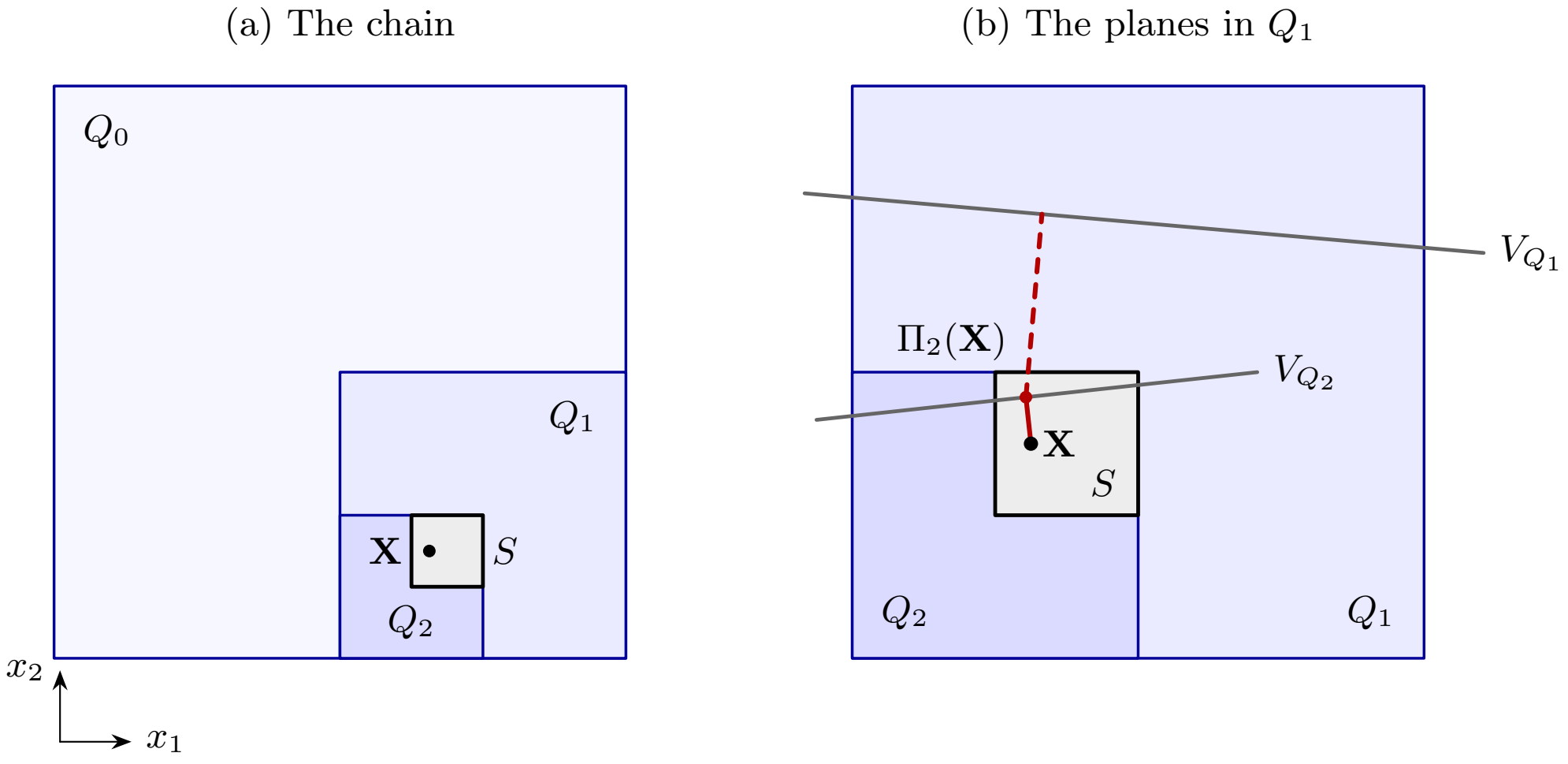
\begin{figure}[H]
\centering
\begin{tikzpicture}[x=1cm,y=1cm,>=Stealth,font=\small,
  line cap=round,line join=round]
  \tikzset{
    good/.style={draw=blue!60!black,line width=.65pt},
    stop/.style={draw=black,fill=gray!14,line width=.9pt},
    plane/.style={draw=gray!80!black,line width=.9pt},
    gap/.style={draw=red!70!black,line width=1.15pt}}

  \begin{scope}
    \node at (2.4,5.3) {(a) The chain};
    \filldraw[good,fill=blue!3] (0,0) rectangle (4.8,4.8);
    \filldraw[good,fill=blue!8] (2.4,0) rectangle (4.8,2.4);
    \filldraw[good,fill=blue!14] (2.4,0) rectangle (3.6,1.2);
    \filldraw[stop] (3.0,.6) rectangle (3.6,1.2);
    \node[anchor=north west] at (.12,4.68) {$Q_0$};
    \node[anchor=north east] at (4.68,2.28) {$Q_1$};
    \node at (3.0,.3) {$Q_2$};
    \node[anchor=west,inner sep=2pt] at (3.6,.9) {$S$};
    \fill (3.15,.9) circle (1.5pt);
    \node[anchor=east,inner sep=2pt] at (3.0,.9) {${\X}$};
    \draw[->,black] (.05,-.7)--(.65,-.7) node[right] {$x_1$};
    \draw[->,black] (.05,-.7)--(.05,-.1) node[left] {$x_2$};
  \end{scope}

  \begin{scope}[shift={(6.7,0)}]
    \node at (2.4,5.3) {(b) The planes in $Q_1$};
    \filldraw[good,fill=blue!8] (0,0) rectangle (4.8,4.8);
    \filldraw[good,fill=blue!14] (0,0) rectangle (2.4,2.4);
    \filldraw[stop] (1.2,1.2) rectangle (2.4,2.4);
    \node[anchor=south east] at (4.68,.12) {$Q_1$};
    \node[anchor=south west] at (.12,.12) {$Q_2$};
    \node[anchor=south east] at (2.34,1.24) {$S$};

    \draw[plane] (-.4,3.9)--(5.3,3.4)
      node[right,text=black] {$V_{Q_1}$};
    \draw[plane] (-.3,2.0)--(3.4,2.4)
      node[right,text=black] {$V_{Q_2}$};

    \coordinate (H) at (1.5,1.8);
    \coordinate (F) at (1.45782,2.19001);
    \coordinate (G) at (1.59246,3.72523);
    \draw[gap] (H)--(F);
    \draw[gap,dashed] (F)--(G);
    \fill (H) circle (1.7pt);
    \fill[red!70!black] (F) circle (1.5pt);
    \node[anchor=west,inner sep=2.5pt] at (H) {${\X}$};
    \node[anchor=south east,inner sep=1.5pt] at (1.36,2.44) {$\Pi_2({\X})$};
  \end{scope}
\end{tikzpicture}
\caption{The proof of Lemma~\ref{lem:alpha-beta-key} when ${\X}$ lies in a
stopping cube $S$, in a fixed-time slice ($n=2$). Here $k=2$, so $Q_2$ is
the parent of $S$. In (b), the solid segment is
$d_{p}({\X},V_{Q_2})$ and the dashed segment is
$d_{p}(\Pi_2({\X}),V_{Q_1})$.}
\label{fig:stopping-plane-chain}
\end{figure}

\begin{proof}
Good cubes satisfy \eqref{eq:alpha_beta:2-9}. Thus, as in the proof of
Lemma~\ref{lem:consecutive-planes}, for every good cube $P$ we have
\[
V_P\cap B_{p}({\bZ}_P,C_*\ell(P))\ne\varnothing,
\qquad\text{and so}\qquad
d_{p}({\bZ}_P,V_P)\le C_*\ell(P).
\]

Fix $Q\in\mathcal G(R)$ and ${\X}\in3Q$. Since
$d_{p}({\X},{\bZ}_Q)\le3D_n\ell(Q)$, we have
\begin{align}\label{e:distance-V_Q}
d_{p}({\X},V_Q)\le(C_*+3D_n)\ell(Q).
\end{align}
Thus, if $\ell_{\X}\ge\ell(Q)$, we have the desired estimate. Now suppose that $\ell_{\X}<\ell(Q)$. We compare the planes
along the chain of dyadic cubes which contain ${\X}$ (see
Figure~\ref{fig:stopping-plane-chain}).

For $j\ge0$, let $Q_j$ be the dyadic cube containing ${\X}$ with
\[
\ell(Q_j)=2^{-j}\ell(Q).
\]
Then $Q_j\subset Q_0\subset3Q\subset3R$, and so $Q_j\in\mathcal D_{n+1}(3R)$.
If a stopping cube contains $Q_j$, then it contains ${\X}$, and so
$\ell(Q_j)\le\ell_{\X}$. Thus,
\[
\ell(Q_j)>\ell_{\X} \implies Q_j\in\mathcal G.
\]
In particular, $Q_0$ is good.

For each good $Q_j$, let $\Pi_j$ be the orthogonal projection onto
$V_{Q_j}$. Since these planes are vertical, these projections fix the time coordinate.
Since ${\X}\in Q_j$, we have by the same reasoning as \eqref{e:distance-V_Q} that 
\[
d_{p}({\X},\Pi_j({\X}))
=d_{p}({\X},V_{Q_j})
\le(C_*+D_n)\ell(Q_j).
\]
Thus,
\[
d_{p}(\Pi_j({\X}),{\bZ}_{Q_j})
\le(C_*+2D_n)\ell(Q_j).
\]
Recalling from Notation \ref{not:constants} that $C_2>C_*+2D_n$, this gives
\[
\Pi_j({\X})\in V_{Q_j}\cap B_{p}({\bZ}_{Q_j},C_2\ell(Q_j)).
\]
Also,
\[
\Pi_0({\X})\in V_{Q_0}\cap B_{p}({\bZ}_{Q_0},(C_*+4D_n+2)\ell(Q)).
\]
In particular, we can use \eqref{eq:alpha_beta:2-13} at $\Pi_j({\X})$ and
\eqref{eq:alpha-beta-neighbour-planes} at $\Pi_0({\X})$.

First, we apply \eqref{eq:alpha-beta-neighbour-planes} with $Q'=Q_0$
to get
\[
\begin{aligned}
d_{p}({\X},V_Q)
&\le d_{p}({\X},\Pi_0({\X}))
+d_{p}(\Pi_0({\X}),V_Q)\\
&\le d_{p}({\X},V_{Q_0})
+C_M\bigl(\alpha_\mu(Q)+\alpha_\mu(Q_0)\bigr)\ell(Q).
\end{aligned}
\]
Now let $j\ge1$ and suppose that $Q_j$ is good. Since its part
$Q_{j-1}$ is also good, we have by \eqref{eq:alpha_beta:2-13} that 
\[
\begin{aligned}
d_{p}({\X},V_{Q_{j-1}})
&\le d_{p}({\X},\Pi_j({\X}))
+d_{p}(\Pi_j({\X}),V_{Q_{j-1}})\\
&\le d_{p}({\X},V_{Q_j})
+C_M\bigl(\alpha_\mu(Q_j)+\alpha_\mu(Q_{j-1})\bigr)
\ell(Q_{j-1}).
\end{aligned}
\]
Since $\ell(Q_{j-1})=2\ell(Q_j)$, summing these inequalities gives
\[
d_{p}({\X},V_Q)
\le d_{p}({\X},V_{Q_k})
+C_M\alpha_\mu(Q)\ell(Q)
+C_M\sum_{j=0}^{k}\alpha_\mu(Q_j)\ell(Q_j)
\]
for any $k$ such that $Q_k$ is good.

If $\ell_{\X}>0$, let $S$ be the stopping cube containing ${\X}$. Since
$\ell(S)=\ell_{\X}<\ell(Q)$, we have $S=Q_m$ for some $m\ge1$. Take $k=m-1$.
Then $\ell(Q_k)=2\ell_{\X}$, so $Q_k$ is good and
\[
d_{p}({\X},V_{Q_k})\le(C_*+D_n)\ell(Q_k)=2(C_*+D_n)\ell_{\X}.
\]
If $\ell_{\X}=0$, every $Q_k$ is good, and
\[
d_{p}({\X},V_{Q_k})
\le(C_*+D_n)2^{-k}\ell(Q)\longrightarrow0.
\]
In this case we let $k\to\infty$. In both cases, the cubes $Q_j$ in the
sum are good, contain ${\X}$, lie in $3Q$ and satisfy $\ell(Q_j)\le\ell(Q)$.
Thus \eqref{eq:alpha_beta:2-20} follows.
\end{proof}

%%%%%%%%%%%%%%%%%%%%%%%%%%%%%%%%%%%%%%%%%%%%%%%%%%%%%%%%%%
%%%%%%%%%%%%%%%%%%%%%%%%%%%%%%%%%%%%%%%%%%%%%%%%%%%%%%%%%%
%%%%%%%%%%%%%%%%%%%%%%%%%%%%%%%%%%%%%%%%%%%%%%%%%%%%%%%%%%
%%%%%%%%%%%%%%%%%%%%%%%%%%%%%%%%%%%%%%%%%%%%%%%%%%%%%%%%%%
%%%%%%%%%%%%%%%%%%%%%%%%%%%%%%%%%%%%%%%%%%%%%%%%%%%%%%%%%%
%%%%%%%%%%%%%%%%%%%%%%%%%%%%%%%%%%%%%%%%%%%%%%%%%%%%%%%%%%

\subsection{Proof of Proposition~\ref{prop:alpha-beta}} \label{s:proof_alpha-beta}

After all the work above, we are ready to show Proposition~\ref{prop:alpha-beta}. Recall that we already showed that this implies the finiteness of the Jones' square function \eqref{e:finite-jones-beta-main-thm} asserted in the corresponding implication of Theorem~\ref{thm:main}, so proving Proposition~\ref{prop:alpha-beta} will effectively finish our proof of Theorem~\ref{thm:main}.

For $Q\in\mathcal D_{n+1}$, we put
\[
\beta_\mu(Q)^2:=\inf_{V}\frac{1}{\ell(Q)^{n+3}}
\int_{3Q}d_{p}({\Y},V)^2\,d\mu({\Y}),
\]
where the infimum is over vertical planes. Notice that
$B_{p}({\X},r)\subset3Q$ for ${\X}\in Q$ and $r\le\ell(Q)$. This is why we
integrate over $3Q$.

\begin{lem}[Summability of beta numbers over good cubes] \label{lem:beta-sum}
With $K$ and $R$ as in the stopping construction (see Subsection~\ref{subsec:preliminaries_partI}), we have
\begin{equation}
\sum_{Q\in\mathcal G(R)}\beta_\mu(Q)^2\mu(Q)
\le C_M(1+N)\mu(3R).
\label{eq:alpha_beta:2-29}
\end{equation}
Here $C_M$ depends only on $M$ and $n$.
\end{lem}

\begin{proof}
Fix $Q\in\mathcal G(R)$ and ${\X}\in3Q$. In the next displays, the
sums run over the cubes $P$ in \eqref{eq:alpha_beta:2-20}, that is,
\[
P\in\mathcal G,\qquad {\X}\in P\subset3Q,\qquad\ell(P)\le\ell(Q).
\]
Since there is at most one such cube of each side length, we have 
\[
\sum_P\ell(P)\le2\ell(Q).
\]
By the Cauchy--Schwarz inequality,
\[
\Bigl(\sum_P\alpha_\mu(P)\ell(P)\Bigr)^2
\le2\ell(Q)\sum_P\alpha_\mu(P)^2\ell(P).
\]
Thus, by applying Lemma~\ref{lem:alpha-beta-key} and the above inequality, we get 
\[
d_{p}({\X},V_Q)^2
\le C_M\alpha_\mu(Q)^2\ell(Q)^2
+C_M\ell(Q)\sum_P\alpha_\mu(P)^2\ell(P)
+C\ell_{\X}^2.
\]

Now we integrate this estimate over $3Q$. Since $L_+>3D_n$, we have
$3Q\subset B_Q^+$. Since $Q$ is good, this gives
\[
\mu(3Q)\le\mu(B_Q^+)\le M\ell(Q)^{n+1}.
\]
By Tonelli's theorem,
\[
\int_{3Q}\sum_P\alpha_\mu(P)^2\ell(P)\,d\mu({\X})
=\sum_{\substack{P\in\mathcal G,\ P\subset3Q\\ \ell(P)\le\ell(Q)}}
\alpha_\mu(P)^2\ell(P)\mu(P).
\]
Also, since the stopping cubes are disjoint,
\[
\int_{3Q}\ell_{\X}^2\,d\mu({\X})
=\sum_{S\in\operatorname{Stop}}\ell(S)^2\mu(S\cap3Q).
\]
We take $V=V_Q$ in the definition of $\beta_\mu(Q)$ and divide by
$\ell(Q)^{n+3}$. We get
\begin{equation}\label{eq:beta-one-cube}
\begin{aligned}
\beta_\mu(Q)^2
\le C_M\alpha_\mu(Q)^2
+C_M\sum_{\substack{P\in\mathcal G,\ P\subset3Q\\ \ell(P)\le\ell(Q)}}
\alpha_\mu(P)^2\,\frac{\ell(P)\mu(P)}{\ell(Q)^{n+2}}
+C\sum_{\substack{S\in\operatorname{Stop}\\ S\cap3Q\ne\varnothing}}
\frac{\ell(S)^2\mu(S)}{\ell(Q)^{n+3}}.
\end{aligned}
\end{equation}

It remains to sum over $Q\in\mathcal G(R)$. Consider for the moment a fixed cubes $P\in\mathcal G$ and $S \in \operatorname{Stop}$. For each side length, there are ony finitely many cubes $Q \in \mathcal G(R)$ such that $P\subset3Q$ and $\ell(Q)\ge\ell(P)$,
then $3Q$ contains the cube of side $\ell(Q)$ which contains $P$. Using this and the fact that
$\mu(Q)\le M\ell(Q)^{n+1}$ (since it is good), we have  
\[
\sum_{\substack{Q\in\mathcal G(R),\ P\subset3Q\\ \ell(Q)\ge\ell(P)}}
\frac{\mu(Q)}{\ell(Q)^{n+2}}
\lesssim\frac{M}{\ell(P)}\sum_{j\ge0}2^{-j}
\lesssim\frac{M}{\ell(P)}.
\]
Recalling from \eqref{eq:alpha-beta-stop-scale} that $\ell(Q)\ge\frac{\ell(S)}{2}$ for any $Q\in\mathcal G(R)$ such that 
$S\cap3Q\ne\varnothing$, as similar argument implies  
\[
\sum_{\substack{Q\in\mathcal G(R)\\ S\cap3Q\ne\varnothing}}
\frac{\mu(Q)}{\ell(Q)^{n+3}}
\lesssim\frac{M}{\ell(S)^2}\sum_{j\ge0}2^{-2j}
\lesssim\frac{M}{\ell(S)^2}.
\]
Thus, multiplying \eqref{eq:beta-one-cube} by $\mu(Q)$, summing over
$Q\in\mathcal G(R)$, changing the order of summation and applying the previous two estimates, we get
\[
\sum_{Q\in\mathcal G(R)}\beta_\mu(Q)^2\mu(Q)
\le C_M\sum_{P\in\mathcal G}\alpha_\mu(P)^2\mu(P)
+C_M\sum_{S\in\operatorname{Stop}}\mu(S).
\]
The first sum is at most $CN\mu(3R)$ by
\eqref{eq:good-alpha-carleson}. Since the stopping cubes are disjoint and lie
in $3R$, the second sum is at most $\mu(3R)$. This gives
\eqref{eq:alpha_beta:2-29}.
\end{proof}

\begin{proof}[\textbf{Proof of Proposition~\ref{prop:alpha-beta}}]
Recalling the start of Section \ref{s:stopping-time}, it suffices to prove the estimate \eqref{e:finite-jones-beta} for $\mu$-a.e. $\X \in K$, for some set $K$ which satisfies \eqref{eq:uniform-density-bounds} and \eqref{eq:uniform-alpha-bounds} with some constants $c_K, C_K, r_0 , N > 0$. We fix such a set and such a collection of constants for the remainder of the proof. We also fix a dyadic side length
$\ell_0<1$ with $2C_1\ell_0<r_0$ (where $C_1$ is as in Notation \ref{not:constants}) and a cube $R \in \cubes_{n+1}$ such that $R \cap K \neq\emptyset$. Note that these choices imply that 
\begin{align}
    2C_1\ell(R) < r_0, 
\end{align}
which is the underlying assumption of the stopping-time construction in Section \ref{s:stopping-time}. Since $K$ is covered, up to a set of $\mu$-measure zero, by the collection of all cubes with a fixed side-length, it suffices to prove that 
\begin{align}\label{e:bounded-Jones-K-cap-R}
\int_0^1\beta_\mu({\X},r)^2\,\frac{dr}{r}<\infty
\qquad\text{for }\mu\text{-a.e. }{\X}\in K\cap R,
\end{align}

We now prove \eqref{e:bounded-Jones-K-cap-R}. Let $\mathcal G(R)$ be the associated collection of good cubes from the above stopping-time construction. Consider a fix ${\X}\in K\cap R$ for the moment, and let $Q\subset R$ be a dyadic cube containing
${\X}$. For each radius
$\frac{\ell(Q)}{2}<r\le\ell(Q)$, we have
$B_{p}({\X},r)\subset3Q$. Thus, 
\[
\beta_\mu({\X},r)^2
\le\Bigl(\frac{\ell(Q)}{r}\Bigr)^{n+3}\beta_\mu(Q)^2
\le2^{n+3}\beta_\mu(Q)^2.
\]
Moreover, since ${\X}\in3Q\cap K$, we have $Q\in\mathcal G(R)$ by \eqref{e:good-intersect-K}. It is therefore easy to see that 
\begin{align}\label{e:continuous-bounded-by-discrete}
\int_0^{\ell(R)}\beta_\mu({\X},r)^2\,\frac{dr}{r}
\le C\sum_{\substack{Q\in\mathcal G(R)\\ {\X}\in Q}}
\beta_\mu(Q)^2.
\end{align}
The sum on the right-hand side of \eqref{e:continuous-bounded-by-discrete} is a measurable function of ${\X}$. Thus, by Tonelli's
theorem and Lemma~\ref{lem:beta-sum},
\begin{align}\label{e:integral-bound}
\int_{K\cap R}
\sum_{\substack{Q\in\mathcal G(R)\\ {\X}\in Q}}\beta_\mu(Q)^2\,d\mu({\X})
\le\sum_{Q\in\mathcal G(R)}\beta_\mu(Q)^2\mu(Q)
\le C_M(1+N)\mu(3R)<\infty.
\end{align}
It now follows from \eqref{e:continuous-bounded-by-discrete} and \eqref{e:integral-bound} that 
\begin{align}\label{e:small-r}
    \int_0^{\ell(R)}\beta_\mu({\X},r)^2\,\frac{dr}{r} < \infty \quad \mbox{ for } \quad  \mbox{$\mu$-a.e. ${\X}\in K\cap R$}. 
\end{align}
For the radii $\ell(R)<r<1$, we use that fact that
\[
\beta_\mu({\X},r)^2
\le\frac{\mu(B_{p}({\X},r))}{r^{n+1}}
\le\frac{\|\mu\|}{r^{n+1}}
\]
to get 
\begin{align}\label{e:large-r}
\int_{\ell(R)}^1\beta_\mu({\X},r)^2\,\frac{dr}{r}
\le\|\mu\|\int_{\ell(R)}^1\frac{dr}{r^{n+2}}<\infty \quad \mbox{ for every } \quad \X \in K \cap R.
\end{align}
The inequality \eqref{e:bounded-Jones-K-cap-R} now follows immediately from \eqref{e:small-r} and \eqref{e:large-r}. 
\end{proof}

\newpage 

\appendix

\section{Measurability of the coefficients}
\label{app:measurability}

\begin{lem}\label{lem:measurability}
Let $\mu$ be a locally finite Borel measure on $\mathbb R^{n+1}$. The functions
\[
({\X},r)\longmapsto\mu(B_{p}({\X},r)),\qquad
({\X},r)\longmapsto\beta_\mu({\X},r),\qquad
({\X},r)\longmapsto\alpha_\mu({\X},r)
\]
are Borel measurable on $\mathbb R^{n+1}\times(0,\infty)$. The lower and
upper $(n+1)$-densities of $\mu$ are Borel measurable. For every $b>0$, the
functions
\[
{\X}\longmapsto\int_0^b\alpha_\mu({\X},r)^2\,\frac{dr}{r},
\qquad
{\X}\longmapsto\int_0^b\beta_\mu({\X},r)^2\,\frac{dr}{r}
\]
are Borel measurable, with values in $[0,\infty]$.
\end{lem}

\begin{proof}
The set $\{({\X},r,{\Y}):{\Y}\in B_{p}({\X},r)\}$ is closed in
$\mathbb R^{n+1}\times(0,\infty)\times\mathbb R^{n+1}$, so its indicator
function is Borel measurable. By Tonelli's theorem, so is
\[
({\X},r)\longmapsto
\int\mathbf 1_{B_{p}({\X},r)}({\Y})\,d\mu({\Y})
=\mu(B_{p}({\X},r)).
\]

Next we consider $\beta_\mu$. Write ${\Y}=(y,t)$. For a vertical plane
$V_{\xi,c}$,
\[
d_{p}({\Y},V_{\xi,c})=|\xi\cdot y-c|.
\]
Thus, for fixed $(\xi,c)$, the function
\[
({\X},r)\longmapsto
\int_{B_{p}({\X},r)}|\xi\cdot y-c|^2\,d\mu({\Y})
\]
is Borel measurable, as above. For fixed $({\X},r)$, it is continuous in
$(\xi,c)$ by dominated convergence, since the ball is bounded. Choose a
countable dense set $\mathcal P_0\subset S^{n-1}\times\mathbb R$. By this
continuity,
\[
\beta_\mu({\X},r)^2
=\frac{1}{r^{n+3}}
\inf_{(\xi,c)\in\mathcal P_0}
\int_{B_{p}({\X},r)}|\xi\cdot y-c|^2\,d\mu({\Y}).
\]
This is a countable infimum of Borel functions, so $\beta_\mu$ is
Borel measurable.

The infimum over all vertical planes is attained. Indeed, if
$\mu(B_{p}({\X},r))=0$, every plane gives
zero. Otherwise, the integral tends to infinity as $|c|\to\infty$,
uniformly in $\xi$. So continuity and compactness give a minimizing plane.

Now we consider $\alpha_\mu$. For $a\ge0$ and
$(\xi,c)\in S^{n-1}\times\mathbb R$, put
\[
D_{a,\xi,c}({\X},r)
=\dist_{B_{p}({\X},r)}(\mu,a\lambda_{V_{\xi,c}}).
\]
First we show that $D_{a,\xi,c}$ is lower semicontinuous in $({\X},r)$.
Suppose $({\X}_k,r_k)\to({\X},r)$, and fix a $1$-Lipschitz function $f$ whose
support is a compact subset of the interior of $B_{p}({\X},r)$. Then
$\supp f\subset B_{p}({\X}_k,r_k)$ for large $k$. Hence
\[
\Bigl|\int f\,d\mu-a\int f\,d\lambda_{V_{\xi,c}}\Bigr|
\le\liminf_{k\to\infty}D_{a,\xi,c}({\X}_k,r_k).
\]
These functions $f$ give the same supremum as all the admissible
functions. Indeed, let $g$ be $1$-Lipschitz with $\supp g\subset B_{p}({\X},r)$.
Then $g$ vanishes on the boundary of $B_{p}({\X},r)$. For $\varepsilon>0$, the
function
\[
g_\varepsilon=\max\{g-\varepsilon,0\}-\max\{-g-\varepsilon,0\}
\]
is $1$-Lipschitz, its support is a compact subset of the interior of
$B_{p}({\X},r)$, and $|g_\varepsilon-g|\le\varepsilon$. Since $\mu$ and
$\lambda_{V_{\xi,c}}$ are finite on $B_{p}({\X},r)$, the integrals of
$g_\varepsilon$ converge to those of $g$ as $\varepsilon\to0$. So,
taking the supremum over $f$ proves the claim. In particular,
$D_{a,\xi,c}$ is Borel measurable.

Next we fix a ball $B$ and show that $\dist_B(\mu,a\lambda_{V_{\xi,c}})$
is continuous in $(a,\xi,c)$. Suppose
\[
a_k\to a,\qquad\xi_k\to\xi,\qquad c_k\to c,
\]
and write $V_k=V_{\xi_k,c_k}$ and $V=V_{\xi,c}$. As in the proof of
Lemma~\ref{lem:flat-coefficients}, there are maps $T_k$ which preserve
$d_{p}$, converge to the identity uniformly on compact sets and
satisfy
\[
\lambda_{V_k}=(T_k)_\#\lambda_V.
\]
Let $f$ be $1$-Lipschitz with $\supp f\subset B$. Then $|f|\le r(B)$, and
$f\circ T_k-f$ vanishes outside $2B$ for large $k$. Since $f$ is
$1$-Lipschitz,
\[
\Bigl|\int f\,d\lambda_{V_k}-\int f\,d\lambda_V\Bigr|
=\Bigl|\int(f\circ T_k-f)\,d\lambda_V\Bigr|
\le\lambda_V(2B)\sup_{{\Y}\in2B}d_{p}(T_k({\Y}),{\Y}).
\]
The right-hand side does not depend on $f$, and it tends to zero. Thus
\[
\dist_B(a_k\lambda_{V_k},a\lambda_V)
\le a_k\lambda_V(2B)\sup_{{\Y}\in2B}d_{p}(T_k({\Y}),{\Y})
+|a_k-a|\,r(B)\lambda_V(B)\longrightarrow0.
\]
By the triangle inequality,
\[
\bigl|\dist_B(\mu,a_k\lambda_{V_k})-\dist_B(\mu,a\lambda_V)\bigr|
\le\dist_B(a_k\lambda_{V_k},a\lambda_V)\longrightarrow0.
\]

It follows that
\[
\alpha_\mu({\X},r)
=\frac{1}{r^{n+2}}
\inf_{\substack{a\in\mathbb Q,\ a\ge0\\ (\xi,c)\in\mathcal P_0}}
D_{a,\xi,c}({\X},r).
\]
This is a countable infimum of Borel functions, so $\alpha_\mu$ is Borel
measurable. By Tonelli's theorem, the two square functions are Borel
measurable.

Finally, recall that
\[
\Theta_\mu({\X},r)=\frac{\mu(B_{p}({\X},r))}{r^{n+1}}.
\]
If $r_j\downarrow r$, then
$B_{p}({\X},r_j)\downarrow B_{p}({\X},r)$, and so
$\Theta_\mu({\X},r_j)\to\Theta_\mu({\X},r)$, since $\mu$ is finite on bounded sets.
Therefore we may use rational radii in the densities:
\[
\Theta^{n+1}_*({\X},\mu)
=\lim_{k\to\infty}
\inf_{\substack{q\in\mathbb Q\\ 0<q<1/k}}\Theta_\mu({\X},q),
\qquad
\Theta^{n+1,*}({\X},\mu)
=\lim_{k\to\infty}
\sup_{\substack{q\in\mathbb Q\\ 0<q<1/k}}\Theta_\mu({\X},q).
\]
Both functions are Borel measurable.
\end{proof}

We use the lemma as follows. The set $F$ in the proof of
Theorem~\ref{thm:main} is Borel, since the graphs are closed and the
densities are Borel measurable. For the sets $K$ satisfying
\eqref{eq:uniform-density-bounds}, we may take $r_0$ rational. Then the
density bounds can be checked at rational radii in $(0,r_0]$, since a
bound at rational $q>r$ passes to $r$ by letting $q\downarrow r$. Together with the measurability of the $\alpha$-square
function, this shows that these sets can be taken Borel.

\newpage

\bibliography{Ref}
\bibliographystyle{alpha}
\end{document}